\documentclass[reqno,fleqn]{amsart}
\usepackage{amsfonts,amsthm,amssymb}
\usepackage{lipsum}
\usepackage[utf8]{inputenc}
\usepackage[english]{babel}
\usepackage{lmodern,textcomp}
\usepackage{amsmath}
\usepackage{mathtools}
\usepackage{latexsym}
\usepackage{tikz}
\usepackage{fancyvrb}
\usepackage{epsfig}
\usepackage{pstricks,slashbox,multirow}
\usepackage{graphicx}
\usepackage{rotating}
\usetikzlibrary{positioning}
\usepackage{subcaption}
\usepackage{datetime}
\newtheorem{theorem}{Theorem}

\newtheorem{lemma}[theorem]{Lemma}

\newtheorem{illu}[theorem]{Illustration}

\newtheorem{cor}[theorem]{Corollary}

\newtheorem{conj}{Conjecture}

\newtheorem{dfn}{Definition}
\newtheorem{prm}{Problem}
\newtheorem{oprm}{Open Problem}

\newtheorem{rem}[theorem]{Remark}
\newtheorem{res}[theorem]{Result}
\newtheorem{note}{Note}

\title[Harary's Sum and Integral Sum Graphs - A Survey in detail]{Harary's Sum and Integral Sum Graphs - A Survey in detail}

\author{\sc  V. Vilfred Kamalappan$^1$ and Lowell W. Beineke$^2$} 

\address{$^1$Department of Mathematics (Formerly), Central University of Kerala, Periye, Kasaragod,   \linebreak Kerala, India - 671 316.} 
\email{vilfredkamal@gmail.com}

\address{$^2$Purdue~ University, Fort~ Wayne, Indiana~ 46805, U.S.A.}
\email{beineke@pfw.edu}

\dedicatory{Dedicated to the memory of \\ Professor Frank Harary [11.03.1921 to 04.01.2005]}

\date{}

\begin{document}

\begin{abstract}
Harary introduced the concepts of sum and integral sum graphs. A graph $G$ is a {\em sum graph} if the vertices of $G$ can be labeled with distinct positive integers so that $e = uv$ is an edge of $G$ if and only if the sum of the labels on vertices $u$ and $v$ is also a label in $G$. An {\em integral sum graph} is also defined just as sum graph, the difference being that the labels may be any distinct integers. This survey presents results obtained by different authors, in detail, on sum and integral sum graphs.
\end{abstract}

\subjclass[2010]{05C78, 05C15, 05C75.}
\keywords{Integral sum graph, sum graph, edge sum class, integral sum edge chromatic number, covering number, independence number, chromatic number, clique number, perfect graph.}

\maketitle


\tableofcontents
	
\section{Introduction}

A {\em graph labeling} is an assignment of integers to the vertices or edges, or both, subject to certain conditions \cite{g26}. 

\subsection{Historical Note} 
 
Interest in graph labeling problems seems to emanate in the mid 1960’s from a long standing conjecture of Ringel \cite{ri64} on graceful labeling that every tree is graceful and a paper by Rosa \cite{ro67}. In the intervening years more than 350 graph labelings techniques have been studied in well over 3600 papers \cite{g26} and by which labeling graphs has become a multidimensional problem. Despite the large number of papers, there are relatively few general results or methods on graph labelings. Whereas in the study of sum and integral sum labeling of graphs, many general results are obtained and we bring out these in this survey. Most of these results are contributed by Beineke, Chen, Frank Harary, Gallian, Julia K. Abraham, Kala, Mary Florida, Nicholas, Rubin Mary, Sajidha, Shine Raj, Somasundaram, Suryakala, Vilfred \cite{c96, c06, g26, h90, h94, js24, ns01, nv02, s21}, \cite{vb14}-\cite{vs14}. 

In 1990, Harary \cite{h90} introduced the notion of a sum graph. A graph $G(V, E)$ is called a {\em sum graph} if there is a bijection $f$ from $V$ to a set of positive integers $S$ such that $xy\in E(G)$ if and only if $f(x) + f(y) \in S$, $x,y\in V(G)$. In 1994, Harary \cite{h94} generalized sum graphs by permitting $S$ to be any set of integers and called these graphs as {\em integral sum graphs}. To distinguish between the two types, we call sum graphs that use only positive integers ${\mathbb N}$-\textit{sum graphs} and those with any integers ${\mathbb Z}$-\textit{sum graphs}. In a sum graph, vertex with the highest label cannot be adjacent to any other vertex and thereby every sum graph must contain isolated vertices. Whereas integral sum graphs need not have isolated vertices.  

Labeled graphs serve as useful models for a broad range of applications such as X-ray crystallography, radar, astronomy, coding theory, circuit design, data base management, secret sharing schemes, cryptology, ultrasound screens, network passwords, communication network addressing, image processing schemes for monitoring air quality \cite{g26}. Slamet et al. \cite{ss06} show an application of sum graph labelings to distribute secret information to a set of people so that only an authorized set of people can reconstruct the secret. Sutton \cite{s01}, in his Ph.D. thesis, introduced two methods of graph labelings that generalize the notion of sum graphs and have applications to storage and manipulation of relational database. 

\subsection{Terminology and Notation}

All terms not defined here can be found in Harary \cite{h69}. Throughout this article, we consider only simple undirected graph $G(V, E)$ with vertex set $V(G)$ and edge set $E(G)$. For additional reading on related graph labeling problems, please refer \cite{g26}. 

For any non-empty set of integers $S$, $G^+(S)$ denotes {\em the integral sum graph on the set $S$}. Harary \cite{h94} also introduced families of sum graphs $G_n = G^+([1, n])$ and integral sum graphs $G_{-n, n} = G^+([-n,n])$ where $n\in \mathbb{N}$. The extension of Harary graphs to all intervals of integers was introduced by Vilfred in \cite{VFS} and \cite{vm12}: for any integers $r$ and $s$ with $-r \leq 0 < s$, let $G_{-r, s} = G^{+}([-r, s])$ where $[r, s]$ = $\{r,r + 1,\dots,s\}$. Here, we present some basic definitions and notations which are used in this survey.	

\begin{dfn}	 \cite{h90} A graph $G$ is a {\em sum graph} or {\em  $\mathbb{N}$-sum graph} if the vertices of $G$ can be labeled with distinct positive integers so that $uv$ is an edge of $G$ if and only if the sum of the labels on vertices $u$ and $v$ is also a label in $V(G).$
\end{dfn}

\begin{dfn}	 \cite{h94} An {\em integral sum graph} or {\em $\mathbb{Z}$-sum graph} is also defined just as sum graph, the difference being that the labels may be any distinct integers. 
\end{dfn} 

\begin{dfn}	 \cite{d05} \quad The \emph{join} of two graphs $A$ and $B$, denoted by $A * B$, is the graph $A \cup B$ together with edges joining each vertex of $A$ with all the vertices of $B$.
\end{dfn} 

\vspace{.1cm}
The following notations are used to keep formulas relatively brief \cite{vb14}:

\begin{enumerate}
	
	\item $n$ will always denote a positive integer, and $G_n$ is the sum graph $G^+([1,n])$.
	
	\item We denote the number of edges of a graph $G$ by $||G||$ = $|E(G)|$. 
	
	\item We denote the sum graph $G^+([1, n])$ by $G^+_n$ when it is labeled and by $G_n$ when it is unlabeled.
	
	\item Two unlabeled graphs are said to be {\em comparable} if one is a subgraph of the other, while two labeled graphs are {\em comparable} if one is a subgraph of the other with the labels preserved.  
	
	Clearly, any two Harary graphs $G_m$ and $G_n$ are comparable, $m,n\in\mathbb{N}$.  In contrast, it is easy to check that labeled graphs $G^+_3$ and $G^+_2 * G^{+}(\{3\})$, given in figures 1 and 2, are not comparable even though as unlabeled graph  $G_3$ is a (spanning) subgraph of unlabeled graph $G_2 * G^{+}(\{3\})$ which are given in figures 3 and 4. 
\end{enumerate}

\medskip
\begin{center}
	\begin{tikzpicture}
		\node (a1) at (9,1.5)  [circle,draw,scale=0.6]{1};
		\node (a2) at (10,0.25)  [circle,draw,scale=0.6]{2};
		\node (a3) at (8.5,-0.25)  [circle,draw,scale=0.6]{3};
		
		\draw (a1) -- (a2);
		
		\node (b4) at (11.5,1.5)  [circle,draw,scale=0.6] {1};
		\node (b5) at (11.5,0)  [circle,draw,scale=0.6]{2};
		\node (b6) at (13,.5)  [circle,draw,scale=0.6] {3};
		
		\draw (b4) -- (b6);
		\draw (b5) -- (b6);
		
		\node (a1) at (15.5,1.5)  [circle,draw,scale=0.6]{~ ~};
		\node (a2) at (16.5,0.25)  [circle,draw,scale=0.6]{~ ~};
		\node (a3) at (15,-0.25)  [circle,draw,scale=0.6]{~ ~};
		
		\draw (a1) -- (a2);
		
		\node (b4) at (18,1.5)  [circle,draw,scale=0.6] {~ ~};
		\node (b5) at (18,0)  [circle,draw,scale=0.6]{~ ~};
		\node (b6) at (19.5,.5)  [circle,draw,scale=0.6] {~ ~};
		
		\draw (b4) -- (b6);
		\draw (b5) -- (b6);
		
	\end{tikzpicture}
	
	\vspace{.1cm}
	{\small	\hspace{.5cm} Fig. 1. $G^+_{3}$ \hspace{.5cm} Fig. 2. $G^+_2 * G^+(\{3\})$ \hspace{.75cm}  Fig. 3. $G_{3}$ \hspace{.3cm} Fig. 4. $G_2 * G^+(\{3\})$ }
	
\end{center}

\medskip
Note that $G^+([-s,-1]) \cong G^+([1,s])$ so that a labeling with only negatives is the same as one with only positives, $s\in\mathbb{N}$.  Furthermore, if every label in a sum graph is replaced by its negative, then the two graphs are isomorphic (see page 17 in \cite{ls21}). 

Vilfred characterized integral sum graphs $G_{-r, s}$ as follows (see \cite{ls21, vm12}):

\begin{theorem} [{\cite[Theorem 3.4(a)]{ls21}}] \label{1.2.4}
	If $r, s \in \mathbb{Z}$ with  $r < 0 < s$, then $G_{r, s}$ = $K_1 * (G_{r} * G_s)$. \hfill $\Box$	
\end{theorem} 

The graph operation of the join, which we denote here by $*$, is both associative and commutative \cite{vm12}. Integral sum graphs $G_{0,6}$ = $K_1(0)*G_6$, $G_{-1,5}$ = $K_1(0)*G_{-1}*G_5$, $G_{-2,4}$ = $K_1(0)*G_{-2}*G_4$, $G_{-3,3}$ = $K_1(0)*G_{-3}*G_3$ are given in figures 5 to 8.

\medskip	
\begin{center}
	\begin{tikzpicture}[scale = 1]
		
		\node (a6) at (8,1.5)  [circle,draw,scale=0.6] {6};
		\node (a0) at (9,2)  [circle,draw,scale=0.6, fill = green!30]{0};
		\node (a1) at (10,1.5)  [circle,draw,scale=0.6]{1};
		\node (a2) at (10.5,.5)  [circle,draw,scale=0.6] {2};
		\node (a3) at (10,-.5)  [circle,draw,scale=0.6] {3};
		\node (a4) at (9,-.5)  [circle,draw,scale=0.6] {4};
		\node (a5) at (8,.5)  [circle,draw,scale=0.6]{5};
		
		\draw (a0)[green] -- (a1);
		\draw (a0)[green] -- (a2);
		\draw (a0)[green] -- (a3);
		\draw (a0)[green] -- (a4);
		\draw (a0)[green] -- (a5);
		\draw (a0)[green] -- (a6);
		
		\draw (a1) -- (a2);
		\draw (a1) -- (a3);
		\draw (a1) -- (a4);
		\draw (a1) -- (a5);
		
		\draw (a2) -- (a3);
		\draw (a2) -- (a4);
		\node (b6) at (12,1.5)  [circle,draw,scale=0.6, fill = blue!30] {-1};
		\node (b0) at (13,2)  [circle,draw,scale=0.6, fill = green!30]{0};
		\node (b1) at (14,1.5)  [circle,draw,scale=0.6]{1};
		\node (b2) at (14.5,.5)  [circle,draw,scale=0.6] {2};
		\node (b3) at (14,-.5)  [circle,draw,scale=0.6] {3};
		\node (b4) at (13,-.5)  [circle,draw,scale=0.6] {4};
		\node (b5) at (12,0)  [circle,draw,scale=0.6]{5};
		
		\draw (b0)[green] -- (b1);
		\draw (b0)[green] -- (b2);
		\draw (b0)[green] -- (b3);
		\draw (b0)[green] -- (b4);
		\draw (b0)[green] -- (b5);
		\draw (b0)[green] -- (b6);
		
		\draw (b6)[red] -- (b1);
		\draw (b6)[red] -- (b2);
		\draw (b6)[red] -- (b3);
		\draw (b6)[red] -- (b4);
		\draw (b6)[red] -- (b5);
		
		\draw (b1) -- (b2);
		\draw (b1) -- (b3);
		\draw (b1) -- (b4);
		
		\draw (b2) -- (b3);
		
	\end{tikzpicture}
	
	\medskip
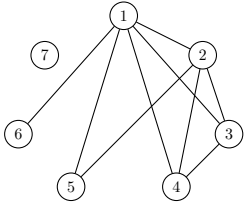
	{\small Fig. 5. $G_{0,6}$ = $K_1(0)*G_6$  \hspace{1cm} Fig. 6. $G_{-1,5}$ = $K_1(0)*G_{-1}*G_5$ }
	
\end{center}

\medskip	
\begin{center}
	\begin{tikzpicture}[scale = 1]
		
		\node (c6) at (16,1.5)  [circle,draw,scale=0.6, fill = blue!30] {-1};
		\node (c0) at (17,2)  [circle,draw,scale=0.6, fill = green!30]{0};
		\node (c1) at (18,1.5)  [circle,draw,scale=0.6]{1};
		\node (c2) at (18.5,.5)  [circle,draw,scale=0.6] {2};
		\node (c3) at (18,-.5)  [circle,draw,scale=0.6] {3};
		\node (c4) at (17,-.5)  [circle,draw,scale=0.6] {4};
		\node (c5) at (16,.5)  [circle,draw,scale=0.6, fill = blue!30]{-2};
		
		\draw (c0)[green] -- (c1);
		\draw (c0)[green] -- (c2);
		\draw (c0)[green] -- (c3);
		\draw (c0)[green] -- (c4);
		\draw (c0)[green] -- (c5);
		\draw (c0)[green] -- (c6);
		
		\draw (c6)[red] -- (c1);
		\draw (c6)[red] -- (c2);
		\draw (c6)[red] -- (c3);
		\draw (c6)[red] -- (c4);
		
		\draw (c5)[red] -- (c1);
		\draw (c5)[red] -- (c2);
		\draw (c5)[red] -- (c3);	
		\draw (c5)[red] -- (c4);
		
		\draw (c1) -- (c2);
		\draw (c1) -- (c3);	
		
		\node (d0) at (21.5,2)  [circle,draw,scale=0.6, fill = green!30]{0};
		\node (d1) at (22.5,1.5)  [circle,draw,scale=0.6]{1};
		\node (d2) at (23,.5)  [circle,draw,scale=0.6] {2};
		\node (d3) at (22.5,-.5)  [circle,draw,scale=0.6] {3};
		\node (d4) at (20.5,-.5)  [circle,draw,scale=0.6, fill = blue!30] {-3};
		\node (d5) at (20,.5)  [circle,draw,scale=0.6, fill = blue!30]{-2};
		\node (d6) at (20.5,1.5)  [circle,draw,scale=0.6, fill = blue!30] {-1};
		
		\draw (d0)[green] -- (d1);
		\draw (d0)[green] -- (d2);
		\draw (d0)[green] -- (d3);
		\draw (d0)[green] -- (d4);
		\draw (d0)[green] -- (d5);
		\draw (d0)[green] -- (d6);
		
		\draw (d6)[red] -- (d1);
		\draw (d6)[red] -- (d2);
		\draw (d6)[red] -- (d3);
		\draw (d5)[red] -- (d1);
		\draw (d5)[red] -- (d2);
		\draw (d5)[red] -- (d3);
		
		\draw (d4)[red] -- (d1);
		\draw (d4)[red] -- (d2);
		\draw (d4)[red] -- (d3);
		
		\draw (d6)[blue] -- (d5);
		
		\draw (d1) -- (d2);
		
	\end{tikzpicture}
	
	\medskip
	{\small Fig. 7. $G_{-2, 4}$ = $K_1(0)*G_{-2}*G_4$ \hspace{.5cm} Fig. 8. $G_{-3,3}$ = $K_1(0)*G_{-3}*G_3$ }
\end{center}

\section{Structural properties of Harary's graphs $G_n$ and $G_{-n,n}$ and graphs $G^c_n$, $G_{0, s}$, $G_{-r, s}$}  
	
Here, we present different structural properties of integral sum graphs $G^+(S)$. In particular, structural properties of 
Harary's graphs $G_n$ and $G_{-n, n}$, and also of graphs $G^c_n$, $G_{0, n}$, $G_{-m, n}$, $m,n\in\mathbb{N}$.  

\subsection{Degree of vertices and the number of edges in $G_n$ and $G_{-r, s}$}

It is observed \cite{ls21} in the Harary graph $G_n$, the vertex labeled 1 has degree $n-2$, and in going from the vertex labeled $i$ to that labeled $i+1$, the degree
goes down by 1, with one exception, when it remains the same. That exception is for $i$ = $k$ when $n$ = $2k$ as well as for $i$ = $k+1$ when $n$ = $2k+1$. Figures 9 and 12 illustrate this for the cases $n$ = 7 and 8. By summing these degrees, we can deduce the number of edges in $G_n$. These facts are summarized in the following results.
	
	\begin{theorem} [{\cite[Theorem 3.2(b)]{ls21}}]  \label{a2}
	{\rm 	The degree of the vertex with label $i$ in $G^+_n$ is 
		\[ \hspace{1cm} \deg (i) = \left\{ \begin{array}{ll} 
		n-i-1 & \mbox{if $1 \leq i \leq \lfloor\frac{n}{2}\rfloor$,} \\
		n-i & \mbox{if $\lfloor\frac{n}{2}\rfloor+1 \leq i \leq n$. \hspace{6.5cm} $\Box$	 } 
	\end{array} \right. 
		\] }
	\end{theorem}	

	Thus, the degree sequence of $G_{2k}$ is $2k-2, 2k-3, \dots, k-1$, $k-1$, $\dots, 1, 0$, while that of $G_{2k+1}$ is $2k-1, 2k-2, \dots, k, k, \dots, 1, 0$. By adding these numbers, we find the number of edges in $G_n$.  The graph $G_{2k}$ has $k(k-1)$ edges and the graph $G_{2k+1}$ has $k^2$ edges.  Combining these results, we get the next theorem.

	\begin{theorem} [{\cite[Theorem 3.2(c)]{ls21}}]  \label{a3} 
	{\rm 	The graph $G_n$	has $\left\lfloor (n-1)^2/4 \right\rfloor$ edges.} \hfill $\Box$
	\end{theorem}	
	
	We now turn to the degrees of vertices in the interval graphs $G_{-r,s}$, $r\in\mathbb{N}_0$ and $s\in\mathbb{N}$. 	

\begin{theorem}[{\cite[Theorem 3.4(c)]{ls21}}]  \label{a4}
	{\rm In $G_{-r,s}$ of order $n$, the degree of the vertex with label $i$ is 
		\[ \deg ~i~ = \left\{ \begin{array}{ll}
		n+i & \mbox{if $-r \leq i \leq -\lfloor\frac{r}{2}\rfloor-1$,} \\
		 n+i-1  & \mbox{if $-\lfloor\frac{r}{2}\rfloor \leq i \leq -1 $,} \\
		n-1  & \mbox{if $i = 0$,}\\
		n-i-1 & \mbox{if $1 \leq i \leq \lfloor\frac{s}{2}\rfloor$,}\\
		n-i & \mbox{if $\lfloor\frac{s}{2}\rfloor+1 \leq i \leq s$ ~\text{where}~ $n$ = $r+s+1$. \hspace{1cm} $\Box$}
		\end{array} \right.	
		\] } 
\end{theorem}	

	By combining Theorem \ref{a3} with the property $G_{-r,s} = K_1(0)*(G_{r}*G_s)$, we find $m(-r, s) = ||G_{-r,s}||$, the number of edges in $G_{-r,s}$ as follows.
	
	\begin{theorem} [{\cite[Theorem 2.4]{vm12}}]  \label{a5}
	{\rm 	For $r,s\in\mathbb{N}$, the number of edges in 	$G_{-r,s}$ is 
		
\hspace{1cm}	$m(-r, s)$ = $||G_{-r,s}||$ = $\frac{1}{4}(r^2+4rs+s^2+3(r+s))$ - $\frac{1}{2}(\left\lfloor\frac{r}{2}\right\rfloor + \left\lfloor\frac{s}{2}\right\rfloor)$

\hspace{4.2cm}	 =	$rs + r + s + \left\lfloor\frac{(r-1)^2}{4}\right\rfloor + \left\lfloor\frac{(s-1)^2}{4}\right\rfloor$. } 
\end{theorem}	
\begin{proof} First note that $|E(G_n)|$ = $n(n-1)/4-\left\lfloor n/2 \right\rfloor /2$, $n\in\mathbb{N}$.
		
Furthermore, $G_{-r,s}$ = $K_1*(G_{-r}*G_s)$ = $G^+(S)$ where $S$ = $[-r, s]$. 
		
		Here, the term $K_1$ is realized by the vertex with label 0 that is adjacent to all other vertices, the term $G_{-r}$ is $G^+([-r, -1])$ which is isomorphic to $G_r$, and the term
		is $G_{s}$ is $G^+([1, s])$. Thus, edges of $G_{-r,s}$
		are: (i) the edges of the term $G_{-r}$, (ii) the edges of the term $G_{s}$, (iii) the $r+s$ edges obtained by joining the $K_1$ term with each vertex of the terms $G_{-r}$ and $G_{s}$, and (iv) the $rs$ edges obtained by joining the vertices of $G_{-r}$ with the vertices of $G_{s}$.
		
		$\therefore$ $|E(G_{-r,s})|$ =	$|E(G_{-r})|$ + $|E(G_{s})| + r+s+rs$ 
		 
		\hspace{2cm} = $\frac{1}{2}(\binom{r}{2}-\left\lfloor\frac{r}{2}\right\rfloor)$ + $\frac{1}{2}(\binom{s}{2}-\left\lfloor\frac{s}{2}\right\rfloor)$ + $rs + r + s$
		 			
				\hspace{2cm} = $\frac{1}{4}(r^2+4rs+s^2+3(r+s))$ - $\frac{1}{2}(\left\lfloor\frac{r}{2}\right\rfloor + \left\lfloor\frac{s}{2}\right\rfloor)$

			\hspace{2cm} =	$rs + r + s + \left\lfloor\frac{(r-1)^2}{4}\right\rfloor + \left\lfloor\frac{(s-1)^2}{4}\right\rfloor$.
	\end{proof}
	
\medskip	
\begin{center}
	\begin{tikzpicture}[scale =0.85]
		
		\node (a1) at (9,3.25)  [circle,draw,scale=0.6]{1};
		\node (a2) at (10.5,2.5)  [circle,draw,scale=0.6]{2};
		\node (a3) at (11,1)  [circle,draw,scale=0.6] {3};
		\node (a4) at (10,0)  [circle,draw,scale=0.6] {4};
		\node (a5) at (8,0)  [circle,draw,scale=0.6] {5};
		\node (a6) at (7,1)  [circle,draw,scale=0.6]{6};
		\node (a7) at (7.5,2.5)  [circle,draw,scale=0.6] {7};
		
		\draw (a1) -- (a6);
		\draw (a1) -- (a5);
		\draw (a1) -- (a4);
		\draw (a1) -- (a3);
		\draw (a1) -- (a2);
		
		\draw (a2) -- (a5);
		\draw (a2) -- (a4);
		\draw (a2) -- (a3);
		
		\draw (a3) -- (a4);
		
		\node (b1) at (15,3.25)  [circle,draw,scale=0.6]{1};
		\node (b2) at (16.5,2.5)  [circle,draw,scale=0.6]{2};
		\node (b3) at (17,1)  [circle,draw,scale=0.6] {3};
		\node (b4) at (16,0)  [circle,draw,scale=0.6] {4};
		\node (b5) at (14,0)  [circle,draw,scale=0.6] {5};
		\node (b6) at (13,1)  [circle,draw,scale=0.6]{6};
		\node (b7) at (13.5,2.5)  [circle,draw,scale=0.6] {7};
		
		{\color{red}
			\draw (b7) -- (b1);
			\draw (b7) -- (b2);
			\draw (b7) -- (b3);
			\draw (b7) -- (b4);
			\draw (b7) -- (b5);
			\draw (b7) -- (b6);
			
			\draw (b6) -- (b2);
			\draw (b6) -- (b3);
			\draw (b6) -- (b4);
			\draw (b6) -- (b5);
			
			\draw (b5) -- (b3);
			\draw (b5) -- (b4);}
		
		
		\node (c1) at (21,3.25)  [circle,draw,scale=0.6]{1};
		\node (c2) at (22.5,2.5)  [circle,draw,scale=0.6]{2};
		\node (c3) at (23,1)  [circle,draw,scale=0.6] {3};
		\node (c4) at (22,0)  [circle,draw,scale=0.6] {4};
		\node (c5) at (20,0)  [circle,draw,scale=0.6] {5};
		\node (c6) at (19,1)  [circle,draw,scale=0.6]{6};
		\node (c7) at (19.5,2.5)  [circle,draw,scale=0.6] {7};
		
		{\color{red}
			\draw (c7) -- (c1);
			\draw (c7) -- (c2);
			\draw (c7) -- (c3);
			\draw (c7) -- (c4);
			\draw (c7) -- (c5);
			\draw (c7) -- (c6);}
		
		\draw (c6) -- (c1);
		{\color{red}
			\draw (c6) -- (c2);
			\draw (c6) -- (c3);
			\draw (c6) -- (c4);
			\draw (c6) -- (c5);}
		
		\draw (c5) -- (c1);
		\draw (c5) -- (c2);
		{\color{red}
			\draw (c5) -- (c3);	
			\draw (c5) -- (c4);}
		
		\draw (c4) -- (c3);
		\draw (c4) -- (c2);
		\draw (c4) -- (c1);
		
		\draw (c3) -- (c2);
		\draw (c3) -- (c1);
		
		\draw (c2) -- (c1);

	\end{tikzpicture}
	
	\medskip	
	\hspace{1cm}	Fig. 9. $G_{7}$ \hspace{2.5cm} Fig. 10. $G^{c}_{7}$
	\hspace{2.5cm} Fig. 11. $G_{7} \cup G^{c}_{7} = K_7$
\end{center}

\medskip	
\begin{center}
	\begin{tikzpicture}[scale =1.1]
		\node (a8) at (8.5,2)  [circle,draw,scale=0.6] {8};
		\node (a7) at (7.5,1.25)  [circle,draw,scale=0.6] {7};
		\node (a1) at (9.75,2)  [circle,draw,scale=0.6]{1};
		\node (a2) at (10.75,1.25)  [circle,draw,scale=0.6]{2};
		\node (a3) at (10.75,.25)  [circle,draw,scale=0.6] {3};
		\node (a4) at (9.75,-.5)  [circle,draw,scale=0.6] {4};
		\node (a5) at (8.5,-.5)  [circle,draw,scale=0.6] {5};
		\node (a6) at (7.5,.25)  [circle,draw,scale=0.6]{6};
		
		\draw (a1) -- (a7);
		\draw (a1) -- (a6);
		\draw (a1) -- (a5);
		\draw (a1) -- (a4);
		\draw (a1) -- (a3);
		\draw (a1) -- (a2);
		
		\draw (a2) -- (a6);
		\draw (a2) -- (a5);
		\draw (a2) -- (a4);
		\draw (a2) -- (a3);
		
		\draw (a3) -- (a5);
		\draw (a3) -- (a4);
		
		\node (b8) at (12.75,2)  [circle,draw,scale=0.6] {8};
		\node (b7) at (11.75,1.25)  [circle,draw,scale=0.6] {7};
		\node (b1) at (14,2)  [circle,draw,scale=0.6]{1};
		\node (b2) at (15,1.25)  [circle,draw,scale=0.6]{2};
		\node (b3) at (15,.25)  [circle,draw,scale=0.6] {3};
		\node (b4) at (14,-.5)  [circle,draw,scale=0.6] {4};
		\node (b5) at (12.75,-.5)  [circle,draw,scale=0.6] {5};
		\node (b6) at (11.75,.25)  [circle,draw,scale=0.6]{6};
		
		{\color{red}
			\draw (b8) -- (b1);
			\draw (b8) -- (b2);
			\draw (b8) -- (b3);
			\draw (b8) -- (b4);
			\draw (b8) -- (b5);
			\draw (b8) -- (b6);
			\draw (b8) -- (b7);
			
			\draw (b7) -- (b2);
			\draw (b7) -- (b3);
			\draw (b7) -- (b4);
			\draw (b7) -- (b5);
			\draw (b7) -- (b6);
			
			\draw (b6) -- (b3);
			\draw (b6) -- (b4);
			\draw (b6) -- (b5);
			
			\draw (b5) -- (b4);}
		
		\node (c8) at (17,2)  [circle,draw,scale=0.6] {8};
		\node (c7) at (16,1.25)  [circle,draw,scale=0.6] {7};
		\node (c1) at (18.25,2)  [circle,draw,scale=0.6]{1};
		\node (c2) at (19.25,1.25)  [circle,draw,scale=0.6]{2};
		\node (c3) at (19.25,.25)  [circle,draw,scale=0.6] {3};
		\node (c4) at (18.25,-.5)  [circle,draw,scale=0.6] {4};
		\node (c5) at (17,-.5)  [circle,draw,scale=0.6] {5};
		\node (c6) at (16,.25)  [circle,draw,scale=0.6]{6};
		
		{\color{red}
			\draw (c8) -- (c1);
			\draw (c8) -- (c2);
			\draw (c8) -- (c3);
			\draw (c8) -- (c4);
			\draw (c8) -- (c5);
			\draw (c8) -- (c6);
			\draw (c8) -- (c7);
			
			\draw (c7) -- (c2);
			\draw (c7) -- (c3);
			\draw (c7) -- (c4);
			\draw (c7) -- (c5);
			\draw (c7) -- (c6);
			
			\draw (c6) -- (c3);
			\draw (c6) -- (c4);
			\draw (c6) -- (c5);
			
			\draw (c5) -- (c4);}
		
		\draw (c7) -- (c1);
		\draw (c6) -- (c1);
		\draw (c5) -- (c1);
		\draw (c4) -- (c1);
		\draw (c3) -- (c1);
		\draw (c2) -- (c1);
		
		\draw (c6) -- (c2);
		\draw (c5) -- (c2);
		\draw (c4) -- (c2);
		\draw (c3) -- (c2);
		
		\draw (c5) -- (c3);
		\draw (c4) -- (c3);
	\end{tikzpicture}
	
	\medskip
	\hspace{.5cm}	Fig. 12. $G_{8}$ \hspace{2.5cm} Fig. 13. $G^{c}_{8}$
	\hspace{2cm} Fig. 14. $G_{8} \cup G^{c}_{8} = K_8$
\end{center}

\subsection{Structural properties of Harary's graph $G_n$} 

Consider sum graphs $G^+([s+1, s+n])$ for $s \in \mathbb{N}_0$ and $n \in \mathbb{N}$. When $s$ = 0, $G^+([s+1, s+n])$ = $G^+([1, n])$ = $G^+_n$, $n \in \mathbb{N}$. Here, we present structural properties of sum graphs $G^+([s+1, s+n])$ for $s \in \mathbb{N}_0$ and $n \in \mathbb{N}$. It is convenient to say that a set of vertices is {\em subscript-labeled} if the label on each of the vertex is same as its subscript.

\begin{theorem} [{\cite[Theorem 2.1]{vb14}}] \label{a6}\label{1} {\rm Let $r \in\mathbb{N}$ and $S \subseteq [r, 2r]$.  $G^+(S)$ is a totally disconnected graph.} 
\end{theorem} 

\begin{proof}  Clearly, the sum of the labels of any two vertices in $G^+(S)$ is at least $r + (r+1)$, which is greater than any of the other labels. Hence all vertices of $G^+(S)$ are isolated vertices only. 
\end{proof}

\begin{theorem} [{\cite[Theorem 2.2]{vb14}}]  \label{2}  {\rm Let $s \in \mathbb{N}_0$, $n \geq 3$ and $S = [s+1, s+n]$, $n\in \mathbb{N}$. 

\rm{(a)} For $s \geq n-2$, $G^+(S)$ is totally disconnected.

\rm{(b)} For $s \leq n-3$, $G^+(S)  \cong G_{n-s} \cup K^c_s$.}

\end{theorem} 

\begin{proof}  Part (a) follows from Theorem \ref{a6}. 
	
For part (b), we note that when $(s+1) + (s+2) \leq s+n$, the graph $G^{+}(S)$ has edges. That is when $s \leq n-3$, $G^{+}(S)$ has edges. Let the vertices of $G^+(S)$ be $u_{s+1}, u_{s+2}, \dots , u_{s+n}$ and be subscript-labeled, and similarly for the vertices $v_1, v_2, \dots, v_{n-s}$ in $G_{n-s}$.  Further, for $k = 1, 2, \dots , s$, let $w_k$ be a vertex with label $n-s+k$, and let $W_s$ be the totally disconnected graph with vertex-set $\{w_1, w_2, \dots , w_s\}$.  
 
Define mapping $f: V(G^+(S))$ $\rightarrow$ $V(G^+_{n-s} \cup W_s)$ such that $f(u_{s+i})$ = $v_i$ for $i$ = $1, 2, \dots, n-s$ and $f(u_{n+j})$ = $w_j$ for $j$ = $1, 2, \dots, s$.  This mapping is clearly a one-to-one correspondence between the sets of vertices of the two graphs.  We now show that it is an isomorphism by looking at the neighborhoods of the vertices $u_i$ in $G^+(S)$.  The neighbors of $u_{s+1}$ are $u_{s+2}, u_{s+3}, \dots, u_{n-1}$, the neighbors of $u_{s+2}$ are $u_{s+1}, u_{s+3}, \dots, u_{n-2}$, and so on.  In addition, for $n-s+1 \leq \l \leq s+n$,  $u_{\l}$ is of degree 0 .  It follows that for $i,j \leq n-s$, $u_{s+i}$ is adjacent to $u_{s+j}$ if and only if $v_i$ is adjacent to $v_j$ in $G^+_{n-s}$.  Furthermore, if $\l > n-s$, then $f(u_{\l})$ is in $W_s$, so isolated vertices are mapped to isolated vertices.  Consequently, $f$ is an isomorphism. 
\end{proof}

For convenience, if graph $F$ is a subgraph of graph $G$ without the vertex labels, this will often be denoted by $F \subseteq _\mathrm{wvl} G$. Similarly, if $F$ is isomorphic to $G$ without the vertex labels, we often write  $F \cong _\mathrm{wvl} G$ and if $F$ is an induced subgraph of $G$, this will sometimes be denoted by $F \leq G$ and by $F \leq_\mathrm{wvl} G$. Following theorem is a corollary to the previous theorem.

\begin{cor} [{\cite[Corollary 2.3]{vb14}}] \label{3}  {\rm If $0 \leq s \leq n-3$ and $S = [s+1, s+n]$, then $G_{n-s} \subseteq _\mathrm{wvl}  G^+(S)$.  In particular, $G_{n-2}  \subseteq _\mathrm{wvl}  G^+([2, n])  \subseteq _\mathrm{wvl} G_{n-1}$.} 
\end{cor}

\begin{proof}\quad It follows from Theorem \ref{2} that, for $n \geq 3$, $G^+([2,n])$ $\cong_{wvl}$ $G_{n-2} \cup W_1$ where $W_1$ consists of a single isolated vertex $w$. This establishes the first inclusion. The second follows from the fact that $G_{n-1}$ is the only maximal sum graph of order $n- 1$. 
\end{proof}

\subsection{Structural properties of $G_n$ and $G^c_n$ using anti-sum labeling}

In \cite{vb14}, Vilfred defined {\em anti-sum labeling} and using it, obtained more structural properties of $G_n$ and $G^c_n$, and are presented in this subsection. The {\it underlying graph} of an integral sum graph is obtained by removing all vertex labels. In \cite{vb14}, we call the vertices $i$ and $n+1-i$ of $G^+_n$ as {\em supplementary}, $1 \leq i \leq n$.  
	
\begin{dfn}	 \cite{vb14} A graph $G$ is an {\em anti-integral sum graph} or {\em anti-$\mathbb{Z}$-sum graph}  if the vertices of $G$ can be labeled with distinct integers so that $e = uv$ is an edge of $G$ if and only if the sum of the labels on vertices $u$ and $v$ is not a vertex label in $G$. 
\end{dfn}

Clearly, $f$ is an integral sum labeling of graph $G$ if and only if $f$ is an anti-integral sum labeling of $G^c$. Indeed, many results on anti-sum graphs are simply analogues to the corresponding results on sum graphs and are stated without proof. A simple property of complements is that if $v$ is a vertex of graph $G$, then ${(G-v)}^c \cong G^c - v$. Therefore, we add results on anti-sum graphs as consequences of results on sum graphs. Sum graph $G_7$, anti-sum graph $G^{c}_{7}$ and $G_{7} \cup G^{c}_{7} = K_7$ are given in figures 9 to 11. In figures 12 to 14 graphs $G_8$, $G^c_8$ and $G_{8} \cup G^{c}_{8} = K_8$ are given. 
	
\begin{dfn} \label{d2.2.2}	 \cite{d05} \quad A graph $G$ is a \emph{split graph} if its vertices can be partitioned into a clique and a stable set. $A$ \emph{ clique} in a graph is a set of pair-wise adjacent vertices and an \emph{ independent set} or \emph{ stable set} in a graph is a set of pair-wise non-adjacent vertices. 
\end{dfn}

Clearly, $G_n$ and $G^c_n$ are split graphs; $[1, m], [1, m+1], [m+1, 2m], [m+2, 2m+1]$ are cliques and $[m+1, 2m], [m+2, 2m+1], [1, m], [1, m+1]$ are stable sets in $G_{2m}, G_{2m+1}, G^c_{2m}, G^c_{2m+1}$, respectively.

\begin{theorem} [{\cite[Theorem 2.4]{s21}}] \label{4.8} \quad {\rm  For $n\in\mathbb{N}$, $G^+_n - \{\left\lfloor n/2 \right\rfloor\}$ $\cong_{wvl}$ $G^+_{n-1}$. }
\end{theorem}

\begin{proof}\quad Proof is based on the principle of mathematical induction on $n$, the order of graph $G^+_n$, $n\in\mathbb{N}$. For $n$ = 1,2, the result is true. Assume the result for all $k < n$, $k\in\mathbb{N}$. That is the result is true for all $G^+_k$ such that $k < n$, $k\in\mathbb{N}$. Now, we consider graph $G^+_n$ with odd and even cases of $n$ separately and present proof only for the even case since proof for the two cases are similar. Let $n$ = $2m$ and the graph be $G^+_{2m}$, $m\in\mathbb{N}$. Let $V(G^+_{2m})$ = $\{v_1, v_2, . . . , v_{2m}\}$ and $V(G^+_{2m-1})$ = $\{w_1, w_2, . . . , w_{2m-1}\}$ and the vertices of the two graphs be subscript-labeled. Our aim is to prove that $G^+_{2m}-\{v_m\}$ $\cong_{wvl}$ $G^+_{2m-1}$, $m\in\mathbb{N}$.  Clearly, set $K$ = $\{v_1, v_2, . . . , v_m\}$ is a clique and set $L$ = $\{v_{m+1}, v_{m+2}, . . . , v_{2m}\}$ is an independent set in $G^+_{2m}$. Also, it is clear that the set $K_1$ = $\{v_1, v_2, . . . , v_{m-1}, v_{m+1}\}$ is a clique and set $L_1$ = $\{v_{m+2}, v_{m+3}, . . . , v_{2m}\}$ is an independent set in $G^+_{2m}-\{v_m\}$ (see figure 12 in the case of $G_8^+$ for more clarity), and set $K_2$ = $\{w_1, w_2, . . . , w_{m}\}$ is a clique and set $L_2$ = $\{w_{m+1}, w_{m+2}, . . . , w_{2m-1}\}$ is an independent set in $G^+_{2m-1}$.
	
	Define mapping $f: V(G^+_{2m}-\{v_m\}) \to V(G^+_{2m-1})$ by $f(v_i)$ = $w_i$, $1 \leq i \leq m-1$ and $f(v_j)$ = $w_{j-1}$, $m+1 \leq j \leq 2m$. Clearly, $f$ preserves adjacency on the clique sets $K_1$ and $K_2$ as well as on the independant sets $L_1$ and $L_2$. We now prove that $f$ also preserves adjacency on every $(v_i, v_j)$ where $v_i\in K_1$ and $v_j\in L_1$. For $v_i\in K_1$ and $v_j\in L_1$, $1 \leq i \leq m-1$ or $i$ = $m+1$ and $m+2 \leq j \leq 2m$. Therefore, for $1 \leq i \leq m-1$ and $m+2 \leq j \leq 2m$, $f((v_i, v_j))$ = $(f(v_i), f(v_j))$ = $(w_i, w_{j-1})$ and $f((v_{m+1}, v_j))$ = $(f(v_{m+1}), f(v_j))$ = $(w_m, w_{j-1})$. This implies, for $1 \leq i \leq m-1$ or $i$ = $m+1$ and $m+2 \leq j \leq 2m$, $f((v_i, v_j))$ = $(w_i, w_{j-1})$ where $(u_i, v_j)\in K_1 \times L_1$ and $(w_i, w_j)\in K_2 \times L_2$. Thus, $f$ preserves adjacency. This implies that $G^+_{2m}-\{v_{m}\})$ $\cong_{wvl}$ $G^+_{2m-1}$, $m\in\mathbb{N}$.
	
	Similarly, we can prove that $G^+_{2m+1}-\{v_{\left\lfloor (2m+1)/2 \right\rfloor}\})$ $\cong_{wvl}$ $G^+_{2m}$, $m\in\mathbb{N}$.
\end{proof}

\begin{theorem} [{\cite[Theorem 2.5]{s21}}]   \label{4.9} \quad {\rm  For $n\in\mathbb{N}$, $G^+_{0,n}-\{\left\lfloor n/2 \right\rfloor\}$ $\cong_{wvl}$ $G^c_{n}$. }
\end{theorem}

\begin{proof}\quad For convenience, we consider odd and even cases of $n$, order of the graph ${(G^+_{n})}^c$, separetely and since the two cases are so similar, here we prove only for the even case. Let $n$ = $2m$ and the graph be $G^+_{0,2m}$, $m\in\mathbb{N}$. Let $V(G^+_{0, 2m})$ = $\{v_0, v_1, v_2, \dots, v_{2m}\}$ and $V({(G^+_{2m})}^c)$ = $\{w_1, w_2, \dots, w_{2m}\}$, both subscript-labeled. Our aim is to prove that $G^+_{0,2m}-\{m\}$ $\cong_{wvl}$ ${(G^+_{2m})}^c$, $m\in\mathbb{N}$. Clearly, in $G^+_{0,2m}-\{m\}$, the set $K_1$ = $\{v_0, v_1, \dots, v_{m-1}, v_{m+1}\}$ is a clique and $L_1$ = $\{v_{m+2}, v_{m+3}, \dots, v_{2m}\}$ is an independent set and in ${(G^+_{2m})}^c$, the set $K_2$ = $\{w_m, w_{m+1}, \dots, w_{2m}\}$ is a clique and $L_2$ = $\{w_1, w_2, \dots, w_{m-1}\}$ is an independent set.
	
	Define the mapping $f: V(G^+_{0, 2m}-\{v_m\}) \to V({(G^+_{2m})}^c)$ by $f(v_i)$ = $w_{2m-i}$ for $i$ = 0 to $m-1$ and $f(v_j)$ = $w_{2m-j+1}$ for $j$ = $m+1$ to $2m$. Then, the proof is similar to Theorem \ref{4.8}.
	
	Similarly, we can prove that $G^+_{0, 2m+1}-\{v_m\}$ $\cong_{wvl}$ ${(G^+_{2m+1})}^c$, $m\in\mathbb{N}$.
\end{proof}

\begin{theorem}[{\cite[Theorems 2.4 and 3.1]{vb14}}] \label{2.11}  {\rm Let $n \geq 3$ and $n\in\mathbb{N}$. 
	\begin{enumerate}
		\item[(i)] If any pair of supplementary vertices are removed from $G_n$, then the result is isomorphic to $G_{n-2}$ without the vertex labels.  
			
		\item[(ii)] If any pair of supplementary vertices are removed from $G^c_n$, then the result is isomorphic to $G^c_{n-2}$ without the vertex labels. 
	\end{enumerate}} 
\end{theorem} 

\begin{proof}  (i)~ For convenience, the odd and even cases are considered separately, and because the two cases are so similar, we prove only the even case here.  Let $n = 2s$ and let $V(G_{2s} ) = \{v_1, v_2, \dots, v_{2s}\}$ and $V(G_{2s-2}) = \{w_1, w_2, \dots, w_{2s-2}\}$, both subscript-labeled. In $G_{2s}$, the set $K=\{v_1, v_2, \dots, v_{s}\}$ induces a clique and the set $L=\{v_{s+1}, v_{s+2}, \dots, v_{2s}\}$ an independent set, and similarly for the vertices in $G_{2s-2}$.
	
	Let $v_j$ and $v_{2s+1-j}$  be a pair of supplementary vertices in $G_{2s}$, and define the mapping $f: V(G_{2s} - \{v_j, v_{2s+1-j}\}) \rightarrow V(G_{2s-2})$ by $f(v_{i}) = w_i$ for $i = 1, 2, \dots, j-1$; $f(v_{i}) = w_{i-1}$ for $i = j+1, j+2, \dots, 2s-j$; and $f(v_{i}) = w_{i-2}$ for $i = 2s-j+2$, $2s-j+3$, $\dots$, $2s$.  This mapping is clearly bijective between the sets of vertices of the two graphs. Furthermore, the images of the clique $K$ and the independent set $L$ are also a clique and independent set.  The preservation of the edges and non-edges between the two sets $K$ and $L$ can be established in a straightforward way using the definition of a sum graph.  From this, it follows that $G_{2s} - \{v_j, v_{2s+1-j}\} \cong_\mathrm{wvl} G_{2s-2}$. 
\end{proof}

This theorem can be extended in the following way:

\begin{theorem}  [{\cite[Corollary 2.5 and Theorem 3.1]{vb14}}] \normalfont\label{4}
	{\rm Let $n > 2k \geq 2$ and $k,n\in\mathbb{N}$. 
		\begin{enumerate}
			\item[(i)] If $k$ pairs of supplementary vertices are removed from  $G_n$, then the result is isomorphic to $G_{n-2k}$ without the vertex labels.  
			\item[(ii)] If $k$ pairs of supplementary vertices are removed from $G^c_n$, then the result is isomorphic to $G^c_{n-2k}$ without the vertex labels. \hfill $\Box$
	\end{enumerate}} 
\end{theorem}

In \cite{vb14}, it is observed that for all $r$ and $s$ with $0 < r < s$, $G^+_r \cup G^{+}([r+1, s])$ is a (spanning) subgraph of $G^+_s$. The following is a related result. Note that Theorems \ref{14} to \ref{16} are stated in \cite{vb14} without proof and here we present proof to these.  
	
\begin{theorem} [{\cite[Theorems 2.6 and 3.2]{vb14}}] \normalfont \label{14} 
		{\rm Let $n,s\in\mathbb{N}$ and $n\geq s+2$.
		\begin{enumerate}
			\item[(a)] $G^+_s \cup G^+([s+1,n]) \leq G^+_n \leq G^+_1 * G^+([2,n])$.
				
			\vspace{.1cm} 
			\item[(b)] $G^+_1\cup (G^+([2,n]))^c \leq {(G^+_n)}^c \leq {(G^+_s)}^c * {(G^+ ([s+1,n]))}^c$. 
		\end{enumerate} }
\end{theorem}
	
	\begin{proof}\quad 
		\begin{enumerate}
			\item[(a)]  Clearly, $G^+_s \cup G^+([s+1, n])$ is a spanning subgraph of $G^+_n$, $n\geq s+2$ and $n,s\in\mathbb{N}$. Now, let us prove the other part. Let $V(G^+_n)$ = $\{u_1,u_2,\dots,u_n\}$ and the vertices be subscript-labeled. Using the definition of sum labeling of graphs, $G^+_n$ = $G^+( [2, n])$ $\cup$ $K_1(u_1)$ $\cup$ $(u_1,  u_2)$ $\cup$ $(u_1, u_3)$ $\cup$ $\dots$ $\cup$ $(u_1, u_{n-1})$. And also the graph $G^+_1 * G^+( [2, n])$ = $G^+( [2, n])$ $\cup$ $K_1(u_1)$ $\cup$ $(u_1, u_2)$ $\cup$ $(u_1, u_3)$ $\cup$ $\dots$ $\cup$ $(u_1, u_{n-1})$ $\cup$ $(u_1, u_n)$. 
			
			\item[(b)]  Taking complement of the graphs in relation (a), we get the result (b).
		\end{enumerate} 
	\end{proof}  
	
	\begin{theorem} [{\cite[Theorems 2.7]{vb14}}] \label{15} {\rm  Let $n$ and $r$ be positive integers.	
			\begin{itemize} 
				\item  [\rm (a)]   Labeled graphs $G^+_n$ and $G^+_s * G^+( [s+1,  n])$ are non-comparable for $n \geq 4$ and $2 \leq s < n$.
				
				\item  [\rm (b)]	The underlying graphs of $G^+_n$ and $G^+_2 * G^+([3,n])$ are non-comparable for $n \geq 5$.  
				
				\noindent
				\vspace{.1cm} 
				For $n \geq 4$,  $||G_2^+ * G^+( [3, n])||$ = $||G^+([3, n] )||$ $+$ $2(n-2)$ = $||G^+_n||+2$. 
			\end{itemize} 
		}
	\end{theorem} 
	
	\begin{proof}\quad 
		\begin{itemize} 
			\item  [\rm (a)]   For $n \geq 4$ and $n > s \geq 2$, vertices with label $1$ and $s$ are adjacent in $G^+_n$ but non-adjacent in $G^+_s$ and thereby they are non-adjacent in $G^+_s * G^+([s+1, n])$. This implies, labeled graph $G^+_n$ cannot be a subgraph of $G^+_s * G^+( [s+1, n])$. Also, vertices with label $s$ and $n$ are adjacent in $G^+_s * G^+( [s+1, n])$ but  non-adjacent in the labeled graph $G^+_n$ and thereby $G^+_s * G^+( [s+1, n])$ cannot be a subgraph of $G^+_n$. Hence the two graphs are non-comparable.
			
			\item  [\rm (b)]	For $n \geq 4$, the sum graph $G_n^+ \setminus \{ 1, 2 \}$ = $G^+( [3, n]);$ $E(G_n^+)$ = $E(G^+( [3, n]))$ $\cup$ $\{(1,n-1)$, $(1,n-2)$, $\cdots$, $(1,3)$, \textbf{(1,2),} $(2,n-2)$, $(2,n-3)$, $\cdots$, $(2,3)\}$. And $E(G_2^+ * G^+([3, n]))$ = $E(G^+( [3, n]))$ $\cup$ $\{ \textbf{(1,n)}$, $(1,n-1)$, $\cdots$, $(1,3)$, $\textbf{(2,n), (2,n-1)}$, $(2,n-2)$, $\cdots$, $(2,3) \}$. Thus, for $n \geq 4$, $||G_n^+||$ + 2 =  $||G_2 ^+ * G^+([3, n] )||$ = $||G^+( [3, n] )||$ + $2(n-2)$. And thereby, the underlying graph of $G_n^+$ cannot be a super graph of the underlying graph of $G_2^+ * G^+( [3, n]).$
			
			\medskip
			\noindent
			\textit{Claim.} For $n\geq 5$, the underlying graph of $G_2 ^+ * G^+([3, n] )$ cannot be a super graph of the underlying graph of $G_n^+$. 
			
			For $n\geq 3$, in $G_2 ^+ * G^+([3, n] )$, vertices 1 and 2 are nonadjacent and each of them makes a cycle of length 3 with each edge of $G^+([3, n])$. When $n$ = $5$ and $n$ = $6$, graphs $G_2 ^+ * G^+([3, 5] )$ and $G_2 ^+ * G^+([3, 6] )$ have no cycle of length 3 since $G_2^+$, $G^+([3, 5])$ and $G^+([3,6])$ are totally disconnected graphs, by Theorem \ref{2} whereas graph $G_5^+$ has one cycle of length 3 and graph $G_6^+$ has 2 cycles, each of length 3. Thus, the underlying graph of $G_2^+ * G^+([3, n])$ cannot be a super graph of the underlying graph of $G_n^+$ for $n$ = 5 and 6. 
			
			For $n\geq 7$, number of cycles, each of length 3 and has 1 or 2 as one of its vertex label in $G_2 ^+ * G^+([3, n] )$ is $2\times ||G^+( [3, n])||$. On the other hand, for $n \geq 7$, in $G_n^+$, vertices 1 and 2 are adjacent, each of them makes a cycle of length 3 with each edge of subgraph $G^+([3, n])$ and edge $(1,2)$ makes cycles of length 3 with each vertex $3, 4, \dots, n-2$. Thus, for $n \geq 7$, $G_n^+$ has $n-2-2$ = $n-4$ (both 1 and 2 cannot be adjacent to $n$ and $n-1$) more number of cycles of length 3 than that of $G_2 ^+ * G^+([3, n] )$ and thereby the underlying graph of  $G_2 ^+ * G^+([3, n] )$ cannot be a super graph of the underlying graph of  $G_n^+$. Hence, the claim is true and thereby the theorem is proved.
		\end{itemize} 
	\end{proof}
	
	\begin{theorem} [{\cite[Theorem 3.3]{vb14}}] \label{16} {\rm   Let $n$ and $s$ be positive integers.
			
			\begin{itemize}
				\item  [\rm (a)] For $2 \leq s < n$ and $n \geq 4$, labeled graphs ${(G^+_n)}^c$ and ${(G^+_s)}^c \cup {(G^+([ s+1, n]))}^c$ are non-comparable.
				
				\item  [\rm (b)] For $n \geq 5$, the underlying graphs of ${(G^+_n)}^c$ and $G^{\rm{c}}_2 \cup (G^+([3,n]))^{\rm{c}}$ are non-comparable.  
		\end{itemize} }
	\end{theorem} 
	
	\begin{proof}\quad Result follows from Theorem \ref{15} by considering complement in each case.
	\end{proof}
	
	The next two theorems give additional connection between Harary's sum graphs, their complements and anti-sum graphs.
	
	\begin{theorem} [{\cite[Theorems 3.4]{vb14}}] \label{2.13} \quad {\rm \small{For $n \geq 2$, 
				$G^+_n \cong_\mathrm{wvl} G^{\rm{c}}_n - \{(1, n), (2, n-1),  \dots,  (\lfloor n/2 \rfloor$, $\lceil n/2 \rceil + 1)\}$.}}  
	\end{theorem} 
	
	\begin{proof}  Let $V(G^+_n) = \{v_1, v_2, \dots, v_n\}$ and $V(G^c_n) = \{w_1, w_2, \dots, w_n\}$, with both graphs subscript-labeled, and let $H_n$ = $G^c_n - \{(1, n), (2, n-1), ...,  (\lfloor n/2 \rfloor,  \lceil n/2 \rceil + 1)\}$.  Define a mapping $f$ from $G^+_n$ to $H_n$ with $f(v_i)$ = $w_{n-i+1}$ for $i = 1, 2, \dots, n$.  Now, for $1 \leq i < j \leq n$, $v_i$ adjacent to $v_j$ in $G^+_n$ implies, $i + j \leq n$ which in turn implies, $(n-i+1) + (n-j+1) > n$, so that $w_{n-i+1}$ is adjacent to $w_{n-j+1}$ in $H_n$.  Hence if two vertices are adjacent in $G^+_n$, their images under $f$ are adjacent in $H_n$.
		
		Now, let us prove the converse part, namely, that if two vertices are not adjacent in $G^+_n$, their images under $f$ are not adjacent in $H_n$.  To this end, let $v_i$ and $v_j$ be non-adjacent in $G^+_n$.  Then $i+j \geq n+1$, so $(n-i+1) + (n-j+1) \leq n+1$.  If $(n-i+1) + (n-j+1) \leq n$, then $w_{n-i+1}$ is not adjacent to $w_{n-j+1}$ in $G^c_n$ and hence not in $H_n$.  That leaves only the case where $(n-i+1) + (n-j+1) = n+1$.  But this implies that $i+j = n+1$, which violates one of the hypotheses and completes the proof.
	\end{proof}
	
	\begin{theorem} [{\cite[Theorem 3.1]{vs14}}]  \label{2.14}
		{\rm  For $n \geq 3$, the underlying graphs of $K_n$, $2 G_n$ $\cup$ $\{(1,n)$, $(2,n-1)$, $\dots$, $(\left\lfloor n/2 \right\rfloor$, $n+1-\left\lfloor n/2 \right\rfloor)\}$ and $G^c_n \cup (G^c_n-\{(1,n), (2,n-1),\dots,(\left\lfloor n/2 \right\rfloor, n+1- \left\lfloor n/2 \right\rfloor)\})$ are isomorphic.  }
	\end{theorem}
	
	\begin{proof} Using Theorem \ref{2.13}, the underlying graphs of $G_n$ and $G^c_n -\{(1,n)$, $(2,n-1)$, $\dots$, $(\left\lfloor n/2 \right\rfloor$, $n+$1 - $\left\lfloor n/2 \right\rfloor)\}$ are isomorphic. This implies that the underlying graphs of $K_n$ = $G_n \cup G^c_n$, $2 G_n$ $\cup$ $\{(1,n)$, $(2,n-1)$, $\dots$, $(\left\lfloor n/2 \right\rfloor, n+1-\left\lfloor n/2 \right\rfloor)\}$, and $G^c_n \cup (G^c_n - \{(1,n)$, $(2,n-1) ,\dots ,(\left\lfloor n/2 \right\rfloor, n+1-\left\lfloor n/2 \right\rfloor)\})$ are isomorphic. 
	\end{proof} 
	
\begin{lemma} [{\cite[Lemma 2.3]{vr14}}] \label{4.7} \quad {\rm  Let $V(G_n)$ = $\{v_1,v_2,\dots,v_n\}$ = $V(G_{n}^c)$ where $v_j$ be the vertex with integral sum labeling $j$ in $G_n$ and anti-integral sum labeling $j$ in $G_{n}^c$, $1 \leq j \leq n$ and $n\in\mathbb{N}$. Then, 
	\begin{enumerate}
		\item [\rm (i)] $G_{0,n}$ $\cong$ $G_{n+2}\setminus \{v_{n+2}\}$;

		\item [\rm (ii)] $G_{n+2}$ $\cong$ $(G_{n} * \{v_{n+1}\})$ $\cup$ $\{v_{n+2}\}$;

		\item [\rm (iii)] $G_{n+2}^c$ $\cong$ $(G_{n}^c$ $\cup$ $\{v_{n+1}\})$ $*$ $\{v_{n+2}\}$ and

		\item [\rm (iv)] $G_{-1,n}$ $\cong$ $G_{n+4}\setminus \{v_{n+3} ,v_{n+4}\}$, without the vertex labels.  
	\end{enumerate}		  }
\end{lemma}
	
\begin{proof}\quad Let $V(G_{0,n})$ = $\{ u_0,u_1,u_2,\dots, u_{n}\}$, $V(G_{-1,n})$ = $\{ u_0,u_1,u_2,\dots, u_{n}, u_{-1}\}$, and $V(G_{j})$ = $\{ v_1,v_2,\dots, v_{j} \}$ and the graphs be subscript-labeled, $j,n\in\mathbb{N}$.
		
	\begin{enumerate}
		\item [\rm (i)]  We have $G_{-m,n}$ = $K_1(0) * ((-G_m) * G_n)$ and $G_{0,n} = K_1(0) * G_n$, $m,n \in\mathbb{N}$. Define $f:$ $V(G_{0,n})$ $\rightarrow$ $V(G_{n+2} \setminus \{ v_{n+2} \})$ such that $f(u_i)$ = $v_{i+1}$ and $f((u,v))$ = $(f(u),f(v))$ for every $(u,v)\in E(G_{0,n})$, $i = 0,1,\dots,n$. Now,  for $x \neq y$, $(u_x,u_y)\in E(G_{0,n})$ if and only if $0 < x+y < n+1$ if and only if $2 < (x+1)+(y+1) < n+3$ if and only if $3 \leq (x+1)+(y+1) \leq n+2$ if and only if $(v_{x+1},v_{y+1})$ = $(f(u_x),f(u_y))\in E(G_{n+2})$ = $E(G_{n+2} \setminus \{ v_{n+2} \})$. This implies, $f$ is a bijective mapping and preserves adjacency. Hence, $G_{0,n}$ $\cong$ $G_{n+2} \setminus \{v_{n+2} \}$, without the vertex labels.  
		
		\item [\rm (ii)]  Using $(i)$, we obtain,  $G_{n+2}$ $\cong$ $G_{0,n}$ $\cup$ $\{u_{n+2}\}$ $\cong$ $(G_n * K_1)$ $\cup$ $\{v_{n+2}\}$ $\cong$  $(G_{n}$ $*$ $\{v_{n+1}\})$ $\cup$ $\{v_{n+2}\}$, without the vertex labels, $n\in\mathbb{N}$. 
		
		\item [\rm (iii)] Using $(ii)$, we get, $G_{n+2}^c$ $\cong$ ${((G_{n}*\{v_{n+1}\})\cup \{v_{n+2}\})}^c$ $\cong$   ${(G_{n}*\{v_{n+1}\})}^c  * \{v_{n+2}\}$ $\cong$ $(G_{n}^c \cup \{v_{n+1}\}) * \{v_{n+2}\}$, without the vertex labels, $n\in\mathbb{N}$. 
		
		\item [\rm (iv)] We have $G_{-1,n} \cong K_1(0) * (K_1(-1) * G_n) \cong K_1(-1) * (K_1(0) * G_{n}) \cong K_1(-1) * G_{0,n} \cong K_1(-1) * (G_{n+2} \setminus \{v_{n+2} \})$, without the vertex labels, using $(i)$, $n \in\mathbb{N}$. Using Theorem \ref{4}, graph $G_{n+4}\setminus \{v_{1},v_{2},v_{n+3},v_{n+4}\}$ is isomorphic to $G_n$, without the vertex labels. And so $((G_{n+4}\setminus \{v_{1},v_{2},v_{n+3},v_{n+4}\})$ $*$ $K_1)$ $*$ $K_1$ $\cong$ $G_{-1,n}$, without the vertex labels. Define $f:$ $V(G_{-1,n})$ $\rightarrow$ $V(G_{n+4} \setminus \{ v_{n+3},v_{n+4} \})$ such that $f(u_0)$ = $v_1$, $f(u_{-1})$ = $v_2$, $f(u_i)$ = $v_{i+2}$ for $i$ = $1,2,\dots,n$ and $f((u,v))$ = $(f(u),f(v))$ for every $(u,v)\in E(G_{-1,n})$. Now, consider images of edges incident at each point $u_0$ and $u_{-1}$, seperately. In $G_{-1,n}$, integral sum labeling of $u_0$ and $u_{-1}$ are 0 and -1, respectively, $u_0$ and $u_{-1}$ are adjacent and each one is adjacent to $u_j$ for $j$ = $1,2,\dots,n$. Now,  $f((K_1(0),u_i))$ = $f((u_0, u_i))$ = $(f(u_0), f(u_i))$ = $(v_1,v_{i+2})\in E(G_{n+4} \setminus \{v_{n+4},v_{n+3} \})$ for every $i$, $i$ = $1,2,\dots,n;$ $f((K_1(0),u_{-1}))$ = $f((u_0,u_{-1}))$ = $(f(u_0),f(u_{-1}))$ = $(v_1,v_2)\in E(G_{n+4} \setminus \{v_{n+3},v_{n+4} \} )$ and $f((u_{-1}$, $u_j))$ = $(f(u_{-1}),f(u_j))$ = $(v_2,v_{j+2})\in E(G_{n+4} \setminus \{v_{n+3},v_{n+4} \})$ for every $j$, $j$ = 1 to $n$. Therefore, $f$ is a bijective mapping preserving adjacency and hence, $G_{-1,n}$ $\cong$ $G_{n+4} \setminus \{v_{n+3},v_{n+4} \}$, without the vertex labels. 
	\end{enumerate}
\end{proof}

\begin{res} [{\cite[Result 2.4]{vr14}}]  \quad \label{8.10b} {\bf [Algorithm to generate $G_{n}$ and $G^c_{n}$]}\\ 
	{\rm	Starting with either $G_0$ or $G_1$ and using results (ii) and (iii) of Lemma \ref{4.7} for $n$ = $2,4,\dots$ or $n$ = $3,5,\dots$, one can generate sum graphs $G_n$ and anti-sum graphs $G^c_{n}$ of any order without using definitions of sum and anti-sum labeling. }
\end{res} 
	
\subsection{Structural properties of integral sum graphs $G_{0, s}$} 
	
	In this section, we present a few structural properties of  $G_{0,s}$ similar to that of $G_n$ as given in Theorems \ref{4.9}, \ref{2.11} and \ref{4}. 
	
	\begin{theorem} [{\cite[Theorem 3.2]{vs14}}] \label{19} {\rm For $n \in \mathbb {N}$, the underlying graphs of $G_{0,n}$ and $G^c_{n+1}$ are isomorphic.}
	\end{theorem}
	
\begin{proof}\quad Using the definition of anti-sum labeling, we obtain, $G^c_{n+1} = (((G^c_n \cup K_1(n+1))-\{(1,n)$, $(2,n-1)$, $\dots$, $(\left\lfloor n/2 \right\rfloor, n+1-\left\lfloor n/2 \right\rfloor)\})$ $\cup$ $\{(1,n+1),(2,n+1)$, $\dots$, $(n,n+1)\})$ $\cong$ $(G^c_n-\{(1,n)$, $(2,n-1)$, $\dots$, $(\left\lfloor n/2 \right\rfloor$, $n+1-\left\lfloor n/2 \right\rfloor)\})$ $\cup$ $K_1(n+1)$ $\cup$ $\{(1,n+1),(2,n+1),\dots,(n,n+1)\}$  which is isomorphic to the underlying graph of $G_n$ $\cup$ $K_1(n+1)$ $\cup$ $\{(1,n+1),(2,n+1),\dots,(n,n+1)\}$ = $G_n * K_1(n+1)$ where $K_1(n+1)$ represents a vertex with vertex label $n+1$ in the graph $G^c_n$ $\cup$ $K_1(n+1)$. This implies, the underlying graphs of $G^c_{n+1}$, $G_n*K_1(n+1)$ and $G_{0,n}$ are isomorphic since the underlying graphs of $G_n * K_1(n+1)$  and $G_{0,n}$ are isomorphic. Hence the result.
\end{proof}

\begin{theorem}[{\cite[Theorem 4.1]{vs14}}] \label{20} { \quad \rm 	For $n \geq 3$, the underlying graphs of $G_{0,n}-\{0,n\}$ and $G_{0,n-2}$ are isomorphic.}
\end{theorem}
\begin{proof}
	Let $V(G_{0,n-2}) = \{u_0,u_1,\dots,u_{n-2}\}$ and $V(G_{0,n}) = \{v_0,v_1,\dots,v_n\}$ where $i$ and $j$ are integral sum labels of $u_i$ and $v_j$, respectively, $0 \leq i \leq n-2$ and $0 \leq j \leq n$. Define mapping $f : V(G_{0,n-2})\rightarrow V(G_{0,n})$ such that $f(u_i) = v_{i+1}$ for $i = 0,1,\dots,n-2$. Now, $u_i$ and $u_j$ are adjacent in $G_{0,n-2}$ if and only if $i \neq j, 0\leq i,j\leq n-2$ and $0+1 = 1 \leq i+j \leq n-2$ if and only if $i+1 \neq j+1, 1\leq i+1,j+1 \leq  n-1$ and $1+2 = 3 \leq(i+1)+(j+1) \leq n = (n-1)+1$ if and only if $v_{i+1}$ and $v_{j+1}$ are adjacent in $G_{0,n}-\{0,n\}$ if and only if $f(u_i)$ and $f(u_j)$ are adjacent in $G_{0,n}-\{0,n\}$. This implies, the mapping $f$ is bijective, preserves adjacency and $f(G_{0,n}-\{0,n\}) = G_{0,n-2}$.
\end{proof}

\begin{theorem} [{\cite[Theorem 4.2]{vs14}}]  \label{2.24} { \quad \rm For $n\geq 5$, the underlying graphs of $G_{0,n}-(\{0,n,n-1$, $n-2\}$ $\cup$ $[n] \cup [n-1])$ and $G_{0,n-4}$ are isomorphic where $[k]$ in $G^+(S)$ denotes the set of all edges of $G^+(S)$ whose edge-sum value is $k, k \in S$.}
\end{theorem}
\begin{proof} Using the definition of integral sum labeling, we obtain isomorphic graphs of the underlying graphs of $G_{0,n}-(\{n,n-1\} \cup [n] \cup [n-1])$ and $G_{0,n-2}$ where $[k]$ in $G^+(S)$ denotes the set of all edges of $G^+(S)$ whose edge-sum value is $k, k \in S$ \cite{vm13c}. Using Theorem \ref{4}, the underlying graphs of $G_{0,n-2}-\{0,n-2\}$ and $G_{0,n-4}$ are isomorphic. Hence the result. 
\end{proof}

Generalizing the above theorem, we get the following result.

\medskip
\begin{theorem} [{\cite[Theorem 4.3]{vs14}}]  \label{2.25} { \quad \rm 
		For $n \geq 3$, the underlying graphs of $G_{0,n}-\{0,n\}$ and $G_{0,n-2}$ are isomorphic and for $n \geq 2r+3$ and $r \in \mathbb {N}$, the underlying graphs of $G_{0,n}- (\{0,n,n-1,n-2,\dots$, $n-2r+1$, $n-2r\}$ $\cup$ $([n] \cup [n-1] \cup \dots \cup [n-2r+1]))$ and $G_{0,n-2r-2}$ are isomorphic.} \hfill $\Box$
\end{theorem}

\begin{theorem} [{\cite[Theorem 4.4]{vs14}}] \label{2.ac} { \quad \rm 
		For $n \geq 5$, the underlying graphs of $G_{0,n}-\{0,1,n-1,n\}$ and $G_{0,n-4}$ are isomorphic where $u_j$ is the vertex of $G_{0,n}$ with integral sum label $j$, $j$ = $0,1,\dots,n$.}
\end{theorem}
\begin{proof} 
	Using Theorem \ref{19}, the underlying graphs of $G_{0,n}$ and $G^c_{n+1}$ are isomorphic and from the structure of these graphs (Graphs $G_{0,6}$ and $G^c_{7}$ are isomorphic without vertex labels and are given in Figures 1 and 6.), vertex with integral sum label $j$ in $G_{0,n}$ and vertex with anti-integral sum label $n-j+1$ in $G^c_{n+1}$ are of same degree and thereby the underlying graphs of $G_{0,n}-\{0,1,n-1,n\}$ and $G^c_{n+1} - \{n+1,n,2,1\}$ are isomorphic, $0 \leq j \leq n$. Using Theorem \ref{4}, the underlying graphs of $G^c_{n+1} - \{1,2,n,n+1\}$ and $G^c_{n-3}$ are isomorphic and using Theorem \ref{19}, the underlying graphs of $G^c_{n-3}$ and $G_{0,n-4}$ are isomorphic. This implies, the underlying graphs of $G_{0,n}-\{0,1,n-1, n\}$ and $G_{0,n-4}$ are isomorphic. Hence the result. 
\end{proof}

\begin{theorem} [{\cite[Theorem 4.5]{vs14}}] \label{24} { \quad \rm 
		For $n \geq 3$, the following pairs of underlying graphs of 
		\begin{enumerate}  
			\item[(i)] $G^c_{0,n-1}$ and $G_{0,n-2} \cup K_1(n-1)$ are isomorphic and 
			\item[(ii)]$K_n$ and $G_{0,n-1} \cup (G_{0,n-2} \cup K_1(n-1))$ are isomorphic \\
			where $K_1(n-1)$ is an isolated vertex with label $n-1$. 
	\end{enumerate}}
\end{theorem}   

\begin{proof} \quad We have the underlying graphs of $(G_{0,n-1})^c = $ the underlying graphs of $(K_1(0)+G_{n-1})^c = $ the underlying graphs of $K_1(0) \cup G^c_{n-1}$ whose underlying graph is isomorphic to the underlying graph of $K_1(0) \cup G_{0,n-2}$ which is isomorphic to the underlying graph of $K_1(n-1) \cup G_{0,n-2}$
	, hence result (i) is proved, where $K_1(0)$ and $K_1(n-1)$are isolated vertices in the underlying graphs of $K_1(0) \cup G_{0,n-2}$ and $K_1(n-1) \cup G_{0,n-2}$, respectively. This implies, the underlying graphs of $(G_{0,n-2} \cup K_1(n-1))^c$ and $G_{0,n-1}$ are isomorphic. This implies, the underlying graphs of $K_n$ and $G_{0,n-1} \cup (G_{0,n-2} \cup K_1(n-1))$ are isomorphic, hence result (ii) is proved, where $K_1(n-1)$ is an isolated vertex with label $n-1$. Hence we get the result. 
\end{proof}

\begin{theorem} [{\cite[Theorem 3.3]{vs14}}] \normalfont\label{25}
	{\rm For all $m,n \in \mathbb {N}_0$, $K_{m+n+1} \cong_{wvl} G_{-m,n} \cup G_{0,m-1} \cup G_{0,n-1}$. }  
\end{theorem}

\begin{proof} For $m,n\in\mathbb {N}_0$, $G_{-m,n}$ $\cong$ $K_1(0) * (G_{-m} * G_n)$ and hence, $G^c_{-m,n}$ $\cong$ $K_1(0)$ $\cup$ ${(G_{-m} * G_n)}^c$ = $K_1(0)$ $\cup$ $(G^c_{-m}$ $\cup$ $G^c_n)$  which is isomorphic to the underlying graph of $K_1(0)$ $\cup$ $G_{0,m-1}$ $\cup$ $G_{0,n-1}$, using Theorem \ref{19}. Hence the result.
\end{proof}

\begin{cor} [{\cite[Corollary 3.4]{vs14}}] \normalfont\label{10}
	{\rm  For $n \geq 2$,  
		
		$K_{2n}$ $\cong_{wvl}$ $G_{-(n-1),n}$ $\cup$ $G_{0,-(n-2)}$ $\cup$ $G_{0,n-1}$ 
		$\cong_{wvl}$ $G_{-(n-1),n}$ $\cup$ $G_{0,n-2}$ $\cup$ $G_{0,n-1}$ $\cong_{wvl}$ $2 G_n$ $\cup$ $n.P_2$	
		
		\medskip
		and for all $n \in \mathbb {N}$, $K_{2n+1}$ $\cong_{wvl}$ $G_{-n,n}$ $\cup$ $2 G_{0,n-1}$ $\cong_{wvl}$ $K_1(0) * (G_{-n} * G_{n})$ $\cup$ $2 G_{0,n-1}$. \hfill $\Box$ }  
\end{cor}

\subsection{On families of integral sum graphs $G$ with $\Delta(G)$ = $|V(G)|-1$}

Study on $\mathbb{Z}$-sum graphs $G^+(S)$ with vertex or vertices of degree = $|S|-1$ is presented here. 
  
\begin{dfn}	 \cite{vn09} \quad Let $G$ be a connected graph with maximum degree $\Delta (G)$ and $V_\Delta (G)$ = $\{ x \in V(G) : deg(x) = \Delta (G) \}.$
\end{dfn} 

\begin{theorem}[{\cite[Fact 4]{c98}}]\quad \label{a0} {\rm Let $f$ be an integral sum labeling of a non-trivial graph $G$ of order $n$. Then, $f(x) \neq 0$ for every vertex $x$ of $G$ if and only if $\Delta(G)$ $<$ $n-1.$} \hfill $\Box$
\end{theorem}

\begin{lemma} [{\cite[Lemma 3]{vn09}}] \label{l29} \quad {\rm  Let $G$ be a connected graph and $x,u,v,w \in V(G)$. If $f$ is an integral sum labeling of $G$ with $f(u)$ = $-f(x)$, then for each $v \in N(u)$ such that $f(u)+f(v)$ = $f(w)$, either $w = x$ or $xw \in E(G)$. } 
\end{lemma} 

\begin{proof}\quad Let $f$ be an integral sum labeling of $G$, $f(u)$ = $-f(x)$ and $v \in N(u)$. If $w \neq x$, then $f(u)+f(v)$ = $f(w)$ implies, $f(x)+f(w)$ = $f(v)$ which implies that $xw \in E(G)$, using the definition of integral sum labeling. Hence we get the result.  
\end{proof}

\begin{theorem}[{\cite[Theorem 4]{vn09}}] \quad\label{2.43} {\rm Let $f$ be an integral sum labeling of a graph $G$ of order $n \geq 4$. If $G$ has at least two vertices of degree $n-1$ each, then
\begin{enumerate}
\item[{\rm (i)}]	there exists a vertex $x$ of degree $n-1$ such that $f(x)$ = $0$ and
\item[{\rm (ii)}]	for every vertex $y \neq x$ and of degree $n-1$, there exists a vertex $y'$ with degree $< n-1$ such that $f(y) + f(y')$ = $0.$
\end{enumerate}}
\end{theorem}

\begin{proof}\quad (i)~ The result follows directly from Theorem \ref{a0}.

\medskip
\noindent
(ii)~ Let $V(G)$ = $\{ x, y, v_1, v_2 , . . . ,v_{n-2} \}$, $d(y)$ = $d(x)$ = $n-1$ and $f(x) = 0$. If for $i$ = $1, 2, \dots, n-2$, $f(y) + f(v_i)$ $\neq$ $0$ for all $i$, then without loss of generality, let us assume that  $f(y) + f(v_i)$ = $f(v_{i+1})$ since $y$ is adjacent to all other vertices of $G$, $1$ $\leq$ $i$ $\leq$ $n-3$. Now, $f(y) + f(v_{n-2})$ $\neq$ $f(y)$ since $f(v_{n-2})$ $\neq$ $0$. This implies, $f(y) + f(v_{n-2})$ = $0$ or $f(v_j)$ for some $j$, $1$ $\leq$ $j$ $\leq$ $n-3$. If $f(y) + f(v_{n-2})$ = $f(v_j)$ for some $j$, $1$ $\leq$ $j$ $\leq$ $n-3$, then by applying the relation $f(y) + f(v_i)$ = $f(v_{i+1})$ repeatedly, we get, $f(v_j)$ = $f(y) + f(v_{n-2})$ = $f(y) + (n-2-j)f(y)+ f(v_j)$, $1$ $\leq$ $j$ $\leq$ $n-2$. This implies, $f(v_{n-2})$ = $f(y)$ + $((n-2-(n-2))f(y) + f(v_{n-2})$ = $f(y)$ + $f(v_{n-2})$ which implies, $f(y)$ = 0. This implies, $f(y)$ = 0, a contradiction to the given condition that $f(y) \neq 0$. Hence, it must be the case that $f(y) + f(v_{n-2})$ = 0, which implies that $f(v_{n-2}) + f(v_i) \neq 0$ for $i$ = $1, 2, \dots, n-3$.  

\medskip
We claim that $d(v_{n-2}) < n-1$. Assume that $d(v_{n-2})$ = $n-1$. Then, $v_{n-2}$ is adjacent to all other vertices in $G$. Now, for $i$ = $1, 2, \dots, n-3$, the relations $f(y) + f(v_{n-2})$ = 0 and $f(y) + f(v_i)$ = $f(v_{i+1})$ imply that $f(v_{n-2}) + f(v_{i+1})$ = $f(v_i)$. Furthermore, $f(v_{n-2}) + f(v_1)$ $\neq$ 0 since $f(v_{n-2}) + f(v_1)$ = 0 implies, $f(v_1)$ = $f(y)$. Also, for $i$ = $1,2,\dots,n-3$, $f(v_{n-2})$ + $f(v_1)$ $\neq$ $f(v_i)$, otherwise, $f(v_1)$ = $f(v_{i+1})$, contradicting the fact that $f$ is an integral sum labeling of $G$. This implies that $f(v_{n-2})$ + $f(v_1)$ = $f(y)$, which implies, $f(v_1)$ = $2f(y)$ since we have $f(y) + f(v_{n-2})$ = 0. Substituting the above relation repeatedly in $f(y) + f(v_i)$ = $f(v_{i+1})$, $1 \leq i \leq n-3$, we obtain $f(v_{n-2})$ = $(n-1)f(y)$ which implies that $n.f(y)$ = 0 since $f(y) + f(v_{n-2})$ = 0. This implies, $f(y)$ = 0, a contradiction. Hence, $d(v_{n-2}) < n-1$. By taking $y'$ = $v_{n-2}$, we get the result. 
\end{proof}  

\begin{theorem} [{\cite[Theorem 5]{vn09}}] \quad \label{2.5} {\rm For $n \geq 3$, let $f$ be an integral sum labeling of a connected graph $G$ of order $n$ with at least two vertices each of degree $n-1$. If $y \in V(G)$ such that $d(y)$ = $n-1$ and $f(y) \neq 0$, then for every vertex $v \in V(G)$, $f(v)$ = $k.f(y)$ where $k \in \{-1,0,1,2,\dots, n-2\}$. }
\end{theorem}

\begin{proof}\quad Using Theorem \ref{2.43}(i), there exists a vertex, say, $x$ of degree $n-1$ with $f(x)$ = $0$. Let $V(G)$ = $\{ x, y, v_1, v_2 , . . . ,v_{n-2} \}$ with $d(y)$ = $n-1$. Without loss of generality, we assume that  $f(y) + f(v_i)$ = $f(v_{i+1})$, $1$ $\leq$ $i$ $\leq$ $n-3$. Suppose there exists a vertex $v_i \in V(G)\setminus \{ x, y \}$, such that $f(v_i) \neq k.f(y)$ for every non-zero integer $k$, $1 \leq i \leq n-2$. Now, for $i$ = $1, 2, \dots, n-3$, we have $f(v_{i+1})$ = $f(y)+f(v_i)$, and thus $f(v_i)$ = $f(y) - f(v_{i+1})$ which implies that the label of every vertex of $V(G)\setminus \{ x, y \}$ is not an integer multiple of $f(y)$, a contradiction to Theorem \ref{2.43}(ii). Hence, $f(v)$ = $k.f(y)$ for every $v \in V(G)$, where $k$ is an integer.
	
	\medskip
	We claim that $-(n-2) \leq k \leq 1$. We already have $f(x)$ = 0 = $0.f(y)$ and $f(y)$ = $1.f(y)$. Let $f(v_1)$ = $k.f(y)$ for some non-zero integer $k$. Then, using the relation $f(v_{i+1})$ = $f(y)+f(v_i)$ for $i$ = $1, 2, \dots, n-3$, we obtain $f(v_j)$ = $(k+j-1)f(y)$, $j$ = $1, 2, \dots, n-2$. This implies that $f(v_{n-2})$ = $(k+n-3)f(y)$, which implies $k+n-3$ = -1. This implies that $k$ = $-(n-2)$, $f(v_1)$ = $-(n-2) f(y)$, $f(v_2)$ = $-(n-3) f(y)$, $\dots$, $f(v_{n-2})$ = $-(n-(n-1)) f(y)$ = $- f(y)$, $f(x)$ = $0.f(y)$ and $f(y)$ = $1.f(y)$. Thus for every vertex $v \in V(G)$, $f(v)$ = $k.f(y)$, where $k \in \{ 0, 1, -1, -2, \dots, -(n-2) \}$. But $G^+(S) \cong_{wvl} G^+(-S)$ for any non-empty $S \subset \mathbb{Z}$. Hence the result.
\end{proof}

\begin{theorem} [{\cite[Theorem 6]{vn09}}] \quad \label{3.10} {\rm Every integral sum graph $G$, except $G_{-1,1} \cong K_3$, has at the most two vertices of degree $\left|V(G)\right| - 1.$}
\end{theorem}

\begin{proof}\quad All the three vertices of integral sum graph $G_{-1,1} \cong K_3$ are of degree 2 = $\left|V(K_3)\right| - 1$. If the theorem is not true, then let $G$ be an integral sum graph of order $n$ with at least 3 of its vertices,  each be of degree $n-1$, $n \geq 4$. Let $f$ be an integral sum labeling of $G$. Then, by Theorem \ref{a0}, one vertex, say, $x$ of $G$ has degree $n-1$ and $f(x)$ = 0. Suppose, $y,z \in V(G)\setminus \{ x \}$ such that $d(y)$ = $d(z)$ = $n-1$, and $f(y),f(z) \in \mathbb{Z} \setminus \{ 0 \}$. Then, by  Theorem \ref{2.5}, every vertex label of $G$ must be an integer multiple of $f(y)$, in particular, let $f(z)$ = $k.f(y)$, where $k \in \mathbb{Z} \setminus \{ 0,1,-1 \}$ since $(i)$ $k$ = 0 implies, $f(z)$ = 0 = $f(y)$; $(ii)$ $k$ = 1 implies, $f(z)$ = 1 = $f(y)$; and $(iii)$ $k$ = -1 implies, $f(y)$ + $f(z)$ = 0 which implies, $d(z) < n-1$ by Theorem \ref{2.43}(ii). This implies, $f(y)$ = $(1/k) f(z)$ where $k\in\mathbb{Z}\setminus\{0,1,-1\}$. This is a contradiction since $f(y)$ must be an integer. Hence the result follows. 
\end{proof}

\begin{theorem} [{\cite[Theorem 7]{vn09}}] \quad \label{3.11} {\rm For $n \geq 4$, there is an integral sum graph (denoted by $G_{\Delta n}$) of order $n$ with exactly two vertices of degree $n-1$. This graph is unique up to isomorphism. } 
\end{theorem}

\begin{proof}\quad Let $G$ be an integral sum graph of order $n \geq 4$ and with exactly two vertices of degree $n-1$. Let $V(G)$ = $\{ v, v_0, v_1, v_2, \cdots, v_{n-2} \}$, where $d(v)$ = $d(v_0)$ = $n-1$. Let $f$ be an integral sum labeling of $G$ with $f(v_0)$ = 0 and $f(v)$ = $k$, $k \in \mathbb{N}$. Then, using Theorem \ref{2.5}, $f(v_i)$ = $-ik$, $i$ = $0,1,2,\dots,n-2$ and for $i \neq j$ and $0 \leq i,j \leq n-2$, vertices $v_i$ and $v_j$ are adjacent in $G$ if and only if $f(v_i)$ + $f(v_j)$ $\geq$ $-(n-2)k$; thus if and only if $i+j \leq n-2$. This condition is independent of choice of the integer $k$. This proves the uniqueness of unlabeled graph $G_{\Delta n}$. In the above labeling by taking $k$ = 1, we get $G_{\Delta n}$ $\cong$ $G_{-1,n-2}$, $n \geq 4$. Hence we get the result.  
\end{proof}

\begin{theorem} [{\cite[Theorem 10]{vn09}}] \quad \label{3.12} {\rm  For $m,n \in \mathbb{N}$, $G_{-1,n}$ $\cong$ $G_{\Delta (n+2)}$ and for $2 \leq m \leq n$, $G_{-m,n}$ has exactly one vertex of degree $m + n$ = $\Delta(G_{-m,n})$. }
\end{theorem}

\begin{proof}\quad For $n \in \mathbb{N}$ and $S$ = $\{ -1, 0, 1, 2, \dots, n \}$, $G^+(S)$ = $G_{-1,n}$. In this graph, vertices with labels $0$ and $-1$ are the only vertices that are adjacent to all other vertices when $n \geq 2$ and thereby each one of these two vertices is of degree $n+1$ = $\Delta(G_{-1,n})$. Clearly, $G^+(S)$ = $G_{-1,n}$ = $G_{\Delta (n+2)}$ for $n+2 \geq 4$. For $m,n \in N$, $2 \leq m \leq n$ and $S$ = $\{ -m$, $-(m-1)$, $\dots$, $-1$, $0$, $1$, $2$, $\dots$, $n \}$, $G^+(S)$ = $G_{-m,n}$ and in this graph vertices with labels $-1$ and $-m$ are not adjacent whereas the vertex with vertex label $0$ is the only vertex which is adjacent to all other vertices and thereby it is of degree $m+n$ = $\Delta(G_{-m,n})$. Hence we get the result. 
\end{proof}

\subsection{Hamiltonian graph property of integral sum graphs $G_{-r, s}$, $r,s\in\mathbb{N}$}

Any sum graph is non-Hamiltonian as it contains isolated vertex or vertices. But integral sum graphs $G_{-r, s}$ need not be non-Hamiltonian for all values of $r$ and $s$, $r \in \mathbb{N}_0$ and $s \in \mathbb{N}$. Here, we present Hamiltonian and and non-Hamiltonian graph properties of integral sum graphs $G_{0, s}$ and $G_{-r, s}$, $r,s\in \mathbb{N}$. 

\begin{theorem} [{\cite[Theorem 11]{vn09}}] \quad \label{3.13} {\rm For $n$ $\geq$ $3$, graph $G_{\Delta n}$ is Hamiltonian.}
\end{theorem}

\begin{proof}\quad For $n \geq 3$, graph $G_{\Delta n}$ $\cong$ $G_{-1,n-2}$ using Theorem \ref{3.11}. In this graph vertices with labels 0 and -1 are adjacent and each is adjacent to all other vertices and hence $\delta (G_{\Delta n} ) \geq 2$ and $G_{\Delta n}$ is a 2-connected graph for $n \geq 3$. This implies, $d(u) \geq 2$ for every vertex $u \in V(G_{\Delta n})$. Let $f$ be an integral sum labeling of $G_{\Delta n};$ $v, v_0, v_1, v_2, \cdots, v_{n-2}$  be its vertices and -1, 0, 1, 2, . . . , $n$-2 be the corresponding vertex labels, $n \geq 3$ and $n\in \mathbb{N}$. Now, when $n$ is even, consider the following sequence of vertices, $v$, $v_{n-2}$, $v_0$, $v_{n-3}$, $v_1$, $v_{n-4}$, $v_2$, $\cdots$, $v_{n/2}$, $v_{(n-4)/2}$, $v_{(n-2)/2}$, $v$. This sequence contains all the $n$ vertices of graph $G_{\Delta n}$ for $n \geq 3$. Moreover, any two consecutive elements of the sequence are adjacent vertices in $G_{\Delta n}$, follows from the definition of integral sum labeling. Hence, it is a Hamiltonian cycle of $G_{\Delta n}$. Similarly, when $n$ is odd, the sequence of vertices $v$, $v_{n-2}$, $v_0$, $v_{n-3}$, $v_1$, $v_{n-4}$, $v_2$, $\cdots$, $v_{(n-5)/2}$, $v_{(n-1)/2}$, $v_{(n-3)/2}$, $v$ is a Hamiltonian cycle of $G_{\Delta n}$. Hence, for $n \geq 3$, $G_{\Delta n}$ is a Hamiltonian graph. 
	
	Figures 15 and 16 represent integral sum graphs $G_{\Delta 8}$ and  $G_{\Delta 9}$. Each of the graph contains a Hamiltonian cycle, as given in the proof of the above theorem, and is shown with dotted lines.
\end{proof}

\begin{center}
	\begin{tikzpicture}[scale =0.8]
		
		\node (a0) at (18,3)  [circle,draw,scale=0.6]{0};
		\node (a1) at (19.5,1.75)  [circle,draw,scale=0.6]{1};
		\node (a2) at (19.5,.5)  [circle,draw,scale=0.6] {2};
		\node (a3) at (18,-.75)  [circle,draw,scale=0.6] {3};
		\node (a4) at (16,-.75)  [circle,draw,scale=0.6] {4};
		\node (a5) at (14.5,0)  [circle,draw,scale=0.6]{5};
		\node (a6) at (14,1.2)  [circle,draw,scale=0.6] {6};
		\node (a7) at (16,3)  [circle,draw,scale=0.6] {-1};
		
		\draw (a0)  -- (a1);
		\draw (a0) -- (a2);
		\draw (a0) -- (a3);
		\draw (a0) -- (a4);
		\draw (a0) [line width=0.2mm] [dashed] -- (a5);
		\draw (a0) [line width=0.2mm] [dashed] -- (a6);
		\draw (a0) -- (a7);
		
		\draw (a7) -- (a0);
		\draw (a7) -- (a1);
		\draw (a7) -- (a2);
		\draw (a7) [line width=0.2mm] [dashed] -- (a3);
		\draw (a7) -- (a4);
		\draw (a7) -- (a5);
		\draw (a7) [line width=0.2mm] [dashed] -- (a6);
		
		\draw (a1) -- (a2);
		\draw (a1) -- (a3);
		\draw (a1) [line width=0.2mm] [dashed] -- (a4);
		\draw (a1) [line width=0.2mm] [dashed] -- (a5);
		
		\draw (a2) [line width=0.2mm] [dashed] -- (a3);
		\draw (a2) [line width=0.2mm] [dashed] -- (a4);
		
		\node (c0) at (25.5,3.5)  [circle,draw,scale=0.6] {0};
		\node (c1) at (27,2.75)  [circle,draw,scale=0.6]{1};
		\node (c2) at (27.5,1.5)  [circle,draw,scale=0.6]{2};
		\node (c3) at (26.5,-0.25)  [circle,draw,scale=0.6] {3};
		\node (c4) at (24,-.5)  [circle,draw,scale=0.6] {4};
		\node (c5) at (22.5,0)  [circle,draw,scale=0.6] {5};
		\node (c6) at (21.75,1)  [circle,draw,scale=0.6]{6};
		\node (c7) at (21.5,2.25)  [circle,draw,scale=0.6] {7};
		\node (c8) at (23,3.5)  [circle,draw,scale=0.6] {-1};
		
		\draw (c8) -- (c0);
		\draw (c8) -- (c1);
		\draw (c8) -- (c2);
		\draw (c8) -- (c4);
		\draw (c8) [line width=0.2mm] [dashed] -- (c3);
		\draw (c8) -- (c4);
		\draw (c8) -- (c5);
		\draw (c8) -- (c6);
		\draw (c8) [line width=0.2mm] [dashed] -- (c7);
		
		\draw (c0) -- (c1);
		\draw (c0) -- (c2);
		\draw (c0) -- (c3);
		\draw (c0) -- (c4);
		\draw (c0) -- (c5);
		\draw (c0) [line width=0.2mm] [dashed] -- (c6);
		\draw (c0) [line width=0.2mm] [dashed] -- (c7);
		
		\draw (c1) -- (c2);
		\draw (c1) -- (c3);
		\draw (c1) -- (c4);
		\draw (c1) [line width=0.2mm] [dashed] -- (c5);
		\draw (c1) [line width=0.2mm] [dashed] -- (c6);
		
		\draw (c2) -- (c3);
		\draw (c2) [line width=0.2mm] [dashed] -- (c4);
		\draw (c2) [line width=0.2mm] [dashed] -- (c5);
		
		\draw (c3) [line width=0.2mm] [dashed] -- (c4);
		
	\end{tikzpicture}
	
	\medskip	
	{\small	Fig. 15. $G_{\Delta 8}$ with Hamiltonian cycle
		\hspace{1cm} Fig. 16. $G_{\Delta 9}$ with Hamiltonian cycle }
\end{center}

\begin{theorem} [{\cite[Theorem 2.3]{vm12}}] \quad \label{3.14} {\rm For $m,n\in \mathbb{N}$, graph $G_{-m, n}$ is Hamiltonian whereas $G_{0, n}$ is non-Hamiltonian.}
\end{theorem}

\begin{proof}\quad For $n\in \mathbb{N}$, graph $G_{n}$ contains an isolated vertex and thereby graph $G_{0, n}$ $\cong$ $K_1(0) * G_n$ contains a vertex of degree 1. This implies, integral sum graph $G_{0, n}$ is non-Hamiltonian for $n\in \mathbb{N}$. 
	
	In Theorem \ref{3.13}, it is proved that graph $G_{\Delta n}$ =  $G_{-1, n-2}$ is Hamiltonian for $n \geq 3$. Consider $G_{-m, n}$ for $m \geq 2$, $n \geq 3$ and $m,n\in \mathbb{N}$. Let $w_0$, $u_1$, $u_2$, $\cdots$, $u_{m}$, $v_1$, $v_2$, $\cdots$, $v_{n}$ be the vertices of $G_{-m, n}$ with integral sum labeling 0,-1,-2,..., $-m$, 1,2,..., $n$, respectively. For $m \geq 2$, $n \geq 3$ and $m,n\in \mathbb{N}$, using the definition of integral sum labeling, $w_0 v_{n-1}$, $v_{n-1} v_1$, $v_1 v_{n-2}$, $v_{n-2} v_2$, . . . , $v_{\left\lfloor \frac{n}{2} \right\rfloor } u_{m}$, $u_{m} v_n$, $v_n u_1$, $u_1 u_{m-1}$, $u_{m-1} u_2$, $u_2 u_{m-2}$, . . . , $u_{\left\lfloor \frac{n+1}{2} \right\rfloor }$ $w_0$ are edges in $G_{-m, n}$ as the integral sum labels of $u_i$ and $v_j$ are $-i$ and $j$, respectively, $1 \leq i \leq m$ and $1 \leq j \leq n$. Thus, $w_0$ $v_{n-1}$ $v_1$ $v_{n-2}$ $v_2$ . . . $v_{\left\lfloor \frac{n}{2} \right\rfloor }$ $u_{-m}$ $v_n$ $u_1$ $u_{m-1}$ $u_2$ $u_{m-2}$ . . . $u_{\left\lfloor \frac{n+1}{2} \right\rfloor }$ $w_0$ is a Hamiltonian cycle in $G_{-m, n}$ for $m \geq 2$, $n \geq 3$ and $m,n\in \mathbb{N}$. Hence the result.
	
	Integral sum graphs $G_{-2,3}$ and $G_{-3,3}$ are shown in Figures 17 and 18. Each graph contains a Hamiltonian cycle which is indicated by dotted lines. 
\end{proof}

\begin{center}
	\begin{tikzpicture}[scale =0.9]
		
		\node (z) at (18,2.5)  [circle,draw,scale=0.6]{0};
		\node (y1) at (20,1.5)  [circle,draw,scale=0.6]{1};
		\node (y2) at (20,0)  [circle,draw,scale=0.6] {2};
		\node (y3) at (18,-.75)  [circle,draw,scale=0.6] {3};
		\node (x1) at (16,1.5)  [circle,draw,scale=0.6] {-1};
		\node (x2) at (16,0)  [circle,draw,scale=0.6]{-2};
		
		\draw (z)  -- (y1);
		\draw (z)  [line width=0.2mm] [dashed] -- (y2);
		\draw (z)  -- (y3);
		\draw (z) [line width=0.2mm] [dashed]  -- (x1);
		\draw (z)  -- (x2);
		\draw (x1)  -- (y1);
		\draw (y1) [line width=0.2mm] [dashed] -- (y2);
		\draw (y1) [line width=0.2mm] [dashed] -- (x2);
		
		\draw (y2) -- (x1);
		\draw (y2)  -- (x2);
		
		\draw (y3) [line width=0.2mm] [dashed] -- (x1);
		\draw (y3) [line width=0.2mm] [dashed]  -- (x2);
		
		\node (w) at (26,2.5)  [circle,draw,scale=0.6]{0};
		\node (v1) at (28,1.5)  [circle,draw,scale=0.6]{1};
		\node (v2) at (28,0)  [circle,draw,scale=0.6] {2};
		\node (v3) at (27,-.75)  [circle,draw,scale=0.6] {3};
		\node (u1) at (24,1.5)  [circle,draw,scale=0.6] {-1};
		\node (u2) at (24,0)  [circle,draw,scale=0.6]{-2};
		\node (u3) at (25,-.75)  [circle,draw,scale=0.6]{-3};
		
		\draw (w)  -- (v1);
		\draw (w) [line width=0.2mm] [dashed] -- (v2);
		\draw (w)  -- (v3); 
		\draw (w)  -- (u1);
		\draw (u2) [line width=0.2mm] [dashed] -- (w);
		\draw (w)  -- (u3);
		\draw (v1) [line width=0.2mm] [dashed] -- (v2);
		\draw (u1) [line width=0.2mm] [dashed] -- (u2);
		\draw (v1) [line width=0.2mm] [dashed] -- (v3);
		\draw (u1) [line width=0.2mm] [dashed] -- (u2);
		\draw (u1) [line width=0.2mm] [dashed] -- (u3); 
		\draw (u1)  -- (v2);
		\draw (u2)  -- (v1); 
		\draw (u1)  -- (v1);
		\draw (u2)  -- (v2);
		\draw (u3) [line width=0.2mm] [dashed] -- (v3);
		
	\end{tikzpicture}
	\medskip
	
	{\small Fig. 17. $G_{-2,3}$ showing Hamiltonian cycle \hspace{1cm}  Fig. 18. $G_{-3,3}$ showing Hamiltonian cycle}
\end{center}

\begin{oprm} [{\cite[Open problem 1]{vm12}}] \quad \label{3.15} {\rm Find the number of distinct Hamiltonian cycles that exist in the integral sum graphs $G_{-m,n}$, $m,n\in \mathbb{N}$. \hfill $\Box$}
\end{oprm}

\section{On Maximal integral sum Graphs}

In \cite{vm12}, Vilfred and Florida investigated maximal integral sum graphs. Their main results are (i) For $m \in  \mathbb{N}$, $G_{-m,m}$ is the maximal integral sum graph of order $2m+1$ and for $m \geq 2$, $G_{-(m-1),m}$ is the maximal integral sum graph of order $2m$; (ii) $G_{0,n}$ is a spanning subgraph of $G_{-1,n-1}$ without vertex labels, $n\in \mathbb{N}$; and (iii) For $n \geq 4$, $G_{0,n}$ is not a maximal integral sum graph of order $n+1$. These are presented in this section.  

\begin{dfn}	 {\rm \cite{vm12}} \quad An integral sum graph or sum graph with underlying graph $G$ is said to be {\em maximal} if $G$ is not a spanning subgraph of the underlying graph of any other integral sum graph or sum graph, respectively. 
	
A {\em sum-maximum integral sum graph} of order $n$ is a maximal integral sum graph of order $n$ with maximum size.
\end{dfn} 

Clearly, if integral sum graph $G^+(S)$ is such that $0 \in S$, then the vertex with label $0$ has the maximum degree $\left|V(G^+(S))\right|$ $-1$ in $G^+(S)$, using Theorem \ref{a0}. Hence, {\em a sum-maximum integral sum graph contains $0$ as a vertex label.} 

\subsection{Comparison of integral sum graphs $G_{-k,n-k}$ and $G_{-r,n-r}$} 

We start with presenting results on isomorphism of complete multipartite graphs of same order. These results are used to compare integral sum graphs $G_{-m,n}$ of same order of $m+n+1$, $m\in \mathbb{N}_0$ and $n\in \mathbb{N}$. Comparing of integral sum graphs of same order means comparing of integral sum graphs of same order without vertex labels until otherwise it is mentioned.

\begin{lemma} [{\cite[Lemma 2.6]{vm12}}] \quad {\rm  Let $1 \leq m_1 \leq m_2 \leq \dots \leq m_k$ and $1 \leq r_1 \leq r_2 \leq \dots \leq r_n$. $K_{m_1,m_2,\dots,m_k} \cong K_{r_1,r_2,\dots,r_n}$ if and only if $k$ = $n$ and $m_i = r_i$ for $i = 1,2,\dots,n$. }
\end{lemma}
\begin{proof} Given, $1 \leq m_1 \leq m_2 \leq \dots \leq m_k$ and $1 \leq r_1 \leq r_2 \leq \dots \leq r_n$. Then, 

 \medskip
 $K_{m_1,m_2,\dots,m_k} \cong K_{r_1,r_2,\dots,r_n}$ if and only if $K^c_{m_1,m_2,\dots,m_k} \cong K^c_{r_1,r_2,\dots,r_n}$ 

\hspace{4.2cm} if and only if $K_{m_1} \cup K_{m_2} \cup \dots \cup K_{m_k} \cong K_{r_1} \cup K_{r_2} \cup \dots \cup K_{r_n}$ 

\hspace{4.2cm} if and only if $k = n$ and $m_1 = r_1$, $m_2 = r_2$, \dots, $m_n = r_n$  since 
\\
 $1 \leq m_1 \leq m_2 \leq \dots \leq m_k$ and $1 \leq r_1 \leq r_2 \leq \dots \leq r_n$. Hence we get the result.
\end{proof}

\begin{lemma} [{\cite[Lemma 2.7]{vm12}}] \quad \label{2.8} {\rm Let $1 \leq m_1 \leq m_2 \leq \dots \leq m_k$, $1 \leq r_1 \leq r_2 \leq \dots \leq r_n$ and $m_1+m_2+\dots+m_k$ = $r_1+r_2+\dots+r_n$. Then, $K_{m_1,m_2,\dots,m_k}$ and $K_{r_1,r_2,\dots,r_n}$ are comparable if and only if $k$ = $n$ and $m_i = r_i$ for $i = 1,2,\dots,n$. In that case, the two graphs are same.  }
\end{lemma}

\begin{proof} Without loss of generality, assume that $k \leq n$. Then, 

\medskip
\noindent
$K_{m_1,m_2,\dots,m_k}$ and $K_{r_1,r_2,\dots,r_n}$ are comparable if and only if $K^c_{m_1,m_2,\dots,m_k}$ and $K^c_{r_1,r_2,\dots,r_n}$ are comparable if and only if $K_{m_1} \cup K_{m_2} \cup \dots \cup K_{m_k}$ and $K_{r_1} \cup K_{r_2} \cup \dots \cup K_{r_n}$ are comparable. If the result is not true, then let $i$ be the smallest integer such that $m_1 = r_1$, $m_2 = r_2$, \dots, $m_{i-1} = r_{i-1}$ and $m_i < r_i$, $1 \leq i \leq k$. Then, $K_{m_i}$ is a proper subgraph of $K_{r_i}$ and since  $1 \leq m_1 \leq m_2 \leq \dots \leq m_k$, ~$1 \leq r_1 \leq r_2 \leq \dots \leq r_n$ and $m_1+m_2+\dots+m_k$ = $r_1+r_2+\dots+r_n$, there exists $j > i$ such that $m_j > r_j$, $2 \leq j \leq k$ and thereby $K_{m_j}$ is a proper super graph of $K_{r_j}$. Now, $1 \leq i < j \leq k$, $m_i < r_i$,~ $m_i \leq m_j$,~ $r_i \leq r_j$ and $r_j < m_j$ which implies, $K_{m_i}$ is a proper subgraph of $K_{r_i}$ which is a subgraph of $K_{r_j}$ which is a proper subgraph of $K_{m_j}$, $1 \leq i < j \leq k$. This implies that $K_{m_i} \cup K_{m_j}$ and $K_{r_i} \cup K_{r_j}$ are non-comparable graphs and thereby $K_{m_1} \cup K_{m_2} \cup \dots \cup K_{m_k}$ and $K_{r_1} \cup K_{r_2} \cup \dots \cup K_{r_n}$ are non-comparable when $m_i < r_i$,~ $m_j > r_j$ and $1 \leq m_1 \leq m_2 \leq \dots \leq m_k$,~ $1 \leq r_1 \leq r_2 \leq \dots \leq r_n$,~ $m_1+m_2+\dots+m_k$ = $r_1+r_2+\dots+r_n$ and $1 \leq i < j \leq k$. Hence, the result follows by the method of contradiction.  
\end{proof}

\begin{lemma} [{\cite[Lemma 2.8]{vm12}}] \quad \label{2.9} {\rm Let $k,n,r\in\mathbb{N}$, $k \leq n-k$, $r \leq n-r$, and $k \neq r$. $K_{k,n-k}$ and $K_{r,n-r}$ are non-comparable graphs.  }
\end{lemma}

\begin{proof} Similar to the proof given in Lemma \ref{2.8}.
\end{proof}

Results on comparison of $\mathbb{Z}$-sum graphs $G_{-r, s}$ using results obtained on complete multi-partite graphs are presented below. 

\begin{theorem} \quad \label{5.1} {\rm Let $k \leq n-k$, $r \leq n-r$ and $k,n,r\in N$. Then, $G_{-k,n-k}$ and $G_{-r,n-r}$, without the vertex labels, are comparable if and only if $k$ = $r$. }
\end{theorem}

\begin{proof} Using the definition of $G_{-m,n}$, integral sum graphs $G_{-k,n-k}$ and $G_{-r,n-r}$ are comparable if and only if $K_1(0) * (G_{-k} * G_{n-k})$ and $K_1(0) * (G_{-r} * G_{n-r})$ are comparable if and only if $G_{-k} * G_{n-k}$ and $G_{-r} * G_{n-r}$ are comparable if and only if $G_{-k} \cup G_{n-k} \cup K_{k,n-k}$ and $G_{-r} \cup G_{n-r} \cup K_{r,n-r}$ are comparable if and only if $K_{k,n-k}$ and $K_{r,n-r}$ are comparable and $G_{-k} \cup G_{n-k}$ and $G_{-r} \cup G_{n-r}$ are comparable, $k \leq n-k$ and $r \leq n-r$ if and only if $k = r$, using Lemma \ref{2.9}. 
\end{proof}

\begin{cor} [{\cite[Corollary 2.13]{vm12}}] \quad \label{5.1a} {\rm Let $k \leq n-k$, $r \leq n-r$ and $k,n,r\in N$. Then, for $k$ $\neq$ $r$, $G_{-k,n-k}$ and $G_{-r,n-r}$ are non-comparable graphs.  \hfill $\Box$}
\end{cor}

\begin{theorem} [{\cite[Theorem 2.17]{vm12}}] \quad \label{5.2} {\rm  Let $n\in \mathbb{N}$ and $n \geq 3$. $G_{0,n}$ is a spanning subgraph of $G_{-1,n-1}$ without vertex labels. }
\end{theorem} 

\begin{proof} For $n \geq 3$, using Corollary \ref{3}, we have $G^+([2,n])$ is a subgraph of $G_{n-1}$. Using Theorem \ref{14}(a), $G_n$ is a spanning subgraph of $G_1 * G^+([2,n])$. This implies that $G_n$ is a spanning subgraph of $G_{1} * G_{n-1}$. This then implies, $K_1 * G_n$ = $G_{0,n}$ is a spanning subgraph of $K_1 * (G_1 * G_{n-1})$ $\cong$ $K_1(-1) * (K_1(0) * G_{n-1})$ $\cong$ $G_{-1,n-1}$. Hence we get the result.
\end{proof}

\begin{theorem} [{\cite[Theorem 3.6]{vm12}}] \quad \label{5.3} {\rm For $n \geq 3$ and $n \in  \mathbb{N}$, $G_{0,n}$ is not a maximal integral sum graph of order $n+1.$} 
\end{theorem}

\begin{proof} Comparing underlying graphs, clearly,
$G_{0,3}$ is a proper subgraph of $G_{-1,2}$ which is a maximal integral sum graph of order 4. By Theorem \ref{5.2}, 
$G_{0,n}$ is a spanning subgraph of $G_{-1, n-1}$, $n\in\mathbb{N}$.
\end{proof}

\subsection{Maximal integral sum graphs $G_{-r, s}$, $r,s\in\mathbb{N}$}

In 2012, Vilfred \cite{vm12} proved that for $r,s\in\mathbb{N}$, integral sum graphs $G_{-r, s}$ are maximal integral sum graphs of order $n$ = $r+s+1$ and in 2013, Tiwari and Tripathi \cite{tt13} independently obtained result on maximum size of integral sum graph of a given order. These are presented in this section. We start with the following two lemmas related to it.

\begin{lemma} [{\cite[Lemma 3.1]{vm12}}] \label{5.10} \quad {\rm {\small Let $k,m \in  \mathbb{N}$ and $S$ = $[-m, m]$.  
\\
$(i)$ \ $E(G^+(S \cup \{2m\}))$ = $E(G^+(S)) \cup \{(0,2m),(-m,2m)\}$ and 

 ~~~$E(G^+(S \cup \{-2m\}))$ = $E(G^+(S)) \cup \{(0,-2m),(m,-2m)\}$. 
\\
$(ii)$  $E(G^+(S \cup \{2m+k\}))$ = $E(G^+(S)) \cup \{(0,2m+k)\}$ and 

~~ $E(G^+(S \cup \{-2m-k\}))$ = $E(G^+(S)) \cup \{(0,-2m-k)\}$. 
\\
$(iii)$ For $1 \leq k < m$, 

~~~$||G^+(S \cup \{m+k\})||$ = $||G^+(S \cup \{-m-k\})||$ = $||G^+(S)||$ + 1 + $\left\lfloor \frac{3(m-k+1)}{2}\right\rfloor$.}  \hfill $\Box$ }
\end{lemma}

\begin{proof} Statements (i) and (ii) follow from the definition of integral sum labeling.
	\\
	(iii) $G^+(S)$ = $G_{-m,m}$. For $1 \leq k < m$, when we include $m+k$ with the elements of $S$, then the additional edges other than those of $G_{-m,m}$ in the resulting integral sum graph $G^+(S \cup \{m+k\})$ are $(0, m+k)$, $(m+k,-k-i)$, and $(k+j, m-j)$ for $i$ = $0,1,\dots,m-k$ and $j$ = $0,1,\dots,\left\lfloor (m-k-1)/2 \right\rfloor$. Thus,
	
	$|E(G^+(S \cup \{m+k\}))|$ = $|E(G^+(S))|$ + 1 + $(m-k+1)$ + $\left\lfloor (m-k-1)/2 \right\rfloor$ + 1
	
\hspace{3.4cm}	= $|E(G^+(S))|$ + 1+  $(m-k+1)$ + $\left\lfloor (m-k+1)/2 \right\rfloor$
	
\hspace{3.4cm}		= $|E(G^+(S))|$ + 1 + $\left\lfloor 3(m-k+1)/2 \right\rfloor$ 
	\\
	since for any positive integer $x$, $x + \left\lfloor x/2 \right\rfloor$ = $\left\lfloor 3x/2 \right\rfloor$.
\end{proof}

\begin{lemma} [{\cite[Lemma 3.2]{vm12}}] \label{5.11} \quad {\rm {\small Let $h,k,m \in  \mathbb{N}$ and $S = [-m, m]$.  
\\
$(1.1)$~ {$||G^+(S \cup \{m+k,-m-k\})||$ = $||G^+(S)|| + 3 + 2 \left\lfloor \frac{3(m-k+1)}{2}\right\rfloor$ if $1 \leq k \leq m-1$.}
 \\    
$(1.2)$~ {$||G^+(S \cup \{m+k,-m-h\})||$ = $||G^+(S)|| + 3 + \left\lfloor \frac{3(m-k+1)}{2}\right\rfloor$ $+$ $\left\lfloor \frac{3(m-h+1)}{2}\right\rfloor$ if $1 \leq k,h \leq m-1$.}
\\   
$(2.1)$~ {$||G^+(S \cup \{2m,-m-k\})||$ = $||G^+(S)|| + 4 + \left\lfloor \frac{3(m-k+1)}{2}\right\rfloor$  if $1 \leq k \leq m-1$.}
\\
$(2.2)$~ {$||G^+(S \cup \{2m,-2m\})||$ = $||G^+(S)|| + 5$.}
 \\   
$(3.1)$~ {$||G^+(S \cup \{m+k,-2m-k\})||$ = $||G^+(S)|| + 3 + \left\lfloor \frac{3(m-k+1)}{2}\right\rfloor$ if $1 \leq k \leq m-1$.}
\\
$(3.2)$~ {$||G^+(S \cup \{m+k,-2m-h\})||$ = $||G^+(S)|| + 2 + \left\lfloor \frac{3(m-k+1)}{2}\right\rfloor$ if $1 \leq k \leq m-1$ and $k < h$.}
\\
$(3.3)$~ {$||G^+(S \cup \{m+k,-2m-h\})||$  = $||G^+(S)|| + 3 + \left\lfloor \frac{3(m-k+1)}{2}\right\rfloor$ if $1 \leq h < k \leq m-1$.}
\\
$(4.1)$~ {$||G^+(S \cup \{2m,-3m\})||$ = $||G^+(S)|| + 4$.}
\\    
$(4.2)$~ {$||G^+(S \cup \{2m,-2m-k\})||$ = $||G^+(S)|| + 4$ if $1 \leq k \leq m-1$.}
\\
$(4.3)$~ {$||G^+(S \cup \{2m,-2m-k\})||$ = $||G^+(S)|| + 3$ if $m < k$.}
\\
$(5.1)$~  {$||G^+(S \cup \{2m+k,-2m-h\})||$ = $||G^+(S)|| + 3$ if $k-h$ = 0 or $1 \leq k-h \leq m$ or $-m \leq k-h \leq -1$. }
\\
$(5.2)$~ {$||G^+(S \cup \{2m+k,-2m-h\})||$  = $||G^+(S)|| + 2$ if $m < k-h$ or $k-h < -m$.}} \hfill $\Box$}
\end{lemma}

\begin{theorem} [{\cite[Theorem 3.3]{vm12}}] \quad \label{5.12} {\rm  Let $S$ = $[-m,m]$, $m \geq 2$, $m,x,y\in  \mathbb{Z}$ and $x,y \notin S$. Then, $||G^+(S \cup \{x,y\})||$ is maximum when $\{x,y\}$ = $\{m+1,-m-1\}.$} 
\end{theorem}

\begin{proof} Let $M$ be the maximum value of $||G^+(S \cup \{x,y\})||$ for all possible values of $x$ and $y$, $x,y\in  \mathbb{Z}$. Using Lemma \ref{5.11}, value of $M$ for different possible values of $x$ and $y$ under different cases is obtained as follows: When 
	
$(1.1)$~ {$x = m+k$ and $y = -m-k$, for $1 \leq k \leq m-1$, 
	
	\hspace{1.5cm} $M$ = $||G^+(S)||$ + 3 + $2 \left\lfloor \frac{3(m-k+1)}{2}\right\rfloor$;}

$(1.2)$~ { $x = m+k$ and $y = -m-h$, for $1 \leq k,h \leq m-1$ and $k \neq h$, 
	
	\hspace{1.5cm} $M$ = $||G^+(S)||$ + 3 + $\left\lfloor \frac{3(m-k+1)}{2}\right\rfloor$ + $\left\lfloor \frac{3(m-h+1)}{2}\right\rfloor$;}

$(2.1)$~ { $x = 2m$ and $y = -m-k$, for $1 \leq k \leq m-1$,
	
	\hspace{1.5cm} $M|$ = $||G^+(S)||$ + 4 + $\left\lfloor \frac{3(m-k+1)}{2}\right\rfloor$;}

$(2.2)$~ {$ x = 2m$ and $y = -2m$,
	
	\hspace{1.5cm} $M$ = $||G^+(S)||$ + 5;}

$(3.1)$~ {$x = m+k$ and $y = -2m-k$, for $1 \leq k \leq m-1$,
	
	\hspace{1.5cm}  $M$ = $||G^+(S)||$ + 3 + $\left\lfloor \frac{3(m-k+1)}{2}\right\rfloor$;}

$(3.2)$~ {$x = m+k$ and $y = -2m-h$, for $1 \leq k \leq m-1$ and $k < h$,
	
	\hspace{1.5cm} $M$ = $||G^+(S)||$ + 2 + $\left\lfloor \frac{3(m-k+1)}{2}\right\rfloor$;}

$(3.3)$~ {$x = m+k$ and $y = -2m-h$, for $1 \leq h < k \leq m-1$,
	
\hspace{1.5cm}  $M$ = $||G^+(S)||$ + 3 + $\left\lfloor \frac{3(m-k+1)}{2}\right\rfloor$;}

$(4.1)$~ {$x = 2m$ and $y = -3m,$
	
	\hspace{1.5cm} $M$ = $||G^+(S)||$ + 4;}

$(4.2)$~ {$x = 2m$ and $y = -2m-k$, for $1 \leq k \leq m-1$,
	
	\hspace{1.5cm} $M$ = $||G^+(S)||$ + 4 for $1 \leq k \leq m-1$;}

$(4.3)$~ {$x = 2m$ and $y = -2m-k$, for $m < k$,
	
	\hspace{1.5cm} $M$ = $||G^+(S)||$ + 3;}

$(5.1)$~  {$x = 2m+k$ and $y = -2m-h$, for $k-h$ = 0 or $1 \leq k-h \leq m$ or $-m \leq k-h \leq -1$,
	
	\hspace{1.5cm} $M$ = $||G^+(S)||$ + 3; and}

$(5.2)$~ {$x = 2m+k$ and $y = -2m-h$, for $m < k-h$ or $k-h < -m$, 
	
	\hspace{1.5cm}  $M$ = $||G^+(S)||$ + 2.}	
\end{proof}

\begin{theorem} [{\cite[Theorem 3.4]{vm12}}] \label{5.7} \quad {\rm \ 
\begin{enumerate}
	\item [\rm (i)] For $m \in  \mathbb{N}$, $G_{-m,m}$ is a maximal integral sum graph of order $2m+1$.  
	\item [\rm (ii)] For $m \geq 2$, $G_{-(m-1),m}$ is a maximal integral sum graph of order $2m$.
\end{enumerate}		
 } 
\end{theorem}

\begin{proof} 
\begin{enumerate}
	\item [\rm (i)] The result is proved by induction on $m$ for integral sum graphs of odd order $2m+1$. When $m$ = 1, $G_{-1,1}$ is a maximal integral sum graph of order 3 and when $m$ = 2, $G_{-2,2}$ is a maximal integral sum graph of order 5. Hence, the statement is true for $m$ = 1 and $m$ = 2. Assume that the statement is true for some $m \geq 2$. That is, that $G_{-m,m}$ is a maximal integral sum graph of order $2m+1$. We have $||G_{-m,m}||$ = $3m(m+1)/2$ - $\left\lfloor \frac{m}{2}\right\rfloor$, using Theorem \ref{a5}. Consider the case for $m+1$. Let $G$ be a maximal integral sum graph of order $2(m+1)+1$ = $2m+3$ such that $G$ = $G^+(S)$ for some non-empty set $S$. We claim that $S$ = $[-m-1,m+1]$. That is, $G$ = $G_{-(m+1),m+1}$. Using Theorem \ref{a5}, $||G_{-(m+1),m+1}||$ = $3(m+1)(m+2)/2$ - $\left\lfloor \frac{m+1}{2} \right\rfloor$. By our assumption, $||G|| \geq$ $||G_{-(m+1),m+1}||$. Since $G$ is a maximal integral sum graph of order $2(m+1)+1$, $G$ contains $G_{-m,m}$, a maximal integral sum graph of order $2m+1$, as a subgraph. Otherwise, $G$ cannot be a maximal integral sum graph of order $2(m+1)+1$. Thus, $–m,-m+1,\dots,0,1,\dots,m-1$, $m\in S$. Let $S$ = $\{x,y,–m,-m+1,\dots,m-1,m\}$ where $x,y\in\mathbb{Z}$ and $x,y\notin [–m,m]$. The graph $G$ is a maximal integral sum graph of order $2m+3$ when $\{x,y\}$ = $\{-m-1,m+1\}$, using Theorem \ref{5.12}. That is, when $\{x,y\}$ = $\{-m-1, m+1\}$, the graph $G$ = $G^+(S\setminus \{x,y\}) \cup \{-m-1, m+1\}$ = $G_{-(m+1),m+1}$ is a maximal integral sum graph of order $2(m+1)+1$. Hence, the claim is true. Therefore, by the principle of mathematical induction, the statement is true for all positive integers $m$.

	\item [\rm (ii)] The proof is similar to (i). 
\end{enumerate}
\end{proof}

The next statement follows immediately from Theorem \ref{5.7}.

\begin{cor} [{\cite[Corollary 3.5]{vm12}}] \label{5.8} {\rm  For $m \in  \mathbb{N}$, a maximal integral sum graph of order $n$ are 
\\
$(i)$ ~$G_{-2m,2m}$ and $G_{-(2m-1),2m+1}$ when $n$ = $4m+1$; 
\\
$(ii)$ $G_{-2m,2m+1}$ when $n$ = $4m+2$; 
\\
$(iii)$ $G_{-(2m+1),2m+1}$ when $n$ = $4m+3$; and 
\\
$(iv)$ $G_{-(m-1),m}$ when $n$ = $2m$.} \hfill $\Box$
\end{cor} 

In Theorem \ref{5.3}, it is shown that  for $n \geq 4$, $G_{0,n}$ is not a maximal integral sum graph of order $n+1$. Now, consider maximal integral sum graph of order $n+1$.

\begin{theorem} [{\cite[Theorem 3.7]{vm12}}] \label{5.6} \quad {\rm For $0 < r \leq n-r$, $G_{-r,n-r}$ is a maximal integral sum graph of order $n+1$. }
\end{theorem}

\begin{proof} \quad By Corollary \ref{5.8}, $G_{-2m,2m}$ and $G_{-(2m-1),2m+1}$, $G_{-2m,2m+1}$, $G_{-(2m+1),2m+1}$, or $G_{-(2m-1),2m}$  are maximal integral sum graph of order $n+1$ when
$n+1$ = $4m+1$, $4m+2$, $4m+3$ or $4m$, respectively. 

By Corollary \ref{5.1a}, for fixed $n$ = $4m$ and various values for $r$ such that $r \leq n-r$, integral sum graphs $G_{-r, n-r}$, each of order $n+1$ are all non-comparable.  This implies that $G_{-r,n-r}$ cannot be a proper spanning subgraph of an integral sum graph. Hence, the result.
\end{proof}

Figures 5 to 8 show integral sum graphs $G_{0,6,}$, $G_{-1,5,}$, $G_{-2,4}$ and $G_{-3,3}$. 

\medskip
In \cite{tt13}, Tiwari and Tripathi  independently obtained the following result on maximum size of integral sum graph of a given order. 

\begin{theorem} [{\cite[Theorem 3.1]{tt13}}] \quad \label{5.9} {\rm Let $M_n =  ||G^{+}(S)||$ denote  the maximum size for integral sum graph $G^+(S)$ of order $n$. Then, $M_n$ = $\left\lceil \frac{3{(n-1)}^2}{8} \right\rceil + \left\lceil \frac{n-1}{2} \right\rceil$. Moreover, there exists a sum graph of order $n$ and size $m$, for each $m \leq M_n$, except for $m = M_n - 1$ when $n~ \equiv 1 ~( mod ~4)$.}  \hfill $\Box$
\end{theorem}

\subsection{Results on spanning subgraph/super graph of an integral sum graph}

So far we could find out maximal integral sum graphs of the form $G_{-r,s}$, having at least one of it’s vertices is of degree $r+s$, $r,s\in\mathbb{N}$. There remain several open questions.

\begin{oprm} [{\cite[Open Problem 1]{vm12}}] \label{5.12a} \quad {\rm Do there exist any other type of maximal integral sum graph of order $n+1$ that have no vertex of degree n? \hfill $\Box$}
\end{oprm}

\begin{oprm} [{\cite[Open Problem 2]{vm12}}] \label{5.13} \quad {\rm Although no proper spanning super graph of a maximal integral sum graph is an integral sum graph, the converse is not known. That is, are there proper spanning subgraphs of a maximal integral sum graph that are integral sum graphs. \hfill $\Box$}
\end{oprm}

\begin{oprm} [{\cite[Open Problem 3]{vm12}}] \label{5.14} \quad {\rm For which graphs $G$ and $H$, are the graphs $K_1 * (G \cup H)$ and $K_1 * (G * H)$ are integral sum graphs, when $G$ and $H$ are either both sum graphs or are not both sum graphs? \hfill $\Box$}
\end{oprm}

It is easy to see that $K_1 * (K_1 \cup K_3)$, $K_1 * (K_2 \cup K_3)$  and $K_1 * (K_3 \cup K_3)$ are not integral sum graphs, though they are spanning sub-graph of integral sum graphs $G_{-1,3}$, $G_{-2,3}$ and $G_{-3,3}$, respectively. These graphs are presented in Figures 19, 20 and 21.  Also, consider the following problem.  

\begin{prm} [{\cite[Problem 1]{vl26}}] \quad \label{5.15} {\rm Show that $K_1 * (K_3 \cup K_3)$ is a spanning subgraph of $G_{-3,3}$ but not an integral sum graph. }
\end{prm}
\medskip
\noindent
{\bf Solution.} \quad Figures 8, 20 and 21 show graphs $G_{-3,3}$, $K_1 * (K_2 \cup K_3)$ and $K_1 * (K_3 \cup K_3)$. Let the vertices of $G_{-3,3}$ be subscript-labeled and $K_1 * (K_3 \cup K_3)$ = $K_1(w) * (K_3(u_1,u_2,u_3) \cup K_3(v_1,v_2,v_3))$. If $K_1 * (K_3 \cup K_3)$ is an integral sum graph, then the integral sum labeling of vertex $w$ must be 0 since $w$ is the only vertex with degree 6 in the graph of order 7. Define an adjacency preserving bijective mapping $f:$ $V(K_1 * (K_3 \cup K_3))$ $\to$ $V(G_{-3,3})$ such that $f(w)$ = 0, $f(u_1)$ = 1, $f(u_2)$ = 2, $f(u_3)$ = -3, $f(v_1)$ = -1, $f(v_2)$ = -2 and $f(v_3)$ = 3. From the mapping, $f(K_1 * (K_3 \cup K_3))$ is a subgraph of $G_{-3,3}$. And by simple algebraic calculation, we can show that $K_1 * (K_3 \cup K_3)$ is not an integral sum graph.

\medskip	
\begin{center}
	\begin{tikzpicture}
\node (a0) at (10.5,2.5)  [circle,draw,scale=0.6]{$w$};
\node (a1) at (11.5,1.75)  [circle,draw,scale=0.6]{$u_1$};
\node (a2) at (11.5,.5)  [circle,draw,scale=0.6] {$u_2$};
\node (a3) at (10.5,-.5)  [circle,draw,scale=0.6] {$u_3$};
\node (a6) at (9,.5)  [circle,draw,scale=0.6]{$v_1$};

\draw (a0) -- (a1);
\draw (a0) -- (a2);
\draw (a0) -- (a3);
\draw (a0) -- (a6);


\draw (a2) -- (a3);

\draw (a1) -- (a2);
\draw (a1) -- (a3);	

	\node (c6) at (13.5,1.75)  [circle,draw,scale=0.6] {$v_1$};
	\node (c0) at (14.75,2.5)  [circle,draw,scale=0.6]{$w$};
	\node (c1) at (16,1.75)  [circle,draw,scale=0.6]{$u_1$};
	\node (c2) at (16,.5)  [circle,draw,scale=0.6] {$u_2$};
	\node (c3) at (15,-.5)  [circle,draw,scale=0.6] {$u_3$};
	\node (c5) at (13.5,.5)  [circle,draw,scale=0.6]{$v_2$};
	
	\draw (c0) -- (c1);
	\draw (c0) -- (c2);
	\draw (c0) -- (c3);
	\draw (c0) -- (c5);
	\draw (c0) -- (c6);
	
	\draw (c6) -- (c5);
	
	\draw (c2) -- (c3);
	
	\draw (c1) -- (c2);
	\draw (c1) -- (c3);	
	
	\node (d6) at (18,1.75)  [circle,draw,scale=0.6] {$v_1$};
	\node (d0) at (19.5,2.5)  [circle,draw,scale=0.6]{$w$};
	\node (d1) at (21,1.75)  [circle,draw,scale=0.6]{$u_1$};
	\node (d2) at (21,.5)  [circle,draw,scale=0.6] {$u_2$};
	\node (d3) at (19,-.5)  [circle,draw,scale=0.6] {$v_3$};
	\node (d4) at (20,-.5)  [circle,draw,scale=0.6] {$u_3$};
	\node (d5) at (18,.5)  [circle,draw,scale=0.6]{$v_2$};
	
	\draw (d0) -- (d1);
	\draw (d0) -- (d2);
	\draw (d0) -- (d3);
	\draw (d0) -- (d4);
	\draw (d0) -- (d5);
	\draw (d0) -- (d6);
	
	\draw (d6) -- (d3);
	\draw (d6) -- (d5);
	
	\draw (d5) -- (d3);
	
	\draw (d4) -- (d1);
	\draw (d4) -- (d2);
	
	\draw (d1) -- (d2);
	
	\end{tikzpicture}

\vspace{.1cm}	
	Fig. 19. $K_1 * (K_1 \cup K_3)$ \hspace{.4cm}   Fig. 20. $K_1 * (K_2 \cup K_3)$ \hspace{.4cm} Fig. 21. $K_1 * (K_3 \cup K_3)$
	
\end{center}

\medskip
By analysing the above examples, the following  result is obtained.

\begin{theorem} [{\cite[Theorem 53]{vl26}}] \label{5.16} \quad {\rm Any proper spanning subgraph or supergraph of an integral sum graph is not an integral sum graph. }
\end{theorem}

\begin{proof} \quad Let $S$ be a non-empty set of integers and $G^+(S)$ be its integral sum graph. For a given $S$, its integral sum graph $G^+(S)$ is unique. Hence, the result.
\end{proof}

On the other hand, for a given positive integer $n$, there may be more than one integral sum graph of order $n$. For example, $G_n$, $G_{0, n-1}$, $G_{-1, n-2}$, $G_{-2, n-3}$, $\dots$, $G_{-r, n-r-1}$ are integral sum graphs, each of order $n$, $r \leq n$ and $r,n\in\mathbb{N}$. Clearly, different integral sum graphs $G^+(S)$ of order $n$ depend on the existence of number of different sets $S$, each of order $n$. 

\section{Number of cycles of length 3 and 4 in $G_{-m,n}$}

Number of cycles of length 3 and 4 of graphs $G_{-m,n}$ are calculated in \cite{vk11, vr14} and these results are presented in this section, $m\in \mathbb{N}_0$ and $n\in \mathbb{N}$. Throughout this section, we denote the number of distinct sub-graphs $H$ in graph $G$ as $|H|_G$ and the size of graph $G$ as $||G||$ = $|E(G)|$. 

In \cite{ls21} and \cite{vm12}, we have $G_{-m,n}$ = $K_1 * (G_{-m}*G_{n})$, $G_{-m,n}^c$ = $K_1(0) \cup (G_{-m}^c \cup G_{n}^c)$, $\left|E(G_n)\right|$ = $\frac{1}{2}(nC_2 - \left\lfloor \frac{n}{2} \right\rfloor)$, $\left|E(G_{n}^c)\right|$ = $\frac{1}{2}(nC_2 + \left\lfloor \frac{n}{2} \right\rfloor)$, $\left|E(G_{2n})\right|$ = $n^2-n$ = $\left|E(G_{2n-1}^c)\right|$ and $\left|E(G_{2n+1})\right|$ = $n^2$ = $\left|E(G_{2n}^c)\right|$ where $\left\lfloor x \right\rfloor$ denotes the floor of $x$, $m,n \in\mathbb{N}_0$. 

\subsection{Number of edges of $G_{0,n}$, and $G_{-m,n}$}

In this subsection, we present results related to the number of edges of $G_{0,n}$, and $G_{-m,n}$, $m,n \in\mathbb{N}_0$. These results are used in the calculations of cycles of length 3 and 4 in $G_{n}$, $G_{0,n}$, $G_{-m,n}$, $G^c_{n}$, $G^c_{0,n}$, and $G^c_{-m,n}$, $m,n\in\mathbb{N}$.

\begin{theorem} [{\cite[Theorem 2.4]{vm12}}] \quad \label{8.1} {\rm  Let $m\in\mathbb{N}_0$, $n\in \mathbb{N}$ and $m+n \geq 3$. 
		
		\hspace{.5cm} $||G_{-m,n}||$ = $\frac{1}{4}(m^2+n^2+3(m+n)+4mn) -\frac{1}{2}(\lfloor \frac{m}{2}\rfloor + \lfloor \frac{n}{2}\rfloor)$. 
		
In particular, 
\begin{enumerate}
	\item [\rm (i)] $||G_{0,n}||$ = $\frac{n(n+3)}{4} -\frac{1}{2}(\lfloor \frac{n}{2})\rfloor$,  		
		
	\item [\rm (ii)]	$||G_{-n,n}||$ = $\frac{3n(n+1)}{2} -\lfloor \frac{n}{2}\rfloor$ and 
		
	\item [\rm (iii)]	$||G_{-(n-1),n}||$ = $\frac{n(3n-1)}{2}$, $n \in\mathbb{N}.$ \hfill $\Box$
\end{enumerate} }
\end{theorem}

\begin{cor} [{\cite[Corollary 2.5]{vm12}}] \quad \label{8.2} {\rm Let $m\in \mathbb{N}$.
\begin{enumerate}
	\item [\rm (i)] $||G_{-1,3}||$ = $||G_{-2,2}||$.
	
	\item [\rm (ii)] $||G_{-1,4m-1}|| < ||G_{-2m,2m}||$ =  $||(G_{-(2m-1),2m+1}||$ for  $m \geq 2$.
	
	\item [\rm (iii)] $||G_{-1,4m}||$ $<$ $||G_{-2m,2m+1}||$, 
	
	\hspace{.5cm} $||G_{-1,4m+1}|| < ||G_{-2m,2m+2}|| < ||G_{-(2m+1),2m+1}||$.
	
	\item [\rm (iv)]  $G_{-2m,2m} \neq G_{-1,4m-1}$, 
	
	~~$G_{-2m,2m} \neq G_{-(2m-1),2m+1}$, 
	
	~	$G_{-(2m+1),2m+1} \neq G_{-1,4m+1}$, and  
	
	~	$G_{-(2m+1),2m+1} \neq G_{-2m,2m+2}$.  \hfill   $\Box$ 
\end{enumerate}		}		 
\end{cor}

\subsection{Number of cycles of length 3 in $G_{n}$, and $G_{-m,n}$}

In \cite{vk11}, number of cycles of length 3 of graphs $G_{2k}$, $G_{2k+1}$, $G_{2k}^c$, $G_{2k+1}^c$, $G_{0,n}$ and $G_{-m,n}$ are calculated and presented here. 

\begin{theorem} [{\cite[Theorem 2.1]{vk11}}] \quad \label{8.3}  {\rm Let $n \geq 3$. 
		
		$\left|C_3\right|_{G_{n}}$ = $\left|C_3\right|_{G_{n-2}} + ||G_{n-2}||$ and $\left|C_3\right|_{G_{n}^c}$ = $\left|C_3\right|_{G_{n-2}^c} + ||G_{n-2}^c||$.   }   
\end{theorem}

\begin{proof}\quad Let $V(G_n)$ = $\{u_1,u_2,\dots,u_n\}$ = $V(G^c_n)$ where $u_j$ is the vertex with sum labeling $j$ in $G_n$ and anti-sum labeling $j$ in $G^c_n$, $1 \leq j \leq n$. Then, using Theorem \ref{2.11}, $G_n \setminus \{u_1,u_n\}$ and unlabeled graph $G_{n-2}$ are isomorphic and $G^c_n \setminus \{u_1,u_n\}$ and unlabeled graph $G^c_{n-2}$ are isomorphic. Hence, the number of cycles of length 3 in $G_n$ = number of cycles of length 3 in $G_{n-2}$ + number of cycles of length 3, each cycle of length 3 with $u_1$ as a vertex in $G_n$ and number of cycles of length 3 in $G^c_n$ = number of cycles of length 3 in $G^c_{n-2}$ + number of cycles of length 3, each cycle of length 3 with $u_n$ as a vertex in $G^c_n$, $n \geq 3$. In $G_n$, vertices $u_1$ and $u_n$ are non-adjacent and end points of each edge of subgraph $G_n \setminus \{u_1,u_n\}$ are adjacent to $u_1$ and non-adjacent to $u_n$ whereas in $G^c_n$, vertices $u_1$ and $u_n$ are adjacent and end points of each edge of subgraph $G^c_n \setminus \{u_1,u_n\}$ are adjacent to $u_n$ and non-adjacent to $u_1$. Thus, each edge of subgraph $G_n \setminus \{u_1,u_n\}$ makes one cycle of length 3 with a pair of edges at $u_1$ in $G_n$ and each edge of subgraph $G^c_n \setminus \{u_1,u_n\}$ makes one cycle of length 3 with a pair of edges at $u_n$ in $G^c_n$. This implies, $\left|C_3\right|_{G_{n}}$ = $\left|C_3\right|_{G_{n-2}} + ||G_{n-2}||$ and $\left|C_3\right|_{G^c_{n}}$ = $\left|C_3\right|_{G^c_{n-2}} + ||G^c_{n-2}||$, $n \geq 3$.
\end{proof}

\begin{cor} [{\cite[Corollary 2.2]{vk11}}] \label{8.4} \quad {\rm Let $m,n \in\mathbb{N}$. 
		
		$\left|C_3\right|_{G_{2n+2}}$ = $\frac{(n-1)n(n+1)}{3}$ = $\left|C_3\right|_{G_{2n+1}^c}$ and $\left|C_3\right|_{G_{2n+3}}$ = $\frac{n(n+1)(2n+1)}{6}$ = $\left|C_3\right|_{G_{2n+2}^c}$.  }
\end{cor}

\begin{proof}\quad Applying Theorem \ref{8.3} successively on $G_{2n+2}$, $G_{2n+3}$, $G_{2n+2}^c$ and $G_{2n+1}^c$, we obtain, 

$\left|C_3\right|_{G_{2n+2}}$ 	= $|E(G_{2n})|$ + $|E(G_{2n-2})|$ + $\dots$ + $|E(G_{6})|$ + $|E(G_{4})|$ + $|C_3|_{G_4}$
	 
\hspace{1.4cm} = $(n^2-n) + ({(n-1)}^2-(n-1)) + \dots + (3^2-3) + (2^2-2) + 0$ = $\frac{(n-1)n(n+1)}{3}$;
	
$\left|C_3\right|_{G_{2n+3}}$ 	= $|E(G_{2n+1})|$ + $|E(G_{2n-1})|$ + $\dots$ + $|E(G_{5})|$ + $|E(G_{3})|$ + $|C_3|_{G_3}$

\hspace{1.4cm} = $n^2 + {(n-1)}^2 + \dots + 2^2 + 1^2 + 0$ = $\frac{n(n+1)(2n+1)}{6}$;

\vspace{.2cm} 
$\left|C_3\right|_{G_{2n+2}^c}$ = $|E(G_{2n}^c)|$ + $|E(G_{2n-2}^c)|$ + $\dots$ + $|E(G_{4}^c)|$ + $|E(G_{2}^c)|$ + $|C_3|_{G_2^c}$

\hspace{1.4cm} = $n^2 + {(n-1)}^2 + \dots + 2^2 + 1^2 + 0$  
= $\frac{n(n+1)(2n+1)}{6}$;

\vspace{.2cm}
$\left|C_3\right|_{G_{2n+1}^c}$ = $|E(G_{2n-1}^c)|$ + $|E(G_{2n-3}^c)|$ + $\dots$ + $|E(G_{5}^c)|$ + $|E(G_{3}^c)|$ + $|C_3|_{G_3^c}$

\hspace{1.4cm} = $(n^2-n) + ({(n-1)}^2-(n-1)) + \dots + (3^2-3) + (2^2-2) + 0$ = $\frac{(n-1)n(n+1)}{3}$.		
\end{proof}

\begin{theorem} [{\cite[Theorem 2.3]{vk11}}] \quad \label{8.5}  {\rm Let $m\in \mathbb{N}_0$ and $n\in \mathbb{N}$. 
		
	~(i)	$\left|C_3\right|_{G_{-m, n}}$ = $\left|C_3\right|_{G_{m}}$ + $\left|C_3\right|_{G_{n}}$ + $(n+1)||G_{m}||$ + $(m+1)||G_{n}||$ + $mn$ and 
		
	(ii)	$\left|C_3\right|_{G_{-m, n}^c}$ = $\left|C_3\right|_{G_{m}^c}$ + $\left|C_3\right|_{G_{n}^c}$.  }   
\end{theorem}

\begin{proof}\quad For $m\in \mathbb{N}_0$ and $n\in \mathbb{N}$, using Theorem \ref{1.2.4}, we get, 
	
	$G_{-m, n}$ = $K_1 * (G_{-m} * G_{n})$ = $K_1 * (G_{-m} \cup G_{n} \cup K_{m,n})$ and $G^c_{-m, n}$ = $K_1(0) \cup (G^c_{-m} \cup G^c_{n})$ where vertices of $K_{m,n}$ are vertices of $G_{-m}$ and $G_{n}$. Here, $K_1$ is the vertex with integral sum label 0 and adjacent to all other vertices in $G_{-m, n}$ and in $G^c_{-m, n}$, it is an isolated vertex. For $m,n \geq 5$, each graph $G_{-m}$ and $G_{n}$ contains cycle(s) of length 3 but the bipartite graph $K_{m,n}$ is a triangle free graph and for $m,n \geq 4$, each graph $G^c_m$ and $G^c_n$ contains cycle(s) of length 3. 
	
	Clearly, $\left|C_3\right|_{G^c_{-m, n}}$ = $\left|C_3\right|_{G^c_{m}}$ + $\left|C_3\right|_{G^c_{n}}$ since $G^c_{-m}$ and $G^c_n$ are disjoint subgraphs in $G^c_{-m,n}$ and also  $\left|C_3\right|_{G^c_{-m}}$ = $\left|C_3\right|_{G^c_{m}}$. 
	
	Now, let us consider all cycles of length 3 in the graph $G_{-m,n}$. We classify all triangles in $G_{-m,n}$ into the following 7 types. 
	\begin{enumerate}
		\item [\rm (i)] Triangles, each of its three edges are edges of $G_{-m}$. 
	
		\item [\rm (ii)] Triangles, each of its three edges are edges of $G_n$. 
	
		\item [\rm (iii)] Triangles, each with an edge of $K_{m,n}$ and the other two edges a pair of edges at $K_1$. 
	
		\item [\rm (iv)] Triangles, each with an edge of $G_{-m}$ and the other two edges a pair of edges at $K_1$. 
	
		\item [\rm (v)] Triangles, each with an edge of $G_n$ and the other two edges a pair of edges at $K_1$. 
	
		\item [\rm (vi)] Triangles, each with a pair of edges at a vertex of $G_n$ and the third, an edge of $G_{-m}$. 
	
		\item [\rm (vii)] Triangles, each with a pair of edges at a vertex of $G_{-m}$ and the third, an edge of $G_n$. 
	\end{enumerate}
	From the above cases, we get, 
	
	 $\left|C_3\right|_{G_{-m, n}}$ = $\left|C_3\right|_{G_{-m}}$ + $\left|C_3\right|_{G_{n}}$ + $|E(K_{m, n})|$ (since every edge of $K_{m,n}$ makes a triangle  
	 
	 \hfill with a pair of edges at $K_1$) 
	 
	\hspace{1.5cm} +  $|E(G_{-m})|$ (since each edge of $G_{-m}$ makes a triangle with a pair of edges at $K_1$) 
	 
	\hspace{1.5cm}  +  $|E(G_{n})|$  (since each edge of $G_n$ makes a triangle with a pair of edges at $K_1$) 
	 
	\hspace{1.5cm}  + $n |E(G_{-m})|$ (since each edge of $G_{-m}$ makes a triangle with a pair 
	
	\hfill of edges at every vertex of $G_n$) 
	 
	\hspace{1.5cm}  + $m |E(G_{n})|$ (since each edge of $G_n$ makes a triangle with a pair of
	
	\hfill  edges at every vertex of $G_{-m}$) 
	 
	\hspace{1.5cm} 	 = $\left|C_3\right|_{G_{m}}$ + $\left|C_3\right|_{G_{n}}$ + $(n+1)||G_{m}||$ + $(m+1)||G_{n}||$ + $mn$. 
	
	Hence the result.                        
\end{proof}

\begin{cor} [{\cite[Corollary 2.4]{vk11}}] \label{8.6} \quad {\rm Let $m,n \in\mathbb{N}$. 
	\begin{enumerate}
	\item [\rm (i)] $\left|C_3\right|_{G_{-2m,2n}}$ = $\frac{1}{3}(m+n)(m^2+5mn+n^2-1)$. 

	\item [\rm (ii)] $\left|C_3\right|_{G_{-2m,2n+1}}$ = $\frac{1}{6}(2(m^3+n^3)+12mn(m+n)+3(2m^2+n^2+4mn) + 4m+n)$. 

	\item [\rm (iii)] $\left|C_3\right|_{G_{-(2m+1),2n}}$ = $\frac{1}{6}(2(m^3+n^3)+12mn(m+n)+3(m^2+2n^2+4mn) + m+4n)$.

	\item [\rm (iv)] $\left|C_3\right|_{G_{-(2m+1),2n+1}}$ = $\frac{1}{6}(m+n)(2{(m+n)}^2+9(m+n)+6mn+13) + mn+1$. 

	\item [\rm (v)] $\left|C_3\right|_{G_{-2m,2n}^c}$ = $\frac{(m-1)m(2m-1)}{6}$ + $\frac{(n-1)n(2n-1)}{6}$. 

	\item [\rm (vi)] $\left|C_3\right|_{G_{-2m,2n+1}^c}$ = $\frac{(m-1)m(2m-1)}{6}$ + $\frac{(n-1)n(n+1)}{3}$. 

	\item [\rm (vii)] $\left|C_3\right|_{G_{-(2m+1),2n}^c}$ = $\frac{(m-1)m(m+1)}{3}$ + $\frac{(n-1)n(2n-1)}{6}$. 

	\item [\rm (viii)] $\left|C_3\right|_{G_{-(2m+1),2n+1}^c}$ = $\frac{(m-1)m(m+1)}{3}$ + $\frac{(n-1)n(n+1)}{3}$.  
	\end{enumerate} }
\end{cor}

\begin{proof} 	Using Corollary \ref{8.4}, $\left|C_3\right|_{G_{2n+2}}$ = $\frac{(n-1)n(n+1)}{3}$ = $|C_3|_{G^c_{2n+1}}$ and $\left|C_3\right|_{G_{2n+1}}$ = $\frac{(n-1)n(2n-1)}{6}$ = $|C_3|_{G^c_{2n}}$. Also, we have, $|E(G_{2m})|$ = $(m-1)m$ = $|E(G^c_{2m-1})|$ and $|E(G_{2m+1})|$ = $m^2$ = $|E(G^c_{2m})|$, $m,n\in\mathbb{N}$. Using the above results and Theorem \ref{8.5}, we get,
	\begin{enumerate}
		\item [\rm (i)] $\left|C_3\right|_{G_{-2m,2n}}$ = $\frac{(m-2)(m-1)m}{3}$ + $\frac{(n-2)(n-1)n}{3}$ + $(2n+1)(m-1)m$ + $(2m+1)(n-1)n$ + $4mn$
		
	\hspace{1.4cm}	= $\frac{1}{3}(m+n)(m^2+5mn+n^2-1)$. 
		
		\item [\rm (ii)] $\left|C_3\right|_{G_{-2m,2n+1}}$ 	
		= $\frac{(m-2)(m-1)m}{3}$ + $\frac{(n-1)n(2n-1)}{6}$ + $2(n+1)(m-1)m$ + $(2m+1)n^2$ + $2m(2n+1)$
		
\hspace{1.5cm} = $\frac{1}{6}(2(m^3+n^3)+12mn(m+n) +3(2m^2+4mn+n^2) + 4m+n)$. 
		
		\item [\rm (iii)] $\left|C_3\right|_{G_{-(2m+1),2n}}$ = $\frac{1}{6}(2(m^3+n^3)+12mn(m+n)+3(m^2+4mn+2n^2) + m+4n)$.
		
		\item [\rm (iv)] $\left|C_3\right|_{G_{-(2m+1),2n+1}}$ 
			= $\frac{1}{6}(m-1)m(2m-1)$ + $\frac{1}{6}(n-1)n(2n-1)$ 
			
			\hfill + $2(n+1)m^2$ + $2(m+1)n^2$ + $(2m+1)(2n+1)$
		
	\hspace{1.4cm}		= $\frac{1}{6}(2(m^3+n^3)$ + $12mn(m+n)$ + $9(m^2+2mn+n^2)$ + $6mn$ + $13(m+n)+6)$
		
	\hspace{1.4cm}	= $\frac{1}{6}(m+n)(2(m^2+n^2-mn)$ + $12mn$ + $9(m+n)+13)$ + $mn+1$
			
\hspace{1.4cm} 	= $\frac{1}{6}(m+n)(2{(m+n)}^2$ + $6mn$ + $9(m+n)+13)$ + $mn+1$.
		\end{enumerate} 
		Results (v) to (viii) follow directly from formulae. Hence the result.
\end{proof}

Any property of numbers is interesting and important. From the above study, the following simple properties of natural numbers are obtained.

\begin{theorem} [{\cite[Theorems 2.5]{vk11}}] \quad \label{8.7}  {\rm Let $m,n\in \mathbb{N}$. 
	\begin{enumerate}
		\item [\rm (i)] $(m+n)(m^2+5mn+n^2-1)$ is divisible by 3. 
		
	\item [\rm (ii)] $2(m^3+n^3)$ + $3m^2+4m+n$ and $(m+n)(2{(m+n)}^2 + 9(m+n)+13)$ are divisible by 6.    
	\end{enumerate} }
\end{theorem}
\begin{proof}\quad The result follows from Corollary \ref{8.6}. 
\end{proof}

\subsection{Number of cycles of length 4 in $G_{n}$,  and $G_{-m,n}$}

The number of cycles of length 4 in $G_{n}$, $G_{0,n}$ and $G_{-m,n}$ are obtained in \cite{vr14} and these results are presented in this subsection. At first, the number of cycles of length 4 in $G_{2k}$, $G_{2k+1}$, $G_{2k}^c$ and $G_{2k+1}^c$ are obtained and then, using these results, the number of cycles of length 4 in $G_{0,n}$, $G_{0,n}^c$, $G_{-m,n}$ and $G_{-m,n}^c$ are derived, $k \geq 2$ and $m,n \in\mathbb{N}$. 

\begin{theorem} [{\cite[Theorems 2.1]{vr14}}] \label{8.8} \quad {\rm Let $n \geq 2$. 
\begin{enumerate}
	\item [\rm (i)] $\left|C_4\right|_{G_{2n+2}}$ = $\left|C_4\right|_{G_{2n}}$ + $\frac{(n-1)n(7n-11)}{6}$ = $\frac{(n-1)n(n+1)(7n-10)}{24}$.

	\item [\rm (ii)]  $\left|C_4\right|_{G_{2n+2}^c}$ = $\left|C_4\right|_{G_{2n}^c}$ + $\frac{(n-1)n(7n+1)}{6}$ = $\frac{(n-1)n(n+1)(7n+6)}{24}$.  
\end{enumerate} }
\end{theorem}

\begin{proof}\quad (i) Let $V(G_{2n+2})$ = $\{u_1,u_2,\dots,u_{2n+2}\}$ = $V(G^c_{2n+2})$ where $u_i$ is the vertex with sum labeling $i$ in $G_{2n+2}$ and anti-sum labeling $i$ in $G^c_{2n+2}$, $1 \leq i \leq 2n+2$ and $n\in\mathbb{N}$. $\{u_1,u_2,\dots,u_{n+1}\}$ is a clique and $\{u_{n+2},u_{n+3},\dots,u_{2n+2}\}$ is a stable set to $G_{2n+2}$.  Using Theorem \ref{2.11}, graph $G_{2n+2}\setminus\{u_1,u_{2n+2}\}$ is isomorphic to $G_{2n}$, without the vertex labels. In $G_{2n+2}$, $u_1$ is adjacent to $u_2,u_3,\dots,u_{2n+1};$ $u_{2n+2}$ is an isolated vertex and $u_{2n+1}$ is a pendant vertex. Therefore, $\left|C_4\right|_{G_{2n+2}}$ = $\left|C_4\right|_{G_{2n}}$ $+$ number of cycles of length four, each with $u_1$ as a vertex in $G_{2n+2}$. Also, none of $u_{2n+1}$ and $u_{2n+2}$ is a vertex of any cycle of length $4$ in $G_{2n+2}.$
	
Let $(u_1 , u_i , u_j , u_k)$ be any cycle of length $4$ (with $u_1$ as a vertex) in $G_{2n+2}$, $1 < i,j,k < 2n+1$ and $i,j,k$ are all different. Under the above conditions, the following three types of $C_4s$ arise in $G_{2n+2}$. Type-1: $u_i,u_j,u_k\in \{u_2,u_3,\dots,u_{n+1}\}$, Type-2: $u_i,u_j\in \{u_2,u_3,\dots,u_{n+1}\}$ and $u_k\in \{u_{n+2}$, $u_{n+3}$, $\dots,u_{2n}\}$, and Type-3: $u_i\in \{u_2,u_3,\dots,u_{n+1}\}$ and $u_j,u_k\in \{u_{n+2}$, $u_{n+3}$, $\dots,u_{2n}\}$. Now, let us calculate the number of $C_4s$ in $G_{2n+2}$ in each type. 

\medskip
   \textbf{Number of $C_4s$ under Type-1:} Here , $u_i,u_j,u_k\in \{u_2,u_3,\dots,u_{n+1}\}$ in $G_{2n+2}$. Number of ways of selecting $3$ vertices $u_i,u_j,u_k$ out of $u_2,u_3,\dots,u_{n+1}$ is $\binom{n}{3}$. There are $3$ different $C_4s$ with $u_1, u_i , u_j , u_k$ as vertices under type-1, namely, $(u_1$, $u_i$, $u_j$, $u_k)$, $(u_1,u_i,u_k,u_j)$ and $(u_1,u_j,u_i,u_k)$. Therefore, total number of $C_4s$ of type-1 in $G_{2n+2}$ = $3 \times \binom{n}{3}$ =  $n(n-1)(n-2)/2$.

\medskip
\textbf{Number of $C_4s$ under Type-2:} Here, $u_i,u_j\in \{u_2,u_3,\dots,u_{n+1}\}$ and $u_k\in \{u_{n+2}$, $u_{n+3}$, $\dots,u_{2n}\}$. Consider all possible cycles, each of length $4$ and with vertices $u_1, u_i , u_j$ and $u_k$ in $G_{2n+2}.$
	
	When $k$ = $2n$, $u_k$ = $u_{2n}$ is adjacent to $u_1$ and $u_2$ only. And under this case, $u_2$ = $u_i$ or $u_2$ = $u_j$. W.l.g., assume $u_2 = u_i$. This implies, $2 = i < 3 \leq j \leq n+1$. And any $C_4$ under this case is of the form $(u_1 , u_{k} , u_i , u_j)$ = $(u_1 , u_{2n} , u_2 , u_j)$, $u_j\in\{u_3,u_4,\dots,u_{n+1}\}$ and number of such $C_4s$ is $\left|\{u_3,u_4,\dots,u_{n+1}\}\right|$ = $n-1$. 
	
	When $k = 2n-1$, $u_k = u_{2n-1}$ is adjacent to $u_1, u_2$ and $u_3$ only and thereby $d(u_k)$ = 3 = $2n+2-(2n-1)$. And any $C_4$ of type-2 is of the form $(u_1, u_{2n-1}, u_2,u_x)$ or $(u_1, u_{2n-1}, u_3, u_y)$ where $u_x\in \{u_3, u_4,\dots,u_{n+1}\}$ and $u_y\in \{u_2, u_4, u_5,\dots,u_{n+1}\}$. Number of such $C_4s$ is $2(n-1)$. 
	
	When $k = 2n-2$, $u_k = u_{2n-2}$ is adjacent to $u_1, u_2 , u_3$ and $u_4$ only and thereby $d(u_k)$ = $4$ = $2n+2-(2n-2)$. Therefore, number of such $C_4s$ is $(4-1)(n-1)$ = $3(n-1)$. 
	
	In general, when $k = 2n+2-x$ and $2 \leq x \leq n$, $u_k = u_{2n+2-x}$ is adjacent to $u_1, u_2, \dots, u_x$ only and thereby $d(u_k)$ = $d(u_{2n+2-x})$ = $x$. And number of $C_4$ of the form $(u_1, u_{2n+2-x}, u_i, u_j)$ is $(x-1)(n-1)$ where $u_i\in \{u_2,u_3,\dots,u_x\}$ and $u_j\in \{u_2,u_3,\dots,u_{n+1}\}\setminus \{u_i\}$. 
	
	$\therefore$ Total number of $C_4s$ of type-2 in $G_{2n+2}$ is $$\sum_{x=2}^n{(x-1)(n-1)} = (n-1)\big(\sum_{x=1}^{n-1}{x}\big) = \frac{n(n-1)^2}{2}.$$ 

\textbf{Number of $C_4s$ under Type-3:} In this type, $u_i\in \{u_2,u_3,\dots,u_{n+1}\}$ and $u_j,u_k\in \{u_{n+2},u_{n+3}$, $\dots, u_{2n}\}$ in $G_{2n+2}$, $j \neq k$. Here, $u_j$ and $u_k$ are adjacent to $u_1$ for every $j,k\in \{n+2,n+3,\dots,2n+1 \}$ in $G_{2n+2}$, $j \neq k$. W.l.g., assume, $j < k$. If $u_j$ and $u_k$ are adjacent to $u_i$, then $j+i \leq 2n+2$ and $k+i \leq 2n+2$ which implies, $j+i < k+i \leq 2n+2.$
	
	For $1 \leq x \leq n$, $u_{n+1+x}$ is adjacent to $u_1,u_2,\dots$, $u_{n+1-x}$ in $G_{2n+2}$ and hence $d(u_{n+1+x})$ = $n+1-x$. In $G_{2n+2}$, $u_{n+1}$ is non-adjacent to $u_{n+2}$ and $u_{2n+1}$ is a pendant vertex and hence neither $u_{n+1}$ nor $u_{2n+1}$ is a vertex of any $C_4$ of type-3 in $G_{2n+2}$. 
	
	When $u_{k}$ = $u_{2n+2-x}$ and $2 \leq x \leq n-1$, different possibilities of $u_i$ in $C_4s$ of type-3 in $G_{2n+2}$ are $u_2,u_3,\dots,u_x$. And corresponding to each pair of $u_i$ and $u_k$, different possible $u_js$ are $u_{k-1}$, $u_{k-2}$, $\dots$, $u_{n+2}$ in $G_{2n+2}$. Therefore, number of $C_4s$ of type-3 in $G_{2n+2}$ with $u_k$ = $u_{2n+2-x}$ is $(x-1)(k-1-(n+1))$ = $(x-1)(n-x)$. Hence, total number of $C_4s$ of type-3 in $G_{2n+2}$ is	
$$\sum_{x=2}^{n-1}{(n-x)(x-1)} = \sum_{x=1}^{n-2}{(n-1-x)x} = (n-1)\big(\sum_{x=1}^{n-2}{x}\big) - \sum_{x=1}^{n-2}{x^2} = \frac{n(n-1)(n-2)}{6}.$$	
	When $u_i,u_j,u_k\in \{u_{n+2},u_{n+3},\dots,u_{2n}\}$, cycle $C_4$ of the form $(u_1 , u_i , u_j , u_k)$ doesn't exist in $G_{2n+2}$ since  $\{u_{n+2},u_{n+3},\dots$, $u_{2n+2}\}$ is a stable set to split graph $G_{2n+2}$. 
	
\medskip
	Adding all $C_4s$ in the three types, we obtain, total number of $C_4s$ in $G_{2n+2}$ with $u_1$ as a vertex is 
		$$\frac{n(n-1)(n-2)}{2} + \frac{n(n-1)^2}{2} + \frac{n(n-1)(n-2)}{6} = \frac{n(n-1)(7n-11)}{6},~ n \geq 2.$$

\medskip
 Therefore, for $n \geq 2$,  
 
 $\left|C_4\right|_{G_{2n+2}}$ = $\left|C_4\right|_{G_{2n}}$ + $(7n^3-18n^2+11n)/6$ 

\hspace{1.4cm}	= $\frac{1}{6}((7n^3-18n^2+11n)$ + $(7{(n-1)}^3 - 18{(n-1)}^2+11(n-1))$ + $\left|C_4\right|_{G_{2n-2}}$ 
	
\hspace{1.4cm}	= $\frac{1}{6}((7n^3-18n^2+11n)$ + $(7{(n-1)}^3 - 18{(n-1)}^2$ + $11(n-1))$

\hfill  + $\dots$ + $(7\times 2^3 -18\times 2^2 + 11\times 2))$ + $\left|C_4\right|_{G_4}$ 
	
\hspace{1.4cm}	= $\frac{1}{6}((7n^3-18n^2+11n)$ + $(7{(n-1)}^3 - 18{(n-1)}^2$ + $11(n-1))$  

\hfill + $\dots$ + $(7\times 2^3 -18\times 2^2 + 11\times 2))$ + 0  
	
\hspace{1.4cm}	= $(n-1)n(n+1)(7n-10)/24$. 

\medskip	
	Now, let us prove the result on $G^c_{2n+2}$. Consider, graph $G^c_{2n+2}$, $n \in\mathbb{N}$. $\{u_1$, $u_2$, $\dots$, $u_{n}\}$ is a stable set and $\{u_{n+1},u_{n+2},\dots,u_{2n+2}\}$ is a clique to split graph $G^c_{2n+2}$.  Using Theorem \ref{2.11}, graph $G^c_{2n+2} \setminus \{u_1,u_{2n+2}\}$ is isomorphic to $G^c_{2n}$, without the vertex labels. In $G^c_{2n+2}$, $u_{2n+2}$ is adjacent to $u_1,u_2,\dots,u_{2n+1}$ and $u_{1}$ is a pendant vertex. Hence, $u_1$ is not a vertex in any cycle of length 4 in $G^c_{2n+2}$. Therefore, $\left|C_4\right|_{G^c_{2n+2}}$ = $\left|C_4\right|_{G^c_{2n}}$ + number of cycles of length four, each with $u_{2n+2}$ as a vertex in $G^c_{2n+2}$. 
	
	Let $(u_{2n+2} , u_k , u_j , u_i)$ be any cycle of length $4$ in $G^c_{2n+2}$, $2 \leq i,j,k \leq 2n+1$ and $i,j,k$ are all distinct. Under the above conditions, the following three types of $C_4s$ arise in $G^c_{2n+2}$. Type-1: $u_i,u_j,u_k\in \{u_{n+1},u_{n+2},\dots$, $u_{2n+1}\}$, Type-2: $u_j,u_k\in \{u_{n+1}$, $u_{n+2}$, $\dots$, $u_{2n+1}\}$ and $u_i\in \{u_{2},u_{3},\dots,u_{n}\}$, and Type-3: $u_k\in \{u_{n+1}$, $u_{n+2}$, $\dots$, $u_{2n+1}\}$ and $u_i,u_j\in \{u_{2},u_{3},\dots$, $u_{n}\}$. Now, let us calculate number of $C_4s$ in $G^c_{2n+2}$ in each type. W.l.g. assume that $i < j < k$. 

\medskip
\textbf{Number of $C_4s$ under Type-1:} Here, $u_i,u_j,u_k\in \{u_{n+1},u_{n+2},\dots$, $u_{2n+1}\}$ in $G^c_{2n+2}$. Number of ways of selecting $3$ vertices $u_i,u_j,u_k$ out of $u_{n+1},u_{n+2},\dots$, $u_{2n+1}$ is $\binom{n+1}{3}$. There are $3$ different $C_4s$ in $G^c_{2n+2}$ with $u_{2n+2}, u_i , u_j , u_k$ as vertices under type-1, namely, $(u_{2n+2},u_k,u_j,u_i)$, $(u_{2n+2},u_k,u_i,u_j)$ and $(u_{2n+2}$, $u_j$, $u_k$, $u_i)$. Hence, total number of $C_4s$ of type-1 in $G^c_{2n+2}$ is $3\times \binom{n+1}{3}$ =  $(n+1)n(n-1)/2$. 

\medskip
\textbf{Number of $C_4s$ under Type-2:} Here, $u_k,u_j\in \{u_{n+1},u_{n+2},\dots,u_{2n+1}\}$ and $u_i\in \{u_{2},u_{3},\dots,u_{n}\}$. Consider all possible cycles, each of length $4$ and with the vertices $u_{2n+2}, u_i , u_j$ and $u_k$ in $G^c_{2n+2}$.
	
	When $i$ = 2, $u_i$ = $u_{2}$ is adjacent to $u_{2n+2}$ and $u_{2n+1}$ only. And under this case, $d(u_i)$ = 2, $u_k$ = $u_{2n+1}$ and $u_j$ = $u_{2n}, u_{2n-1}, \dots, u_{n+1}$. Number of such $C_4s$ is $\left|\{u_{2n}, u_{2n-1}, \dots, u_{n+1}\}\right|$ = $n$. 
	
	When $i$ = 3, $u_i$ = $u_{3}$ is adjacent to $u_{2n+2}$, $u_{2n+1}$ and $u_{2n}$ only and thereby $d(u_i)$ = 3. And any $C_4$ of type-2 is of the form $(u_{2n+2}$, $u_{3}$, $u_{2n+1}$, $u_x)$ or $(u_{2n+2}$, $u_{3}$, $u_{2n}$, $u_y)$ where $u_x\in \{u_{2n}$, $u_{2n-1}$, $\dots$, $u_{n+1}\}$ and $u_y\in \{u_{2n+1}$, $u_{2n-1}$, $u_{2n-2}$, $\dots$, $u_{n+1}\}$. Number of such $C_4s$ is $2n$. 
	
	In general, when $i$ = $x$ and $2 \leq x \leq n$, $u_i$ = $u_{x}$ is adjacent to $u_{2n+2}$, $u_{2n+1}$, $\dots$, $u_{2n+2-(x-1)}$ only and thereby $d(u_i)$ = $x$ and number of $C_4s$ of the form $(u_{2n+2},u_i,u_y,u_z)$ is $(x-1)n$ where $u_y\in \{u_{2n+1}$, $u_{2n}$, $\dots$, $u_{n+1}\}$ and $u_z\in \{u_{2n+1}$, $u_{2n}$, $\dots$, $u_{n+1}\}\setminus \{u_y\}$.

\medskip	
	$\therefore$ Total number of $C_4s$ of type-2 in $G^c_{2n+2}$ is 	$$\sum_{x=2}^n{(x-1)n} = n\big(\sum_{x=1}^{n-1}{x}\big) = \frac{(n-1)n^2}{2}.$$  

	\textbf{Number of $C_4s$ under Type-3:} Here, $u_k\in \{u_{n+1},u_{n+2},\dots,u_{2n+1}\}$ and $u_i,u_j\in \{u_{2},u_{3},\dots, u_{n}\}$, $i \neq j$. Consider all possible cycles, each of length 4 and with the vertices $u_{2n+2}$, $u_k$, $u_j$ and $u_i$ in $\overline{G_{2n+2}}$. For a given $i$, $2 \leq i \leq n-1$, $j$ takes values $i+1,i+2,\dots,n$ and possible values of $k$ are $2n+2-1,2n+2-2,\dots,2n+2-(i-1)$.
	
	$\therefore$ Total number of $C_4s$ of type-3 in $G^c_{2n+2}$ is 
	$$\sum_{i=2}^{n-1}{(n-i)(i-1)} = \sum_{i=1}^{n-2}{i(n-1-i)} = (n-1)\big(\sum_{i=1}^{n-2}{i}\big) - \sum_{i=1}^{n-2}{i^2} = n(n-1)(n-2)/6.$$	
	When $u_i,u_j,u_k\in \{u_{2},u_{3},\dots,u_{n}\}$, cycle $C_4$ of the form $(u_{2n+2}, u_k , u_j, u_i)$ doesn't exist in $G^c_{2n+2}$ since  $\{u_{2},u_{3}$, $\dots$, $u_{n}\}$ is a stable set to split graph $G^c_{2n+2}$. 

\medskip
	Adding all $C_4s$ in the three types, we obtain: 

\medskip
	Total number of $C_4s$ with $u_{2n+2}$ as a vertex in $G^c_{2n+2}$ is 
	$$\frac{(n-1)n(n+1)}{2} + \frac{(n-1)n^2}{2} + \frac{(n-2)(n-1)n}{6} = \frac{(n-1)n(7n+1)}{6}, ~n \geq 2.$$ 
	Therefore, for $n \geq 2$, 
	
	$\left|C_4\right|_{G^c_{2n+2}}$ = $\left|C_4\right|_{G^c_{2n}}$ + $\frac{7n^3-6n^2-n}{6}$ 
	
\hspace{1.3cm}	= $\frac{1}{6}((7n^3-6n^2-n)$ + $(7{(n-1)}^3 - 6{(n-1)}^2-(n-1)))$ + $\left|C_4\right|_{G^c_{2n-2}}$ 
	
\hspace{1.3cm}	= $\frac{1}{6}((7n^3-6n^2-n)$ + $(7{(n-1)}^3 - 6{(n-1)}^2 -(n-1))$ + $\dots$ 

\hfill + $(7\times 2^3 -6\times 2^2 -2))$ + $\left|C_4\right|_{G^c_{4}}$ 
	
\hspace{1.3cm}	= $\frac{1}{6}((7n^3-6n^2-n)$ + $(7{(n-1)}^3 - 6{(n-1)}^2 -(n-1))$ + $\dots$ + $(7\times 1^3 -6\times 1^2 -1))$ 
	
\hspace{1.3cm}	= $\frac{(n-1)n(n+1)(7n+6)}{24}$.  
	
	Hence the result is true in this case. 

\medskip
\noindent
(ii) Proof is similar to case (i).  
\end{proof}

\begin{theorem} [{\cite[Theorems 2.2]{vr14}}] \label{8.9} \quad {\rm Let $n \geq 2$. 
\begin{enumerate}
	\item [\rm (i)] $\left|C_4\right|_{G_{2n+3}}$ = $\left|C_4\right|_{G_{2n+1}}$ + $\frac{(n-1)n(7n+1)}{6}$ = $\frac{1}{24}(n-1)n(n+1)(7n+6)$ = $\left|C_4\right|_{G_{2n+2}^c}$. 

	\item [\rm (ii)] $\left|C_4\right|_{G_{2n+1}^c}$ = $\left|C_4\right|_{G_{2n-1}^c}$ + $\frac{(n-1)n(7n-11)}{6}$ = $\frac{(n-1)n(n+1)(7n-10)}{24}$ = $\left|C_4\right|_{G_{2n+2}}$.   
\end{enumerate}  }
\end{theorem}

\begin{proof}\quad (i) Let $V(G_{2n+3})$ = $\{u_1,u_2,\dots,u_{2n+3}\}$ = $V(G^c_{2n+3})$ where $u_i$ is the vertex with sum labeling $i$ in $G_{2n+3}$ and anti-sum labeling $i$ in $G^c_{2n+3}$, $1 \leq i \leq 2n+3$ and $n \in\mathbb{N}$. At first, let us prove the result on $G_{2n+3}$, $n \in\mathbb{N}$. $\{u_1,u_2,\dots,u_{n+2}\}$ is a clique and $\{u_{n+3},u_{n+4},\dots,u_{2n+3}\}$ is a stable set to $G_{2n+3}$.  Using Theorem \ref{2.11} graph $G_{2n+3}\setminus\{u_1,u_{2n+3}\}$ is isomorphic to $G_{2n+1}$, without the vertex labels. Also, in $G_{2n+3}$, $u_1$ is adjacent to $u_2,u_3,\dots,u_{2n+2};$ $u_{2n+3}$ is an isolated vertex and $u_{2n+2}$ is a pendant vertex. Therefore, $\left|C_4\right|_{G_{2n+3}}$ = $\left|C_4\right|_{G_{2n+1}}$ + number of cycles of length four, each with $u_1$ as a vertex in $G_{2n+3}$. Also, none of $u_{2n+2}$ and $u_{2n+3}$ is a vertex of any cycle of length $4$ in $G_{2n+3}.$
	
	Let $(u_1 , u_i , u_j , u_k)$ be any cycle of length $4$ in $G_{2n+3}$, $1 < i,j,k < 2n+2$ and $i,j,k$ are all distinct. Under the above conditions, the following three types of $C_4s$ with $u_1$ as a vertex arise in $G_{2n+3}$. Type-1: $u_i,u_j,u_k\in \{u_2,u_3,\dots,u_{n+2}\}$, Type-2: $u_i,u_j\in \{u_2,u_3,\dots,u_{n+2}\}$ and $u_k\in \{u_{n+3},u_{n+4},\dots,u_{2n+1}\}$ and Type-3: $u_i\in \{u_2,u_3,\dots,u_{n+2}\}$ and $u_j,u_k\in \{u_{n+3},u_{n+4},\dots,u_{2n+1}\}$. Now, let us calculate number of $C_4s$ in $G_{2n+3}$ under each type. 
		
\medskip
\textbf{Number of $C_4s$ under Type-1:} Here , $u_i,u_j,u_k\in \{u_2,u_3,\dots,u_{n+2}\}$ in $G_{2n+3}$. Number of ways of selecting $3$ vertices $u_i,u_j,u_k$ out of $u_2,u_3,\dots,u_{n+2}$ is $\binom{n+1}{3}$. There are $3$ different $C_4s$ with $u_1, u_i , u_j , u_k$ as vertices under type-1, namely, $(u_1,u_i,u_j,u_k)$, $(u_1,u_i,u_k,u_j)$ and $(u_1,u_j,u_i,u_k)$. Therefore, total number of $C_4s$ under type-1 in $G_{2n+3}$ is $3 \times \binom{n+1}{3}$ =  $(n+1)n(n-1)/2$. \\

\medskip
\textbf{Number of $C_4s$ under Type-2:} Here, $u_i,u_j\in \{u_2,u_3,\dots,u_{n+2}\}$ and $u_k\in \{u_{n+3}$, $u_{n+4},\dots,u_{2n+1}\}$. Consider all possible cycles, each of length $4$ and with vertices $u_1, u_i , u_j$ and $u_k$ in $G_{2n+3}.$
	
	When $k$ = $2n+1$, $u_k$ = $u_{2n+1}$ is adjacent to $u_1$ and $u_2$ only. And under this case, $u_2$ = $u_i$ or $u_2$ = $u_j$. W.l.g., assume $u_2$ = $u_i$. This implies, $2$ = $i < 3$ $\leq$ $j$ $\leq$ $n+2$ and any $C_4$ under this case is of the form $(u_1 , u_{k} , u_i , u_j)$ = $(u_1 , u_{2n+1} , u_2 , u_j)$, $u_j\in\{u_3$, $u_4$, $\dots$, $u_{n+2}\}$ and number of such $C_4s$ is $\left|\{u_3,u_4,\dots,u_{n+2}\}\right|$ = $n$. 
	
	When $k = 2n$, $u_k = u_{2n}$ is adjacent to $u_1, u_2$ and $u_3$ only. And any $C_4$ of type-2 is of the form $(u_1, u_{2n}, u_2,u_x)$ or $(u_1, u_{2n}, u_3, u_y)$, $u_x\in \{u_3$, $u_4,\dots,u_{n+2}\}$ and $u_y\in \{u_2,u_4,u_5,\dots,u_{n+2}\}$. Number of such $C_4s$ is $2n$. 
	
	When $k = 2n-1$, $u_k = u_{2n-1}$ is adjacent to $u_1, u_2 , u_3$ and $u_4$ only and thereby $d(u_k)$ = $4$. Therefore, number of such $C_4s$ is $(4-1)n$ = $3n$. 
	
	In general, when $k = 2n+3-x$ and $2 \leq x \leq n$, $u_k = u_{2n+3-x}$ is adjacent to $u_1, u_2, \dots, u_x$ only and thereby $d(u_k)$ = $d(u_{2n+3-x})$ = $x$ and number of $C_4s$ of the form $(u_1, u_{2n+3-x}, u_i, u_j)$ in $G_{2n+3}$ is $(x-1)n$ where $u_i\in \{u_2,u_3,\dots,u_x\}$ and $u_j\in \{u_2,u_3,\dots,u_{n+2}\}\setminus \{u_i\}$. Therefore, total number of $C_4s$ of type-2 in $G_{2n+3}$ is $\sum_{x=2}^n{(x-1)n}$ = $n\big(\sum_{x=1}^{n-1}{x}\big)$ = $(n-1)n^2/2$. \\ 

\medskip
\textbf{Number of $C_4s$ under Type-3:} Here, $u_i\in \{u_2,u_3,\dots,u_{n+2}\}$ and $u_j,u_k\in \{u_{n+3}$, $u_{n+4}$, $\dots,u_{2n+1}\}$ and $u_j$ and $u_k$ are adjacent to $u_1$ for every $j,k\in \{n+3,n+4$, $\dots$, $2n+2\}$ in $G_{2n+3}$, $j \neq k$. W.l.g., assume that $j < k$. If $u_j$ and $u_k$ are adjacent to $u_i$, then $j+i < k+i \leq 2n+3.$
	
	In $G_{2n+3}$, $u_{n+2+x}$ is adjacent to $u_1,u_2,\dots$, $u_{n+1-x}$, $1 \leq x \leq n$ and thereby $d(u_{n+2+x})$ = $n+1-x$. Also, $u_{n+1}$ and $u_{n+2}$ are non-adjacent to $u_{n+3}$ and $u_{2n+2}$ is a pendant vertex. Hence, none of $u_{n+1}$, $u_{n+2}$ and $u_{2n+2}$ is a vertex of any $C_4$ of type-3 in $G_{2n+3}$. 
	
	When $u_{k}$ = $u_{2n+3-x}$ and $2 \leq x \leq n-1$, different possibilities of $u_i$ in $C_4s$ of type-3 in $G_{2n+3}$ are $u_2,u_3,\dots,u_x$. And corresponding to each pair of $u_i$ and $u_k$, different possibilities of $u_js$ are $u_{k-1}$, $u_{k-2}$, $\dots$, $u_{n+3}$ in $G_{2n+3}$. Therefore, number of $C_4s$ of type-3 in $G_{2n+3}$ with $u_k$ = $u_{2n+3-x}$ is $(x-1)(k-1-(n+2))$ = $(x-1)(n-x)$. Hence, total number of $C_4s$ of type-3 in $G_{2n+3}$ is $\sum_{x=2}^{n-1}{(n-x)(x-1)}$ = $\sum_{x=1}^{n-2}{(n-1-x)x}$ = $n(n-1)(n-2)/6$.
	
	Cycle $C_4$ of the form $(u_1 , u_i , u_j , u_k)$ with $u_i,u_j,u_k\in \{u_{n+3},u_{n+4},\dots,u_{2n+3}\}$ doesn't exist in $G_{2n+3}$ since  $\{u_{n+3},u_{n+4}$, $\dots$, $u_{2n+3}\}$ is a stable set to split graph $G_{2n+3}$. 
	
	Adding all $C_4s$ in the three types, we obtain, total number of $C_4s$ in $G_{2n+3}$ with $u_1$ as a vertex is $$\frac{(n-1)n(n+1)}{2} + \frac{(n-1)n^2}{2} + \frac{(n-2)(n-1)n}{6} = \frac{(n-1)n(7n+1)}{6},~ n \geq 2.$$ 
	$\therefore$ For $n \geq 2$,  

\medskip	
	$\left|C_4\right|_{G_{2n+3}}$ = $\left|C_4\right|_{G_{2n+1}}$ + $\frac{7n^3-6n^2-n}{6}$ 
	
\hspace{1.3cm}	= $\frac{1}{6}((7n^3-6n^2-n)$ + $(7{(n-1)}^3 - 6{(n-1)}^2-(n-1)))$ + $\left|C_4\right|_{G_{2n-1}}$ 
	
\hspace{1.3cm}	= $\frac{1}{6}((7n^3-6n^2-n)$ + $(7{(n-1)}^3 - 6{(n-1)}^2 -(n-1))$ + $\dots$ 
	
	\hfill + $(7\times 2^3 -6\times 2^2 -2))$ + $\left|C_4\right|_{G_5}$ 
	
\hspace{1cm}		= $\frac{1}{6}((7n^3-6n^2-n)$ + $(7{(n-1)}^3 - 6{(n-1)}^2 -(n-1))$ + $\dots$  + $(7\times 2^3 -6\times 2^2 -2))$ + 0 
	
\hspace{1cm}		= $(n-1)n(n+1)(7n+6)/24$. 
	
\medskip
Now, let us to prove the result on $G^c_{2n+1}$. Consider, graph $G^c_{2n+1}$, $n \in\mathbb{N}$. $\{u_1,u_2,\dots,u_{n}\}$ is a stable set and $\{u_{n+1},u_{n+2},\dots,u_{2n+1}\}$ is a clique to $G^c_{2n+1}$.  Using Theorem \ref{2.11}, graph $G^c_{2n+1} \setminus \{u_1,u_{2n+1}\}$ is isomorphic to $G^c_{2n-1}$, without the vertex labels. In $G^c_{2n+1}$, $u_{2n+1}$ is adjacent to $u_1,u_2,\dots,u_{2n}$ and $u_{1}$ is a pendant vertex. Hence, $u_1$ is not a vertex in any cycle of length 4 in $G^c_{2n+1}$. Therefore, $\left|C_4\right|_{G^c_{2n+1}}$ = $\left|C_4\right|_{G^c_{2n-1}}$ $+$ number of cycles of length four, each with $u_{2n+1}$ as a vertex in $G^c_{2n+1}$. 
	
	Let $(u_{2n+1} , u_k , u_j , u_i)$ be any cycle of length $4$ with $u_{2n+1}$ as a vertex in $G^c_{2n+1}$,  $2 \leq i,j,k \leq 2n$ and $i,j,k$ are all distinct. Under the above conditions, the following three types of $C_4s$ with $u_{2n+1}$ as a vertex arise in $G^c_{2n+1}$. Type-1: $u_i,u_j,u_k\in \{u_{n+1},u_{n+2},\dots,u_{2n}\}$, Type-2: $u_k,u_j\in \{u_{n+1},u_{n+2},\dots,u_{2n}\}$ and $u_i\in \{u_{2},u_{3},\dots,u_{n}\}$ and Type-3: $u_k\in \{u_{n+1},u_{n+2},\dots,u_{2n}\}$ and $u_i,u_j\in \{u_{2}$, $u_{3},\dots,u_{n}\}$. Now, let us calculate number of $C_4s$ in $G^c_{2n+1}$ in each type. W.l.g., assume that $i < j < k$. 
		
\medskip
\textbf{Number of $C_4s$ under Type-1:} Here , $u_i,u_j,u_k\in \{u_{n+1},u_{n+2},\dots,u_{2n}\}$ in $G^c_{2n+1}$. Number of ways of selecting $3$ vertices $u_i,u_j,u_k$ out of $u_{n+1},u_{n+2},\dots$, $u_{2n}$ is $\binom{n}{3}$. There are $3$ different $C_4s$ with $u_{2n+1}, u_i , u_j , u_k$ as vertices under type-1, namely, $(u_{2n+1},u_k,u_j,u_i)$, $(u_{2n+1},u_k,u_i,u_j)$ and $(u_{2n+1},u_j,u_k,u_i)$. Hence, total number of $C_4s$ of type-1 in $G^c_{2n+1}$ is $3\times \binom{n}{3}$ =  $n(n-1)(n-2)/2$. 

\medskip
\textbf{Number of $C_4s$ under Type-2:} Here, $u_k,u_j\in \{u_{n+1},u_{n+2},\dots,u_{2n}\}$ and $u_i\in \{u_{2},u_{3},\dots,u_{n}\}$. Consider all possible cycles $(u_{2n+1}, u_i , u_j , u_k)$ in $G^c_{2n+1}$.
	
	When $i$ = $2$, $u_i$ = $u_{2}$ is adjacent to $u_{2n+1}$ and $u_{2n}$ only. And under this case, $d(u_i)$ = $2$, $u_k$ = $u_{2n}$ and $u_j$ = $u_{2n-1}, u_{2n-2}, \dots, u_{n+1}$. Number of such $C_4s$ is $\left|\{u_{2n-1}, u_{2n-2}, \dots, u_{n+1}\}\right|$ = $n-1$. 
	
	When $i = 3$, $u_i = u_{3}$ is adjacent to $u_{2n+1}, u_{2n}$ and $u_{2n-1}$ only and $d(u_i) = 3$. And any $C_4$ of type-2 is of the form $(u_{2n+1}, u_{3}, u_{2n},u_x)$ or $(u_{2n+1}, u_{3}, u_{2n-1}, u_y)$ where $u_x\in \{u_{2n-1}$, $u_{2n-2},\dots,u_{n+1}\}$ and $u_y\in \{u_{2n},u_{2n-2},u_{2n-3},\dots,u_{n+1}\}$. Number of such $C_4s$ is $2(n-1)$. 
	
	In general, when $i$ = $x$ and $2 \leq x \leq n$, $u_i$ = $u_{x}$ is adjacent to $u_{2n+1}, u_{2n},\dots$, $u_{2n+1-(x-1)}$ only and thereby $d(u_i)$ = $x$ and $u_k$ takes values $u_{2n}, u_{2n-1},\dots$, $u_{2n+1-(x-1)}$ and $u_j\in \{u_{2n},u_{2n-1},\dots$, $u_{n+1} \}$ $\setminus \{u_k \}$. Therefore, number of such $C_4s$ of the form $(u_{2n+1},u_i,u_k,u_j)$ is $(x-1)(n-1)$, $2 \leq x \leq n$. Here, $j$ need not be less than $k.$
	
\medskip
	Total number of $C_4s$ of type-2 in $G^c_{2n+1}$ is 
	$$\sum_{x=2}^n{(x-1)(n-1)} = (n-1)\big(\sum_{x=1}^{n-1}{x}\big) = (n-1)^2n/2.$$ 
\textbf{Number of $C_4s$ under Type-3:} Here, $u_k\in \{u_{n+1},u_{n+2},\dots,u_{2n}\}$ and $u_i,u_j\in \{u_{2}$, $u_{3}$, $\dots, u_{n}\}$, $i \neq j$. Consider all possible cycles $(u_{2n+1}, u_i, u_j, u_k)$ in $G^c_{2n+1}$. For a given $i$, $2 \leq i \leq n-1$, $j$ takes values $i+1,i+2,\dots,n$ and possible values of $k$ are $2n,2n-1,\dots,2n-(i-2)$. Hence, total number of $C_4s$ of type-3 in $G^c_{2n+1}$ is 
	$$\sum_{i=2}^{n-1}{(n-i)(i-1)} = \sum_{i=1}^{n-2}{i(n-1-i)} = \frac{(n-2)(n-1)^2}{2} - \frac{(n-2)(n-1)(2n-3)}{6} = \frac{(n-2)(n-1)n}{6}.$$	
	Cycle $C_4$ of the form $(u_{2n+1}, u_k , u_j, u_i)$ with $u_i,u_j,u_k\in \{u_{2},u_{3},\dots,u_{n}\}$ doesn't exist in $G_{2n+1}^c$ since  $\{u_{2},u_{3}$, $\dots$, $u_{n}\}$ is a stable set to split graph $G^c_{2n+1}$. 
	
\medskip
Adding all $C_4s$ in the three types, we obtain, total number of $C_4s$ with $u_{2n+1}$ as a vertex in $G^c_{2n+1}$ is 
	$$\frac{(n-2)(n-1)n}{2} + \frac{(n-1)^2n}{2} + \frac{(n-2)(n-1)n}{6} = \frac{(n-1)n(7n-11)}{6},~ n \geq 2.$$ 
	$\therefore$ For $n \geq 2$, 
	
\medskip
	$\left|C_4\right|_{G^c_{2n+1}}$ = $\left|C_4\right|_{G^c_{2n-1}}$ + $\frac{7n^3-18n^2+11n}{6}$ 
	
	\hspace{1.4cm} = $\frac{1}{6}((7n^3-18n^2+11n)$ + $(7{(n-1)}^3 - 18{(n-1)}^2+11(n-1)))$ + $\left|C_4\right|_{G^c_{2n-3}}$ 
	
	\hspace{1.4cm} = $\frac{1}{6}((7n^3-18n^2+11n)$ + $(7{(n-1)}^3 - 18{(n-1)}^2 +11(n-1))$ + $\dots$  
	
	\hfill + $(7\times 2^3 -18\times 2^2 +11\times 2))$ + $\left|C_4\right|_{G^c_{3}}$ 
	
	\hspace{1.4cm} = $\frac{1}{6}((7n^3-18n^2+11n)$ + $(7{(n-1)}^3 - 18{(n-1)}^2 +11(n-1))$ + $\dots$ 
	
	\hfill + $(7\times 2^3 -18\times 2^2 +11\times 2))$ + 0 
	
	\hspace{1.4cm} = $\frac{(n-1)n(n+1)(7n-10)}{24}$. 
	
\medskip
Rest of the result follows from Theorem \ref{8.3}.  
\end{proof}

\begin{theorem} [{\cite[Theorem 2.5]{vr14}}] \label{8.10c} {\rm \quad Let $n\in\mathbb{N}$. 
		\begin{enumerate}
			\item [\rm (i)]  $|E(G_{0,n})|$ = $|E(G_{n+2})|$ = $n$ + $|E(G_{n})|$, and 
			\\
			$|E(G_{-1,n})|$ = $|E(G_{n+4})|-1$ = $2n+1$ + $|E(G_{n})|$. 
			
			\item [\rm (ii)]  $\left|C_3\right|_{G_{0,n}}$ = $\left|C_3\right|_{G_{n+2}}$, and $\left|C_3\right|_{G_{-1,n}}$ = $\left|C_3\right|_{G_{n+4}}$. 
			
			\item [\rm (iii)] $\left|C_4\right|_{G_{0,n}}$ =  $\left|C_4\right|_{G_{n+2}}$, and $\left|C_4\right|_{G_{-1,n}}$ =  $\left|C_4\right|_{G_{n+4}}$.
	\end{enumerate}}
\end{theorem}

\begin{proof}\quad Result follows from Lemma \ref{4.7}. 
\end{proof} 

\begin{theorem} [{\cite[Theorem 2.6]{vr14}}] \label{8.11} {\rm \quad Let $n\in\mathbb{N}$. 
\begin{enumerate}
	\item [\rm (i)]  $\left|C_4\right|_{G_{0,2n}}$ = $\left|C_4\right|_{G_{2n+2}}$ = $\frac{(n-1)n(n+1)(7n-10)}{24}$. 

	\item [\rm (ii)]   $\left|C_4\right|_{G_{0,2n+1}}$ = $\left|C_4\right|_{G_{2n+3}}$ = $\frac{(n-1)n(n+1)(7n+6)}{24}$. 

	\item [\rm (iii)]  $\left|C_4\right|_{G_{-1,2n}}$ = $\left|C_4\right|_{G_{2n+4}}$ = $\frac{n(n+1)(n+2)(7n-3)}{24}$. 

	\item [\rm (iv)] $\left|C_4\right|_{G_{-1,2n+1}}$ =  $\left|C_4\right|_{G_{2n+5}}$ = $\frac{n(n+1)(n+2)(7n+13)}{24}$.
\end{enumerate}}
\end{theorem}

\begin{proof}\quad Result follows from Theorems \ref{8.9} and \ref{8.10c}. 
\end{proof}

\begin{theorem} [{\cite[Theorem 2.7]{vr14}}] \label{8.12} {\rm  \quad Let $n\in\mathbb{N}$.
\begin{enumerate}
	\item [\rm (i)] Number of $P_3s$ in $G_{2n}$ such that each $P_3$ = $uvw$ with $uw\notin E(G_{2n})$ 
	
	is $\frac{(n-1)n(2n-1)}{6}$, $u,v,w\in V(G_{2n})$. 
		
	\item [\rm (ii)] Number of $P_3s$ in $G_{2n+1}$ such that each $P_3$ = $uvw$ with $uw \notin E(G_{2n+1})$ 
	
	is $\frac{(n-1)n(n+1)}{3}$, $u,v,w\in V(G_{2n+1})$.  
\end{enumerate} }
\end{theorem}

\begin{proof} \quad Let $V(G_{2n})$ = $\{u_1,u_2,\dots,u_{2n}\}$ where $u_j$ is the vertex of $G_{2n}$ with sum labeling $j$, $j$ = $1,2,\dots,2n$. $\{u_1,u_2,\dots$,  $u_n\}$ is a clique and $\{u_{n+1},u_{n+2},\dots,u_{2n}\}$ is a stable set to split graph $G_{2n}$ and vertex $u_n$ is non-adjacent to $u_{n+1}$, $u_{n+2}$, $\dots$, $u_{2n}$. Each required $P_3$ in $G_{2n}$ contains at least one element of $\{u_{n+1}$, $u_{n+2}$, $\dots$, $u_{2n-1}\}$. In $G_{2n}$, counting of $P_3s$ such that each $P_3$ = $uvw$ and $uw \notin E(G_{2n})$ is done as follows, $u,v,w\in V(G_{2n})$. W.l.g., assume that $1 \leq i < j < 2n-k \leq 2n-1$. For $1 \leq k \leq n-1$, vertex $u_{2n-k}$ is adjacent to $v_i$ for $i$ = $1,2,\dots,k$ and $P_3$ = $u_{2n-k}u_iu_j$ is a required path on the 3 vertices for $j$ = $k+1,k+2,\dots,2n-k-1$. Therefore, in $G_{2n}$, number of $P_3s$ such that each $P_3$ = $uvw$ with $uw \notin E(G_{2n})$ and $u,v,w\in V(G_{2n})$ is
	 $$\sum_{k=1}^{n-1}\big(\sum_{i=1}^{k}(2n-2k-1)\big) = \sum_{k=1}^{n-1}k(2n-1-2k) = \frac{(n-1)n(2n-1)}{2} - \frac{(n-1)n(2n-1)}{3} = \frac{(n-1)n(2n-1)}{6}.$$	
	Similarly, let $V(G_{2n+1})$ = $\{u_1,u_2,\dots,u_{2n+1}\}$ where $u_j$ is the vertex of $G_{2n+1}$ with sum labeling $j$, $j$ = $1,2,\dots,2n+1$. $\{u_1,u_2,\dots$,  $u_{n+1}\}$ is a clique and $\{u_{n+2},u_{n+3},\dots,u_{2n+1}\}$ is a stable set to split graph $G_{2n+1}$ and vertex $u_{n+1}$ is non-adjacent to $u_{n+2}$, $u_{n+3}$, $\dots$, $u_{2n+1}$. Each required $P_3$ in $G_{2n+1}$ contains at least one element of $\{u_{n+2}$, $u_{n+3}$, $\dots$, $u_{2n}\}$. In $G_{2n+1}$, counting of $P_3s$ such that each $P_3$ = $uvw$ and $uw \notin E(G_{2n+1})$ is done as follows, $u,v,w\in V(G_{2n+1})$. W.l.g., assume that $1 \leq i < j < 2n+1-k \leq 2n$. For $1 \leq k \leq n-1$, vertex $u_{2n+1-k}$ is adjacent to $v_i$ for $i$ = $1,2,\dots,k$ and $P_3$ = $u_{2n+1-k}u_iu_j$ is a required path on the 3 vertices for $j$ = $k+1,k+2,\dots,2n+1-k-1$. Therefore, in $G_{2n+1}$, number of $P_3s$ such that each $P_3$ = $uvw$ with $uw \notin E(G_{2n+1})$ and $u,v,w\in V(G_{2n+1})$ is $$\sum_{k=1}^{n-1}\big(\sum_{i=1}^{k}(2n-2k)\big) = \sum_{k=1}^{n-1}k(2n-2k) = \frac{2n(n-1)n}{2} - \frac{2(n-1)n(2n-1)}{6} = \frac{(n-1)n(n+1)}{3}.$$ 
		Hence the result. 
\end{proof}

\begin{theorem} [{\cite[Theorem 2.8]{vr14}}] \label{8.13} {\rm  \quad Let  $m,n \geq 2$ and $m,n\in\mathbb{N}$. 
\begin{enumerate}
	\item [\rm (i)] $\left|C_4\right|_{G_{-m,n}}$ = $\left|C_4\right|_{G_{m}}$ + $\left|C_4\right|_{G_{n}}$ + $mC_2.nC_2$ + number of $C_4s$ with $K_1$ as a vertex in $G_{-m,n}$ 
	
\hspace{1.1cm}	= $\left|C_4\right|_{G_{m}}$ + $\left|C_4\right|_{G_{n}}$ + $3(\left|C_3\right|_{G_{m}}$ + $\left|C_3\right|_{G_{n}})$ 
	
\hfill	+ $2(n.||G_m|| + m.||G_n||)$ + $mC_2.nC_2 + n.mC_2 + m.nC_2$ 
	
\hfill		+ (number of $P_3s$ in $G_{-m}$, each $P_3 = uvw$ with $uw\notin E(G_{-m}))$ 
	
\hfill	+ (number of $P_3s$ in $G_n$, each $P_3$ = $uvw$ with $uw\notin E(G_n))$. 
	
	\item [\rm (ii)] $\left|C_4\right|_{G_{-m,n}^c}$ = $\left|C_4\right|_{G_{m}^c}$ + $\left|C_4\right|_{G_{n}^c}$.  
\end{enumerate} }   
\end{theorem}

\begin{proof}  Using Theorem \ref{1.2.4}, we have $G_{-m,n}$ = $K_1  * (G_{-m} * G_n)$ = $K_1  * (G_{-m} \cup G_n \cup K_{m,n})$, and $G^c_{-m,n}$ = $K_1(0) \cup (G^c_{-m} \cup G^c_n)$ where the vertices of $K_{m,n}$ are vertices of $G_{-m} \cup G_n$, $m,n \in\mathbb{N}_0$. Here, $K_1$ is the vertex with integral sum label 0 and adjacent to all other vertices in $G_{-m,n}$ and an isolated vertex in $G^c_{-m,n}$. Clearly, $\left|C_4\right|_{G^c_{-m,n}}$ = $\left|C_4\right|_{G^c_{m}}$ + $\left|C_4\right|_{G^c_{n}}$ since $G^c_{m}$ and $G^c_{n}$ are disjoint subgraphs in $G^c_{-m,n}$ and $G^c_{m}$ and $G^c_{-m}$ are isomorphic graphs without vertex labels. Now, 
	
\medskip
	$C_4$s in $G_{-m} \cup G_n \cup K_{m,n}$ is ($C_4$s in $G_{-m}$) $\cup$ ($C_4$s in $G_n$) $\cup$ ($C_4$s in $K_{m,n}$) and 

\medskip
\noindent
$\left|C_4\right|_{G_{-m,n}}$ = Number of $C_4$s, each $C_4$ with $K_1$ as a vertex in $G_{-m,n}$

\hfill + number of $C_4$s, each $C_4$ without $K_1$ as a vertex in $G_{-m,n}$ 

\hspace{1.1cm} = $\left|C_4\right|_{G_{m}}$ + $\left|C_4\right|_{G_{n}}$ + $\left|C_4\right|_{K_{m,n}}$ + number of $C_4$s with $K_1$ as a vertex in $G_{-m,n}$ 

\hspace{1.1cm} = $\left|C_4\right|_{G_{m}}$ + $\left|C_4\right|_{G_{n}}$ + $\binom{m}{2} \times \binom{n}{2}$ + number of $C_4$s with $K_1$ as a vertex in $G_{-m,n}$ since

\hfill $K_{m,n}$ is a complete bipartite graph and number of $C_4$s in $K_{m,n}$ is $\binom{m}{2} \times \binom{n}{2}$, $m,n \geq 2$.

\medskip
Let $V(G_{-m,n})$ = $\{u_0,u_1,u_2,\dots,u_{m+n}\}$ where, in $G_{-m,n}$ = $K_1 * (G_{-m} + G_n)$, $u_0$ is the vertex $K_1$ with integral sum labeling $0$, $u_i$ is the vertex of $G_{-m}$ with integral sum labeling $-i$ for $i$ = $1,2,\dots,m$ and $u_{m+j}$ is the vertex of $G_n$ with integral sum labeling $j$, $j$ = $1,2,\dots,n$. Let $1 \leq \left|i\right| < \left|j\right| < \left|k\right| \leq m+n$ and $(u_0,u_i,u_j,u_k)$ be any cycle of length $4$ with $u_0$ as a vertex in $G_{-m,n}$. The following types of $C_4s$ with $u_0$ as a vertex arise. Type-1: $u_i,u_j,u_k\in V(G_{-m});$ Type-2: $u_i,u_j,u_k\in V(G_n);$ Type-3: $u_i,u_j\in V(G_{-m})$ and $u_k\in V(G_n)$ and Type-4: $u_i\in V(G_{-m})$ and $u_j,u_k\in V(G_n)$. Let us calculate number of $C_4s$ with $K_1$ as a vertex in $G_{-m,n}$ in each type. 

\medskip   
\textbf{Number of $C_4s$ under Type-1:} Here, $u_i,u_j,u_k \in V(G_{-m})$. In this case, $C_4$ is formed in $G_{-m,n}$ with vertices $u_0$, $u_i$, $u_j$ and $u_k$, either $(u_i,u_j,u_k)$ is a cycle of length $3$ in $G_{-m}$ or $u_iu_ju_k$ is a path of length $2$ in $G_{-m}$ with $u_iu_k \notin E(G_{-m})$. When $(u_i,u_j,u_k)$ is a cycle of length 3 in $G_{-m}$, possible type-1 $C_4s$ in $G_{-m,n}$ with vertices $u_0,u_i,u_j,u_k$ are $(u_0,u_i,u_j,u_k)$, $(u_0,u_i,u_k,u_j)$ and $(u_0,u_j,u_i,u_k)$. Hence, number of $C_4s$ in $G_{-m,n}$ with vertices $u_0,u_i,u_j,u_k$ when $(u_i,u_j,u_k)$ is a cycle of length 3 in $G_{-m}$ is $3.\left|C_3\right|_{G_{m}}$. Similarly, when $u_iu_ju_k$ is a path of length 2 in $G_{-m}$ and $u_iu_k \notin E(G_{-m})$, then the only possible type-1 $C_4$ in $G_{-m,n}$ with vertices $u_0,u_i,u_j,u_k$ is $(u_0,u_i,u_j,u_k)$. Thus, number of $C_4s$ in $G_{-m,n}$ with vertices $u_0,u_i,u_j,u_k$ when $u_iu_ju_k$ is a path of length 2 in $G_{-m}$ but $u_iu_k$ is not an edge of $G_{-m}$ = Number of $P_3s$ in $G_{-m}$, each $P_3$ is not a subgraph of any $C_3$ of $G_{-m}$. Hence, number of $C_4s$ of type-1 in $G_{-m,n}$ = $3.\left|C_3\right|_{G_{m}}$ + number of $P_3s$ in $G_{-m}$ such that each $P_3$ is not a subgraph of any $C_3$ of $G_{-m}$. 

\medskip
\textbf{Number of $C_4s$ under Type-2:} Here, $u_i,u_j,u_k \in V(G_{n})$. Similar to type-1, we obtain, number of $C_4s$ of type-2 in $G_{-m,n}$ = $3.\left|C_3\right|_{G_{n}}$ + number of $P_3s$ in $G_n$ such that each $P_3$ is not a subgraph of any $C_3$ of $G_n$.  

\medskip
\textbf{Number of $C_4s$ under Type-3:} Here, $u_i,u_j \in V(G_{-m})$ and $u_k\in V(G_n)$. In this case, $C_4$ is formed in $G_{-m,n}$ with vertices $u_o,u_i,u_j,u_k$ such that either $u_i$ and $u_j$ are adjacent or $u_i$ and $u_j$ are non-adjacent whereas $u_k$ takes all vertices of $G_n$. When $u_i$ and $u_j$ are adjacent, possible $C_4s$ of type-3 in $G_{-m,n}$ with vertices $u_0,u_i,u_j,u_k$ are $(u_0,u_i,u_j,u_k)$, $(u_0,u_i,u_k,u_j)$ and $(u_0,u_j,u_i,u_k)$. 

$\therefore$ Number of $C_4s$ of type-3 in $G_{-m,n}$ with vertices $u_0,u_i,u_j,u_k$ when $u_i$ and $u_j$ are adjacent 

\hfill = $3\times||G_{-m}|| \times ($ number of vertices of $G_n)$ = $3n\times ||G_{-m}||$. 

Similarly, when $u_i$ and $u_j$ are non-adjacent, the only possible $C_4$ of type-3 in $G_{-m,n}$ with vertices $u_0,u_i,u_j,u_k$ is $(u_0,u_i,u_k,u_j)$. Number of non-adjacent pair of vertices in $G_{-m}$ = $\binom{m}{2}$ - number of adjacent pair of vertices in $G_{-m}$ = $\binom{m}{2}$ - $||G_{-m}||$. Hence, number of $C_4s$ of type-3 in $G_{-m,n}$ with vertices $u_0,u_i,u_j,u_k$ when $u_i$ and $u_j$ are non-adjacent = $n(\binom{m}{2} - ||G_{-m}||)$. 

\medskip
$\therefore$ Number of $C_4s$ of type-3 in $G_{-m,n}$ 

\hspace{.5cm} = Number of $C_4s$ of type-3 in $G_{-m,n}$ with vertices $u_0,u_i,u_j,u_k$ when $u_i$ and $u_j$ are adjacent 

\hspace{.5cm}  + number of $C_4s$ of type-3 in $G_{-m,n}$ with vertices $u_0,u_i,u_j,u_k$ when $u_i$ and $u_j$ are non-adjacent 

\hspace{.5cm} = $n(\binom{m}{2} + 2\times ||G_{-m}||)$. 

\medskip
\textbf{Number of $C_4s$ under Type-4:} Here, $u_i \in V(G_{-m})$ and $u_j,u_k\in V(G_n)$. Similar to type-3, we obtain, number of $C_4s$ of type-4 in $G_{-m,n}$ = $m(\binom{n}{2} + 2 \times ||G_{n}||)$. 

\medskip	
Therefore, for $m,n \geq 2$, 

$\left|C_4\right|_{G_{-m,n}}$ = Number of $C_4s$ of type-1 in $G_{-m,n}$ + number of $C_4s$ of type-2 in $G_{-m,n}$  

\hfill + number of $C_4s$ of type-3 in $G_{-m,n}$ + number of $C_4s$ of type-4 in $G_{-m,n}$ 

\hspace{1.5cm} = $\left|C_4\right|_{G_{m}}$ + $\left|C_4\right|_{G_{n}}$ + $\binom{m}{2} \times \binom{n}{2}$ + $3 \times \left|C_3\right|_{G_{m}}$ + $3 \times \left|C_3\right|_{G_{n}}$  

\hspace{1.5cm}  + (number of $P_3s$ in $G_{-m}$ $\ni$ each $P_3$ is not a subgraph of any $C_3$ of $G_{-m}$) 

\hspace{1.5cm}  + (number of $P_3s$ in $G_n$ $\ni$ each $P_3$ is not a subgraph of any $C_3$ of $G_n$) 

\hfill + $n(\binom{m}{2} + 2\times ||G_{-m}||)$ + $m(\binom{n}{2} + 2\times ||G_n||)$. 

\medskip
Hence the result.  
\end{proof}

\begin{cor} [{\cite[Corollary 2.9]{vr14}}] \quad \label{8.14} {\rm  \quad Let  $m,n\in\mathbb{N}$. 
\begin{enumerate}
	\item [\rm (i)] $\left|C_4\right|_{G_{-2m,2n}}$ = $\frac{(m-1)m(7m^2+m-18)}{24}$ + $\frac{(n-1)n(7n^2+n-18)}{24}$ 
		+ $mn(4mn+6(m+n)-11)$.  

	\item [\rm (ii)] $\left|C_4\right|_{G_{-2m,2n+1}}$ = $\frac{(m-1)m(7m^2+m-18)}{24}$ + $\frac{(n-1)n(7n^2+17n-2)}{24}$
	\begin{flushright}
		+ $m(4m-3)(2n+1)$ + $mn(4mn+2m+6n+1)$.  
	\end{flushright}

	\item [\rm (iii)] $\left|C_4\right|_{G_{-(2m+1),2n}}$ = $\frac{(m-1)m(7m^2+17m-2)}{24}$ + $\frac{(n-1)n(7n^2+n-18)}{24}$
	\begin{flushright}
		+ $(2m+1)n(4n-3)$ + $mn(4mn+6m+2n+1)$. 
	\end{flushright}

	\item [\rm (iv)] $\left|C_4\right|_{G_{-(2m+1),2n+1}}$ = $\frac{(m-1)m(7m^2+17m-2)}{24}$ + $\frac{(n-1)n(7n^2+17n-2)}{24}$ 
	\begin{flushright}
		+ $(mn+m+n)(2m+1)(2n+1)$ + $4mn(m+n)$ + $2(m^2+n^2)$.   
	\end{flushright}

	\item [\rm (v)] $\left|C_4\right|_{G_{-2m,2n}^c}$ = $\frac{(m-2)(m-1)m(7m-1)}{24}$ + $\frac{(n-2)(n-1)n(7n-1)}{24}$. 

	\item [\rm (vi)] $\left|C_4\right|_{G_{-2m,2n+1}^c}$ = $\frac{(m-2)(m-1)m(7m-1)}{24}$ + $\frac{(n-1)n(n+1)(7n-10)}{24}$. 

	\item [\rm (vii)] $\left|C_4\right|_{G_{-(2m+1),2n}^c}$ = $\frac{(m-1)m(m+1)(7m-10)}{24}$ + $\frac{(n-2)(n-1)n(7n-1)}{24}$. 
	
	\item [\rm (viii)] $\left|C_4\right|_{G_{-(2m+1),2n+1}^c}$ = $\frac{(m-1)m(m+1)(7m-10)}{24}$ + $\frac{(n-1)n(n+1)(7n-10)}{24}$.  
	\end{enumerate} }
\end{cor}

\begin{proof}\quad For $m,n \in\mathbb{N}$, using Theorems \ref{a5}, \ref{8.8},  \ref{8.9}, \ref{8.12},  \ref{8.13},  and Corollary \ref{8.6}, we obtain the following.
\begin{enumerate}
	\item [\rm (i)]  $\left|C_4\right|_{G_{-2m,2n}}$ = $\left|C_4\right|_{G_{2m}}$ + $\left|C_4\right|_{G_{2n}}$ + $3(\left|C_3\right|_{G_{2m}}$ + $\left|C_3\right|_{G_{2n}})$
	\begin{flushright}
		+ $4(n.||G_{2m}|| + m.||G_{2n}||)$ + $\binom{2m}{2}\binom{2n}{2}$ + $2n\binom{2m}{2}$ + $2m \binom{2n}{2}$
	\end{flushright}
	\begin{flushright} 
		+ number of $P_3s$ in $G_{-2m}$ $\ni$ each $P_3$ is not a subgraph of any $C_3$ of $G_{-2m}$ 
	\end{flushright}
	\begin{flushright}
		+ number of $P_3s$ in $G_{2n}$ $\ni$ each $P_3$ is not a subgraph of any $C_3$ of $G_{2n}$
	\end{flushright}
\hspace{.3cm}	= $\frac{(m-1)m(7m^2-31m+34)}{24}$ + $\frac{(n-1)n(7n^2-31n+34)}{24}$ 
		+ $3(\frac{(m-2)(m-1)m}{3}$ + $\frac{(n-2)(n-1)n}{3})$  

	\begin{flushright}
		+ $4mn(m-1)$ + $4mn(n-1)$ + $mn(2m-1)(2n-1)$ + $2mn(2m-1)$  
	\end{flushright}
	\begin{flushright}
		+ $2mn(2n-1)$ + $\frac{(m-1)m(2m-1)}{6}$ + $\frac{(n-1)n(2n-1)}{6}$ 
	\end{flushright} 
	\begin{flushright}  
		= $\frac{(m-1)m(7m^2+m-18)}{24}$ + $\frac{(n-1)n(7n^2+n-18)}{24}$ + $mn(4mn+6(m+n)-11)$. 
	\end{flushright}

	\item [\rm (ii)] $\left|C_4\right|_{G_{-2m,2n+1}}$
	= $\left|C_4\right|_{G_{2m}}$ + $\left|C_4\right|_{G_{2n+1}}$ + $3(\left|C_3\right|_{G_{2m}}$ + $\left|C_3\right|_{G_{2n+1}})$ 
	\begin{flushright}
		+ $2((2n+1)||G_{2m}|| + 2m.||G_{2n+1}||)$ + $\binom{2m}{2} \binom{2n+1}{2}$ + $(2n+1)\binom{2m}{2}$ + $2m\binom{2n+1}{2}$ 
	\end{flushright}
	\begin{flushright}
	+ number of $P_3s$ 	in $G_{-2m}$ such that each $P_3$ is not a subgraph of any $C_3$ of $G_{-2m}$ 
	\end{flushright}
	\begin{flushright}
	+ number of $P_3s$ in $G_{2n+1}$ such that each $P_3$ is not a subgraph of any $C_3$ of $G_{2n+1}$
	\end{flushright}
	= $\frac{(m-2)(m-1)m(7m-17)}{24}$ + $\frac{(n-2)(n-1)n(7n-1)}{24}$ + $(m-2)(m-1)m$ + $\frac{(n-1)n(2n-1)}{2}$ 
	\begin{flushright} 
		+ $2(2n+1)(m-1)m$ + $4mn^2$ + $m(2m-1)(2n+1)n$ + $(2n+1)m(2m-1)$ 
	\end{flushright}
	\begin{flushright}
		+ $2m(2n+1)n$ + $\frac{(m-1)m(2m-1)}{6}$ + $\frac{(n-1)n(n+1)}{3}$
	\end{flushright}
	=  $\frac{(m-1)m(7m^2+m-18)}{24}$ + $\frac{(n-1)n(7n^2+17n-2)}{24}$ 
		+ $m(4m-3)(2n+1)$ + $mn(4mn+2m+6n+1)$. 

	Similarly, we obtain the following. 

\medskip
	\item [\rm (iii)] $\left|C_4\right|_{G_{-(2m+1),2n}}$ 
	= $\frac{(m-1)m(7m^2+17m-2)}{24}$ + $\frac{(n-1)n(7n^2+n-18)}{24}$  
	\begin{flushright}
		+ $(2m+1)n(4n-3)$ + $mn(4mn+6m+2n+1)$.
	\end{flushright}

	\item [\rm (iv)] $\left|C_4\right|_{G_{-(2m+1),2n+1}}$ 
	= $\left|C_4\right|_{G_{2m+1}}$ + $\left|C_4\right|_{G_{2n+1}}$ + $3(\left|C_3\right|_{G_{2m+1}}$ + $\left|C_3\right|_{G_{2n+1}})$ 
	\begin{flushright}
		+ $2((2n+1)||G_{2m+1}|| + (2m+1)||G_{2n+1}||)$ + $\binom{2m+1}{2}\binom{2n+1}{2}$ + $(2n+1)\binom{2m+1}{2}$ + $(2m+1)\binom{2n+1}{2}$	
	\end{flushright}
	\begin{flushright}
		 + number of $P_3s$ in $G_{-2m-1}$
such that each $P_3$ is not a subgraph of any $C_3$ of $G_{-2m-1}$ 
	\end{flushright}
	\begin{flushright}
		+ number of $P_3s$ in $G_{2n+1}$ such that each $P_3$ is not a subgraph of any $C_3$ of $G_{2n+1}$ 
	\end{flushright}
	= $\frac{(m-2)(m-1)m(7m-1)}{24}$ + $\frac{(n-2)(n-1)n(7n-1)}{24}$ + $\frac{(m-1)m(2m-1)}{2}$ + $\frac{(n-1)n(2n-1)}{2}$ 
	\begin{flushright}
		+ $2(2n+1)m^2$ + $2(2m+1)n^2$ + $(2m+1)m(2n+1)n$ + $(2n+1)(2m+1)m$ 
	\end{flushright}
	\begin{flushright}
		+ $(2m+1)(2n+1)n$ + $\frac{(m-1)m(m+1)}{3}$ + $\frac{(n-1)n(n+1)}{3}$
	\end{flushright}   
	= $\frac{(m-1)m(7m^2+17m-2)}{24}$ + $\frac{(n-1)n(7n^2+17n-2)}{24}$
	\begin{flushright}
		+ $2(m^2(2n+1)$ + $(2m+1)n^2)$ + $(mn+m+n)(2m+1)(2n+1)$. 
	\end{flushright}
\end{enumerate}	

\medskip
Results $(v)$ to $(viii)$ follow from $G^c_{-m,n}$ = $K_1(0)$ $\cup$ $G^c_{-m}$ $\cup$ $G^c_{n}$ and using Theorem \ref{8.9}.    
\end{proof}

\subsection{Algebraic expression of natural numbers divisible by 6 and 24}

From Theorems \ref{8.8} and Corollary \ref{8.14}, the following simple properties of natural numbers are obtained.

\begin{theorem} [{ \cite[Theorem 2.10]{vr14}}] \quad \label{8.15} {\rm Let $n\in\mathbb{N}$. \begin{enumerate}
	\item [\rm (i)(a)]  $n(n+1)(7n-4)$ = $n(7n^2+3n-4)$ is divisible by 6.
	\item [\rm (b)]  $n(n+1)(7n+8)$ = $n(7n^2+15n+8)$ is divisible by 6.

	\item [\rm (ii)(a)] $n(n+1)(n+2)(7n-3)$ is divisible by 24.
	\item [\rm (b)] $n(n+1)(n+2)(7n+13)$ is divisible by 24.
	\item [\rm (c)] $n(n+1)(7n^2+15n-10)$ is divisible by 24. 
	\item [\rm (d)] $n(n+1)(7n^2+31n+22)$ is divisible by 24.
\end{enumerate}	}
\end{theorem}

\begin{proof}\quad We use Theorem \ref{8.8} and Corollary \ref{8.14} to prove these results.
\begin{enumerate}
\item [\rm (i)(a)] Using (i) in Theorems \ref{8.8}, we get,
	
	$\frac{(n-1)n(7n-11)}{6}\in\mathbb{N}$, $n \geq 2$ and $n\in\mathbb{N}$. 	
	\\
	$\Rightarrow$ $\frac{n(n+1)(7n-4)}{6}\in\mathbb{N}$,  $n\in\mathbb{N}$. 
	
	This implies, $n(n+1)(7n-4)$ = $n(7n^2+3n-4)$ is divisible by 6, $n\in\mathbb{N}$.
	
\item [\rm (i)(b)] Using (ii) in Theorems \ref{8.8}, we get,

$\frac{(n-1)n(7n+1)}{6}\in\mathbb{N}$, $n \geq 2$ and $n\in\mathbb{N}$. 
\\
$\Rightarrow$ $\frac{n(n+1)(7n+8)}{6}\in\mathbb{N}$,  $n\in\mathbb{N}$. 

This implies, $n(n+1)(7n+8)$ = $n(7n^2+15n+8)$  is divisible by 6, $n\in\mathbb{N}$.

\item [\rm (ii)(a)] Using (i) in Theorems \ref{8.8}, we get,

$\frac{(n-1)n(n+1)(7n-10)}{24}\in\mathbb{N}$, $n \geq 2$ and $n\in\mathbb{N}$. 
\\
$\Rightarrow$ $\frac{n(n+1)(n+2)(7n-3)}{24}\in\mathbb{N}$, $n\in\mathbb{N}$. 

 This implies, $n(n+1)(n+2)(7n-3)$ is divisible by 24, $n\in\mathbb{N}$.

\item [\rm (ii)(b)] Using (ii) in Theorems \ref{8.8}, we get,

$\frac{(n-1)n(n+1)(7n+6)}{24}\in\mathbb{N}$, $n \geq 2$ and $n\in\mathbb{N}$. 
\\
$\Rightarrow$ $\frac{n(n+1)(n+2)(7n+13)}{24}\in\mathbb{N}$, $n\in\mathbb{N}$. 

This implies, $n(n+1)(n+2)(7n+13)$ is divisible by 24, $n\in\mathbb{N}$.

\item [\rm (ii)(c)] Using (i) in Corollary \ref{8.14}, we get,

$\frac{(n-1)n(7n^2+n-18)}{24}\in\mathbb{N}$, $n \geq 2$ and $n\in\mathbb{N}$. 
\\
$\Rightarrow$ $\frac{n(n+1)(7n^2+14n+7+(n+1)-18)}{24}\in\mathbb{N}$, $n\in\mathbb{N}$. 

This implies, $n(n+1)(7n^2+15n-10)$ is divisible by 24, $n\in\mathbb{N}$.

\item [\rm (ii)(d)] Using (ii) in Corollary \ref{8.14}, we get,

$\frac{(n-1)n(7n^2+17n-2)}{24}\in\mathbb{N}$, $n \geq 2$ and $n\in\mathbb{N}$. 
\\
$\Rightarrow$ $\frac{n(n+1)(7n^2+14n+7+(17n+17)-2)}{24}\in\mathbb{N}$, $n\in\mathbb{N}$. 

This implies, $n(n+1)(7n^2+31n+22)$ is divisible by 24, $n\in\mathbb{N}$.
\end{enumerate}
\end{proof}

\section{On new families of integral sum graphs}

In this section, we present results on new families of integral sum graphs. It is proved that Banana trees, star, union of stars, Dutch windmill, triangular book, fan with a handle, and triangular book with bookmark are integral sum graphs. 

\subsection{Star and union of stars are integral sum graphs}

In this subsection, we present results that star graphs and union of stars are integral sum graphs. 

\medskip
{\em Star graph} $S_n$ on $n$ vertices is defined as $S_n$ = $K_1 * ((n - 1)K_1)$ = $K_{1,n-1}$, $n \geq 2$ and $n\in\mathbb{N}$. 

In the following theorem, we present a special type of integral sum labeling of Star graph $S_n$ such that, for each value of $t$ and $d$, the sequence $\{ t{(d + 1)}^{i-1}\}^{n-1}_{i=1}$ of vertex labeling $f(v_i)$s is a GP and strictly monotonic increasing sequence, $n \geq 2$ and $d,i,n,t\in\mathbb{N}$. 

\begin{theorem} \label{sg} {\rm Let $n \geq 2$, $S_n$ be a star graph, $V(S_n)$ = $\{u_0, v_1, v_2, \dots, v_{n-1}\}$, $d(u_0)$ = $n - 1$ and $d(v_1)$ = 1 = $d(v_2)$ = $\dots$ = $d(v_{n-1})$. Define labeling 
		
\medskip
		$f:$ $V(S_n) \to \mathbb{Z}$ $\ni$ 	$f(u_0)$ = 0, $f(v_1)$ = $t$,
		
		$f(v_{j+1})$ = 	$d\big( \sum^{j}_{i=1} f(v_i) \big) + t$, $j$ = 1, 2, $\dots$, $n - 2$.

\medskip
\noindent
	Then, $f$ is an integral sum labeling of $S_n$ and the sequence of vertex labeling is $\{f(u_0) = 0, f(v_i) = t{(d + 1)}^{i-1}\}^{n-1}_{i=1}$, $d,t\in\mathbb{N}$.} 
\end{theorem}

\begin{proof}\quad It is easy to prove $f(v_i)$ = $t{(d + 1)}^{i-1}$ by induction on $i$, $i$ = 1,2,...,$n - 1$ where $d,t\in\mathbb{N}$. Also, 
	
	$(i)$ for $i$ = 1, 2, ..., $n$, $f(u_0)$ + $f(v_i)$ = 0 + $t{(d + 1)}^{i-1}$  = $f(v_i)$ and
	
	$(ii)$ for $1 \leq i < j < k \leq n$ and $i,j,k,d,t\in\mathbb{N}$,  
	$f(v_i)$ + $f(v_j)$ $\neq$ $f(v_k)$ since 
	$t{(d + 1)}^i$ + $t{(d + 1)}^j$ = $t{(d + 1)}^i(1+{(d + 1)}^{j-i})$ $\neq$ $t{(d + 1)}^k$, 	$1 \leq j - i < j < k$. This implies that $u_0$ and $v_i$ are adjacent $\forall$  $i$ whereas $v_j$ and $v_k$ are non-adjacent $\forall$ $j,k$, $j \neq k$ and $1 \leq j,k \leq n$. Thus, $f$ is an integral sum labeling of $S_n$ for $n \geq 2$.
\end{proof} 

\medskip
Xu \cite{x99} proved that the union of any three stars is an integral sum graph. An extended result for any number of stars is given by Vilfred and Nicholas in \cite{vn11} and is presented below.

\begin{theorem} [{ \cite[Theorem 2.6]{vn11}}] \label{su} {\rm \quad  The union of stars is an integral sum graph.} 
\end{theorem}

\begin{proof}\quad Let $G$ = $K_{1, n_1} \cup K_{1, n_2} \cup \dots \cup K_{1, n_t}$ be the given graph, $t \geq 1$. Let $v_i$ denote the central vertex and $u_{i,j}$ denote the end vertices of the $i^{th}$ star $K_{1, n_i}$ where $i$ = 1, 2, . . . , $t$ and $j$ = 1,2,...,$n_i$. 

If $n_i$ = 1 for all $i$, $i$ = 1, 2, . . . , $t$, then $G$ is a perfect matching. In this case, the following labeling $f$ is an integral sum labeling for $G$.

$f(u_{1, 1})$ = 1; $f(v_1)$ = 3;

$f(u_{2, 1})$ = $f(u_{1, 1}) + f(v_1)$ = 4; $f(v_2)$ = $f(v_1) + f(u_{2, 1})$ + 1 = 8;

$f(v_3)$ = $f(v_2) + f(u_{2, 1})$ = 12; $f(u_{3, 1})$ = $f(u_{1, 1})$ - $f(v_3)$ = -11;

$f(v_i)$ = $2f(v_{i-1})$, $i$ = 4,5,...;

$f(u_{i, 1})$ = $f(u_{i-1, 1})$ - $f(v_i)$, $i$ = 4,5,.... 

\medskip
Figure 22 shows integral sum labeling of graph $G$ = $5 K_2$ along with the $i^{th}$ star $K_{1,1}$ = $u_{i,1} v_i$. 

\begin{center}
\begin{tikzpicture}
	
	\node (a1) at (0,0)  [circle,draw,scale=0.8]{$3$};
	\node (b1) at (0,1.5)  [circle,draw,scale=0.8] {$1$};
	\node (a2) at (1.5,0)  [circle,draw,scale=0.8]{$8$};
	\node (b2) at (1.5,1.5)  [circle,draw,scale=0.8] {$4$};
	\node (a3) at (3,0)  [circle,draw,scale=0.65]{$12$};
	\node (b3) at (3,1.5)  [circle,draw,scale=0.5] {$-11$};
	\node (a4) at (4.5,0)  [circle,draw,scale=0.65]{$24$};
	\node (b4) at (4.5, 1.5)  [circle,draw,scale=0.5] {$-35$};
	\node (a5) at (6,0)  [circle,draw,scale=0.65]{$48$};
	\node (b5) at (6,1.5)  [circle,draw,scale=0.5] {$-83$};
	
	\node (a6) at (10,0)  [circle,draw,scale=0.75]{$v_i$};
	\node (b6) at (10,1.5)  [circle,draw,scale=0.55] {$u_{i,1}$};
	
	\draw (a1) -- (b1);
	\draw (a2) -- (b2);
	\draw (a3) -- (b3);
	\draw (a4) -- (b4);
	\draw (a5) -- (b5);
	\draw (a6) -- (b6);
	
\end{tikzpicture}

\medskip
Fig. 22. Integral sum labeling of $5 K_2$ and $i^{th}$ star $K_{1,1}$ = $u_{i,1} v_i$	
\end{center}

Now without loss of generality let us assume, $n_1 \geq 2$ and $t \geq 2$. Define a vertex labeling $f$ on $G$ as follows:

$f(v_1)$ = 2; $f(u_{1, j})$ = $2j$-1 for $j$ = 1,2,...,$n_1$; $f(u_{2, 1})$ = $2n_1$ + 1;
\\	
For $i$ = 2,3,...,$t$-1,

$f(v_i)$ = $2f(u_{i, 1})$,

$f(u_{i, j})$ = $f(u_{i, 1}) + (j-1)f(v_i)$ = $(2j-1)f(u_{i, 1})$, $j$ = 2,3,...,$n_i$;

$f(u_{i+1, 1})$ = $f(v_i) + f(u_{i, n_i})$ = $2 f(u_i, 1) + f(u_i, n_i)$ = $(2n_i + 1) f(u_i, 1)$. 

And $f(v_t)$ = 1 - $f(u_t, 1)$ = $1 - (f(v_{t-1}) + f(u_{t-1}, n_{t-1}))$ = $1-(2n_{t-1}+1) f(u_{t-1}, 1)$,

$f(u_{t, j})$ = $f(u_{t, 1}) - (j-1) f(v_t)$, $j$ = 2,3,...,$n_t$.

\medskip
Figure 23 shows an integral sum labeling of union of stars $G$ = $K_{1, 1} \cup$ $K_{1, 2} \cup K_{1, 3} \cup K_{1, 4} \cup K_{1, 5}$. 

\begin{center}
\begin{tikzpicture}
	
	\node (v1) at (0,0)  [circle,draw,scale=0.8]{$2$};
	\node (a1) at (-.75,1.5)  [circle,draw,scale=0.8]{$1$};
	\node (b1) at (0,1.5)  [circle,draw,scale=0.8] {$3$};
	\node (c1) at (.75,1.5)  [circle,draw,scale=0.8]{$5$};
	\node (v2) at (1.75,0)  [circle,draw,scale=0.7]{$14$};
	\node (a2) at (1.75,1.5)  [circle,draw,scale=0.8]{$7$};
	\node (v3) at (3.25,0)  [circle,draw,scale=0.7]{$42$};
	\node (a3) at (2.75,1.5)  [circle,draw,scale=0.7]{$21$};
	\node (b3) at (3.75,1.5)  [circle,draw,scale=0.7] {$63$};
	\node (v4) at (5.5,0)  [circle,draw,scale=0.6]{$210$};
	\node (a4) at (4.5,1.5)  [circle,draw,scale=0.6]{$105$};
	\node (b4) at (5.25,1.5)  [circle,draw,scale=0.6] {$315$};
	\node (c4) at (6,1.5)  [circle,draw,scale=0.6]{$525$};
	\node (d4) at (6.75,1.5)  [circle,draw,scale=0.6]{$735$};
	\node (v5) at (9,0)  [circle,draw,scale=0.5]{$-944$};
	\node (a5) at (7.5,1.5)  [circle,draw,scale=0.6]{$945$};
	\node (b5) at (8.25,1.5)  [circle,draw,scale=0.5] {$1889$};
	\node (c5) at (9,1.5)  [circle,draw,scale=0.5]{$2833$};
	\node (d5) at (9.75,1.5)  [circle,draw,scale=0.5]{$3777$};
	\node (e5) at (10.5,1.5)  [circle,draw,scale=0.5]{$4721$};
	
	\draw (v1) -- (a1);
	\draw (v1) -- (b1);
	\draw (v1) -- (c1);
	\draw (v2) -- (a2);
	\draw (v3) -- (a3);
	\draw (v3) -- (b3);
	\draw (v4) -- (a4);
	\draw (v4) -- (b4);
	\draw (v4) -- (c4);
	\draw (v4) -- (d4);
	
	\draw (v5) -- (a5);
	\draw (v5) -- (b5);
	\draw (v5) -- (c5);
	\draw (v5) -- (d5);
	\draw (v5) -- (e5);
\end{tikzpicture}

\medskip
Fig. 23. integral sum labeling of $G$ = $K_{1, 3} \cup K_{1, 1} \cup K_{1, 2} \cup K_{1, 4} \cup K_{1, 5}$	
\end{center}

In this labeling all the end vertices receive odd integers whereas all the central vertices receive even integers and using similar argument given in the proof of Theorem \ref{6.2.7} (Here $A$ = $\emptyset$, $B$ = $\{v_i:$ $i$ = 1,2, . . . , $t-1\}$ and $C$ = $\{u_{,j} :$ $i$ = 1,2,...,$t$ and $j$ = 1,2,...,$n_i\}$. While considering $e$ = $uv\in E(G^c)$, the following three cases arise. Case-1. $u,v\in B$; Case-2. $u,v\in C$ and Case-3. $u\in B$ and $v\in C$. In each case, we can show that $f(u) + f(v)\notin f(V)$, $uv\in E(G^c)$.), it is easy to verify that the labeling $f$ is an integral sum labeling of $G$. Hence we get the result.	
\end{proof} 

Chen \cite{c98} defined a {\em generalized star} as a tree obtained from a star by extending each edge to a path and proved the following result.

\begin{theorem} [{ \cite[Theorem 4]{c98}}]  \label{gs} \quad {\rm  Every generalized star is an integral sum graph. \hfill $\Box$} 
\end{theorem}

\subsection{Banana trees are integral sum graphs}

Here, we present the result that Banana trees are integral sum graphs that is proved in \cite{vn11} by Vilfred and Nicholas.

\begin{dfn}{\rm\cite{vn11}}\quad A {\em banana tree} is a family of stars with a new vertex adjoined to one end vertex of each star.
\end{dfn}

\begin{dfn}{\rm\cite{c96}}\quad Let $G(V,E)$ be any graph. An edge $u v$ of $G$ is said to be {\em $f$-proper} if $f(u) + f(v)$ = $f(w)$ for some $w \in V(G).$
\end{dfn}

It is easy to prove that given labeling $f$ of a graph $G$ is an integral sum labeling of $G$ if and only if all edges of $G$ are $f$-proper and all edges of $G^c$ are not $f$-proper. 

\begin{theorem} [{ \cite[Theorem 2.4]{vn11}}] \label{6.2.7} \normalfont\quad 
	{\rm Every banana tree $T$ is an integral sum graph. } 
\end{theorem}

\begin{proof}\quad Let $T$ be a banana tree corresponding to the family of stars $\{K_{1,n_1}$, $K_{1,n_2}$, ..., $K_{1,n_t}\}$, $t\in\mathbb{N}$. Let $v_i$ denote the central vertex and $u_{i,j}$, $j$ = $1,2,...,n_i$ denote the end vertices of the $i^{th}$ star $K_{1,n_i}$ and $w$ be the new vertex joining one vertex $u_{i, 1}$ of each star, $i$ = $1, 2, ... , t$. 
	
	If $n_i \leq 2$ for all $i$ = $1, 2, . . . , t$, then $T$ is a generalized star and hence is an integral sum graph using Theorem \ref{gs}. Let $n_i > 2$ for at least one $i$, $1 \leq i \leq t$. Without loss of generality, let $n_1 \leq n_2 \leq \cdots \leq n_t$.  This implies, $n_t \geq 3$. Also, a banana tree with $t$ = 1 is actually a general star which is an integral sum graph, using Theorem \ref{gs}. 
	
	Let $y > x > 0$, $t \geq 2$, $n_0$ = 0 and $t,x,y\in\mathbb{N}$. Define a vertex labeling $f$ on $T$ as follows:
	
	\vspace{.1cm}	
	$f(w) = x$; $f(u_{1, 1})$ = $y$; $f(v_1)$ = $f(w) + f(u_{1, 1})$;
	
	\vspace{.1cm}		
	$f(u_{1, j+1})$ = $f(u_{1, 1}) + j f(v_1)$, ~$j$ = 1, 2, . . . , $n_1$-1.
	
	\vspace{.1cm}
	\noindent		
	For $i$ = 1,2,...,$t$-2,
	
	\vspace{.1cm}		
	$f(u_{i+1, 1})$ = $f(v_i) + f(u_{i, n_i})$ and 
	$f(v_{i+1})$ = $f(w) + f(u_{i+1, 1})$, 
	
	\vspace{.1cm}		
	$f(u_{i+1, j+1})$ = $f(u_{i+1, 1}) + j f(v_{i+1})$, $j$ = 1, 2, . . . , $n_{i+1}$-1.
	
	\vspace{.1cm}		
	$f(u_{t, 1})$ = $f(v_{t-1}) + f(u_{t-1, n_{t-1}})$; 
	
	\vspace{.1cm}		
	$f(v_t)$ = $f(u_{1, 1}) - f(u_{t, 1})$; 
	
	\vspace{.1cm}		
	$f(u_{t, 2})$ = $f(w) + f(u_{t, 1})$;
	
	\vspace{.1cm}		
	$f(u_{t, j+2})$ = $f(u_{t, j+1}) - f(v_t)$, $j$ = 1, 2, . . . , $n_{t}$-2.
	
	\vspace{.1cm}		
	Now, only $v_t$ has negative value and all other vertices have positive values.
	
	\vspace{.1cm}	
	Figure 24 shows the banana tree corresponding to the family of stars $\{K_{1,3}, K_{1,5}$, $K_{1,6}\}$ with an integral sum labeling $f$ given in the proof of the above theorem with $x$ = 1 and $y$ = 2. Here, $t$ = 3, $n_0$ = 0, $n_1$ = 3, $n_2$ = 5 and $n_3$ = 6.
	
\medskip

\begin{center}
	\definecolor{myblue}{RGB}{80,80,160}
	\begin{tikzpicture}[scale =1.2]
		
		\node (v) at (-4,-10.5) [circle,draw,scale=0.7] {1};
		
		\node (a7) at (-7.6,-12)  [circle,draw,scale=0.7] {2};
		\node (a8) at (-6.7,-12)  [circle,draw,scale=0.7]{5};
		\node (a9) at (-5.8,-12)  [circle,draw,scale=0.7]{8};
		\node (h13) at (-6.7,-13.5)  [circle,draw,scale=0.7] {3};
		
		\draw (v) -- (a7);
		\draw (h13) -- (a9);
		\draw (h13) -- (a8);
		\draw (h13) -- (a7);
		
		\node (h13) at (-7.2,-14.2) [label=0:$K_{1, 3}$] {};
		
		\node (a10) at (-4.9,-12)  [circle,draw,scale=0.6] {11};
		\node (a11) at (-4,-12)  [circle,draw,scale=0.6]{23};
		\node (a12) at (-3.1,-12)  [circle,draw,scale=0.6]{35};
		\node (a13) at (-2.2,-12)  [circle,draw,scale=0.6]{47};
		\node (a14) at (-1.3,-12)  [circle,draw,scale=0.6]{59};
		\node (h21) at (-3.1,-13.5)  [circle,draw,scale=0.6] {12};
		
		\draw (h21) -- (a10);
		\draw (h21) -- (a11);
		\draw (h21) -- (a12);
		\draw (h21) -- (a13);
		\draw (h21) -- (a14);
		\draw (v) -- (a10);
		
		\node (h21) at (-3.7,-14.2) [label=0:$K_{1,5}$] {};
		
		\node (a15) at (-0.4,-12)  [circle,draw,scale=0.6] {71};
		\node (a16) at (.5,-12)  [circle,draw,scale=0.6] {72};
		\node (a17) at (1.4,-12)  [circle,draw,scale=0.5] {141};
		\node (a18) at (2.3,-12)  [circle,draw,scale=0.5] {210};
		\node (a19) at (3.2,-12)  [circle,draw,scale=0.5] {279};
		\node (a20) at (4.1,-12)  [circle,draw,scale=0.5] {348};
		\node (h22) at (1.8,-13.5)  [circle,draw,scale=0.5]  {-69};
		
		\draw (h22) -- (a15);
		\draw (h22) -- (a16);
		\draw (h22) -- (a17);
		\draw (h22) -- (a18);
		\draw (h22) -- (a19);
		\draw (h22) -- (a20);
		\draw (v) -- (a15);
		
		\node (h22) at (1.3,-14.2) [label=0:$K_{1,6}$] {};	
	\end{tikzpicture}
	
	\medskip
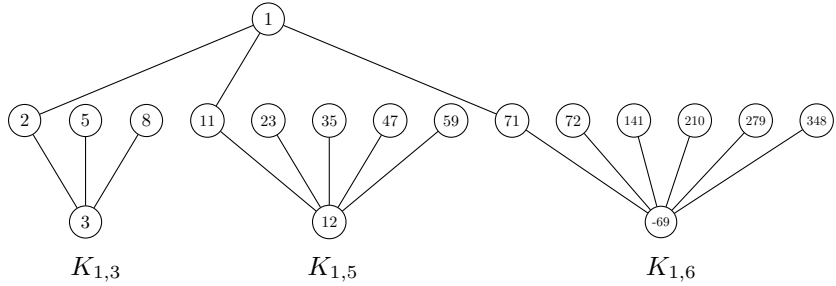
	Fig. 24. Banana tree corresponding to $K_{1, 3}, K_{1, 5}, K_{1, 6}$ with integral sum labeling  \hspace{1cm}
\end{center}
	
	\vspace{.1cm}	
	\noindent
	{\it Claim.} The labels of all the vertices in $T-w$ are of the form $ax+by$ where $a\in \{b,b-1\}$.
	
	We have, $f(w)$ = $x > 0$; $f(u_{1, 1})$ = $y > x$. Now let us calculate the values of $f(u_{i, j})$ and $f(v_i)$, in terms of $x$ and $y$, for all possible values of $i$ and $j$,  $1 \leq j \leq n_i$,~ $1 \leq i \leq t$ and $t \geq 2$.
	
	All the values are calculated from the recurrence relations of $f(u_{i, j})$ and $f(v_i)$ defined earlier.
	
	\noindent
	When $i$ = 1, $f(u_{1, 1})$ = $y$ and $f(v_1)$ = $x+f(u_{1, 1})$ = $x+y$ and 
	
	$f(u_{1, j+1})$ = $f(u_{1, 1})+ j f(v_1)$ = $y + j(x+y) = (j+1)(x+y)–x$, $j$ = 1,2,...,$n_1$-1.
	
	\noindent
	When $i$ = 2 (and $t > 2$),
	
	$f(u_{2, 1})$ = $f(v_1) + f(u_{1, n_1})$ = $(n_1+1)(x+y)-x$,
	
	$f(v_2)$ = $x+f(u_{2, 1})$ = $(n_1+1)(x+y)$ and 
	
	$f(u_{2, j+1})$ = $f(u_{2, 1})+ j f(v_2)$ = $(n_1+1)(j+1)(x+y)-x$, ~ $j$ = 1,2,...,$n_2$-1.
	
	\noindent
	When $i$ = 3 (and $t > 3$),
	
	$f(u_{3, 1})$ = $f(v_2) + f(u_{2, n_2})$ = $(n_1+1)(n_2+1)(x+y)-x$,
	
	$f(v_3)$ = $x+f(u_{3, 1})$ = $(n_1+1)(n_2+1)(x+y)$ and 
	
	$f(u_{3, j+1})$ = $f(u_{3, 1})+ j f(v_3)$ = $(n_1+1)(n_2+1)(j+1)(x+y)-x$, ~ $j$ = 1,2,...,$n_3$-1.

	\noindent
	When $i$ = k (and $t > k \geq 2$),
	
	$f(u_{k, 1})$ = $(n_1+1)(n_2+1) . . . (n_{k-1}+1)(x+y)-x$,
	
	$f(v_k)$ = $(n_1+1)(n_2+1) . . . (n_{k-1}+1)(x+y)$ and 
	
	$f(u_{k, j+1})$ = $(n_1+1)(n_2+1) . . . (n_{k-1}+1)(j+1)(x+y)-x$, $j$ = 1,2,...,$n_k$-1. Here $k$ = 1,2,...,$t$-1.
	
	\noindent
	In particular when $k$ = $t-1$ and $t \geq 3$, we get,
	
	$f(u_{t-1, 1})$ = $(n_1+1)(n_2+1) . . . (n_{t-2}+1)(x+y)-x$,
	
	$f(v_{t-1})$ = $(n_1+1)(n_2+1) . . . (n_{t-2}+1)(x+y)$ and 
	
	$f(u_{t-1, j+1})$ = $(n_1+1)(n_2+1) . . . (n_{t-2}+1)(j+1)(x+y)-x$, ~ $j$ = 1,2,...,$n_{t-1}$-1. 
	
	\noindent
	When $i$ = t and $t \geq 2$,
	
	$f(u_{t, 1})$ = $f(v_{t-1}) + f(u_{t-1, n_{t-1}})$ = $(n_1+1)(n_2+1) . . . (n_{t-1}+1)(x+y)-x$,
	
	$f(v_t)$ = $f(u_{1, 1}) - f(u_{t, 1})$ = $-((n_1+1)(n_2+1) . . . (n_{t-1}+1)-1)(x+y)$ and 
	
	$f(u_{t, 2})$ = $x + f(u_{t, 1})$ = $(n_1+1)(n_2+1) . . . (n_{t-1}+1)(x+y)$ and
	\\
	for $j$ = 2,3,...,$n_t$-1, $f(u_{t, j+1})$ = $f(u_{t, j}) - f(v_{t})$ = $f(u_{t, 2}) - (j-1) f(v_{t})$ 
	
	\hfill = $((n_1+1)(n_2+1) . . . (n_{t-1}+1)j - (j-1))(x+y)$. 
	
	Thus from the above, it is clear that the vertex labels of $T$, except $w$, are of the form either $a(x+y)$ or $a(x+y)-x$, $a\in\mathbb{Z}$. Hence the claim is true.
	
	\medskip
	Now, partition the vertex set $V(T)$ in to $A$, $B$ and $C$ such that $V(T)$ = $A \cup B \cup C$ where $A$ = $\{w\}$, $B$ = $\{u\in V(T): f(u) = (b-1)x+by \}$ = $\{u_{t, 1}\} \cup \{u_{i, j}:$ $i = 1,2,\cdots,t-1$ and $j$ = $1,2,\cdots,n_i\}$ and $C$ = $\{u\in V(T): f(u) = bx+by \}$ = $\{v_i: i = 1,2,\cdots,t\} \cup \{u_{t, j}:$ $j = 2,3,\cdots,n_t$ and $t \geq 2\}$.
	
	Thus, the set $B$ contains all the end vertices of the stars $K_{1,n_1}$, $K_{1, n_2}$, . . . , $K_{1, n_{t-1}}$ and the first end vertex $u_{t,1}$ of the $t^{th}$ star $K_{1,n_t}$ while the set $C$ contains the central vertices of all the stars and end vertices $u_{t,j}$'s of the $t^{th}$ star $K_{1,n_t}$, except $u_{t,1}$. Clearly, all the edges of $T$ are $f$-proper. Let $f(V)$ = $\{ f(u): u\in V(T) \}$. 
	
	It remains to show that every edge $e$ = $uv\in E(T^c)$ is not $f$-proper. Observe that $v_t$ is the only vertex of $T$ with negative label.
	
	If $e$ = $u v_t\in E(T^c)$, then $u\notin V(K_{1, n_t})$ and hence $f(u) < |f(v_t)|$. This implies that $f(u)+f(v_t) < 0$ and hence $e$ = $u v_t\in E(T^c)$ is not $f$-proper. 
	
	Now, consider any edge $e$ = $uv\in E(T^c)$ where both $u$ and $v$ have positive labels. The following four cases arise.
	
	\medskip	
	\noindent
	{\bf Case 1} ~ $u\in A$ and $v\in B \cup C$.
	
	In this case, $u = w$ and if $v\in C$, then $f(u)+f(v)$ = $(a+1)x+ay$ for some integer $a$ and hence $f(u)+f(v)\notin f(V)$.
	
	If $v\in B$, then $f(v) = (b-1)x+by$ for some integer $b$. This implies, $f(u)+f(v)$ = $bx+by > 0$. If $bx+by$ = $f(z)$, then $z\in C \setminus \{v_t\}$ since $f(v_t) < 0$ and so $z = v_i$ or $u_{t, j}$, $i$ = 1, 2, $\dots$, $t-1$ and $j$ = 2, 3, $\dots$, $n_t$.
	
	If $z$ is the central vertex of any star, say $z$ = $v_i$, $1 \leq i \leq t-1$, then $v$ = $u_{i, 1}$ which is a contradiction to $e = uv\in E(T)$. On the other hand, if $z$ = $u_{t, j}$ for any $j$, $2 \leq j \leq n_t$, then $v$ = $u_{t, 1}$ is the only possibility since $f(w)+f(u_{t, 1})$ = $f(u_{t, 2})$. This implies that $uv\in E(T)$ which is a contradiction to our assumption that $uv\in E(T^c)$. Thus, in both cases we get, $f(u)+f(v)\notin f(V)$.
	
	\medskip	
	\noindent
	{\bf Case 2} ~ $u,v\in B.$
	
	In this case, $f(u)+f(v)$ = $(a-2)x+ay$ for some integer $a$ and hence $f(u)+f(v)\notin f(V)$.
	
	\vspace{.1cm}	
	\noindent
	{\bf Case 3} ~ $u,v\in C$.
	
	Here we have to consider the following three sub cases.
	
	\vspace{.1cm}	
	\noindent
	{\bf Case 3.1} ~  $u$ and $v$ are end vertices ($\neq u_{t, 1}$) of the $t^{th}$ star.	
	
	In this case, obviously $f(u)+f(v)$ is not a label of the end vertex of the $t^{th}$ star. This implies that $f(u)+f(v)\notin f(V)$.
	
	\medskip	
	\noindent
	{\bf Case 3.2} ~  $u$ and $v$ are central vertices of two different stars.
	
	In this case, let $u = v_r$ and $v = v_s$, $r < s$, say. Now, by construction $f(v_{r+1})$ = $f(u_{r+1, 1}) + f(w)$ = $f(u_{r,n_r}) + f(v_r) + f(w)$ = $f(u_{r, 1}) + n_r f(v_r) + f(w)$ = $(n_{r}+1)f(v_r)$ since $f(u_{r, 1}) + f(w)$ = $f(v_r)$, $r$ = 1, 2, . . . , $t-1$. 
	
	If $f(v_r)+f(v_s)$ = $f(v_k)$ for some $k$, $r < s < k \leq t-1$, then by dividing on both sides by $f(v_s)$, we see that the left side is not an integer whereas the right side is an integer which is not possible and hence $f(u)+f(v)\notin f(V)$.
	
	Also, $f(v_r)+f(v_s) \neq f(u_{t, j})$ for any $j$, $2 \leq j \leq n_t$ and $r < s < t-1$. Thus, in this case also $f(u)+f(v)\notin f(V)$.
	
	\medskip	
	\noindent
	{\bf Case 3.3} ~ $u = v_i$ and $v = u_{t, j}$, $1 \leq i \leq t-1$ and $j$ = 2,3,...,$n_t$.
	
	In this case, $f(v_i) < |f(v_t)| < f(u_{t, j})$ for every $i$ = 1,2,...,$t$-1 and $j$ = $1,2,...,n_t$. This implies that $f(u)+f(v)\notin f(V)$.
	
	\medskip	
	\noindent
	{\bf Case 4} ~ $u\in B$ and $v\in C$.
	
	In this case either $u$ and $v$ are vertices of different stars or of the same star. Here, it is easy to see that $f(u)+f(v)\notin f(V)$ in the following sub cases:
	
	i)~~ $u = u_{t, 1}$ and $v = u_{t, j}$, $j$ = 2, 3, . . . , $n_t$;
	
	ii)~ $u = u_{t, 1}$ and $v = v_i$, $i$ = 1, 2, . . . , $t-1$;
	
	iii) $u = u_{i, j}$ and $v = v_t$, $i$ = 1, 2, . . . , $t-1$ and $j$ = 1, 2, . . . , $n_i$ and
	
	iv) $u = u_{i, j}$ and $v = u_{t, k}$, $i$ = $1,2,...,t-1$, $j$ = $1,2,...,n_i$ and $k$ = $1,2,...,n_t$.
	
	Thus, in this case also, $e = uv\in E(T^c)$ is not $f$-proper.	
	
	Hence, in all possible cases, every edge $e = uv\in E(T^c)$ is not $f$-proper. This implies that $T$ is an integral sum graph. 	 
\end{proof}

\subsection{Dutch windmill, triangular book, fan with a handle are integral sum}

In \cite{vn10}, it is proved that amalgamation of graphs as given in Theorem \ref{amal} is an integral sum graph. A sufficient condition is obtained for the graph $G * v$ to be an integral sum graph when $G$ is an integral sum graph when $v \notin V(G)$. And using this condition, it is proved that for $n \geq 5$, Dutch $m$-windmill $K_3^{(m)}$ = $m.K_2 * v$, triangular book with bookmark, fan graph $P_n * K_1$ and $G * v$ are integral sum graphs where graph $G$ is the union of stars. These results are presented in this subsection. 

\begin{dfn} \cite{lw88}\quad Given a collection of graphs $G_1$, $G_2$, $\cdots$, $G_k$ and some fixed vertex $v_i\in V(G_i)$, $i$ = $1,2,...,k$, define ${\it Amal \{ (G_i^{(k)}, v_i) \}}$, the {\em amalgamation of graphs} $G_1,G_2,\cdots,G_k$, as the graph obtained by taking the union of the $G_i$'s and identifying $v_1,v_2,\cdots,v_k.$
\end{dfn}

\begin{theorem} [{ \cite[Theorem 6]{vn10}}] {\rm\quad \cite{vn10} \label{amal} Let $f_i$ be an integral sum labeling of $G_i$ with $\Delta (G_i)$ = $|V(G_i)| -1$ and $f_i(x_i)$ = 0, $x_i\in V(G_i)$ for every $i$ = $1, 2, \dots, n$. Then, $Amal{(G_i^{(n)}, x_i)}$ is an integral sum graph, $i$ = $1, 2, \dots, n$.}
\end{theorem}

\begin{proof}\quad Let $G$ = $Amal{(G_i^{(n)}, x_i)}$, $1 \leq i \leq n$. Here the amalgamation of the graphs $G_1$, $G_2$, $\cdots$, $G_n$ is done by identifying the vertices $x_1$, $x_2$, $\cdots$, $x_n$ where $f_i(x_i)$ = $0$, $i$ = $1, 2, \dots, n$. Given, $f_i$ is an integral sum labeling of $G_i$ with $\Delta (G_i)$ = $|V(G_i)|$-1 and $f_i(x_i)$ = 0, $x_i\in V(G_i)$ and so each $G_i$ is a connected graph, $i$ = $1, 2, \dots, n$. In $G$, let $x$ = $x_1$ = $x_2$ = $\cdots$ = $x_n$. Define a vertex labeling function $h$ on $G$ as follows: $h(x)$ = 0; $h(v)$ = $f_1(v)$ if $v \in V(G_1)$ and for $i$ = $2, 3, \dots, n$, $h(v)$ = $(M_{i-1} + 1)f_i(v)$ if $v\in V(G_i)$ where $M_{i-1}$ = $max ~\{ 2  |h(v)| : v \in V(G_{i-1}) \}$. 
	
	If $u$ and $v$ are two different vertices belonging to one $G_i$, then $u v \in E(G_i)$ if and only if there exists $w \in V(G_i)$ such that $f_i(u) + f_i(v)$ = $f_i(w)$ if and only if $(M_{i-1} + 1)f_i(u)$ + $(M_{i-1} + 1)f_i(v)$ = $(M_{i-1} + 1)f_i(w)$ if and only if $h(u)$ + $h(v)$ = $h(w)$. On the other hand, let $u$ and $v$ be two vertices belonging to two different $G_i$'s, say $G_i$ and $G_j$, $1 \leq i < j \leq n$. It is easy to prove that $h(u)$ + $h(v) \neq h(w)$ for any $w\in V(G_k) \setminus \{ x \}$ by considering the following cases of $k:$  (i)  $1 \leq k < i$, (ii)  $k$ = $i$, (iii)  $i < k < j$, (iv) $k$ = $j$ and (v) $j < k \leq n$. Thus $h$ is an integral sum labeling on $G$. Hence the result is proved.
\end{proof}

\subsection{On $\mathbb{Z}$-sum graphs $G * v$ when $G$ is an $\mathbb{Z}$-sum graph and $v \notin V(G)$}

A sufficient condition is obtained for the graph $G * v$ to be an integral sum graph when $G$ is an integral sum graph and $v \notin V(G)$. Using the above condition, it is proved that for $n \geq 5$, fan graph $F_n$ = $P_n * K_1$, Dutch $m$-windmill $K_3^{(m)}$ = $m.K_2 * v$ and $G * v$ are integral sum graphs where graph $G$ is the union of stars.

\begin{dfn} \cite{lw88, vn10}\quad A {\em fan graph} $F_{n-1}$ is the graph obtained by taking $n-3$ concurrent chords at a vertex in a cycle $C_n$, $n \geq 3$. And it can be described as $F_n$ = $P_n * K_1$ where $P_n$ is a path on $n$ vertices, $n \geq 2$. The vertex at which all the $n-3$ chords are concurrent is called the {\em apex vertex}. 
\end{dfn}

\begin{theorem} [{ \cite[Theorem 7]{vn10}}] \label{sc} {\rm [{\em Sufficient Condition for $G$ $*$ $v$ to be an integral sum graph}]\quad Let $G$ be an integral sum graph and $v \notin V(G)$ be a new vertex. Suppose (i)  $\Delta (G)$ $<$ $|V(G)|$ $-$ $1$ and (ii) there exists at least one integral sum labeling $f$ on $G$ with $f(x) \neq -f(y)$ for all vertices $x,y$ of $G$. Then, $G * v$ is also an integral sum graph.} 
\end{theorem}

\begin{proof}\quad Let $f$ be an integral sum labeling of $G$ with $f(x) \neq -f(y)$ for every $x,y\in V(G)$. Since $\Delta (G) < |V(G)| - 1$, using Theorem \ref{a0}, $f(x) \neq 0$ for every $x \in V(G)$. Define a vertex labeling $g$ on $G * v$ such that $g(v)$ = 0 and $g(x)$ = $f(x)$ for every $x\in V(G)$. Given that for every $x,y\in V(G)$, $f(x) \neq -f(y)$ and hence $g(x) + g(y) \neq g(v) = 0$. Also, for every $x,y\in V(G*v \setminus \{v\})$, $xy\in E(G)$ if and only if there exists $w\in V(G)$ such that $f(x) + f(y)$ = $f(w) \neq 0$ if and only if $g(x) + g(y)$ = $g(w) \neq g(v)$, $x,y,w\in V(G*v \setminus \{v\})$ if and only if $xy\in E(G * v)$, $x,y \neq v$. Thus $g$ is an integral sum labeling on $G * v$ and hence the theorem is proved.
\end{proof} 

\begin{note} [{ \cite[Note]{vn10}}] \quad {\rm For $n$ = 2,3,4, the possible integral sum labeling of $P_n$ are $S$ = $\{ 0, a \}$, $\{ 0, a, b:$ $a+b \neq  0 \}$ and $\{a, -a, -2a, -3a\}$, $a \in\mathbb{N}$. Thus, in these cases $S$ contains either $0$ or both $a$ and $-a,a\in\mathbb{N}$. Next two theorems are proved using Theorem \ref{sc}. }
\end{note}

\begin{cor}  [{ \cite[Corollary 7.1]{vn10}}] \label{fan} {\rm \quad For $n \geq 5$, fan graph $F_n$ = $P_n * K_1$ is an integral sum graph.} 
\end{cor}

\begin{proof}\quad Let $G$ = $F_n$ = $P_n * K_1$ and $P_n$ = $u_1 u_2 \cdots u_n$ be the path on $n$ vertices, $n \geq 2$. Define a vertex  labeling $f$ on $P_n$ such that 
	
	$f(u_1)$ = 1, $f(u_2)$ = -1 and 
	
	$f(u_i)$ = $f(u_{i-2}) - f(u_{i-1})$, $i$ = $3, 4, \dots, n$. 
	
	\medskip
	Clearly, $f$ is an integral sum labeling on $P_n$. Also, for $n \geq 5$, $\Delta (P_n)$ = $2 < |V(P_n)|-1$, $f(x) \neq 0$, $f(x) \neq -f(y)$ for every $x,y \in V(P_n)$ and $f$ is an integral sum labeling on $P_n$ and hence applying Theorem \ref{sc}, fan graph $F_n$ = $P_n * K_1$ is an integral sum graph with the integral sum labeling $g$ defined on $P_n * K_1$ as $g(K_1)$ = 0 and $g(u_i)$ = $f(u_i)$, $i$ = $1, 2, \dots, n$. Figures 25 and 26 show graph $F_6$ and its integral sum labeling.
\end{proof}

\begin{center}
	\begin{tikzpicture}[scale =1.2]
		
		\node (c0) at (10,1.5)  [circle,draw,scale=0.7]{$w$};
		\node (c1) at (8,0)  [circle,draw,scale=0.6]{$u_1$};
		\node (c2) at (9,0)  [circle,draw,scale=0.6] {$u_2$};
		\node (c3) at (10,0)  [circle,draw,scale=0.6] {$u_3$};
		\node (c4) at (11,0)  [circle,draw,scale=0.6] {$u_4$};
		\node (c5) at (12,0)  [circle,draw,scale=0.6]{$u_5$};
		\node (c6) at (13,0)  [circle,draw,scale=0.6]{$u_6$};
		
		\draw (c0) -- (c1);
		\draw (c0) -- (c2);
		\draw (c0) -- (c3);
		\draw (c0) -- (c4);
		\draw (c0) -- (c5);
		\draw (c0) -- (c6);
		
		\draw (c1) -- (c2);
		\draw (c2) -- (c3);
		\draw (c3) -- (c4);
		\draw (c4) -- (c5);
		\draw (c5) -- (c6);
		
		\node (c0) at (17,1.5)  [circle,draw,scale=0.7]{$0$};
		\node (c1) at (15,0)  [circle,draw,scale=0.7]{$1$};
		\node (c2) at (16,0)  [circle,draw,scale=0.5] {$-1$};
		\node (c3) at (17,0)  [circle,draw,scale=0.7] {$2$};
		\node (c4) at (18,0)  [circle,draw,scale=0.5] {$-3$};
		\node (c5) at (19,0)  [circle,draw,scale=0.7]{$5$};
		\node (c6) at (20,0)  [circle,draw,scale=0.5]{$-8$};
		
		\draw (c0) -- (c1);
		\draw (c0) -- (c2);
		\draw (c0) -- (c3);
		\draw (c0) -- (c4);
		\draw (c0) -- (c5);
		\draw (c0) -- (c6);
		
		\draw (c1) -- (c2);
		\draw (c2) -- (c3);
		\draw (c3) -- (c4);
		\draw (c4) -- (c5);
		\draw (c5) -- (c6);
		
	\end{tikzpicture}
	
	\medskip
	Fig. 25. Fan graph $F_6$ = $P_6 * K_1$ \hspace{3cm} Fig. 26. Integral sum labeling of $F_6$
	
\end{center}

\begin{cor}  [{ \cite[Corollary 7.2]{vn10}}] \label{wind} {\rm \quad Dutch $m$-windmill $K_3^{(m)}$ = $mK_2 * v$ is an integral sum graph where $v \notin V(mK_2)$. }
\end{cor}
\begin{proof}\quad Harary \cite{h94} proved that matching are integral sum graph. Let matching $mK_2$ = $G(V,E)$ where $E(G)$ = $\{ u_i v_i : i = 1, 2, \dots, m \}$ and $f$ be the vertex labeling defined on $G$ as $f(u_1)$ = 1, $f(v_1)$ = 3; $f(u_2)$ = 4, $f(v_2)$ = 8; $f(u_3)$ = -11, $f(v_3)$ = 12 and $f(v_i)$ = $2f(v_{i-1});$ $f(u_i)$ = $f(u_{i-1}) - f(v_i)$, $i$ = $4, 5, \dots, m$. Clearly, $f$ is an integral sum labeling of $mK_2$ and satisfies the condition (ii) of Theorem \ref{sc}. Hence, using Theorem \ref{sc}, graph $mK_2 * v$ = $G * v$ is an integral sum graph with the integral sum labeling $g$ defined on $G * v$ as $g(v)$ = 0 and $g(u)$ = $f(u)$, for all $u \in V(G)$. Figures 27 and 28 show $5$-windmill $K_3^{(5)}$ = $5K_2 * v$ and its integral sum labeling.
\end{proof}

\begin{cor}  [{ \cite[Corollary 7.3]{vn10}}]  \label{st} {\rm \quad  Let $G$ = $K_{1, n_1} \cup K_{1, n_2} \cup \cdots \cup K_{1, n_t}$, $n_1,n_2,...,n_t,t\in\mathbb{N}$. Then $G * v$ is an integral sum graph.}
\end{cor}
\begin{proof}\quad In Theorem \ref{su}, it is proved that $G$ is an integral sum graph. Moreover, the integral sum labeling of $G$ given in the proof of the theorem satisfies all the conditions of Theorem \ref{sc} and thereby graph $G * v$ is an integral sum graph.
	
	Figure 29 shows graph $G$ = $K_{1, 3} \cup K_{1, 5} \cup K_{1, 6}$ with their vertices and Figure 30 shows graph $G * v$ and its integral sum labeling as given in the proof.
\end{proof}

\medskip	
\begin{center}
	\begin{tikzpicture}[scale=1.2]
	
	\node (w0) at (10,0.4)  [circle,draw,scale=0.8]{$w$};
	\node (a1) at (9.5,2)  [circle,draw,scale=0.7]{$u_1$};
	\node (b1) at (10.5,2)  [circle,draw,scale=0.7] {$v_1$};
	\node (a2) at (11.5,1.5)  [circle,draw,scale=0.7]{$u_2$};
	\node (b2) at (12,0.75)  [circle,draw,scale=0.7] {$v_2$};
	\node (a3) at (11.5,-0.5)  [circle,draw,scale=0.7]{$u_3$};
	\node (b3) at (10.5,-1)  [circle,draw,scale=0.7] {$v_3$};
	\node (a4) at (9.5,-1)  [circle,draw,scale=0.7]{$u_4$};
	\node (b4) at (8.5,-0.5)  [circle,draw,scale=0.7] {$v_4$};
	\node (a5) at (8,0.75)  [circle,draw,scale=0.7]{$u_5$};
	\node (b5) at (8.5,1.5)  [circle,draw,scale=0.7] {$v_5$};
	
	\draw (a1) -- (b1);
	\draw (a2) -- (b2);
	\draw (a3) -- (b3);
	\draw (a4) -- (b4);
	\draw (a5) -- (b5);
	
	\draw (w0) -- (a1);
	\draw (w0) -- (b1);
	\draw (w0) -- (a2);
	\draw (w0) -- (b2);
	\draw (w0) -- (a3);
	\draw (w0) -- (b3);
	\draw (w0) -- (a4);
	\draw (w0) -- (b4);
	\draw (w0) -- (a5);
	\draw (w0) -- (b5);
	
	\node (w0) at (16,0.4)  [circle,draw,scale=0.8]{$0$};
	\node (a1) at (15.5,2)  [circle,draw,scale=0.8]{$1$};
	\node (b1) at (16.5,2)  [circle,draw,scale=0.8] {$3$};
	\node (a2) at (17.5,1.5)  [circle,draw,scale=0.8]{$4$};
	\node (b2) at (18,0.75)  [circle,draw,scale=0.8] {$8$};
	\node (a3) at (17.5,-0.5)  [circle,draw,scale=0.55]{$-11$};
	\node (b3) at (16.5,-1)  [circle,draw,scale=0.7] {$12$};
	\node (a4) at (15.5,-1)  [circle,draw,scale=0.55]{$-35$};
	\node (b4) at (14.5,-0.5)  [circle,draw,scale=0.7] {$24$};
	\node (a5) at (14,0.75)  [circle,draw,scale=0.55]{$-83$};
	\node (b5) at (14.5,1.5)  [circle,draw,scale=0.7] {$48$};
	
	\draw (a1) -- (b1);
	\draw (a2) -- (b2);
	\draw (a3) -- (b3);
	\draw (a4) -- (b4);
	\draw (a5) -- (b5);
	
	\draw (w0) -- (a1);
	\draw (w0) -- (b1);
	\draw (w0) -- (a2);
	\draw (w0) -- (b2);
	\draw (w0) -- (a3);
	\draw (w0) -- (b3);
	\draw (w0) -- (a4);
	\draw (w0) -- (b4);
	\draw (w0) -- (a5);
	\draw (w0) -- (b5);
	
	\end{tikzpicture}
	
	\medskip
	 Fig. 27. Windmill $K_3^{(5)}$ = $5K_2 * K_1$ \hspace{1cm} Fig. 28. Integral sum labeling of $K_3^{(5)}$  	
\end{center}

\medskip
\begin{center}
	\definecolor{myblue}{RGB}{80,80,160}
	\begin{tikzpicture}[scale =1.2]
	
	\node (a7) at (-7.6,-12)  [circle,draw,scale=0.6] {$u_{1,1}$};
	\node (a8) at (-6.7,-12)  [circle,draw,scale=0.6]{$u_{1,2}$};
	\node (a9) at (-5.8,-12)  [circle,draw,scale=0.6]{$u_{1,3}$};
	\node (h13) at (-6.7,-14)  [circle,draw,scale=0.7] {$v_{1}$};
	
	\draw (h13) -- (a7);
	\draw (h13) -- (a8);
	\draw (h13) -- (a9);
	
	\node (h13) at (-7.1,-14.5) [label=0:$K_{1, 3}$] {};
	
	\node (a10) at (-4.9,-12)  [circle,draw,scale=0.6] {$u_{2,1}$};
	\node (a11) at (-4,-12)  [circle,draw,scale=0.6]{$u_{2,2}$};
	\node (a12) at (-3.1,-12)  [circle,draw,scale=0.6]{$u_{2,3}$};
	\node (a13) at (-2.2,-12)  [circle,draw,scale=0.6]{$u_{2,4}$};
	\node (a14) at (-1.3,-12)  [circle,draw,scale=0.6]{$u_{2,5}$};
	\node (h21) at (-3.1,-14)  [circle,draw,scale=0.7] {$v_{2}$};
	
	\%draw (v) -- (a13);
	\draw (h21) -- (a10);
	\draw (h21) -- (a11);
	\draw (h21) -- (a12);
	\draw (h21) -- (a13);
	\draw (h21) -- (a14);
	
	\node (h21) at (-3.5,-14.5) [label=0:$K_{1,5}$] {};
	
	\node (a15) at (-0.4,-12)  [circle,draw,scale=0.6] {$u_{3,1}$};
	\node (a16) at (.5,-12)  [circle,draw,scale=0.6] {$u_{3,2}$};
	\node (a17) at (1.4,-12)  [circle,draw,scale=0.6] {$u_{3,3}$};
	\node (a18) at (2.3,-12)  [circle,draw,scale=0.6] {$u_{3,4}$};
	\node (a19) at (3.2,-12)  [circle,draw,scale=0.6] {$u_{3,5}$};
	\node (a20) at (4.,-12)  [circle,draw,scale=0.6] {$u_{3,6}$};
	\node (h22) at (1.8,-14)  [circle,draw,scale=0.7]  {$v_{3}$};
	
	\draw (h22) -- (a15);
	\draw (h22) -- (a16);
	\draw (h22) -- (a17);
	\draw (h22) -- (a18);
	\draw (h22) -- (a19);
	\draw (h22) -- (a20);
	
	\node (h22) at (1.3,-14.5) [label=0:$K_{1,6}$] {};
	
	\end{tikzpicture}
	
	\medskip
	Fig. 29. Graph $G = K_{1, 3} \cup K_{1, 5} \cup K_{1, 6}$  
	
\end{center}

\medskip

\begin{center}
	\definecolor{myblue}{RGB}{80,80,160}
	\begin{tikzpicture}[scale =1.2]
	
	\node (v) at (-1,-10) [circle,draw,scale=0.7] {0};
	
	\node (a7) at (-8.6,-12)  [circle,draw,scale=0.7] {1};
	\node (a8) at (-7.7,-12)  [circle,draw,scale=0.7]{3};
	\node (a9) at (-6.8,-12)  [circle,draw,scale=0.7]{5};
	\node (h13) at (-6.7,-14)  [circle,draw,scale=0.7] {2};
	
\draw (v) -- (a9);
\draw (v) -- (a8);
\draw (v) -- (a7);
\draw (h13) -- (a9);
\draw (h13) -- (a8);
\draw (h13) -- (a7);
\draw[] [blue](-1.2,-10) to [out=160,in=65] (-6.6,-13.88);
	
\node (h13) at (-7.2,-14.5) [label=0:$K_{1, 3}$] {};
	
	\node (a10) at (-5.9,-12)  [circle,draw,scale=0.7] {7};
	\node (a11) at (-5,-12)  [circle,draw,scale=0.6]{21};
	\node (a12) at (-4.1,-12)  [circle,draw,scale=0.6]{35};
	\node (a13) at (-3.2,-12)  [circle,draw,scale=0.6]{49};
	\node (a14) at (-2.3,-12)  [circle,draw,scale=0.6]{63};
	\node (h21) at (-3.1,-14)  [circle,draw,scale=0.6] {14};
	
	\draw (h21) -- (a10);
	\draw (h21) -- (a11);
	\draw (h21) -- (a12);
	\draw (h21) -- (a13);
	\draw (h21) -- (a14);
	\draw (v) -- (a10);
	\draw (v) -- (a11);
	\draw (v) -- (a12);
	\draw (v) -- (a13);
	\draw (v) -- (a14);
\draw[] [blue](-0.8,-10) to [out=300,in=40] (-2.9,-14);
	
\node (h21) at (-3.5,-14.5) [label=0:$K_{1,5}$] {};
	
	\node (a15) at (-1.4,-12)  [circle,draw,scale=0.6] {77};
	\node (a16) at (-0.5,-12)  [circle,draw,scale=0.5] {153};
	\node (a17) at (0.4,-12)  [circle,draw,scale=0.5] {229};
	\node (a18) at (1.3,-12)  [circle,draw,scale=0.5] {305};
	\node (a19) at (2.2,-12)  [circle,draw,scale=0.5] {381};
	\node (a20) at (3.1,-12)  [circle,draw,scale=0.5] {457};
	\node (h22) at (0.5,-14)  [circle,draw,scale=0.5]  {-76};
	
	\draw (h22) -- (a15);
	\draw (h22) -- (a16);
	\draw (h22) -- (a17);
	\draw (h22) -- (a18);
	\draw (h22) -- (a19);
	\draw (h22) -- (a20);
	\draw (v) -- (a15);
	\draw (v) -- (a16);
	\draw (v) -- (a17);
	\draw (v) -- (a18);
	\draw (v) -- (a19);
	\draw (v) -- (a20);
\draw[] [blue](-0.82,-9.9) to [out=15,in=140] (0.3,-14);
	
\node (h22) at (-0.1,-14.5) [label=0:$K_{1,6}$] {};	
	\end{tikzpicture}
	
\medskip
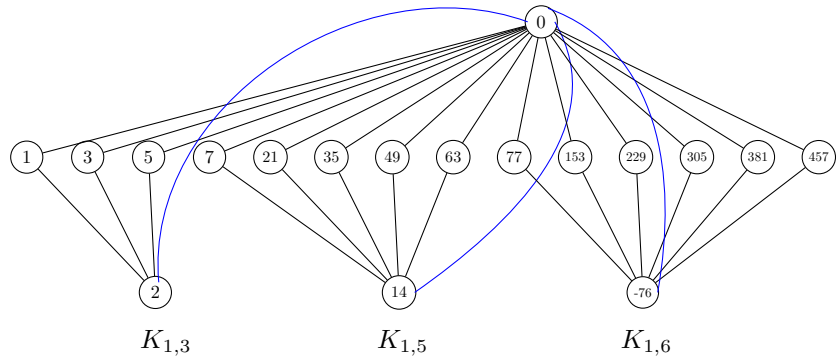
	Fig. 30. Integral sum graph $G * v$ with integral sum labeling
\end{center}

\subsection{Triangular book with bookmark and fan with handle are $\mathbb{Z}$-sum}

Here, we present the result that triangular book with bookmark and fan with handle are integral sum graphs.

\begin{dfn} \cite{vs15}\quad When $k$ copies of $C_n$ share a common edge, it will form an {\em $n$-gon book of $k$ pages} and is denoted by $B(n, k)$. The common edge is called the {\em spine} or {\em base of the book}. 
	
\medskip
	A {\em triangular book} $B(3,n)$ consists of $n$ triangles with a common edge and can be described as $B(3,n)$ = $TB_n$ = $P_2*nK_1$. $TB_n(u,v)$ = $P_2(u,v) * nK_1$ denotes the triangular book $B(3,n)$ with the spine $(u,v)$. Clearly, $TB_0$ = $K_2$ {\em represents a book without pages} or {\em the trivial book}.
\end{dfn}

\begin{dfn} \cite{vs15}\quad An $n$-gon book of $k$ pages $B(n,k)$ with a pendant edge terminating from any one of the end vertices of the spine is called an {\em $n$-gon book with a bookmark}. Triangular book $TB_n(u,v)$ with bookmark $(u,w)$ is denoted by $TB_n(u,v)(u,w)$ where $w$ is the pendant vertex adjacent to $u$.   
\end{dfn}

\begin{dfn} \cite{vs15}\quad If a fan graph $F_n$ has a pendant edge attached with the apex vertex, then the graph is called a {\em fan with a handle} or {\em a palm fan} and is denoted by $F_n^*$.   
\end{dfn}

Triangular book with bookmark are used to decompose complete graphs. Theorem \ref{fan} proves that fan graphs are integral sum. Next theorems show that triangular book with bookmark and fan graph with a handle are integral sum graphs.

\begin{theorem} [{ \cite[Theorem 5.1(i)]{vs15}}] \label{6.12} \normalfont \quad
	{\rm Triangular book with bookmark $TB_n(u_0,v_0)(u_0,w_0)$  is an integral sum graph, $n\in\mathbb{N}$. } 
\end{theorem}

\begin{proof}\normalfont \quad $TB_n(u_0,v_0)(u_0,w_0)$ is of order $n+3$, size $2n+2$ and $(u_0,w_0)$ is the pendant edge terminating at $u_0$ and let $V(TB_n(u_0,v_0)(u_0,w_0))$ = $\{w_0,u_0,v_0,v_1,\dots,v_n\}$. Define mapping $f: V(TB_n(u_0,v_0)(u_0,w_0))$ $\rightarrow$ $\mathbb{N}_0$ $\ni$ $f(u_0) = 0$, $f(v_0) = 2m$, $f(v_i) = 2mi+1$ for $i$ = $1,2,\dots,n$ and $f(w_0) = 2m(n+1)+1$, $m\in\mathbb{N}$. 
	
	Consider the integral sum graph $G^+(S)$ where $S = \{0, 2m, 2m+1, 4m+1$, $6m+1$, $. . . $, $2mn+1, 2m(n+1)+1 : m \in\mathbb{N}\}$ = $f(V(TB_n(u_0,v_0)(u_0,w_0)))$. Our aim is to prove that $f$ is an integral sum labeling of $TB_n(u_0,v_0)(u_0,w_0)$ and $G^+(S)$ = $TB_n(u_0,v_0)(u_0,w_0)$.
	
	$f(u_0) = 0$ implies, $f(u_0) + f(v_i) = f(v_i)$ and $f(u_0) + f(w_0) = f(w_0)$ for $i = 0,1,2,\dots,n$. This implies, $u_0$ is adjacent to $w_0, v_0$ and $v_i$ for $i = 1,2,\dots,n$. For $i = 1,2,\dots,n-1$, $f(v_0)+f(v_i)$ = $f(v_{i+1})$, $f(v_0) + f(v_n)$ = $f(w_0)$, $f(v_0) + f(u_0)$ = $f(v_0)$, $f(v_0) + f(w_0)$ $\neq$ $f(u_0),f(v_0),f(w_0),f(v_j)$ for $j = 1,2,\dots,n$. This implies, $v_0$ is adjacent to $u_0$ and $v_i$ and non-adjacent to $w_0$ for $i = 1,2,\dots,n$. Also, $f(w_0) + f(u_0) = f(w_0)$ and $f(w_0) + f(v_j)$ $\neq$ $f(w_0),f(u_0),f(v_j)$ for $j = 0,1,\dots,n$. This implies, $w_0$ is a pendant vertex adjacent only to $u_0$. For $i,j = 0,1,2,\dots,n$, $f(v_i)+f(w_0) \neq f(u_0),f(v_j)$. Also for $1\leq i,j,k \leq n$, $f(v_i) + f(v_j)\neq f(v_k)$ since $f(v_i) + f(v_j)$ is an even number and $f(v_k)$ is an odd number. This implies that $v_i$ and $v_j$ are non-adjacent in $TB_n(u_0,v_0)(u_0,w_0)$ when $i \neq j$ and $1 \leq i,j \leq n$. Thus, $v_j$ is adjacent only to $u_0$ and $v_0$ for $j = 1,2,\dots,n$.
	
	From all the above conditions integral sum graph $G^+(S)$ is same as $TB_n(u_0,v_0)$ $(u_0,w_0)$ and $f$ is an integral sum labeling of $TB_n(u_0,v_0)(u_0,w_0)$ where $S$ = $\{0, 2m$, $2m+1, 4m+1, \dots, 2mn+1, 2m(n+1)+1: m \in\mathbb{N}\}$.  
	
	In Figure 31, graph $TB_6(u_0,v_0)(u_0,w_0)$ which is a triangular book $B(3, 6)$ with 6 pages, spine $(u_0,v_0)$ and bookmark $(u_0,w_0)$ is presented where $w_0$ is the pendant vertex adjacent to $u_0$ and in Figure 32, the graph with an integral sum labeling is shown.
	\end{proof}  
	
\begin{center}
	\begin{tikzpicture}	[scale =0.9]
\node (b0) at (13,4)  [circle,draw,scale=0.6]{$w_0$};
\node (a0) at (9,3)  [circle,draw,scale=0.6]{$u_0$};
\node (c0) at (9,-0.75)  [circle,draw,scale=0.6] {$v_0$};
\node (b1) at (13,3)  [circle,draw,scale=0.6] {$v_1$};
\node (b2) at (13,2.25)  [circle,draw,scale=0.6] {$v_2$};
\node (b3) at (13,1.5)  [circle,draw,scale=0.6] {$v_3$};
\node (b4) at (13,0.75)  [circle,draw,scale=0.6] {$v_4$};
\node (b5) at (13,0)  [circle,draw,scale=0.6] {$v_5$};
\node (b6) at (13,-0.75)  [circle,draw,scale=0.6] {$v_6$};
		
		\draw (a0) -- (b0);
		\draw (a0) -- (c0);
		
		\draw (a0) -- (b1);
		\draw (a0) -- (b2);
		\draw (a0) -- (b3);
		\draw (a0) -- (b4);
		\draw (a0) -- (b5);
		\draw (a0) -- (b6);
		
		\draw (c0) -- (b1);
		\draw (c0) -- (b2);
		\draw (c0) -- (b3);
		\draw (c0) -- (b4);
		\draw (c0) -- (b5);
		\draw (c0) -- (b6);
		
		\node (b0) at (21,4)  [circle,draw,scale=0.6]{85};
		\node (a0) at (17,3)  [circle,draw,scale=0.7]{0};
		\node (c0) at (17,-0.75)  [circle,draw,scale=0.6] {12};
		\node (b1) at (21,3)  [circle,draw,scale=0.6] {13};
		\node (b2) at (21,2.25)  [circle,draw,scale=0.6] {25};
		\node (b3) at (21,1.5)  [circle,draw,scale=0.6] {37};
		\node (b4) at (21,0.75)  [circle,draw,scale=0.6] {49};
		\node (b5) at (21,0)  [circle,draw,scale=0.6] {61};
		\node (b6) at (21,-0.75)  [circle,draw,scale=0.6] {73};
		
		\draw (a0) -- (b0);
		\draw (a0) -- (c0);
		
		\draw (a0) -- (b1);
		\draw (a0) -- (b2);
		\draw (a0) -- (b3);
		\draw (a0) -- (b4);
		\draw (a0) -- (b5);
		\draw (a0) -- (b6);
		
		\draw (c0) -- (b1);
		\draw (c0) -- (b2);
		\draw (c0) -- (b3);
		\draw (c0) -- (b4);
		\draw (c0) -- (b5);
		\draw (c0) -- (b6);
		\end{tikzpicture}
		
		\medskip
		\hspace{.75cm}	Fig. 31.  $TB_7(u_0,v_0)(u_0,w_0)$ \hspace{1cm} Fig. 32. $TB_7(u_0,v_0)(u_0,w_0)$ with an integral sum labeling 	
	\end{center}
	
\begin{theorem} [{ \cite[Theorem 5.1(ii)]{vs15}}] \label{6.13} \normalfont \quad
	{\rm Fan graph with a handle $F^*_n$ is an integral sum graph, $n\in\mathbb{N}$. } 
\end{theorem}
\begin{proof}\normalfont \quad $F_n = P_n + K_1$ and $F_n^*$ is of order $n+2$ and size $2n$ where $P_n$ is a path on $n$ vertices. Let $V(F_n^*)$ = $\{u_0,v_0,v_1,\dots,v_n\}$ where $u_0$ is the pendant vertex, $v_0$ is the apex vertex and $d(v_0)$ = $n+1$ = $\Delta(F_n^*)$. Define mapping $f : V(F_n^*)\rightarrow \mathbb{N}_0$ such that $f(v_0) = 0$, $f(v_1) = p_m$, the $m^{th}$ Fibonacci number, $m\geq 2$, $f(v_i) = p_{m+i-1}$ for $i = 2,\dots,n$ and  $f(u_0) = p_{m+n}$. Here, $f(v_0) = 0$ $<$ $f(v_1) = p_m < f(v_2) = p_{m+1} < \dots < f(v_n) = p_{m+n-1} < f(u_0) = p_{m+n}$ and for $i-j \neq 1$ and $1 \leq i,j,k \leq n$, $f(v_i) + f(v_j)\neq f(v_k)$. Also $f(v_i) + f(v_{i+1}) = f(v_{i+2})$ for $i = 1,2,\dots,n-2$ and $f(v_{n-1}) + f(v_n) = f(u_0)$, $m \geq 2$. Hence the labeling $f$ is an integral sum labeling of graph $F_n^*$ and thereby $F_n^*$ is an integral sum graph. In Figure 33, graph $F_6^*$, fan graph $F_6 = P_6 * K_1(v_0)$ with handle $v_0 u_0$ is shown and in Figure 34, the graph with an integral sum labeling is shown. 
\end{proof}

\medskip
\begin{center}
	\begin{tikzpicture}[scale =1.2]
	
	\node (c0) at (10,1.5)  [circle,draw,scale=0.7]{$v_0$};
	\node (c1) at (8,0)  [circle,draw,scale=0.7]{$u_1$};
	\node (c2) at (9,0)  [circle,draw,scale=0.7] {$u_2$};
	\node (c3) at (10,0)  [circle,draw,scale=0.7] {$u_3$};
	\node (c4) at (11,0)  [circle,draw,scale=0.7] {$u_4$};
	\node (c5) at (12,0)  [circle,draw,scale=0.7]{$u_5$};
	\node (c6) at (13,0)  [circle,draw,scale=0.7]{$u_6$};
	\node (c7) at (13,1.5)  [circle,draw,scale=0.7]{$u_0$};
	
	\draw (c0) -- (c1);
	\draw (c0) -- (c2);
	\draw (c0) -- (c3);
	\draw (c0) -- (c4);
	\draw (c0) -- (c5);
	\draw (c0) -- (c6);
	\draw (c0) -- (c7);
	
	\draw (c1) -- (c2);
	\draw (c2) -- (c3);
	\draw (c3) -- (c4);
	\draw (c4) -- (c5);
	\draw (c5) -- (c6);
	
	\node (c0) at (17,1.5)  [circle,draw,scale=0.7]{0};
	\node (c1) at (15,0)  [circle,draw,scale=0.7]{2};
	\node (c2) at (16,0)  [circle,draw,scale=0.7] {3};
	\node (c3) at (17,0)  [circle,draw,scale=0.7] {5};
	\node (c4) at (18,0)  [circle,draw,scale=0.7] {8};
	\node (c5) at (19,0)  [circle,draw,scale=0.6]{13};
	\node (c6) at (20,0)  [circle,draw,scale=0.6]{21};
	\node (c7) at (20,1.5)  [circle,draw,scale=0.6]{34};
	
	\draw (c0) -- (c1);
	\draw (c0) -- (c2);
	\draw (c0) -- (c3);
	\draw (c0) -- (c4);
	\draw (c0) -- (c5);
	\draw (c0) -- (c6);
	\draw (c0) -- (c7);
	
	\draw (c1) -- (c2);
	\draw (c2) -- (c3);
	\draw (c3) -- (c4);
	\draw (c4) -- (c5);
	\draw (c5) -- (c6);	
	\end{tikzpicture}
	
	\medskip
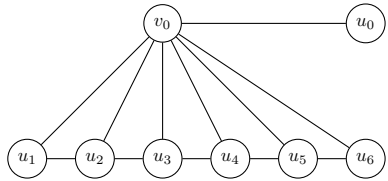
	\hspace{.5cm} 	Fig. 33. Fan with a handle $F^*_6$ \hspace{2.5cm} Fig. 34. $F^*_6$ with an integral sum labeling	
\end{center}

\subsection{On integral sum labeling of graph $P_k * G_n$}

In this subsection, we present our study on integral sum labeling of graph $P_k * G_n$, $k,n \in \mathbb{N}$. We prove that for $k$ = 1 to 3, $P_k * G_n$ is an integral sum graph and is isomorphic to $G_{-(k-1), n}$, and conjecture that $P_k * G_n$ is not an integral sum graph for $k,n\geq 4$, $k,n \in \mathbb{N}$. 

\begin{theorem} {\rm
	Let $k$ = 1 to 3, $n \in \mathbb{N}$ and $*$ be the graph operation of the join. Graph $P_k * G_n$ is an integral sum graph and $P_k * G_n$ is isomorphic to $G_{-(k-1), n}$. }
\end{theorem}
\begin{proof} Clearly, for $n\in\mathbb{N}$, $P_1 * G_n$ = $K_1 * G_n$ $\cong$ $K_1(0) * G_n$ $\cong$ $G_{0, n}$ and $P_2 * G_n$ = $K_2 * G_n$ $\cong$ $K_1(0) * (K_1(-1) * G_n)$ $\cong$ $G_{-1, n}$. 
	
	Also, we prove here graph $P_3 * G_n$ is an integral sum graph $\cong$ $G_{-2, n}$, $n \in \mathbb{N}$.
	Let $V(G_n)$ = $\{v_1, v_2,\ldots,v_m\}$ and $P_3$ = $u_1  u_2 u_3$. Define labeling $f:$ $V(P_3 * G_n)$ $\to\{-2,-1,0,1,2,\ldots,n\}$  such that$f(u_1)$ = -1, $f(u_2)$ = 0, $f(u_3)$ = -2 and $f(v_i)$ = $i$, $i$ = 1 to $n$. Clearly, $f$ is an integral sum labeling on $P_3 * G_n$ and $P_3 * G_n$ $\cong$ $K_1(u_2(0))$ $*$ $(G_2(u_1(-1),u_3(-2)) * G_{n})$. Integral sum labeling of graph $P_3 * G_7$ is given in Figure 35 where $P_3$ = $u_1(-1)  u_2(0) u_3(-2)$. See Figure 35.
	
	Hence the result. 
\end{proof}

\begin{center}
\begin{tikzpicture}[scale =0.9]
				
    \node (a6) at (12.3,4.4) [circle,fill=magenta!30,draw,scale=0.8] {1};
	\node (a0) at (10.5,4.5)  [circle,fill=green,draw,scale=0.8]{0};
	\node (a1) at (8,2)  [circle,fill=blue!30,draw,scale=0.8]{-2};
	\node (a2) at (11.5,-1) [circle,fill=magenta!30,draw,scale=0.8] {5};
	\node (a3) at (13,-.3) [circle,fill=magenta!30,draw,scale=0.8] {4};
	\node (a4) at (13.8,1)  [circle,fill=magenta!30,draw,scale=0.8] {3};
	\node (a5) at (13.5,3.2) [circle,fill=magenta!30,draw,scale=0.8]{2};
	\node (a7) at (10,-1)  [circle,fill=magenta!30,draw,scale=0.8]{6};
	\node (a8) at (8.7,3.7)  [circle,fill=blue!30,draw,scale=0.8]{-1};
	\node (a9) at (8.7,-.25) [circle,fill=magenta!30,draw,scale=0.8]{7};
				
				\draw (a0)[blue, thick] -- (a1);
				\draw (a0) -- (a2);
				\draw (a0) -- (a3);
				\draw (a0) -- (a4);
				\draw (a0) -- (a5);
				\draw (a0) -- (a6);
				\draw (a0) -- (a7);
				\draw (a0) -- (a9);
				
				\draw (a1) -- (a2);
				\draw (a1) -- (a3);
				\draw (a1) -- (a4);
				\draw (a1) -- (a5);
				\draw (a1) -- (a6);
				\draw (a1) -- (a7);
				\draw (a1) -- (a9);
				
				\draw (a8)[blue, thick] -- (a0);
				\draw (a8) -- (a2);
				\draw (a8) -- (a3);
				\draw (a8) -- (a4);
				\draw (a8) -- (a5);
				\draw (a8) -- (a6);
				\draw (a8) -- (a7);
				\draw (a8) -- (a9);
				
				\draw (a6)[magenta, thick] -- (a5);
				\draw (a6)[magenta, thick] -- (a4);
				\draw (a6)[magenta, thick] -- (a3);
				\draw (a6)[magenta, thick] -- (a2);
				\draw (a6)[magenta, thick] -- (a7);
				
				\draw (a5)[magenta, thick] -- (a4);
				\draw (a5)[magenta, thick] -- (a3);
				\draw (a5)[magenta, thick] -- (a2);
				
				\draw (a4) -- (a3);
		\end{tikzpicture}
		
		\vspace{.15cm}
		Fig. 35. $P_3 * G_7$ $\cong$ $G_{-2,7}$
		\label{P_3 * G_7}
\end{center}

\medskip
We propose the following conjecture on $P_k * G_n$ for $k,n \geq 4$ and $k,n\in\mathbb{N}$.

\begin{conj}
	For $k,n \geq 4$, $P_k * G_n$ is not an integral sum graph, $k,n\in\mathbb{N}$. \hfill $\Box$
\end{conj}

\subsection{Integral sum graphs by deleting 1 or 2 elements of $S$ of $G^+(S)$}

In this subsection, we focus on integral sum graphs whose labels are intervals of integers - the integral sum graphs on the bounded sets of consecutive integers with positive and negative entries. Next stage in the theory of sum graphs is the analysis of other sets, and here we consider the sets in which one value has been deleted from an interval. We focus on the deleted value not being one of the endpoints of the interval since those are among the interval sum graphs. We continue to identify each vertex with its value, and the value of an edge is the sum of the values of its end vertices. Then, we go for the analysis of integral sum graphs after removal of two vertices from $G_{-m,n}$, $m,n\in\mathbb{N}$. 

\begin{dfn} \cite{vm13c}  In an integral sum graph, the set of all edges each with same edge sum number, say $i$, is called the {\em edge-sum class} and is denoted by $E_i$. 
\end{dfn} 

We now define (in stages) families of subgraphs of $G_{-m,n}$ and we assume that $m \leq n$ and $m,n\in\mathbb{N}$. Using the definition of the edge-sum class $E_{i}$, set $F_{i} = \{v_{i}\} \cup E_{i}$. We start by defining $H_{-m,n}(i)$ with $-m \leq i \leq n$ and $m,n \in \mathbb{N}$ as the integral sum graph on the set $[-m, n] - \{{i}\}$; that is, $$H_{-m,n}(i) \cong G_{-m,n} - F_{i}.$$

Clearly, we then have 
$H_{-m,n}(-m) \cong G_{-(m-1), n}$ and  $H_{-m,n}(n) \cong G_{-m,n-1}$. 

By analogy, we can of course extend our definition to the integral sum graph when two numbers are removed, $[-m, n] - \{-i, j\}$ with $-m \leq -i <$ $0 < j \leq n$: 

$H_{-m,n}(-i, j) \cong G_{-m,n} - (F_{-i} \cup F_{j})$. 

Also, 	$H_{-m,n}(0) = H_{-m,n} - (E_0 \cup \{ 0 \}) \cong   (G_{-m} * G_n)- E_0$.

With the removal of the restriction on the relative magnitudes of the endpoints of the interval, we have `isomorphism of negatives': 

$H_{-m,n}(-i,j) \cong H_{-m,n}(-j,i)$.  Also, $H_{-m,n}(-m,n) \cong G_{-(m-1),n-1}$. 

The general case is to take any proper subset $W$ of the interval $[-m, n]$ of integers, and  let $\displaystyle E_W =  \bigcup_{i \in W} E_{i}$ and $F_W = W \cup E_W$.  Then $H_{-m,n}(W)$ is defined as $H_{-m,n}(W)$ $\cong$ $G_{-m,n} - F_W$.  Obviously, properties that are analogous to those considered above (such as the negativity property) also hold. Now, let us see a few results on the above.

\begin{theorem}  [{ \cite[Theorem 2.2]{vm13c}}] \label{thm2.1} {\rm Let $i, j, m$, and $n$ be integers with $-m \leq -i < 0 < j \leq n$ and $m \leq n$. Then, the following holds. 
\small 
\begin{enumerate}
	\item  $\|H_{-m,n}(0)\| = \|G_{-m,n}\| - 2m - n$.
			
	\item $\|H_{-m,n}(-i)\|$  = \small $\begin{cases}
		\text{$\|G_{-m,n}\| -2m - n + \left\lfloor \frac{3i+2}{2}\right\rfloor$}\!\! &\!\!\!\!\!\quad \text{if $ 1 \leq i \leq \left\lfloor \frac{m}{2}\right\rfloor$ }\\
		\text{$\|G_{-m,n}\| -2m - n + \left\lfloor \frac{3i+2}{2}\right\rfloor - 1 $}\!\!\!\! &\!\!\!\!\!\quad \text{if $ \left\lfloor \frac{m}{2}\right\rfloor < i \leq m$. }\\
	\end{cases}$
			
	\item 
\begin{enumerate}
	\item [\rm (a)] When  $\left\lfloor \frac{n}{2} \right\rfloor < n-m$,	
				
		$\|H_{-m,n}(j)\|$ = $\begin{cases}
		\text{$\|G_{-m,n}\| - 2m -n + \left\lfloor  \frac{j+2}{2}\right\rfloor$} & \quad \text{if $1 \leq j \leq \left\lfloor \frac{n}{2} \right\rfloor$}\\
		\text{$\|G_{-m,n}\| - 2m - n +  \left\lfloor \frac{j+2}{2}\right\rfloor - 1$} & \quad \text{if $ \left\lfloor \frac{s}{2} \right\rfloor < j \leq m-n$} \\
		\text{$\|G_{-m,n}\| - m  -  2n  + \left\lfloor \frac{3j+2}{2}\right\rfloor - 1 $} & \quad \text{if $ n-m < j \leq n$.} \\
	\end{cases}$
				
	\item [\rm (b)] When $\left\lfloor \frac{n}{2} \right\rfloor \geq n-m$,	
				
		\noindent $\|H_{-m,n}(j)\|$ = $\begin{cases}
					
		\text{$\|G_{-m,n}\| -  2m  - n  + \left\lfloor \frac{j+2}{2}\right\rfloor$ } & \quad \text{if $1 \leq j \leq n-m$}\\
					
		\text{$\|G_{-m,n}\| -  m  -  2n  + \left\lfloor \frac{3j+2}{2}\right\rfloor$ } & \quad \text{if $n-m < j \leq \left\lfloor \frac{n}{2}\right\rfloor $}\\
					
		\text{$\|G_{-m,n}\|  -  m  -  2n  +  \left\lfloor \frac{3j+2}{2}\right\rfloor - 1 $} & \quad \text{if $ \left\lfloor \frac{n}{2} \right\rfloor < j \leq n$.}\\
	\end{cases}$  
\end{enumerate} 	 	
\end{enumerate} }
\end{theorem}		

\begin{proof} 
	\begin{enumerate} 
		\item [(1)] We have  $ H_{-m,n}(0)  =  G_{-m,n}-F_0 $. 
		
		This implies,  $\|H_{-m,n}(0)\| = \|G_{-m,n}-F_0\|$ 
		
		\hfill = $\|G_{-m,n} - \{v_0,~ v_0v_{-i},~ v_{0}v_{j},~ v_{-l}v_{l}/$ $i,l = 1,2,\dots,m;~ j = 1,2,\dots, n\}\|$ 
		
		\hspace{.2cm}	= $ \|G_{-m,n}\| - (2m + n)$ ~since we have taken $m \leq n$.
		
		\item [(2)]  We have $H_{-m,n}(-i)= G_{-m,n}-F_{-i}$,  $1 \leq i \leq m$. Thus, $H_{-m,n}(-i)$ is obtained from $G_{-m,n}$ by removing vertex $v_{-i}$ and edge sum class $E_{-i}$. In $G_{-m,n}$, vertices adjacent to $v_{-i}$ are $v_0$, all the vertices of $G_n$ and a few vertices of $G_{-m}$. Hence, the number of elements of $E_{-i}$ in $G_{-m,n}$ is $|E_{-i}|$ in $G_{-m,n}$ =  $|E_{-i}|$ in $G_{-m}$ + $|E_{-i}|$ in $K_{m,n}$ + 1 (corresponds to $v_0v_{-i}$). 
		
		Also, $d(v_{-i})$ in $G_{-m,n}$ =  $d(v_{-i})$ in $G_{-m}$ + $d(v_{-i})$ in $K_{m,n}$ + $1$. 
		
		$\Rightarrow$	$\|H_{-m,n}(-i)\|$ =  $\|G_{-m,n}\| - |E_{-i}|$ in $G_{-m,n} -1 - d(v_{-i})$ in $G_{-m,n}$ 
		
		\hspace{2.35cm} = $\|G_{-m,n}\| - |E_{-i}|$ in $G_{-m} - |E_{-i}|$ in $K_{m,n}$ 
		
		\hspace{5cm}$ - d(v_{-i})$ in $G_{-m}- d(v_{-i})$ in $K_{-m,n}$ $ -1 $ \\
		Hence,	
		for $ 1 \leq i \leq \left\lfloor \frac{m}{2}\right\rfloor$,
		
		$\|H_{-m,n}(-i)\|$ = $\|G_{-m,n}\| - \left\lfloor \frac{i-1}{2}\right\rfloor - (m-i) - (m-i-1) - n -1$
		
		\hspace{2.1cm}	= $\|G_{-m,n}\|- 2m - n + \left\lfloor \frac{3i+2}{2} \right\rfloor$.\\ 
		And for $\left\lfloor \frac{m}{2}\right\rfloor < i \leq m$, 
		
		$\|H_{-m,n}(-i)\| $ = $\|G_{-m,n}\| - \left\lfloor \frac{i-1}{2}\right\rfloor - (m-i) - (m-i) - n -1$
		
		\hspace{2.1cm}	= $\|G_{-m,n}\|- 2m - n + \left\lfloor \frac{3i+2}{2} \right\rfloor - 1$. 
		
		\item [(3)]  We have $H_{-m,n}(j)= G_{-m,n}-F_{j}$,  $1 \leq j \leq n$. Thus, $H_{-m,n}(j)$ is obtained from $G_{-m,n}$ by removing vertex $v_{j}$ and edge sum class $E_{j}$. Vertices adjacent to $v_{j}$ in $G_{-m,n}$ are $v_0$, all the vertices of $G_{-m}$ and a few vertices of $G_{n}$.\\
		And $|E_{j}|$ in $G_{-m,n}$ =  $|E_{j}|$ in $G_{n}$ + $|E_{j}|$ in $K_{m,n}$ + 1.\\
		Also, $d(v_{j})$ in $G_{-m,n} =  d(v_{j})$ in $G_{n} + d(v_{j})$ in $K_{m,n}+ 1 $. 
		
		$\Rightarrow$	$\|H_{-m,n}(j)\|$ = $\|G_{-m,n}\| - |E_{j}|$ in $G_{-m,n}- d(v_{j})$ in $G_{-m,n}  -1 $
		
		\hspace{2.1cm}	= $\|G_{-m,n}\|$ - $|E_{j}|$ in $G_{n} - |E_{j}|$ in $K_{m,n} - d(v_{j})$ in $G_{n}$ $- d(v_{j})$ in $K_{m,n} -1$ 
		
		\hspace{2.1cm} = $\|G_{-m,n}\| - \left\lfloor \frac{j-1}{2}\right\rfloor - m -1 - |E_{j}|$ in $K_{m,n} - d(v_{j})$ in $G_{n}$ $- d(v_{j})$ in $K_{m,n}$. 
		
		Here, $|E_{j}|$ in $K_{m,n}$ = $\begin{cases}
			\text{$ m $} &\quad \text{if $ 1 \leq j \leq n-m. $}\\
			\text{$ n-j $} &\quad \text{if $ n-m \leq j \leq n$.}	\\	
		\end{cases}$\\
		For further simplification of 	$\|H_{-m,n}(j)\|$, we consider the following two cases. 
		
		\noindent \textbf{Case 1.} $\left\lfloor \frac{n}{2} \right\rfloor < n-m$. 
		
		In this case, $1 \leq j \leq n$ and the following three subcases arise.
		
		\noindent  \textbf{Subcase 1.1.} $1 \leq j \leq \left\lfloor \frac{n}{2} \right\rfloor$.
		
		In this case, {\small
			\begin{eqnarray*}
				\|H_{-m,n}(j)\|  &=& \|G_{-m,n}\| - \left\lfloor \frac{j-1}{2}\right\rfloor - m -1 - m - (n-j-1)  \\
				&=&  \|G_{-m,n}\| - 2m -n + j - \left\lfloor  \frac{j-1}{2}\right\rfloor\\
				&=& \|G_{-m,n}\| - 2m -n + \left\lfloor  \frac{j+2}{2}\right\rfloor
		\end{eqnarray*} }
		\textbf{Subcase 1.2.} $ \left\lfloor \frac{n}{2} \right\rfloor < j \leq n-m $.
		
		In this case, {\small 	
			\begin{eqnarray*}
				\|H_{-m,n}(j)\| &=& \|G_{-m,n}\| - \left\lfloor \frac{j-1}{2}\right\rfloor -m -1 - m - (n-j) \\
				&=& \|G_{-m,n}\| - 2m - n+\left\lfloor  \frac{j+2}{2}\right\rfloor - 1.  
		\end{eqnarray*} }
		
		\noindent
		\textbf{Subcase 1.3.}  $n-m < j \leq n$.
		
		In this case, {\small 	
			\begin{eqnarray*}
				\|H_{-m,n}(j)\| &=& \|G_{-m,n}\| - \left\lfloor \frac{j-1}{2}\right\rfloor - m -1 - (n-j) - (n-j)  \\
				&=&\|G_{-m,n}\| - 2n - m + 2j - \left\lfloor \frac{j-1}{2}\right\rfloor -1 \\
				&=& \|G_{-m,n}\| - 2n - m + \left\lfloor  \frac{3j+2}{2}\right\rfloor -1. 
		\end{eqnarray*} }
		\noindent 
		\textbf{Case 2.} $\left\lfloor \frac{n}{2} \right\rfloor \geq n-m$. \\ 
		
		In this case, $1 \leq j \leq n$ and the following three subcases arise.\\
		
		\noindent  
		\textbf{Subcase 2.1.} $1 \leq j \leq n-m$.\\
		{\small 		\begin{eqnarray*}
				\|H_{-m,n}(j)\| &=& \|G_{-m,n}\| - \left\lfloor \frac{j-1}{2}\right\rfloor - m -1 - m -(n-j-1)  \\
				&=&\|G_{-m,n}\| - 2m - n + \left \lfloor  \frac{j+2}{2}\right\rfloor.  
		\end{eqnarray*} }
		
		\noindent
		\textbf{Subcase 2.2.} $n-n < j \leq\left\lfloor \frac{n}{2} \right\rfloor $.
		{\small 		\begin{eqnarray*}
				\|H_{-m,n}(j)\|&=& \|G_{-m,n}\| -\left\lfloor \frac{j-1}{2}\right\rfloor - m -1 - (n-j) - (n-j-1) \\
				&=& \|G_{-m,n}\|- 2n - m + \left\lfloor  \frac{3j+2}{2}\right\rfloor.\\ 
		\end{eqnarray*} }
		\noindent
		\textbf{Subcase 2.3.} $ \left\lfloor \frac{n}{2} \right\rfloor < j \leq n$.
		{\small 		\begin{eqnarray*}
				\|H_{-m,n}(j)\|  &=& \|G_{-m,n}\| - \left\lfloor \frac{j-1}{2}\right\rfloor - m  -1 - (n-j) - (n-j) \\
				&=& \|G_{-m,n}\|  -  m  -  2n  +  \left\lfloor  \frac{3j+2}{2}\right\rfloor - 1.
		\end{eqnarray*} }
		Combining all the above cases, we get, 
		\begin{enumerate}
			\item [\rm (a)] for  $\left\lfloor \frac{n}{2} \right\rfloor < n-m$,	
			
			{\small 		$\|H_{-m,n}(j)\|$ = $\begin{cases}
					\text{$\|G_{-m,n}\| - 2m -n + \left\lfloor  \frac{j+2}{2}\right\rfloor$} & \quad \text{if $1 \leq j \leq \left\lfloor \frac{n}{2} \right\rfloor$}\\
					\text{$\|G_{-m,n}\| - 2m - n +  \left\lfloor \frac{j+2}{2}\right\rfloor - 1 $} & \quad \text{if $ \left\lfloor \frac{n}{2} \right\rfloor < j \leq n-m $} \\
					\text{$\|G_{-m,n}\| - m  -  2n  + \left\lfloor \frac{3j+2}{2}\right\rfloor - 1 $} & \quad \text{if $n-m < j \leq n$; } \\
				\end{cases}$  }
			
			\item [\rm (b)] for $\left\lfloor \frac{n}{2} \right\rfloor \geq n-m$,	
			
			{\small 		\noindent $\|H_{-m,n}(j)\|$ = $\begin{cases}
					
					\text{$\|G_{-m,n}\| -  2m  - n  + \left\lfloor \frac{j+2}{2}\right\rfloor$ } & \quad \text{if $1 \leq j \leq n-m$}\\
					
					\text{$\|G_{-m,n}\| -  m  -  2n  + \left\lfloor \frac{3j+2}{2}\right\rfloor$ } & \quad \text{if $n-m < j \leq \left\lfloor \frac{n}{2}\right\rfloor $}\\
					
					\text{$\|G_{-m,n}\|  -  m  -  2n  +  \left\lfloor \frac{3j+2}{2}\right\rfloor - 1 $} & \quad \text{if $ \left\lfloor \frac{n}{2} \right\rfloor < j \leq n$.}\\
				\end{cases}$ }
		\end{enumerate}
	\end{enumerate}
	Hence the result.
\end{proof}	

\begin{cor} [{ \cite[Corollary 2.3]{vm13c}}] \label{b2}   
	For  $-n \leq -i < 0 < j \leq n$, {\small 	
		\begin{enumerate}
			\item \noindent $\|H_{-n,n}(0)\| = \|G_{-n,n}\| - 3n$;
			\item  \noindent $\|H_{-n,n}(-j)\| =  \|H_{-n,n}(j) \|$ = \small $ \begin{cases}
				
				\text{$\|G_{-n,n}\| - 3s + \left\lfloor \frac{3j+2}{2}\right\rfloor $}  & ~ \text{if $ 1 \leq j \leq \left\lfloor \frac{n}{2}\right\rfloor$} \\
				\text{$ \|G_{-n,n}\| - 3n + \left\lfloor \frac{3j+2}{2} \right\rfloor - 1 $} & ~ \text{if $\left\lfloor \frac{n}{2}\right\rfloor < j \leq n$;}
			\end{cases} $
			
			\item \noindent  $\|H_{-n,n+1}(-i)\|$ =  $ \begin{cases}
				
				\text{$\|G_{-n,n+1}\| - 3n -1 +  \left\lfloor  \frac{3i+2}{2} \right\rfloor $} & \quad \text {if $1 \leq i \leq \left\lfloor \frac{n}{2}\right\rfloor$} \\
				\text {$\|G_{-n,n+1}\| - 3n -2 + \left\lfloor  \frac{3i+2}{2} \right\rfloor $} & \quad \text {if $\left\lfloor \frac{n}{2}\right\rfloor < i \leq n$;}
			\end{cases} $\\
			
			\item \noindent  $ \|H_{-n,n+1}(j) \| $ =  $ \begin{cases}
				\text {$\|G_{-n,n+1}\| - 3n -2 +  \left\lfloor  \frac{3j+2}{2} \right\rfloor$} & \quad \text{if $1\leq  j \leq \left\lfloor \frac{n+1}{2} \right\rfloor$} \\
				\text{$ \|G_{-n,n+1}\| - 3n -3 +  \left\lfloor  \frac{3j+2}{2} \right\rfloor $} & \quad \text {if $\left\lfloor \frac{n+1}{2}\right\rfloor < j \leq n+1$. \hspace{3cm} $\Box$} 
			\end{cases} $ 		
	\end{enumerate}  }
\end{cor}

For various values of $i > 0$ and $j \geq 0$, we calculate $\|H_{-4,5}(-i)\|$,  $\|H_{-4,5}(j)\|$ and $\|H_{-5,5}(-i)\|$ with respect to integral sum graphs $G_{-4,5}$ and $G_{-5,5}$. integral sum graphs $G_{-4,4}$, $G_{-4,5}$ and $G_{-5,5}$ are presented in Figures 36, 37 and 38.

$\|G_{-4,4}\| = 28, \|G_{-4,5}\| = 35, \|G_{-5,5}\| = 43$; 

$\|H_{-4,5}(-1)\| = 24$,  

$\|H_{-4,5}(-2)\| = 26 = \|H_{-4,5}(-3)\|$, 

$\|H_{-4,5}(-4)\|$ = $\|G_{-3,5}\|= 28$; 

$\|H_{-4,5}(1)\| = 23$, 

$\|H_{-4,5}(2)\| = 25 = \|H_{-4,5}(3)\|$,

$\|H_{-4,5}(4)\| = 27$, 

$\|H_{-4,5}(5)\| = 28$; 

$\|H_{-5,5}(-1)\| = 30$,  

$\|H_{-5,5}(-2)\| = 32 = \|H_{-5,5}(-3)\|$, 

$\|H_{-5,5}(-4)\| = 34$, 

$\|H_{-5,5}(-5)\|= 35$; 

$\|H_{-4,4}(0)\|$ = 16, 

$\|H_{-4,5}(0)\|$ = 22, 

$\|H_{-5,5}(0)\|$ = 28. 

\begin{center}
	\begin{tikzpicture}[scale = 0.4]
		
		\node (a0) at (0,7)  [circle,draw,scale=.6,  ]{0};
		\node (a1) at (3.5,5.5) [circle,draw,scale=.6,  ]{1};
		\node (a2) at (4.5,3) [circle,draw,scale=.6,  ] {2};
		\node (a3) at (3.5,.5) [circle,draw,scale=.6,  ] {3};
		\node (a4) at (1.5,-1) [circle,draw,scale=.6,  ] {4};
		
		\node (a9) at (-3.5,5.5) [circle,draw,scale=.5,  ] {-1};
		\node (a10) at (-4.5,3) [circle,draw,scale=.5,  ] {-2};
		\node (a11) at (-3.5,.5) [circle,draw,scale=.5, ] {-3};
		\node (a12) at (-1.5,-1) [circle,draw,scale=.5,  ] {-4};
		
		\draw (a0)[ gray, thick] --node[][above]{  } (a1);
		\draw (a0)[ gray, thick] --node[near start] {  }  (a2);
		\draw (a0)[ gray, thick] --node[near start] { }  (a3);
		\draw (a0)[ gray, thick] --node[near start] { }  (a4);
		
		\draw (a0)[ gray, thick] --node[near start] {}  (a9);
		\draw (a0)[gray, thick] --node[near start] { }  (a10);
		\draw (a0)[gray, thick] --node[near start] { }  (a11);
		\draw (a0)[gray, thick] --node[near start] { } (a12);
		
		\draw (a1)[ gray, thick] --node[near start] { }  (a9);
		\draw (a1)[gray, thick] --node[near start] { } (a10);
		\draw (a1)[ gray, thick] --node[near start] {}(a11);
		\draw (a1)[ gray, thick] --node[near start] {} (a12);
		
		\draw (a2)[ gray, thick] --node[near start] { } (a9);
		\draw (a2)[gray, thick] --node[near start] {}(a10);
		\draw (a2)[ gray, thick] --node[near start] { } (a11);
		\draw (a2)[ gray, thick] --node[near start] { } (a12);
		
		\draw (a3)[ gray, thick] --node[near start] { } (a9);
		\draw (a3)[gray, thick] --node[near start] { } (a10);
		\draw (a3)[gray, thick, dashed] --node[near start] { } (a10);
		\draw (a3)[ gray, thick] --node[near start] { }(a11);
		\draw (a3)[ gray, thick] --node[near start] { } (a12);
		
		\draw (a4)[ gray, thick] --node[near start] { } (a9);
		\draw (a4)[gray, thick] --node[near start] { }  (a10);
		\draw (a4)[  gray, thick] --node[near start] { }(a11);
		\draw (a4)[ gray , thick] --node[near start] {}(a12);
		
		\draw (a1)[ gray, thick] --node[above] { }  (a2);
		\draw (a1)[ gray, thick] --node[near start] { }(a3);
		
		\draw (a9)[gray, thick] --node[left] { }(a10);
		\draw (a9)[ gray, thick] --node[near end] { }(a11);
		
		\node (b0) at (11.5,7)  [circle,draw,scale=.6,  ]{0};
		\node (b1) at (15,7) [circle,draw,scale=.6,  ]{1};
		\node (b2) at (16.5,5.5) [circle,draw,scale=.6,  ] {2};
		\node (b3) at (17,3) [circle,draw,scale=.6,  ] {3};
		\node (b4) at (16.5,.5) [circle,draw,scale=.6,  ] {4};
		\node (b5) at (14,-1.5)  [circle,draw,scale=.6,  ] {5};
		
		\node (b9) at (9,5.5) [circle,draw,scale=.5,  ] {-1};
		\node (b10) at (8,3) [circle,draw,scale=.5,  ] {-2};
		\node (b11) at (9,.5) [circle,draw,scale=.5, ] {-3};
		\node (b12) at (10.5,-1) [circle,draw,scale=.5,  ] {-4};
		
		\draw (b0)[ gray, thick] --node[][above]{  } (b1);
		\draw (b0)[ gray, thick] --node[near start] {  }  (b2);
		\draw (b0)[ gray, thick] --node[near start] { }  (b3);
		\draw (b0)[ gray, thick] --node[near start] { }  (b4);
		\draw (b0)[ gray, thick] --node[near start] { }  (b5);
		\draw (b0)[ gray, thick] --node[near start] {}  (b9);
		\draw (b0)[gray, thick] --node[near start] { }  (b10);
		\draw (b0)[gray, thick] --node[near start] { }  (b11);
		\draw (b0)[gray, thick] --node[near start] { } (b12);
		
		\draw (b1)[ gray, thick] --node[near start] { }  (b9);
		\draw (b1)[gray, thick] --node[near start] { } (b10);
		\draw (b1)[ gray, thick] --node[near start] {}(b11);
		\draw (b1)[ gray, thick] --node[near start] {} (b12);
		
		\draw (b2)[ gray, thick] --node[near start] { } (b9);
		\draw (b2)[gray, thick] --node[near start] {}(b10);
		\draw (b2)[ gray, thick] --node[near start] { } (b11);
		\draw (b2)[ gray, thick] --node[near start] { } (b12);
		
		\draw (b3)[ gray, thick] --node[near start] { } (b9);
		\draw (b3)[gray, thick] --node[near start] { } (b10);
		\draw (b3)[gray, thick, dashed] --node[near start] { } (b10);
		\draw (b3)[ gray, thick] --node[near start] { }(b11);
		\draw (b3)[ gray, thick] --node[near start] { } (b12);
		
		\draw (b4)[ gray, thick] --node[near start] { } (b9);
		\draw (b4)[gray, thick] --node[near start] { }  (b10);
		\draw (b4)[  gray, thick] --node[near start] { }(b11);
		\draw (b4)[ gray , thick] --node[near start] {}(b12);
		
		\draw (b5)[ gray, thick] --node[near start] { } (b9);
		\draw (b5)[gray, thick] --node[near start] { }(b10);
		\draw (b5)[ gray, thick] --node[near start] { } (b11);
		\draw (b5)[ gray, thick] --node[near start] { } (b12);
		\draw (b5)[ gray, thick, dashed] --node[near start] { } (b12);
		
		\draw (b1)[ gray, thick] --node[above] { }  (b2);
		\draw (b1)[ gray, thick] --node[near start] { }(b3);
		\draw (b1)[ gray, thick] --node[near start] { } (b4);
		\draw (b2)[ gray, thick] --node[near start] { }(b3);
		
		\draw (b9)[gray, thick] --node[left] { }(b10);
		\draw (b9)[ gray, thick] --node[near end] { }(b11);
		
		\node (c0) at (24,7.5)  [circle,draw,scale=.6,  ]{0};
		\node (c1) at (27,7) [circle,draw,scale=.6,  ]{1};
		\node (c2) at (29.5,5) [circle,draw,scale=.6,  ] {2};
		\node (c3) at (29.5,3) [circle,draw,scale=.6,  ] {3};
		\node (c4) at (29,1) [circle,draw,scale=.6,  ] {4};
		\node (c5) at (28,-1)  [circle,draw,scale=.6,  ] {5};
		
		\node (c9) at (21.5,6) [circle,draw,scale=.5,  ] {-1};
		\node (c10) at (20.5,4) [circle,draw,scale=.5,  ] {-2};
		\node (c11) at (20.5,1.75) [circle,draw,scale=.5, ] {-3};
		\node (c12) at (22,0) [circle,draw,scale=.5,  ] {-4};
		\node (c13) at (24,-1.5) [circle,draw,scale=.5,  ] {-5};
		
		\draw (c0)[ gray, thick] --node[][above]{  } (c1);
		\draw (c0)[ gray, thick] --node[near start] {  }  (c2);
		\draw (c0)[ gray, thick] --node[near start] { }  (c3);
		\draw (c0)[ gray, thick] --node[near start] { }  (c4);
		\draw (c0)[ gray, thick] --node[near start] { }  (c5);
		\draw (c0)[ gray, thick] --node[near start] {}  (c9);
		\draw (c0)[gray, thick] --node[near start] { }  (c10);
		\draw (c0)[gray, thick] --node[near start] { }  (c11);
		\draw (c0)[gray, thick] --node[near start] { } (c12);
		\draw (c0)[gray, thick] --node[near start] {} (c13);

		\draw (c1)[ gray, thick] --node[near start] { }  (c9);
		\draw (c1)[gray, thick] --node[near start] { } (c10);
		\draw (c1)[ gray, thick] --node[near start] {}(c11);
		\draw (c1)[ gray, thick] --node[near start] {} (c12);
		\draw (c1)[gray, thick] --node[near start] {} (c13);
		
		\draw (c2)[ gray, thick] --node[near start] { } (c9);
		\draw (c2)[gray, thick] --node[near start] {}(c10);
		\draw (c2)[ gray, thick] --node[near start] { } (c11);
		\draw (c2)[ gray, thick] --node[near start] { } (c12);
		\draw (c2)[gray, thick] --node[near start] {} (c13);
		
		\draw (c3)[ gray, thick] --node[near start] { } (c9);
		\draw (c3)[gray, thick] --node[near start] { } (c10);
		\draw (c3)[gray, thick, dashed] --node[near start] { } (c10);
		\draw (c3)[ gray, thick] --node[near start] { }(c11);
		\draw (c3)[ gray, thick] --node[near start] { } (c12);
		\draw (c3)[gray, thick] --node[near start] {} (c13);
		
		\draw (c4)[ gray, thick] --node[near start] { } (c9);
		\draw (c4)[gray, thick] --node[near start] { }  (c10);
		\draw (c4)[  gray, thick] --node[near start] { }(c11);
		\draw (c4)[ gray , thick] --node[near start] {}(c12);
		\draw (c4)[ gray, thick] --node[near start] {} (c13);
		
		\draw (c5)[ gray, thick] --node[near start] { } (c9);
		\draw (c5)[gray, thick] --node[near start] { }(c10);
		\draw (c5)[ gray, thick] --node[near start] { } (c11);
		\draw (c5)[ gray, thick] --node[near start] { } (c12);
		\draw (c5)[ gray, thick, dashed] --node[near start] { } (c12);
		\draw (c5)[gray, thick] --node[near start] {} (c13);
		
		\draw (c1)[ gray, thick] --node[above] { }  (c2);
		\draw (c1)[ gray, thick] --node[near start] { }(c3);
		\draw (c1)[ gray, thick] --node[near start] { } (c4);
		\draw (c2)[ gray, thick] --node[near start] {} (c3);
		
		\draw (c9)[gray, thick] --node[left] { }(c10);
		\draw (c9)[ gray, thick] --node[near end] { }(c11);
		\draw (c9)[ gray, thick] --node[near end] { }(c12);
		
		\draw (c10)[ gray, thick] --node[left] { } (c11);				
	\end{tikzpicture}				
	
	\medskip
	Fig. 36. $G_{-4, 4}$  \hspace{2.5cm} Fig. 37. $G_{-4, 5}$ \hspace{2.5cm} Fig. 38. $G_{-5, 5}$
\end{center}

\medskip
Our next step is, for a given natural number $n$, to find maximal integral sum subgraph(s) of order $n$ of $G_{-m,n-m}$ (of order $n+1$) where $-m < 0 < n-m$. 

$H_{-m,n-m}(0)$, $H_{-m,n-m}(-i)$ and $H_{-m,n-m}(j)$ are integral sum graphs of order $n$ each, $-i < 0 < j$. They are obtained by simply removing one vertex with label $0, -i$ or $j$ and the corresponding edges, each with edge sum $0, -i$ or $j$ from the respective graph. The following result shows that by calculating values of $\|H_{-m,n-m}(0)\|$, $ \|H_{-m,n-m}(-i)\|$ and $\|H_{-m,n-m}(j)\|$ for all possible values of $i,j > 0$, it is possible to find out maximal integral sum subgraph(s) of order $n$ of $G_{-m,n-m}$, $-m < 0 < n-m$. 

\begin{theorem} [{ \cite[Theorem 2.4]{vm13c}}] \label{b3}  {\rm Let $m$ and $n$ be integers with $-m \leq -i < 0 < n$ and $n-m \geq 0$. Then for all possible values of $i$, integral sum graph $H_{-m,n-m}(-i)$ has maximum number of edges when $i$ = $m$. }
\end{theorem}

\begin{proof}
	For simplicity, let $m \leq n-m$, and $m,n-m\in\mathbb{N}$. Then, for $1 \leq i \leq m$, using Theorem \ref{thm2.1}, 
	\begin{equation*}
		\|H_{-m,n-m}(-i)\| = 
		\begin{cases}
			\text{$\|G_{-m,n}\| - m -n + \left\lfloor \frac{3i+2}{2}\right\rfloor$ } &\quad \text{if $1 \leq i \leq \left\lfloor \frac{m}{2} \right\rfloor$}\\
			
			\text{$\|G_{-m,n}\| - m -n + \left\lfloor \frac{3i+2}{2}\right\rfloor - 1$} &\quad \text{if $\left\lfloor \frac{n}{2} \right\rfloor < i \leq m$.}\\
		\end{cases}	
	\end{equation*}	 
	
	Using Theorem \ref{a5}, we get, \\
	
	\hspace{.2cm}	$\|G_{-m,n-m}\|$ = $\frac{1}{4}(m^2+(n-m)^2+3m+3(n-m)+4m(n-m))-\frac{1}{2}\left( \left\lfloor \frac{m}{2}\right\rfloor + \left\lfloor \frac{n-m}{2}\right\rfloor\right)$, 
	
	\hfill $-m < 0 < n-m$.
	
	\hspace{1cm}	 = $\frac{1}{4}(m^2+(n-m)^2+3n+4m(n-m))-\frac{1}{2}\left( \left\lfloor \frac{m}{2}\right\rfloor + \left\lfloor \frac{n-m}{2}\right\rfloor\right)$,  $-m < 0 < n-m$.
	
	\begin{eqnarray*}
\therefore	\|H_{-m,n-m}(-i)\|
		&=& \!\!\begin{cases}
			\tiny \text{\!$\frac{1}{4}(m^2+(n-m)^2+4m(n-m)+3n)	\!-\! \frac{1}{2}\left( \left\lfloor \frac{m}{2}\right\rfloor \!+ \!\left\lfloor \frac{n-m}{2}\right\rfloor\right)\! -\! m\! -n + \left\lfloor \frac{3i+2}{2}\right\rfloor$} \!\! & \!\!{\tiny\text{if $1 \leq i \leq \left\lfloor \frac{m}{2} \right\rfloor$}}\\
			
			\tiny \text{$\frac{1}{4}(m^2+(n-m)^2+4m(n-m)+3n)	\!-\! \frac{1}{2}\left( \left\lfloor \frac{m}{2}\right\rfloor \!+ \!\left\lfloor \frac{n-m}{2}\right\rfloor\right)\! - m -n + \left\lfloor \frac{3i+2}{2}\right\rfloor - 1$}\!\! &\!\! {\tiny \text{if $\left\lfloor \frac{m}{2} \right\rfloor < i \leq m$}}\\
		\end{cases}\\ 
		&=&	\!\!\begin{cases}
			\tiny \text{\!$\frac{1}{4}(m^2+(n-m)^2+4m(n-m)-n) -\! m\!-\! \frac{1}{2}\left(\left\lfloor \frac{m}{2}\right\rfloor\!+ \!\left\lfloor \frac{n-m}{2}\right\rfloor\right)\!\! + \left\lfloor \frac{3i+2}{2}\right\rfloor$} \!\! & \!\!{\tiny \text{if $1 \leq i \leq \left\lfloor \frac{m}{2} \right\rfloor$}}\\
			
			\tiny \text{$\frac{1}{4}(m^2+(n-m)^2+4m(n-m)-n) -\! m\!-\! \frac{1}{2}\left( \left\lfloor \frac{m}{2}\right\rfloor \!+ \!\left\lfloor \frac{n-m}{2}\right\rfloor\right)\! + \left\lfloor \frac{3i+2}{2}\right\rfloor - 1$} \!\!& \!\!{\tiny \text{if $\left\lfloor \frac{m}{2} \right\rfloor < i \leq m$} }
		\end{cases}\\
		&=&	\!\! \begin{cases}
			\tiny \text{\!$\frac{1}{4}((n+m)^2- 3m^2  - 4m -\! n)\!-\! \frac{1}{2}\left(\left\lfloor \frac{m}{2}\right\rfloor\!+ \!\left\lfloor \frac{n-m}{2}\right\rfloor\right)\!\! + \left\lfloor \frac{3i+2}{2}\right\rfloor$}  & \text{if $1 \leq i \leq \left\lfloor \frac{m}{2} \right\rfloor$}\\
			
			\tiny \text{$\frac{1}{4}((n+m)^2- 3m^2  - 4m -\! n)\!-\! \frac{1}{2}\left(\left\lfloor \frac{m}{2}\right\rfloor\!+ \!\left\lfloor \frac{n-m}{2}\right\rfloor\right)\!\! + \left\lfloor \frac{3i+2}{2}\right\rfloor - 1$} & \text{if $\left\lfloor \frac{m}{2} \right\rfloor < i \leq m$.}\\
		\end{cases}
	\end{eqnarray*}
	
	For further simplification of $\|H_{-m,n-m}(-i)\|$, we consider the following 4 cases.\\
	
	\noindent
	\textbf{Case 1}. $n$ and $r$ are even. \\
	In this case, $\left\lfloor \frac{m}{2}\right\rfloor + \left\lfloor \frac{n-m}{2} \right\rfloor$ = $\frac{n}{2}$.  Then, 
		 
	\begin{eqnarray*}
	\|H_{-m,n-m}(-i)\|
		&=& \begin{cases}
			\text{\!$\frac{1}{4}((n+m)^2- 3m^2  - 4m -\! n)\!-\! \frac{n}{4}\! + \left\lfloor \frac{3i+2}{2}\right\rfloor$}  &\quad \text{if $1 \leq i \leq \left\lfloor \frac{m}{2} \right\rfloor$}\\
			
			\text{$\frac{1}{4}((n+r)^2- 3r^2  - 4r -\! n)\!-\! \frac{n}{4}\!\! + \left\lfloor \frac{3i+2}{2}\right\rfloor - 1$} &\quad \text{if $\left\lfloor \frac{m}{2} \right\rfloor < i \leq m$}\\
		\end{cases}	\\
		&=& \begin{cases}
			\text{\!$\frac{1}{4}((n+m)^2- 3m^2  - 4m -\! 2n) + \left\lfloor \frac{3i+2}{2}\right\rfloor$}  &\quad \text{if $1 \leq i \leq \left\lfloor \frac{m}{2} \right\rfloor$}\\
			
			\text{$\frac{1}{4}((n+m)^2- 3m^2  - 4m -\! 2n) + \left\lfloor \frac{3i+2}{2}\right\rfloor - 1$} &\quad \text{if $\left\lfloor \frac{m}{2} \right\rfloor < i \leq m$.}\\
		\end{cases}	\\
	\end{eqnarray*}
	
	\noindent \textbf{Case 2.}~ $n$ is even and $m$ is odd. \\
	In this case, $\left\lfloor \frac{m}{2}\right\rfloor +$ $\left\lfloor \frac{n-m}{2} \right\rfloor$ = $ \frac{n-2}{2} $. 
	
	\begin{eqnarray*} 
	\therefore	\|H_{-m,n-m}(-i)\|
	&=& \begin{cases}
			\small \text{\!$\frac{1}{4}((n+m)^2- 3m^2  - 4m -\! n)- \frac{n-2}{4} + \left\lfloor \frac{3i+2}{2}\right\rfloor$}  &\quad \text{if $1 \leq i \leq \left\lfloor \frac{m}{2} \right\rfloor$}\\
			
			\small \text{$\frac{1}{4}((n+m)^2- 3m^2  - 4m -\! n) - \frac{n-2}{4} + \left\lfloor \frac{3i+2}{2}\right\rfloor - 1$} &\quad \text{if $\left\lfloor \frac{m}{2} \right\rfloor < i \leq m$}\\
		\end{cases}	\\
		&=& \begin{cases}
			\small \text{\!$\frac{1}{4}((n+m)^2- 3m^2  - 4m -\! 2n + 2) + \left\lfloor \frac{3i+2}{2}\right\rfloor$}  &\quad \text{if $1 \leq i \leq \left\lfloor \frac{m}{2} \right\rfloor$}\\
			
			\small \text{$\frac{1}{4}((n+m)^2- 3m^2  - 4m -\! 2n + 2) + \left\lfloor \frac{3i+2}{2}\right\rfloor - 1$} &\quad \text{if $\left\lfloor \frac{m}{2} \right\rfloor < i \leq r$.}\\
		\end{cases}	\\
	\end{eqnarray*}
	\noindent \textbf{Case 3}. $n$ is odd and $m$ is even. 
	
	In this case, $\left\lfloor \frac{m}{2}\right\rfloor + \left\lfloor \frac{n-m}{2}\right\rfloor$ = $ \frac{n-1}{2}$. 
	
	\begin{eqnarray*}
\therefore	 \|H_{-m,n-m}(-i)\|
		&=& \begin{cases}
			\small \text{\!$\frac{1}{4}((n+m)^2- 3m^2  - 4m -\! n)- \frac{n-1}{4} + \left\lfloor \frac{3i+2}{2}\right\rfloor$}  &\quad \text{if $1 \leq i \leq \left\lfloor \frac{m}{2} \right\rfloor$}\\
			
			\small \text{$\frac{1}{4}((n+m)^2- 3m^2  - 4m -\! n) - \frac{n-1}{4} + \left\lfloor \frac{3i+2}{2}\right\rfloor - 1$} &\quad \text{if $\left\lfloor \frac{m}{2} \right\rfloor < i \leq m$}\\
		\end{cases}	\\
		&=& \begin{cases}
			\small \text{\!$\frac{1}{4}((n+m)^2- 3m^2  - 4m -\! 2n + 1) + \left\lfloor \frac{3i+2}{2}\right\rfloor$}  &\quad \text{if $1 \leq i \leq \left\lfloor \frac{m}{2} \right\rfloor$}\\
			
			\small \text{$\frac{1}{4}((n+m)^2- 3m^2  - 4m -\! 2n + 1) + \left\lfloor \frac{3i+2}{2}\right\rfloor - 1$} &\quad \text{if $\left\lfloor \frac{m}{2} \right\rfloor < i \leq m$.}\\
		\end{cases}	\\
	\end{eqnarray*}
	\normalsize		
	\noindent \textbf{Case 4}.~ $n$ and $m$ are odd.\\
	In this case, $\left\lfloor \frac{m}{2}\right\rfloor + \left\lfloor \frac{n-m}{2} \right\rfloor$ = $ \frac{n-1}{2}$. Thus, as in case 3, we obtain,
	\begin{equation*}
		\tiny{ \|H_{-m,n-m}(-i)\| =} 
		\begin{cases}
			\tiny \text{$\frac{1}{4}((n+m)^2- 3m^2  - 4m -\! 2n + 1) + \left\lfloor \frac{3i+2}{2}\right\rfloor$} &\quad \text{if $1 \leq i \leq \left\lfloor \frac{m}{2} \right\rfloor$}\\
			
			\tiny \text{$\frac{1}{4}((n+m)^2- 3m^2  - 4m -\! 2n + 1) + \left\lfloor \frac{3i+2}{2}\right\rfloor - 1$} &\quad \text{if $\left\lfloor \frac{m}{2} \right\rfloor < i \leq m$.}\\
		\end{cases}	
	\end{equation*}	 
	
	The values of $\|H_{-m,n-m}(-i)\|$ calculated in all the above cases show that for given values of $n$ and $m$ and for all possible values of $i$, the maximum value of $\|H_{-m,n-m}(-i)\|$ occurs when $\left\lfloor \frac{3i+2}{2}\right\rfloor$ is maximum.  And $\left\lfloor \frac{3i+2}{2}\right\rfloor$ attains its maximum value when $i$ is maximum, $ 1 \leq i \leq m$. The maximum possible value of $i$ is $m$ and the maximum value is $\|H_{-m,n-m}(-m)\|$ = $\|G_{-(m-1),n-m}\|$, $n-m \geq m$, $m,n-m\in\mathbb{N}$.
	
	Similarly, we can prove the result when $m \geq n-m$, $m,n-m,n\in \mathbb{N}$.
	
	Thus, for given integers $r$ and $n$ with $-r \leq -i < 0 < n$ and $n \geq r$, and for all possible values of $i$, the value of $\|H_{-m,n-m}(-i)\|$ is maximum when $i = m$ and the maximum value is $\|H_{-m,n-m}(-m)\| = \|G_{-(m-1),n-m}\|$. Hence the result.
\end{proof}

\begin{theorem} [{ \cite[Theorems 2.2, 2.4, and 2.5]{vm13c}}] \label{b4} {\rm Let $i,j,m$ and $n$ be integers with $m \leq n$, $1 \leq i \leq m$ and $1 \leq j \leq n$. Then, for given $m$ and $n$, the following statements are ture.
	\begin{enumerate}
		\item  The maximum number of edges in $H_{-m,n}(-i)$ occurs when $i = m$.
		\item  The maximum number of edges in $H_{-m,n}(j)$ occurs when $j = n$.
		\item  The maximum number of edges in $H_{-m,n}(-i, j)$ occurs when $i = m$ and $j = n$. \hfill $\Box$ 
	\end{enumerate}   }
\end{theorem}

\section{On sum number $\sigma (G)$ and integral sum number $\zeta (G)$ of graph $G$}

In this section, we present the investigation done in \cite{ns01} on some open problems posed by Harary (1994) in \cite{h94} on sum and integral sum numbers of graphs. In a sum graph, vertex with the highest label cannot be adjacent to any other vertex and thereby every sum graph must contain isolated vertices. {\em The sum number of a graph} $G$ is denoted by  $\sigma (G)$ and is defined as the minimum number of isolated vertices that must be added to $G$ so that the resulting graph is a sum graph. A labeling that makes $G$ together with $\sigma (G)$ isolated points a sum graph is called an {\em optimal sum graph labeling}. {\em The integral sum number} $\zeta (G)$ of a graph is the smallest non-negative integer such that $G$ $\cup$ $\zeta (G) K_1$ is an integral sum graph.

In the sum labeling of a graph, vertices whose labels correspond to an edge $uv$ are said to be {\em working vertices}. It has been realized that certain graphs can only be labelled in such a way that all the working vertices are also isolated; such graphs are called {\em exclusive}. It has been shown that $K_n$ \cite{bh89}, $W_n$ \cite{hs95}, and the complete $n$-partite graph $H_{m, n}$ \cite{mr98} are exclusive, while the forest $F$ and $K_{m,n}$ are not. Miller et al. \cite{mr98} opined that the exclusive graphs are more likely to have large sum numbers. In \cite{ns01}, it is proved that $P_n \times P_2$ is exclusive and $\sigma (P_n \times P_2)$ = 3, though the problem of characterization of exclusive graphs remains an open problem.

Determining $\sigma (G)$ and $\zeta (G)$, for a given graph $G$, is an interesting work in the study of sum and integral sum graphs. In 1989, Hao \cite{h89} established a lower bound on $\sigma (G)$ that a sum graph $G$ of order $p$ and size $q$ exists if and only if $q \leq (\binom{p}{2}-\left\lfloor \frac{p}{2} \right\rfloor)/2$ where $\left\lfloor x \right\rfloor$ denotes the integer part of $x$. For some special classes of graphs, such as cycles, trees, complete graphs and complete bipartite graphs, the sum number has been derived, and in the following, we list those useful results.

\begin{theorem}{\rm (Bergstrand et al. \cite{bh89})} \label{6.0.1} \quad {\rm For all complete graph $K_n$ of order $n$ with $n \geq 4$, $\sigma (K_n)$ = $2n-3$. \hfill $\Box$} 
\end{theorem}

\begin{theorem}{\rm (Harary \cite{h90})} \label{6.0.2} \quad {\rm For all cycles $C_n$ of order $n$, $\sigma (C_n)$ = 3 when $n$ = 4 and $\sigma (C_n)$ = 2 when $n \neq 4$. \hfill $\Box$} 
\end{theorem}

\begin{theorem}{\rm (Hartsfield and Smyth \cite{hs95})} \label{6.0.3} \quad {\rm For complete bipartite graph $K_{r,s}$ with $r \leq s$, $\sigma (K_{r,s})$ = $\left\lceil (3r+s-3)/2 \right\rceil$. \hfill $\Box$} 
\end{theorem}

Harary \cite{h94} posed an open problem of characterizing the classes of graphs $G$ for which $\sigma (G)$ = $\zeta (G)$. B. Xu \cite{x99} proved the following theorem in this connection:

\begin{theorem}{\rm (B. Xu \cite{x99})} \label{6.0.4} \quad {\rm For all complete graphs $K_{n}$ with $n \geq 4$, $\zeta (K_{n})$ = $\sigma (K_n)$ = $2n-3$. \hfill $\Box$} 
\end{theorem}

\begin{theorem}{\rm (B. Xu \cite{x99})} \label{6.0.5} \quad {\rm For an arbitary integer $m \geq 1$, $m K_{3}$ is an integral sum graph.  \hfill $\Box$} 
\end{theorem}

\begin{theorem}{\rm (A. Sharary \cite{s96})} \label{6.0.6} {\rm $\zeta (C_{n})$ = 3 if $n$ = 4 and $\zeta (C_n)$ = 0 if $n \neq 4$. \hfill $\Box$} 
\end{theorem}

In \cite{ns01}, it is proved that (i) $\sigma (P_n \times P_2)$ = $\zeta (P_n \times P_2)$ = 3 for any $n \geq 2$ where $P_n$ is a path on $n$ vertices; (ii) $P_n \times P_2$ is exclusive; (iii) $\sigma (mC_n)$ = 2 = $\sigma (C_n)$ when $n \neq 4$ and $m \geq 1$; and (iv) $\zeta (mC_n)$ $\neq$ $\sigma (mC_n)$. Also, the result of Xu \cite{x99} is extended to achieve integral sum labelling for a special class of caterpillars $T$ = $S(x_1,x_2,\dots,x_n)$, where $x_2 > 1$.

Ellingham \cite{e93} proved the conjecture of Harary \cite{h90} that $\sigma (T)$ = 1 for every tree  $T \neq K_1$. Smyth \cite{s91} conjectures that the disjoint union of graphs with sum number 1 has sum number 1. And its more general conjecture \cite{k01} is $\sigma (G \cup H) \leq \sigma (G) + \sigma (H) - 1$. For more details on sum number see \cite{g26}.   

\subsection{Graph $P_n \times P_2$ and its $\sigma (P_n \times P_2)$, $\zeta (P_n \times P_2)$, and exclusiveness   }

The graph $L_n$ = $P_n \times P_2$, $n \geq 2$, known as the {\em ladder graph}, has $2n$ vertices, viz., $u_1, u_2, \dots, u_n$ and $v_1, v_2, \dots, v_n$ and edge set 

$E(L_n)$ = $\{u_iu_{i+1}, v_iv_{i+1}/ i = 1,2,\dots,n-1 \}$ $\cup$ $\{u_i v_i/ i = 1,2,\dots,n\}$.

\begin{theorem} [{\cite[Theorem 2.1]{ns01}}]   \label{6.1.1} 
	{\rm Let $n \geq 2$ and the graph $L_n$ = $P_n \times P_2$. $\sigma (P_n \times P_2)$ = 3.  }
\end{theorem}	

\begin{proof} When $n$ = 2, $L_2$ = $C_4$ and hence $\sigma (L_2)$ = 3 by Theorem \ref{6.0.2}  Suppose $n$ = 3. To get the best possible sum labeling, we have to arrange the labels to get minimum number of isolated vertices by maximizing the number of duplicate sums. We choose distinct non-zero positive integers as the labels of the vertices of the graph $L_3$ such that $u_1+v_1$ = $u_2+v_2$ = $u_3+v_3$ = $k_1$. Then there must be an isolated vertex with label $k_1$. Let $u_1+u_2$ = $k_2$ ($\neq k_1$). Then $v_1+v_2$ $\neq k_2$, since $k_1$ $\neq k_2$. Also, $u_2+u_3\neq k_2$, since $u_1 \neq u_3$. Let $u_2+u_3$ = $k_3$. It is possible to have $v_1+v_2$ = $k_3$ and $v_2+v_3$ = $k_2$. If $k_2 < k_1$, then $k_3 > k_1$ and vice versa. Assume that $k_2 < k_1$. Then $k_3$ must be a label of an isolated vertex. Now, if $k_2$ is a label of any of the vertices of $L_3$, then it may hold either the vertex $u_3$ or $v_1$. Assume $u_3$ = $k_2$ (the case when $v_1$ = $k_2$ becomes redundant). Then, we get 
	
	$u_1$ = $k_2-u_2$ = $k_2-(k_1-v_2)$ = $k_2-k_1+(k_2-v_3)$ = $2k_2-k_1-v_3$ 
	
\hfill	= $2k_2-k_1-(k_1-u_3)$ = $2k_2-2k_1+u_3$ = $3k_2- 2k_1$, 
	
	$u_2$ = $k_1-v_2$ = $k_1-(k_2-v_3)$ = $k_1-k_2+v_3$ = $k_1-k_2+(k_1-u_3)$ = $2k_1-k_2-k_2$ = $2(k_1 - k_2)$, 
	
	$u_3$ = $k_2$, $v_1$ = $k_1-u_1$ = $k_1-(3k_2-2k_1)$ =
	$3(k_1 - k_2)$, 
	
	$v_2$ = $k_2-v_3$ = $k_2-(k_1-u_3)$ = $2k_2 - k_1$, 	
	$v_3$ = $k_1-u_3$ = $k_1 - k_2$. 
	
\medskip
But $u_1+v_3$ = $2k_2 - k_1$ which is the label of $v_2$, whereas, $u_1 v_3$ is not an edge of $L_3$. Hence, the $k_2$ also holds an isolated vertex and so $\sigma (L_3) \geq 3$. In general, 
	
	\hspace{2cm}	$\sigma (L_n) \geq 3$, for any $n \geq 2$. \hfill (1)
	
\medskip
Now, for $L_n$ = $P_n \times P_2$, we shall establish a labeling such that $\sigma (L_n)$ = 3. Consider the integer $4(n+1)$. It can be written as the sum of two distinct positive odd integers in $n+1$ ways, viz., 1, $4n+3$; 3, $4n+1$; $\dots$; $2n+1$, $2n+3$. Neglecting the first pair, we can label the vertices of $L_n$ with the remaining $n$ pairs as $u_i$ = $2i+1$; $v_i$ = $4n-2i+3$ for odd $i$ and $u_i$ = $4n-2i+3$; $v_i$ = $2i+1$ for even $i$, where $i$ = $1,2,\dots,n$. Clearly, for each column edge $u_i v_i$, we get $u_i + v_i$ = $4(n+1)$ for $i$ = $1,2,\dots,n$. Now, $u_i u_{i+1}$ = $4n+2$ if $i$ is odd and $4n+6$ if $i$ is even. Similarly, by construction, $v_i v_{i+1}$ = $4n+6$ if $i$ is odd and $4n+2$ if $i$ is even where $i$ = $1,2,\dots,n-1$. Let $w_1, w_2, w_3$ be three isolated vertices with labels $4n+2$, $4(n+1)$ and $4n+6$.
	
	It remains to be shown that no vertex has a label which is the sum of labels of two non-adjacent vertices. For any vertex $x$ of $L_n$, and for $i$ = 1,2,3, $x+w_i$ would yield an odd number greater than $4n+3$ which is not at all used anywhere in the labeling scheme. Also, for any two non-adjacent vertices $x$ and $y$ of $L_n$, it can be easily verified that $x+y$ is even and $x+y\notin \{w_1, w_2, w_3\}$. Now we could give a labeling by which $\sigma (L_n)$ = 3 for all $n \geq 2$. Hence by (1), we have proved the result.		
\end{proof}	

In the sum labeling of a graph, vertices whose labels correspond to an edge $uv$ are said to be {\em working vertices}. It has been realized that certain graphs can only be labelled in such a way that all the working vertices are also isolated; such graphs are called {\em exclusive}. It has been shown that $K_n$ \cite{bh89}, $W_n$ \cite{hs95}, and the complete $n$-partite graph $H_{m, n}$ \cite{mr98} are exclusive, while the forest $F$ and $K_{m,n}$ are not. Miller et al. \cite{mr98} opined that the exclusive graphs are more likely to have large sum numbers. In \cite{ns01}, it is proved that $P_n \times P_2$ is exclusive and $\sigma (P_n \times P_2)$ = 3, though characterization of exclusive graphs remains an open problem.

\begin{theorem} [{\cite[Theorem 2.2]{ns01}}] \label{6.1.2}
	{\rm Let $n \geq 2$. The graph $L_n$ = $P_n \times P_2$ is exclusive.  }
\end{theorem}	
\begin{proof} For the graph $L_n$, $\sigma (L_n)$ = 3 by Theorem \ref{6.1.1}. Let $I$ = $\{1,2,\dots,n\}$ be the index set. To have an optimal sum labeling of $L_n$, we must label the vertices of $L_n$ in such a way that $u_i+v_i$ = $p$ for all $i\in I$. Then, $p$ must be a label of an isolated vertex, say $w_1$. Also, for any edge $u_i u_{i+1}$, if $u_i + u_{i+1}$ = $q$ for some $i\in I\setminus \{n\}$, and if $v_j$ = $q$ for some $j\in I$, then clearly, $j\notin [i, i+1]$ since otherwise, $u_{i+1} + p$ = $v_i v_{i+1}$, a contradiction (since $p$ is an isolated vertex). Moreover, if $j$ = $i-1$, then $q$ = $v_j$ = $v_{i-1}$ = $u_i + u_{i+1}$; $p$ = $u_{i} + v_i$ = $u_{i-1} + v_{i-1}$ = $u_i + u_{i+1} + u_{i-1}$ which implies, $u_i + v_{i}$ = $u_{i-1} + u_i + u_{i+1}$. This implies $u_{i-1} + u_{i+1}$ = $v_i$, a contradiction since $u_{i-1} u_{i+1}\notin E(L_n)$. Therefore $j \neq i-1$ also. Now, $u_j$ = $p-v_j$ = $p-q$, $p$ = $u_i+v_i$ = $u_{i+1} + v_{i+1}$. Then, $q+v_i+v_{i+1}$ = $u_i + u_{i+1}$ + $v_i+v_{i+1}$ = $2p$; $v_i+v_{i+1}$ = $2p-q$ = $p+(p-q)$ = $p+u_j$ = $u_i+v_{i}+u_j$, which implies that $u_i+u_{j}$ = $v_{j+1}$, a contradiction since $u_iu_j\notin E(L_n)$. Hence, $P_n \times P_2$ is exclusive. 
\end{proof}

Now, we consider the integral sum number of $L_n$ = $P_n \times P_2$.

\begin{theorem} [{\cite[Theorem 2.3]{ns01}}] \label{6.1.3}
	{\rm For the graph $L_n$ = $P_n \times P_2$, $\zeta (L_n)$ = 3.   }
\end{theorem}	

\begin{proof} Consider the $n$ pairs, -3,5; -7,9; -11,13; $\dots$; $-4n+1$, $4n+1$ such that sum of each pair is 2. Label the vertices $u_i$, $v_i$ of $L_n$ for $i$ = $1,2,\dots,n$ such that $u_i$ = $-4i+1$; $v_i$ = $4i+1$ if $i$ is odd and $u_i$ = $4i+1$; $v_i$ = $-4i+1$ if $i$ is even, respectively. Then, for each edge $u_i v_i$, $u_i + v_i$ = 2 for $i$ = $1,2,\dots,n$. For each edge $u_i u_{i+1}$, $u_i + u_{i+1}$ = 6 or -2 as $i$ is odd or even, respectively, where $i$ = $1,2,\dots,n-1$. Similarly, for each edge $v_i v_{i+1}$, $v_i + v_{i+1}$ = -2 or 6 as $i$ is odd or even, respectively, where $i$ = $1,2,\dots,n-1$. Therefore, by similar argument as in Theorem \ref{6.1.2}, and by Theorem \ref{6.0.6}, we conclude that $\zeta (L_n)$ = 3.
\end{proof}

\begin{cor} [{\cite[Corollary 2.4]{ns01}}] \label{6.1.4}
	{\rm Let $n \geq 2$ and graph $L_n$ = $P_n \times P_2$.  $\sigma (L_n)$ = $\zeta (L_n)$ = 3. \hfill $\Box$ }
\end{cor}	

\subsection{Sum number of $C_{n_1} \cup C_{n_1} \cup \dots \cup C_{n_k}$}

In this subsection, we present the result that $\sigma (mC_n)$ = $\sigma (C_n)$ if $n \neq 4$ and $m \geq 1$. And the main result is $\sigma (mC_n)$ $\neq$ $\zeta (mC_n)$, when $n \neq 4$ and $m \geq 1$.

\begin{theorem} [{\cite[Theorem 3.1]{ns01}}] \label{6.3.5}
	{\rm Let $k\in\mathbb{N}$ and $n_i \neq 4$ for $i$ = $1,2,\dots,k$. $\sigma (C_{n_1} \cup C_{n_1} \cup \dots \cup C_{n_k})$ = 2.  }
\end{theorem}

\begin{proof} When $k$ = 1, the result is true by Theorem \ref{6.0.2}. Let $k \geq 2$. Let $x_{i_j}$ denote the $j^{th}$ vertex of the $i^{th}$ cycle. For convenience, we use the same symbol $x_{i_j}$ for the label of the $(i,j)^{th}$ vertex, where $i$ = $1,2,\dots,k$ and $j$ = $1,2,\dots,n_i$. Let $x_{1_1}$ = $a > 1$, $x_{1_2}$ = $b > a$ be any two positive integers. For $j$ = $3,4,\dots,n_1$, let $x_{1_j}$ = $x_{1_{j-1}}$ + $x_{1_{j-2}}$. For $i$ = $2,3,\dots,k$, let $x_{i_{1}}$ = $x_{{(i-1)}_{{n_{i-1}-1}}}$ + $x_{{(i-1)}_{n_{i-1}}}$; $x_{i_{2}}$ = $x_{{(i-1)}_{n_{i-1}-2}}$ + $x_{{(i-1)}_{{n_{i-1}-1}}}$ and $x_{i_{j}}$ = $x_{i_{j-2}}$ + $x_{i_{j-1}}$ for $j$ = $3,4,\dots,n_i$. Since, $n \neq 4$ for any $i$, by construction, the first cycle determines two labels for two consecutive vertices of the second cycle and the second cycle does the same to the next cycle and so on. It can be easily verified that for any edge $uv\in E(C_{n_1} \cup C_{n_1} \cup \dots \cup C_{n_k})$, $u+v$ = $w$, where $w$ is the next consecutive vertex of $u$ or $v$ in the same cycle, or, first or second vertex of the next cycle. Also, for $1 \leq i \leq m \leq r \leq k$ and $1 \leq j \leq p \leq s \leq n_i$, if $x_{i_{j}}$ + $x_{m_{p}}$ = $x_{r_{s}}$, then $x_{i_{j}}$ + $x_{m_{p}}$ = $x_{i_{j}}$ + $x_{m_{p}}$ + $K$, where $K$ is the sum of the positive multiples of the vertex labels of the intermediate vertices between $x_{i_{j}}$ and $x_{m_{p}}$. This implies that $K$ = 0, which is a contradiction. Therefore, $\sigma (C_{n_1} \cup C_{n_1} \cup \dots \cup C_{n_k})$ = 2 when $n_i \neq 4$ for $i$ = $1,2,\dots,k$. 
\end{proof} 	

\begin{cor} [{\cite[Corollary 3.2]{ns01}}] \label{6.3.6}
	{\rm $\sigma (mC_n)$ = $\sigma (C_n)$ = 2 when $n \neq 4$ and $m \geq 1$.  }
\end{cor}	

\begin{proof} Put $n_1$ = $n_2$ = $\dots$ = $n_k$ = $n$ and $k$ = $m$ in Theorem \ref{6.3.5}, we get the result.
\end{proof} 

In \cite{x99}, Xu obtained the following result.	

\begin{cor} [{\cite[Corollary 3.3]{x99}}] \label{6.3.7}
	{\rm For an arbitrary integer $m \geq 1$, $m K_3$ is an integral sum graph. \hfill $\Box$ }
\end{cor}	

\begin{theorem} [{\cite[Theorem 3.4]{ns01}}] \label{6.3.8}
	{\rm If $G_1$ and $G_2$ are two disjoint graphs and  $\Delta (G_i) < |V(G_i)|-1$ ($i$ = 1,2), then $\zeta (G_1 \cup G_2)$ $\leq$ $\zeta (G_1)$ + $\zeta (G_2)$. \hfill $\Box$ }
\end{theorem}	

\begin{cor} [{\cite[Corollary 3.5]{ns01}}] \label{6.3.9}
	{\rm When $n \neq 4$, $\zeta (mC_n)$ = 0 for any $m \geq 1$.   }
\end{cor}	

\begin{proof} By Theorem \ref{6.0.6}, $\zeta (C_n)$ = 0 if $n \neq 4$. Then, $\zeta (mC_n)$ = 0 if $n > 4$ by Theorem \ref{6.3.8}. Also, $\zeta (mC_3)$ = $\zeta (mK_3)$ = 0 by Theorem \ref{6.3.7}. Hence the result is proved.
\end{proof} 	

The following result is an immediate consequence of corollaries \ref{6.3.6} and \ref{6.3.9}.

\begin{theorem} [{\cite[Theorem 3.6]{ns01}}] \label{6.3.10}
	{\rm When $n \neq 4$, 2 = $\sigma (mC_n)$ $\neq$ $\zeta (mC_n)$ = 0, $m \in\mathbb{N}$.  \hfill $\Box$ }
\end{theorem}	

\subsection{On the integral sum number of caterpillar $S(x_1, x_2, \dots, x_n)$, $x_2 \geq 2$}

Following the result of Ellingham \cite{e93} that $\sigma (T)$ = 1 for any tree $T$, it is clear that $\zeta (T)$ = 0 or 1. integral sum labelling of trees are of most interest since Harary's \cite{h94} open problem to characterize the trees $T$ satisfying $\zeta (T)$ = 0. Harary \cite{h94} also proved that all paths and matchings are integral sum graphs and conjectured that any tree $T$ for which $\zeta (T)$ = 0 is a caterpillar. But, this conjecture was disproved by Xu \cite{x99} by showing that all three path trees are integral sum graphs. Xu also considers a {\em caterpillar} as a tree $T$ = $S(x_1, x_2, \dots, x_n)$, $n \geq 2$, where $T$ has a path $P_n$ = $(a_1, a_2, \dots, a_n)$ of $n$ vertices and $x_i$ leaves (end-edges) are attached to each $a_i$ of $P_n$, $i$ = $1,2,\dots,n$, such that removal of the end edges would give a path $P_n$. Then clearly $x_1 \neq 0$. Here, $|V(T)|$ = $n+ \sum_{i=1}^n{x_i}$. Xu \cite{x99} gives an integral sum labelling to a special class of caterpillar $T$ = $S(1, x_2, x_3, \dots, x_n)$. In \cite{ns01}, this result is extended to more caterpillars $S(x_1, x_2, \dots, x_n)$ for $x_2 \geq 2$. Also, it is noticed that the integral sum labelling given by Xu \cite{x99} fails in a few cases and an alternative labelling to a failed case is provided. Here, we use Lemma \ref{l29}. 

\medskip
We briefly reproduce three properties of $T$ = $S(1, x_2, x_3, \dots, x_n)$ given by Xu \cite{x99} here for further observation and extension.

All nodes of $T$ are arranged in a sequence $\overline{V}$ = $\{c_1, c_2, \dots, c_m\}$ where $m$ = $|V(T)|$ satisfying the following three properties:

(1) Let $c_2$ be such an end vertex of $P_n$ as is adjacent to exactly one pendant-vertex of $T$ and this pendant node is written as $c_1$.

(2) For each $i$ = $1,2,\dots,m-1$, the distance between any two consecutive nodes $c_i$ and $c_{i+1}$ is not more that 2.

(3) For each $k\in \{1,2,\dots,m\}$, the subgraph induced by $A_k$ = $\{c_1, c_2, \dots, c_k\}$ is a subtree of $T$. Then for each $i\in \{2,3,\dots,m\}$, there exists a unique integer (node) $j(i) < i$ such that $c_i$ is adjacent to $c_{j(i)}$ in $T$. Then Xu \cite{x99} proved the following theorem:

\begin{theorem} [{\cite[Theorem 2.2]{x99}}] \label{1.3.12}
	{\rm For any integer $x_i \geq 0$ ($i$ = $1,2,\dots,n$) and $n \geq 2$, $S(1, x_2, x_3, \dots, x_n)$ is an integral sum graph. } \hfill $\Box$  
\end{theorem}	

It is stated that the above result was achieved by the labeling scheme $c_1$ = -3, $c_2$ = 1, $c_3$ = -2, $c_4$ = -1. From $i$ = 5 to $i$ = $m$, $c_i$ = $m_{i-1}-c_{j(i)}$ when $c_{j(i)} > 0$ and $c_i$ = $M_{i-1}-c_{j(i)}$ when $c_{j(i)} < 0$, where $m_{i-1}$ = $\min ~\{c_k/ k \leq i-1\}$, $M_{i-1}$ = $\max ~\{c_k/ k \leq i-1\}$. But it is observed in \cite{ns01} by Lemma \ref{l29} that the above scheme of labeling fails for $S(1, 1, x_3, x_4, \dots, x_n)$ and $S(1, 0, x_2, x_3, \dots, x_n)$. For example, in the tree $S(1,1,2)$ the vertices $c_1,c_2,c_3,c_4,c_5,c_6,c_7$ would receive the labels -3, 1, -2, -1, 3, -6, -9, respectively, according to Xu \cite{x99}. Here, $c_5$ = $-c_1$ and $c_1 + c_6$ = $c_7$; but $c_1 c_6$ is not an edge in $S(1,1,2)$. Therefore, by Lemma \ref{l29}, $S(1,1,2)$ is not an integral sum graph. Now, for $S(1, 0, x_2, x_3, \dots, x_n)$, the following integral sum labeling is given in \cite{ns01}: $c_1$ = 1, $c_2$ = 2, $c_3$ = -1, $c_4$ = 3 and from $i$ = 5 to $i$ = $m$, $c_i$ is defined as in \cite{x99}.

In the following, it is proved that $T$ = $S(x_1, x_2, x_3, \dots, x_n)$ with $x_2 > 1$ is an integral sum graph by exhibiting a labeling scheme which is different from that shown in \cite{x99} and thereby extend the result of Xu \cite{x99}.

\begin{theorem} [{\cite[Theorem 4.3]{ns01}}] \label{1.3.13}
	{\rm Let $n,x_2 \geq 2$ and $x_i \geq 0$ for $i$ = $1,2,\dots,n$. The catapillar $T$ = $S(x_1, x_2, x_3, \dots, x_n)$ is an integral sum graph.  }
\end{theorem}	

\begin{proof} The vertex set $V(T)$ is partitioned into two sets, namely, $V_1$ = $\{b_i/ b_i$ is an end-vertex adjacent to $a_1$, $1 \leq i \leq x_1\}$ and $V_2$ = $V(T) \setminus V_1$. Let $T'$ be the tree induced by the set $V_2$. Label the vertices in $V_1$ as $b_j$ = $-2 \deg (a_1) +j$, where $j$ = $1,2,\dots,x_1$. Arrange the vertices in $V_2$ into a sequence $\overline{V}$ = $\{c_1, c_2, \dots, c_k\}$ satisfying the three properties given in the proof of Theorem \ref{1.3.12} in Xu \cite{x99} with respect to the tree $T'$ with $a_1$ = $c_1$, and $k$ = $|V(T)|-x_1$. Now the vertices of $V_2$ are labelled as $c_1$ = $a_1$ = 1, $c_2$ = $a_2$ = $-\deg (a_1)$, $c_3$ = $c_1+c_2$ and for $i$ = 4 to $i$ = $k$, $c_i$ = $m_{i-1}-c_{j(i)}$ when $c_{j(i)} > 0$ and $c_{j(i)}$ = $M_{i-1}-c_{j(i)}$ when $c_{j(i)} < 0$, where $m_{i-1}$ = $\min ~(\{c_k/ k \leq i-1\} \cup V_1)$ and $M_{i-1}$ = $\max ~(\{c_k/ k \leq i-1\} \cup V_1)$. It can be easily checked that the caterpillar $T$ is an integral sum graph, by the labelling given above. 
\end{proof} 

\section{Perfect graph property of integral sum graphs $G_n$, $G_{0,n}$ and $G_{-r,n}$}

In \cite{js24}, it is proved that the integral sum graphs $G_n$, $G_{0,n}$ and $ G_{-r,n}$ are perfect graphs for $r,n \in \mathbb{N}$. We present these results in this section. 

\subsection{Preliminary results on perfect graphs}

We start with Vizing's theorem on the edge chromatic number of a graph.

\begin{theorem} [{\cite[Theorem 1.1]{v65}}] {\rm (Vizing's Theorem)\label{1.1}\quad 
		For any graph $G$, the edge chromatic number satisfies the inequalities, 
		$\Delta(G)\leq \chi^{'}\left(G \right)\leq \Delta(G)+1$.  \hfill $\Box$  } 
\end{theorem} 

\indent A simple graph $G$ is \emph{class 1} if $\chi^{'}\left(G \right) $ = $\Delta(G)$. It is \emph{class 2} if $\chi^{'}\left(G \right) $ = $\Delta(G)+1$.

\medskip
\indent The \emph{clique} of a graph $G$ is the maximal complete subgraph in $G$ and its order is the \emph{clique number} of $G$, denoted by $\omega(G)$ \cite{h69}. Clearly $\omega\left(K_n \right) $ = $n$. Here we denote the clique of the graph $G$ by $C_G$.

\begin{dfn}	 {\rm \cite{d05}} \quad  A graph $G$ is \emph{perfect} if for every induced subgraph of $G$, the clique number and the chromatic number have the same value. Equivalently, for every $A \subseteq V(G)$, $\chi\left(G[A] \right)$ = $\omega(G[A])$.  
	
\medskip
A graph $G$ is {\em $1$-perfect} if $\chi\left(G \right)$ = $\omega(G)$. It is proved that graph $G$ is perfect if and only if every induced subgraph of $G$ is 1-perfect.
\end{dfn}	 

\begin{theorem} [{\cite[Perfect graph theorem]{l72}}]  \label{thm2.5}{\rm The complement of a perfect graph is perfect as well. \hfill $\Box$} 
\end{theorem}	

Several graphs are proved as perfect \cite{b60, d05, hs58}. They include bipartite graphs and their line graphs, chordal graphs, comparability graphs, triangulated graphs, etc. Perfect graphs arise in the statistical competition of block designs and in graph coloring problems. Another application of perfect graphs is the optimal routing of garbage trucks related to urban science problems. The perfect graph is also closely related to perfect channels in communication theory \cite{t73}.

\indent A vertex and an edge of a graph are said to cover each other if they are incident \cite{h69}. A set of vertices (edges) which covers all the edges (vertices) of a graph $G$ is called a \emph{vertex cover} (\emph{edge cover}) for $G$. \emph{Vertex covering number} (\emph{edge covering number}) is the minimum cardinality among all the vertex covers (edge covers) for $G$ and is denoted by $\alpha_0(G)$ or $\alpha_0$ ($\alpha_1(G)$ or $\alpha_1$) \cite{h69}.

\indent A set of vertices (edges) in $G$ is independent if none of them are adjacent. The maximum cardinality among all the vertex (edge) independent set is called its \emph{vertex independence number} (\emph{edge independence number}) and is denoted by $\beta_0(G)$ or $\beta_0$ ($\beta_1(G)$ or $\beta_1$) \cite{h69}. We denote the independent set of the graph $G$ by $I_G$. 

\begin{theorem} [{\cite[Gallai Theorem]{g59}}] {\rm 	For any connected nontrivial graph $G$ of order $n$, $\alpha_0(G) + \beta_0(G)$ = $n$ = $\alpha_1(G) + \beta_1(G)$.   \hfill $\Box$  }  
\end{theorem}	

\begin{lemma}[{\cite[Lemma 2.3]{js24}}] \label{resGn} \quad {\rm 	For every $n \in \mathbb{N}$, the integral sum graph $G_n$ has 
		\begin{enumerate}
			\item [\rm (i)] Clique of $G_n$, $C_{G_n}$ = $\left\lbrace v_i \in V(G_n):{i} = 1, 2, 3, \dots,\left \lceil{\frac{n}{2}} \right\rceil \right\rbrace$  and $\omega(G_n)$ = $\left \lceil{\frac{n}{2}} \right\rceil $.
			
			\item [\rm (ii)] Minimum vertex cover of $G_{n}$ = $ \left\lbrace v_i \in V(G_n): i = 1, 2, 3, \dots,\left \lceil{\frac{n}{2}} \right\rceil - 1 \right\rbrace$
			
			\indent and ${\alpha_0}({G_n})$ = $\left \lceil{\frac{n}{2}} \right\rceil - 1$.
			\item [\rm (iii)] Maximum independent set of $G_{n}$, 
			
			\hspace{1cm} $I_{G_n}$ = $\left \lbrace v_i \in V(G_n):{i} = \left \lceil{\frac{n}{2}}\right \rceil, \left \lceil{\frac{n}{2}}\right \rceil + 1, \left \lceil{\frac{n}{2}}\right \rceil+ 2	,\dots, n \right \rbrace$ 
			
			\indent and ${\beta_0}{(G_n)} $ = 	$ \left \lfloor{\frac{n}{2}}\right \rfloor + 1. $ 
			\item [\rm (iv)] Maximum matching of $G_{n}$ 
			
			\hfill	= $ \{ v_iv_j \in E(G_n):(i,j)= \left( 1,n-1\right), \left(2,n-2\right), \dots$, $\left(\left \lceil{\frac{n}{2}} \right\rceil - 1, \left \lceil{\frac{n}{2}}\right \rceil\right) \} $ 			
			
			\indent and $\beta_1(G_{n})$ = 	$\left \lceil{\frac{n}{2}} \right\rceil - 1$.
			\item [\rm (v)] Any sum graph $G^+(S)$ has no edge cover since the vertex corresponding to the biggest label which is a natural number is an isolated vertex in $G^+(S)$. In particular, $G_n$ has no edge cover, $n\in\mathbb{N}$. \hfill $\Box$
	\end{enumerate}	}
\end{lemma}

\begin{lemma} [{\cite[Lemma 2.4]{js24}}] \label{resGrs} \quad {\rm 	For every integral sum graph $G_{-r,n}$ and $r,n \in \mathbb{N}$,  
		\begin{enumerate}
			\item [\rm (i)] Clique of $G_{-r,n}$,  
			
			\hspace{1cm} $C_{G_{-r,n}}$	= $\left\lbrace v_i \in V({G_{-r,n}}):{i} = 0, 1, 2, \dots,\left \lceil{\frac{n}{2}} \right\rceil, -1, -2, \dots, -\left \lceil{\frac{r}{2}} \right\rceil\right\rbrace$
			
			\noindent and $\omega({G_{-r,n}})$ = $ 1 + \left \lceil{\frac{r}{2}}  \right\rceil + \left \lceil{\frac{n}{2}} \right\rceil $.
			
			\item [\rm (ii)]	Minimum vertex cover of $G_{-r,n}$  
			
			\small \hspace{2cm} =   $\begin{cases} \text{$\left\lbrace v_i \in V(G_{-r,n}):i = -r,\dots, -1, 0, 1, \dots,\left\lceil{\frac{n}{2}}\right \rceil - 1 \right\rbrace $} & \text{if $r\leq n$}\\
				\text{$\left\lbrace v_i \in V(G_{-r,n}):i = -\left\lceil{\frac{r}{2}}\right \rceil+1, \dots, -1, 0, 1, \dots,n\right\rbrace $} & \text{if $r\geq n$,}
			\end{cases}$
			
			\indent and ${\alpha_0} {(G_{-r,n})}$ = $ \min \{r,n\} + \left \lceil{\frac{\max ~\{r,n\}}{2}}\right \rceil $.
			
			\item [\rm (iii)] Maximum independent set of $G_{-r,n}$ 
			
			\hspace{2cm}	= $\begin{cases} \text{$ \left\lbrace v_i \in V({G_{-r,n}}):i =\left\lceil{\frac{n}{2}}\right \rceil, \left\lceil{\frac{n}{2}}\right \rceil + 1, \dots, n-1, n\right\rbrace$} & \text{if $r\leq n$}\\					
				\text{$\left\lbrace v_i \in V({G_{-r,n}}):i = -\left\lceil{\frac{r}{2}}\right \rceil, -\left(\left\lceil{\frac{r}{2}}\right \rceil +1\right), \dots, -r \right\rbrace$} & \text{if $r\geq n$,}\\
			\end{cases}$ 
			
			\indent and $\beta_0(G_{-r,n})$ = $\left\lfloor {\frac{\max~ \{r,n\}}{2}} \right\rfloor + 1$.
			
			\item [\rm (iv)] Minimum edge cover of $G_{-r,n}$
			
		\hspace{2cm}	= \small$\begin{cases} \text{$\left\lbrace v_iv_j \in E(G_{-r,n}) : i + j = n-r \right\rbrace$} & \text{if $k$ is even}\\
				\text{$\left\lbrace v_iv_j \in E(G_{-r,n}) : i + j = n-r ~\&~ \left(i, j\right) = \left(0, \frac{n-r}{2}\right) \right\rbrace$} & \text{if $k$ is odd}.
			\end{cases}$ 
			\\	 	
			\indent
			\hspace{2cm} {\tiny $= \begin{cases}  \text{$\left\lbrace v_iv_j \in E(G_{-r,n}) : (i,j)= (-r,n),(-r+1,n-1),\dots, (0,n-r)\right\rbrace$}  \hfill  \text{if $k$ is even}\\
				{\tiny\text{$\left\lbrace v_iv_j \in E(G_{-r,n}) :(i,j)= (-r,n),(-r+1,n-1),\dots, (0,n-r)\& \left(0, \frac{n-r}{2}\right)\right\rbrace$} } \hfill \text{ if $k$ is odd},
				\end{cases} $}\\
			
			\indent and ${\alpha_1} {(G_{-r,n})}$ = $\left\lceil {\frac{r + n + 1}{2}} \right\rceil $ where $k$ = $r+n+1$.			
			\item [\rm (v)] \small Maximum matching of $G_{-r,n}$ 
			
			\hspace{2cm}	= $\left\lbrace v_iv_j \in E(G_{-r,n}): i + j = n-r \right\rbrace$,
			
			\indent and	$\beta_1(G_{-r,n})$ = $\left\lfloor {\frac{r + n + 1}{2}} \right\rfloor $.			\hfill   $\Box$ 	
	\end{enumerate}}	
\end{lemma}

\begin{note} [{\cite[Note 2.5]{js24}}] \quad {\rm By putting $r=0$ in Lemma \ref{resGrs}, we obtain corresponding properties of $G_{0,n}$, $n \in \mathbb{N}$. }
\end{note}

\begin{note} [{\cite[Note 2.6]{js24}}] \quad {\rm In the integral sum graph $G_{-r,n}$ of even order, every maximum matching is a perfect matching for $r,n \in \mathbb{N}_0$. }
\end{note}

\subsection{$G_n$ and $G^c_n$ are perfect graphs}

Here, we present the result that integral sum graphs $G_n$ and $G^c_n$ are perfect graphs. Let us start with preliminaries that are required in this section. For any integral sum graph $G_n$, $G_{0,n}$ or $G_{-r,n}$, maximum independent set, minimum vertex cover, minimum edge cover, and maximum matching need not be unique but clique is unique, $r,n \in \mathbb{N}$. 

\begin{theorem} [{\cite[Theorem 3.1]{js24}}] \label{Gnpf} \quad {\rm 	For a positive integer $n$, the sum graph $G_n$ is perfect.}
\end{theorem}

\begin{proof} 	Let $G_n$ be the sum graph of order $n$ obtained by the labeling $f : V(G_n) \rightarrow \mathbb{N}$ defined by $f(v_i) = i, 1\leq i\leq n$. Let $H$ be any induced subgraph of $G_n$. Our aim is to prove that $\chi\left(H \right) $= $\omega\left(H \right)$ for all subgraphs $H$ of $G_n$ where $\chi\left(H \right) $ is the chromatic number of $H$ and $\omega\left(H \right) $ is its clique number.
	
	$G_n$ is a split graph with clique ${C_{G_n}}$ = $\left\lbrace v_i \in V(G_n):{i} = 1, 2, \dots,\left \lceil{\frac{n}{2}} \right\rceil \right\rbrace$ and the independent set ${I_{G_n}} = \left \lbrace v_i \in V(G_n):{i} = \left \lceil{\frac{n}{2}}\right \rceil, \left \lceil{\frac{n}{2}}\right \rceil + 1, \left \lceil{\frac{n}{2}}\right \rceil+ 2	,\dots, n \right \rbrace$. To prove the theorem we consider the following three cases of $H$.
	
	\medskip
	\noindent\textit{Case 1.}\quad $V(H) \subseteq V({C_{G_n}})$.
	
	Every induced subgraph of a clique is another clique of order fewer than or equal to the maximal clique. And thereby $H$ itself a clique and $\chi\left(H \right) $= $\omega\left(H \right)$.
	
\medskip
	\noindent\textit{Case 2.}\quad $V(H) \subseteq V({I_{G_n}})$.\\
	\indent Induced subgraph of an independent set is also an independent set and the clique number of any induced subgraph of an independent set is $1$. This implies, $\chi\left(H \right) $= $\omega\left(H \right) $ = $1$ for every $V(H) \subseteq V({I_{G_n}})$.
	
\medskip
	\noindent\textit{Case 3.}\quad $V(H) \subseteq V({C_{G_n}}) \cup V({I_{G_n}})$.
	
	\indent In this case, vertices of $H$ belong to either ${C_{G_n}}$ or ${I_{G_n}}$. Let $V(H)$ = $V({C_H})\cup V({I_H})$ where $V({C_H})$ $\subseteq$ ${C_{G_n}}$ and $V({I_H})$$\subseteq$ ${I_{G_n}}$, $V({C_H})$, $V({I_H})$ $\neq \emptyset $. Let $\omega\left(H \right) = k$. This implies that $|C_H|$ = $k$ = $\chi\left(C_H \right) $ where $C_H$ is the clique in $H$. If $H$ itself is a clique, then $H$ = $C_H$ and thereby $\chi\left(H \right) $ = $\omega\left(H \right)$ = $k$. Hence the result is true in this case. If $H$ is not a clique, let $\omega\left(H \right)$ = $k$. This implies,  $\chi\left(C_H \right) $ = $|C_H|$ = $k$. 
	
	For each vertex $v_j \in V({I_H})$, there exists atleast one $v_i \in V({C_H})$ such that $v_{i}v_{j}$ $\notin E(H)$ (since $i+j > n$) so that there will not be a clique of higher order than that of ${C_H}$ and thus $v_{i}$ and $v_{j}$ can be assigned with the same color. That is, the colors used in ${C_H}$ are enough for coloring $V({I_H})$. This implies, $\chi\left(H \right) $= $\chi\left(C_H \right) $ = $k$. Hence we get $\chi\left(H \right) $ = $k$ = $\omega\left(H \right)$. 
	
	\medskip
	Thus in all the three cases we obtain $\chi\left(H \right) $ = $\omega\left(H \right)$ for every induced subgraph $H$ of $G_n$  and thereby graph $G_n$ is perfect for $n\in\mathbb{N}$. 		
\end{proof}	

\begin{illu} [{\cite[Illustration 3.2]{js24}}] \label{Gnillu} \nonumber\rm
	Consider the sum graph $G_{13}$ as given in Figure 39. Vertices of $G_{13}$ are subscript-labeled. Vertices $v_1, v_2, v_3, v_4, v_5, v_6$ and $v_7$ (vertices joined by edges with blue color) form a clique in $G_{13}$ and $\omega(G_{13})= 7$. These vertices are colored by the colors $c_1, c_2,\dots, c_7$. All the remaining vertices are nonadjacent to $v_7$ and so, they can be colored by $c_7$. Or, the vertices $v_8$,$v_9$,$v_{10}$, $v_{11}$,$v_{12}$ and $v_{13}$ are  nonadjacent to the vertices $v_6$,$v_5$,$v_4$,$v_3$,$v_2$ and $v_1$ respectively, and can be colored by the corresponding colors  $c_1, c_2,\dots, c_7$. Therefore, the chromatic number of the sum graph $G_{13}$ is 7. i.e., $\chi\left(G_{13} \right)$ = 7 = $\omega(G_{13})$. Thus $G_{13}$ is 1-perfect. Clearly every pair of nonadjacent vertices in $G_{13}$ are also nonadjacent in any of its induced subgraphs. This shows that every induced subgraph of $G_{13}$ is 1-perfect. This implies, graph $G_{13}$ is 1-perfect and thereby $G_{13}$ is perfect.
	
\end{illu}

\begin{cor} [{\cite[Corollary 3.3]{js24}}] \quad {\rm For every positive integer $n$, graph $G_{n}^c$ is perfect.}	
\end{cor}  
\begin{proof} The proof follows from Theorems \ref{Gnpf} and \ref{thm2.5}.	
\end{proof}

\begin{theorem} [{\cite[Theorem 3.4]{js24}}] \label{thm Gn chrom} \quad {\rm Let $n \in \mathbb{N}$.
		\begin{enumerate}
			\item [\rm (i)]  The chromatic number of $G_n$ is $\left \lceil{n/2} \right\rceil$.
			
			\item [\rm (ii)] The edge chromatic number of $G_1$ is $0$ and for $n \geq 2$, the edge chromatic number of $G_n$ is $n-2$.
			
			i.e., $\chi^{'}\left( G_1 \right)$ = 0 and for $n \geq 2$, ~$\chi^{'}\left( G_n \right)$ = $n-2$.
			\end{enumerate}}
\end{theorem}	

\begin{proof}\quad 		The sum graph $G_n$ is obtained by the labeling $f : V(G_n) \rightarrow \mathbb{N}$ defined by $f(v_i) = i$ for $i$ = 1 to $n$.
	\\
	(i)~ From Theorem \ref{Gnpf}, $G_n$ is perfect and thereby the chromatic number of $G_n$ is equal to its clique number. This implies, $\chi\left( G_n \right)$ = $\omega(G_n)$ = $\left \lceil{\frac{n}{2}} \right\rceil$. Hence we get the result (i). Vertex coloring and clique of $G_{13}$ are shown in  Figure  39.
	\\
	(ii)~ For $n = 1$, $G_n$ = $G_1$ and so the edge chromatic number $\chi^{'}\left(G_1 \right)$ = $0$. 
	
	For $n \geq 2$, the maximum degree of graph $G_n$ is $\Delta(G_n)$ = $n-2$. Using Theorem \ref{1.1}, we have, $\Delta (G_n)$ $\leq$ ${\chi}^{'}(G_n) \leq \Delta(G_n)+1$. Here, we present the edge coloring of the sum graph $G_n$ with exactly $n-2$ colors so that ${\chi}^{'}(G_n) = n-2$ for $n \geq 2$.
	
	\indent Let $c_k$ denote $k^{th}$ color assigned to an edge and $C_k$ denote the color class of edges each with color $c_k$ in $G_n$, $n\geq 2$. Color the set of edges $\{v_iv_j\in E(G_n) : i+j = k+2,1\leq i,j\leq n, i\neq j\}$ of $G_n$ with the color $c_k$. 
	
	\indent It is clear from the above edge coloring that colors $c_1$ to $c_{n-2}$ are assigned to the edges of  $G_n$ and no more edge colors are required. And all the colors of edges incident at each vertex $v_i$ of sum graph $G_n$ are all distinct since there are exactly $n-2$ possibilities of $j$, $i\neq j$, $1 \leq j \leq n$ for which $i+j = k+2$ in $G_n$, $1 \leq k \leq n-2$. This implies, $\chi^{'}\left( G_n \right)$ = $n-2$ for $n \geq 2$. Hence we get the result (ii).
\end{proof}

\begin{illu} [{\cite[Illustration 3.5]{js24}}] \label{illuGn}\nonumber\rm
	Consider the sum graph $G_{13}$ as given in Figure 40. Vertices of $G_{13}$ are subscript-labeled. The vertices $v_1, v_2, v_3, v_4, v_5, v_6$ and $v_7$ in $G_{13}$ (vertices joined by edges with blue color) form a clique, $\omega(G_{13})$ = 7 and $\Delta(G_{13})$ = 11. Since $G_{13}$ is perfect, $\chi\left(G_n \right)$ = $\omega(G_{13})$ which implies,  $\chi\left(G_n \right)$ = 7 = $\left \lceil{\frac{13}{2}} \right\rceil $.
	As in the proof, color the set of edges $\{v_i,v_j\in E(G_n) : i+j = k+2$, $1 \leq i,j\leq n, ~i\neq j\}$ of $G_n$ with the color $c_k$. The colors of each edge of $G_{13}$ is given as follows.  \\
	${\color{blue}c_1} \rightarrow v_1v_2$;\\
	${\color{black}c_2} \rightarrow v_1v_3$;\\
	${\color{red}c_3} \rightarrow v_1v_4, v_2v_3$;\\
	${\color{gray}c_4} \rightarrow v_1v_5, v_2v_4$;\\
	${\color{orange}c_5} \rightarrow v_1v_6, v_2v_5, v_3v_4$;\\
	${\color{olive}c_6} \rightarrow v_1v_7, v_2v_6, v_3v_5$;\\
	${\color{magenta}c_7} \rightarrow v_1v_8 ,v_2v_7, v_3v_6, v_4v_5$;\\
	${\color{violet}c_8} \rightarrow v_1v_9, v_2v_8, v_3v_7, v_4v_6$;\\
	${\color{brown}c_9} \rightarrow v_1v_{10}, v_2v_9, v_3v_8, v_4v_7, v_5v_6$;\\
	${\color{purple}c_{10}} \rightarrow v_1v_{11}, v_2v_{10}, v_3v_9, v_4v_8, v_5v_7$;\\
	${\color{green}c_{11}} \rightarrow v_1v_{12}, v_2v_{11}, v_3v_{10}, v_4v_9, v_5v_8, v_6v_7$.\\
	Edge coloring of $G_{13}$ with 11 colors, using the method given in the proof of Theorem \ref{thm Gn chrom}, is given in Figure 40. Here, $\chi^{'}\left(G_{13} \right) = 11 = 13-2$.
\end{illu}

\begin{figure}
	\centering
	\begin{minipage}{0.42\textwidth}
		\centering
		\resizebox{1\textwidth}{!}{%
			\begin{tikzpicture}[scale =0.5]
				\node (c1) at (7,1.8) [circle,draw,scale=0.6,fill=magenta] {1};
				\node (c2) at (7.5,4.1) [circle,draw,scale=0.6,fill=orange] {2};
				\node (c3) at (9.2,5.9) [circle,draw,scale=0.6,fill=cyan] {3};
				\node (c4) at (11.5,6.6) [circle,draw,scale=0.6,fill=pink] {4};
				\node (c5) at (13.8,5.9) [circle,draw,scale=0.6,fill=yellow] {5};
				\node (c6) at (15.5,4.1) [circle,draw,scale=0.6, fill=lightgray] {6};
				\node (c7) at (16,1.8) [circle,draw,scale=0.6,fill=green] {7};
				\node (c8) at (15.6,0.1) [circle,draw,scale=0.6,fill=green] {8};
				\node (c9) at (14.5,-1.3) [circle,draw,scale=0.6,fill=green] {9};
				\node (c10) at (12.9,-2.1) [circle,draw,scale=0.5,fill=green] {10};
				\node (c11) at (11,-2.3) [circle,draw,scale=0.5,fill=green] {11};
				\node (c12) at (9.1,-1.7) [circle,draw,scale=0.5,fill=green] {12};
				\node (c13) at (7.5,-0.1) [circle,draw,scale=0.5,fill=green] {13};
				
				\draw (c1)[blue, thick] --node[near start] {$ $} (c2);
				\draw (c1)[blue, thick] --node[near start] {$ $} (c3);
				\draw (c1)[blue, thick] --node[near start] {$ $} (c4);
				\draw (c1)[blue, thick] --node[near start] {$ $} (c5);
				\draw (c1)[blue, thick] --node[near start] {$ $} (c6);
				\draw (c1)[blue, thick] --node[near start] {$ $} (c7);
				\draw (c1)[, thick] --node[near start] {$ $} (c8);
				\draw (c1)[, thick] --node[near start] {$ $} (c9);
				\draw (c1)[, thick] --node[near start] {$ $} (c10);
				\draw (c1)[, thick] --node[near start] {$ $} (c11);
				\draw (c1)[, thick] --node[near start] {$ $} (c12);
				
				\draw (c2)[blue, thick] --node[near start] {$ $} (c3);
				\draw (c2)[blue, thick] --node[near start] {$ $} (c4);
				\draw (c2)[blue, thick] --node[near start] {$ $} (c5);
				\draw (c2)[blue, thick] --node[near start] {$ $}(c6);
				\draw (c2)[blue, thick] --node[near start] {$ $} (c7);
				\draw (c2)[, thick] --node[near start] {$ $} (c8);
				\draw (c2)[, thick] --node[near start] {$ $} (c9);
				\draw (c2)[, thick] --node[near start] {$ $}(c10);
				\draw (c2)[, thick] --node[near start] {$ $} (c11);
				
				\draw (c3)[blue, thick] --node[near start] {$ $}(c4);
				\draw (c3)[blue, thick] --node[near start] {$ $}(c5);
				\draw (c3)[blue, thick] --node[near start] {$ $} (c6);
				\draw (c3)[blue, thick] --node[near start] {$ $} (c7);
				\draw (c3)[, thick] --node[near start] {$ $} (c8);
				\draw (c3)[, thick] --node[near start] {$ $} (c9);
				\draw (c3)[, thick] --node[near start] {$ $} (c10);
				
				\draw (c4)[blue, thick] --node[near start] {$ $} (c5);
				\draw (c4)[blue, thick] --node[near start] {$ $}(c6);
				\draw (c4)[blue, thick] --node[near start] {$ $}(c7);
				\draw (c4)[, thick] --node[near start] {$ $} (c8);
				\draw (c4)[, thick] --node[near start] {$ $} (c9);
				
				\draw (c5)[blue, thick] --node[near start] {$ $} (c6);
				\draw (c5)[blue, thick] --node[near start] {$ $} (c7);
				\draw (c5)[, thick] --node[near start] {$ $} (c8);
				
				\draw (c6)[blue, thick] --node[near start] {$ $} (c7);
		\end{tikzpicture} }%
		
		{\small	Fig. 39. $G_{13}:$ vertex coloring $\&$ clique}
		\label{$G_{13}:$ vertex coloring and clique}
	\end{minipage}\hfill
	\begin{minipage}{0.42\textwidth}
		
		\centering
		\resizebox{1\textwidth}{!}{%
			\begin{tikzpicture}[scale =0.7]
				
				\node (c1) at (7,1.8) [circle,draw,scale=0.8,fill=pink] {1};
				\node (c2) at (7.5,4.1) [circle,draw,scale=0.8,fill=pink] {2};
				\node (c3) at (9.2,5.9) [circle,draw,scale=0.8,fill=pink] {3};
				\node (c4) at (11.5,6.6) [circle,draw,scale=0.8,fill=pink]{4};
				\node (c5) at (13.8,5.9) [circle,draw,scale=0.8,fill=pink]{5};
				\node (c6) at (15.5,4.1) [circle,draw,scale=0.8,fill=pink] {6};
				\node (c7) at (16,1.8) [circle,draw,scale=0.8,fill=pink] {7};
				\node (c8) at (15.6,0.1) [circle,draw,scale=0.8,fill=pink] {8};
				\node (c9) at (14.5,-1.3) [circle,draw,scale=0.8,fill=pink] {9};
				\node (c10) at (12.9,-2.1) [circle,draw,scale=0.7,fill=pink] {10};
				\node (c11) at (11,-2.3) [circle,draw,scale=0.7,fill=pink] {11};
				\node (c12) at (9.1,-1.7) [circle,draw,scale=0.7,fill=pink] {12};
				\node (c13) at (7.5,-0.1) [circle,draw,scale=0.7,fill=pink] {13};
				
				\draw (c1)[blue, thick] --node[near start] {$ $} (c2);
				\draw (c1)[, thick] --node[near start] {$ $} (c3);
				\draw (c1)[red, thick] --node[near start] {$ $} (c4);
				\draw (c1)[gray, thick] --node[near start] {$ $} (c5);
				\draw (c1)[orange, thick] --node[near start] {$ $} (c6);
				\draw (c1)[olive, thick] --node[near start] {$ $} (c7);
				\draw (c1)[ magenta, thick] --node[near start] {$ $} (c8);
				\draw (c1)[violet, thick] --node[near start] {$ $} (c9);
				\draw (c1)[ brown, thick] --node[near start] {$ $} (c10);
				\draw (c1)[purple, thick] --node[near start] {$ $} (c11);
				\draw (c1)[green, thick] --node[near start] {$ $} (c12);
				
				\draw (c2)[red, thick] --node[near start] {$ $} (c3);
				\draw (c2)[gray, thick] --node[near start] {$ $} (c4);
				\draw (c2)[orange, thick] --node[near start] {$ $} (c5);
				\draw (c2)[olive, thick] --node[near start] {$ $}(c6);
				\draw (c2)[magenta, thick] --node[near start] {$ $} (c7);
				\draw (c2)[ violet, thick] --node[near start] {$ $} (c8);
				\draw (c2)[brown, thick] --node[near start] {$ $} (c9);
				\draw (c2)[purple, thick] --node[near start] {$ $}(c10);
				\draw (c2)[green, thick] --node[near start] {$ $} (c11);
				
				\draw (c3)[orange, thick] --node[near start] {$ $}(c4);
				\draw (c3)[olive, thick] --node[near start] {$ $}(c5);
				\draw (c3)[magenta, thick] --node[near start] {$ $} (c6);
				\draw (c3)[violet, thick] --node[near start] {$ $} (c7);
				\draw (c3)[brown, thick] --node[near start] {$ $} (c8);
				\draw (c3)[purple, thick] --node[near start] {$ $} (c9);
				\draw (c3)[ green, thick] --node[near start] {$ $} (c10);
				
				\draw (c4)[magenta, thick] --node[near start] {$ $} (c5);
				\draw (c4)[violet, thick] --node[near start] {$ $}(c6);
				\draw (c4)[brown, thick] --node[near start] {$ $}(c7);
				\draw (c4)[ purple, thick] --node[near start] {$ $} (c8);
				\draw (c4)[green, thick] --node[near start] {$ $} (c9);
				
				\draw (c5)[brown, thick] --node[near start] {$ $} (c6);
				\draw (c5)[purple, thick] --node[near start] {$ $} (c7);
				\draw (c5)[green, thick] --node[near start] {$ $} (c8);
				
				\draw (c6)[green, thick] --node[near start] {$ $} (c7);
			\end{tikzpicture} 		
		}%
		
		{\small	Fig. 40. $G_{13}:$ edge coloring}
		\label{$G_{13}:$ edge coloring.}
	\end{minipage}
\end{figure}

\subsection{Join of two perfect graphs is a perfect graph}

Results are available for the composition of perfect graphs which allows different types of graph operations  \cite{cc85}. But the graph operation join is not available yet. Following is the result corresponding to join of two perfect graphs. 

\begin{theorem} [{\cite[Theorem 3.6]{js24}}] \label{thm join} \quad {\rm Join of two perfect graphs is also perfect.}
\end{theorem}

\begin{proof}
	Let $G$ and $F$ be two perfect graphs with $\chi\left(G \right) $ = $\omega(G)$ = $p$ and $\chi\left(F \right) $ = $\omega(F)$ = $q$. Also let $J$ = $G * F$, the join of the graphs $G$ and $F$.
	When we take the join of two graphs $G$ and $F$, the clique in $G$ together with clique in $F$ form a higher clique of order $p+q$. i.e., $\omega(J)$ = $p+q$. Let $H$ be any induced subgraph of $J$. We have to prove that $\chi\left(H \right) $= $\omega\left(H \right)$.\\ 
	\indent Let $\omega(H)$ = $k$ and let $C_H$ be the clique in $H$ of order $k$. Then $\omega(H)$ = $|C_H|$ = $k$ = $\chi(C_H)$. That is $C_H$ can be colored using $k$ colors.\\ 
	\indent Let $D = V(H)\setminus V(C_H)$. For each $v_j$ $\in$ $D \cap V(G)$
	and $v_y$ $\in$ $D \cap V(F)$, there exist $v_i$ $\in$ $ V(C_H)\cap V(G)$ and $v_x$ $\in$ $ V(C_H)\cap V(F)$ such that $ v_iv_j,v_xv_y \notin E(H)$ so that there does not exist a clique of higher order than that of ${C_H}$. Hence, the colors of $v_i$ and $v_x$ can be assigned to $v_j$ and $v_y$, respectively. And thereby the vertices of $D$ can be colored  using the colors of vertices of $C_H$. This implies, $\chi(H)$ = $\chi(C_H)$ = $\omega(H)$. Thus, for every induced subgraph $H$ of $J$, $\chi\left(H \right) $ = $\omega\left(H \right)$. Thus the graph $J$ is perfect.
\end{proof}	

\subsection{Integral sum graphs $G_{0, n}$ and $G_{-m, n}$ are perfect graphs}

Results related to integral sum graphs $G_{0, n}$ and $G_{-m, n}$ are perfect is presented here, $m,n\in\mathbb{N}$.

\begin{theorem} [{\cite[Theorem 3.7]{js24}}] \label{thm G0n} \quad {\rm 		For every $n \in \mathbb{N}$, the integral sum graph $G_{0,n}$ is perfect.}
\end{theorem}

\begin{proof}
	The integral sum graph $G_{0,n}$ = $K_1 * G_n$ where $K_1$ = $G_1$ is labeled with 0. Then using Theorems \ref{Gnpf} and \ref{thm join}, the graph $G_{0,n}$ is perfect, $n \in \mathbb{N}$. 
\end{proof} 
		
	\begin{theorem} [{\cite[Theorem 3.8]{js24}}] \label{thm G0n chrom}\quad {\rm 
			For every $n \in \mathbb{N}$, the integral sum graph $G_{0,n}$ has
			\\
			(i)~ ~the chromatic number $\chi\left(G_{0,n} \right)$ = $ \left \lceil{\frac{n}{2}} \right\rceil +1 $;
			\\
			(ii)~ the edge chromatic number $\chi^{'}\left(G_{0,n} \right)$ = $ n $.}
	\end{theorem}

	\begin{proof}
		The integral sum graph $G_{0,n}$ is obtained by the labeling $f : V(G_{0,n}) \rightarrow \mathbb{N} \cup \{0\}$ defined by $f(v_i) = i$, $i$ = 0 to $n$.
		\\
		(i)~ We have, $G_{0,n} \cong K_1 * G_n$ which implies, $\chi(G_{0,n})$ = $\chi(K_1)$ + $\chi(G_{n})$ = 1+$ \left \lceil{\frac{n}{2}} \right\rceil$ using Theorem \ref{thm Gn chrom}. Vertex coloring and clique of $G_{0, 12}$ are shown in Figure  41.  
		\\
		(ii)~	For $n \in \mathbb{N}$, the maximum degree of integral sum graph $G_{0,n}$ is $\Delta(G_{0,n})$ = $n$.	Using Theorem \ref{1.1}, we get, $\Delta (G_{0,n}) \leq {\chi}^{'}(G_{0,n}) \leq \Delta(G_{0,n})+1$. Now, we present a proper edge coloring of $G_{0,n}$ with $n$ colors so that  ${\chi}^{'}(G_{0,n})$ = $n$.
		
		\indent Let $c_k$ denote $k^{th}$ color assigned to an edge and $C_k$ denote the color class of edges, each with color $c_k$ in  $G_{0,n}$, $n \in \mathbb{N}$. Color all the edges in the set $\{v_iv_j\in E(G_{0,n}) : i+j = k,~0\leq i,j\leq n, ~i\neq j, 1\leq k\leq n\}$ with the same color $c_k$. Clearly, the edges of $G_{0,n}$ take at the most $n$ colors since $1\leq i+j\leq n$. Also, there are $n$ edges incident at the vertex $v_0$, and all these $n$ edges take $n$ distinct colors. That is, $ \chi^{'} (G_{0,n})$ = $n$, $n \in \mathbb{N}$. Edge coloring of $G_{0,12}$ with 12 colors, using the method given in the proof of Theorem \ref{thm G0n chrom}, is given in Figure 42.
	\end{proof} 
	
	\begin{figure}
		\centering
		\begin{minipage}{0.41\textwidth}
			\centering
			\resizebox{1\textwidth}{!}{%
				\begin{tikzpicture}[scale =0.6]
					\node (c0) at (7.5,-0.1) [circle,draw,scale=0.9,fill=white] {0};
					\node (c1) at (7,1.8) [circle,draw,scale=0.9,fill=purple!50] {1};
					\node (c2) at (7.5,4.1) [circle,draw,scale=0.9,fill=orange] {2};
					\node (c3) at (9.2,5.9) [circle,draw,scale=0.9,fill=cyan!80] {3};
					\node (c4) at (11.5,6.6) [circle,draw,scale=0.9,fill=pink]{4};
					\node (c5) at (13.8,5.9) [circle,draw,scale=0.9,fill=yellow]{5};
					\node (c6) at (15.5,4.1) [circle,draw,scale=0.9,fill=lightgray] {6};
					\node (c7) at (16,1.8) [circle,draw,scale=0.9,fill=lightgray] {7};
					\node (c8) at (15.6,0.1) [circle,draw,scale=0.9,fill=lightgray] {8};
					\node (c9) at (14.5,-1.3) [circle,draw,scale=0.9,fill=lightgray] {9};
					\node (c10) at (12.9,-2.1) [circle,draw,scale=0.7,fill=lightgray] {10};
					\node (c11) at (11,-2.3) [circle,draw,scale=0.7,fill=lightgray] {11};
					\node (c12) at (9.1,-1.7) [circle,draw,scale=0.7,fill=lightgray] {12};
					
					\draw (c0)[magenta, thick] --node[near start] {$ $} (c1);
					\draw (c0)[magenta, thick] --node[near start] {$ $} (c2);
					\draw (c0)[magenta, thick] --node[near start] {$ $} (c3);
					\draw (c0)[magenta, thick] --node[near start] {$ $} (c4);
					\draw (c0)[magenta, thick] --node[near start] {$ $} (c5);
					\draw (c0)[magenta, thick] --node[near start] {$ $} (c6);
					\draw (c0)[, thick] --node[near start] {$ $} (c7);
					\draw (c0)[, thick] --node[near start] {$ $} (c8);
					\draw (c0)[, thick] --node[near start] {$ $} (c9);
					\draw (c0)[, thick] --node[near start] {$ $} (c10);
					\draw (c0)[, thick] --node[near start] {$ $} (c11);
					\draw (c0)[, thick] --node[near start] {$ $} (c12);
					
					\draw (c1)[magenta, thick] --node[near start] {$ $} (c2);
					\draw (c1)[magenta, thick] --node[near start] {$ $} (c3);
					\draw (c1)[magenta, thick] --node[near start] {$ $} (c4);
					\draw (c1)[magenta, thick] --node[near start] {$ $} (c5);
					\draw (c1)[magenta, thick] --node[near start] {$ $} (c6);
					\draw (c1)[, thick] --node[near start] {$ $} (c7);
					\draw (c1)[, thick] --node[near start] {$ $} (c8);
					\draw (c1)[, thick] --node[near start] {$ $} (c9);
					\draw (c1)[, thick] --node[near start] {$ $} (c10);
					\draw (c1)[, thick] --node[near start] {$ $} (c11);
					
					\draw (c2)[magenta, thick] --node[near start] {$ $} (c3);
					\draw (c2)[magenta, thick] --node[near start] {$ $} (c4);
					\draw (c2)[magenta, thick] --node[near start] {$ $} (c5);
					\draw (c2)[magenta, thick] --node[near start] {$ $}(c6);
					\draw (c2)[, thick] --node[near start] {$ $} (c7);
					\draw (c2)[, thick] --node[near start] {$ $} (c8);
					\draw (c2)[, thick] --node[near start] {$ $} (c9);
					\draw (c2)[, thick] --node[near start] {$ $}(c10);
					
					\draw (c3)[magenta, thick] --node[near start] {$ $}(c4);
					\draw (c3)[magenta, thick] --node[near start] {$ $}(c5);
					\draw (c3)[magenta, thick] --node[near start] {$ $} (c6);
					\draw (c3)[, thick] --node[near start] {$ $} (c7);
					\draw (c3)[, thick] --node[near start] {$ $} (c8);
					\draw (c3)[, thick] --node[near start] {$ $} (c9);
					
					\draw (c4)[magenta, thick] --node[near start] {$ $} (c5);
					\draw (c4)[magenta, thick] --node[near start] {$ $}(c6);
					\draw (c4)[, thick] --node[near start] {$ $}(c7);
					\draw (c4)[, thick] --node[near start] {$ $} (c8);
					
					\draw (c5)[magenta, thick] --node[near start] {$ $} (c6);
					\draw (c5)[, thick] --node[near start] {$ $} (c7);
					
			\end{tikzpicture} }%
			
			{\small Fig. 41.	$G_{0,12}:$ vertex coloring $\&$ clique}
			\label{$G_{0,12}$ - vertex coloring and clique}
		\end{minipage}\hfill
		\begin{minipage}{0.41\textwidth}
			\centering
			\resizebox{1\textwidth}{!}{%
				\begin{tikzpicture}[scale =0.7]
					
					\node (c0) at (7.5,-0.1) [circle,draw,scale=0.9,fill=cyan!40] {0};
					\node (c1) at (7,1.8) [circle,draw,scale=0.9,fill=cyan!40] {1};
					\node (c2) at (7.5,4.1) [circle,draw,scale=0.9,fill=cyan!40] {2};
					\node (c3) at (9.2,5.9) [circle,draw,scale=0.9,fill=cyan!40] {3};
					\node (c4) at (11.5,6.6) [circle,draw,scale=0.9,fill=cyan!40]{4};
					\node (c5) at (13.8,5.9) [circle,draw,scale=0.9,fill=cyan!40]{5};
					\node (c6) at (15.5,4.1) [circle,draw,scale=0.9,fill=cyan!40] {6};
					\node (c7) at (16,1.8) [circle,draw,scale=0.9,fill=cyan!40] {7};
					\node (c8) at (15.6,0.1) [circle,draw,scale=0.9,fill=cyan!40] {8};
					\node (c9) at (14.5,-1.3) [circle,draw,scale=0.9,fill=cyan!40] {9};
					\node (c10) at (12.9,-2.1) [circle,draw,scale=0.8,fill=cyan!40]{10};
					\node (c11) at (11,-2.3) [circle,draw,scale=0.8,fill=cyan!40] {11};
					\node (c12) at (9.1,-1.7) [circle,draw,scale=0.8,fill=cyan!40] {12};
					
	\draw (c0)[yellow, thick] --node[near start] {$ $} (c1);
	\draw (c0)[green, thick] --node[near start] {$ $} (c2);
	\draw (c0)[blue, thick] --node[near start] {$ $} (c3);
	\draw (c0)[black, thick] --node[near start] {$ $} (c4);
	\draw (c0)[red, thick] --node[near start] {$ $} (c5);
	\draw (c0)[gray, thick] --node[near start] {$ $} (c6);
	\draw (c0)[orange, thick] --node[near start] {$ $} (c7);
	\draw (c0)[ olive, thick] --node[near start] {$ $} (c8);
	\draw (c0)[ magenta, thick] --node[near start] {$ $} (c9);
	\draw (c0)[violet, thick] --node[near start] {$ $} (c10);
	\draw (c0)[ brown, thick] --node[near start] {$ $} (c11);
	\draw (c0)[purple, thick] --node[near start] {$ $} (c12);
					
	\draw (c1)[blue, thick] --node[near start] {$ $} (c2);
	\draw (c1)[, thick] --node[near start] {$ $} (c3);
	\draw (c1)[red, thick] --node[near start] {$ $} (c4);
	\draw (c1)[gray, thick] --node[near start] {$ $} (c5);
	\draw (c1)[orange, thick] --node[near start] {$ $} (c6);
	\draw (c1)[olive, thick] --node[near start] {$ $} (c7);
	\draw (c1)[ magenta, thick] --node[near start] {$ $} (c8);
	\draw (c1)[violet, thick] --node[near start] {$ $} (c9);
	\draw (c1)[ brown, thick] --node[near start] {$ $} (c10);
	\draw (c1)[purple, thick] --node[near start] {$ $} (c11);
					
	\draw (c1)[blue, thick] --node[near start] {$ $} (c2);
	\draw (c1)[, thick] --node[near start] {$ $} (c3);
	\draw (c1)[red, thick] --node[near start] {$ $} (c4);
	\draw (c1)[gray, thick] --node[near start] {$ $} (c5);
	\draw (c1)[orange, thick] --node[near start] {$ $} (c6);
	\draw (c1)[olive, thick] --node[near start] {$ $} (c7);
	\draw (c1)[ magenta, thick] --node[near start] {$ $} (c8);
	\draw (c1)[violet, thick] --node[near start] {$ $} (c9);
	\draw (c1)[ brown, thick] --node[near start] {$ $} (c10);
	\draw (c1)[purple, thick] --node[near start] {$ $} (c11);
					
	\draw (c2)[red, thick] --node[near start] {$ $} (c3);
	\draw (c2)[gray, thick] --node[near start] {$ $} (c4);
	\draw (c2)[orange, thick] --node[near start] {$ $} (c5);
	\draw (c2)[olive, thick] --node[near start] {$ $}(c6);
	\draw (c2)[magenta, thick] --node[near start] {$ $} (c7);
	\draw (c2)[ violet, thick] --node[near start] {$ $} (c8);
	\draw (c2)[brown, thick] --node[near start] {$ $} (c9);
	\draw (c2)[purple, thick] --node[near start] {$ $}(c10);
					
	\draw (c3)[orange, thick] --node[near start] {$ $}(c4);
	\draw (c3)[olive, thick] --node[near start] {$ $}(c5);
	\draw (c3)[magenta, thick] --node[near start] {$ $} (c6);
	\draw (c3)[violet, thick] --node[near start] {$ $} (c7);
	\draw (c3)[brown, thick] --node[near start] {$ $} (c8);
	\draw (c3)[purple, thick] --node[near start] {$ $} (c9);
					
	\draw (c4)[magenta, thick] --node[near start] {$ $} (c5);
	\draw (c4)[violet, thick] --node[near start] {$ $}(c6);
	\draw (c4)[brown, thick] --node[near start] {$ $}(c7);
	\draw (c4)[ purple, thick] --node[near start] {$ $} (c8);
					
	\draw (c5)[brown, thick] --node[near start] {$ $} (c6);
	\draw (c5)[purple, thick] --node[near start] {$ $} (c7);
						
\end{tikzpicture} 		
				}%
				
	~~~ {\small	Fig. 42. $G_{0,12}:$ edge coloring}\hfill
				\label{$G_{0,12}$ - edge coloring.}
	\end{minipage}
\end{figure}
	
\begin{theorem} [{\cite[Theorem 4.1]{js24}}] \label{thm Grs pf}\quad {\rm 			For $r,n \in \mathbb{N}$, the integral sum graph $G_{-r,n}$ is perfect.}
	\end{theorem}	
	
	\begin{proof} We have, $G_{-r,n} \cong K_1 * (G_{-r} * G_n)$ $\cong$ $K_1 * (G_{n} * G_{-r})$ $\cong$ $(K_1 * G_{n}) * G_{-r}$ $\cong$ $G_{0,n} * G_{-r}$ since join operation is associative as well as commutative among undirected graphs. Graphs $G_r$ and $G_{-r}$ are isomorphic without vertex labeling and so both are perfect by Theorems \ref{Gnpf} whereas $G_{0,n}$ is perfect by Theorem \ref{thm G0n} for $r,n\in\mathbb{N}$. Therefore  using Theorem \ref{thm join}, graph $G_{o,n} * G_{-r}$ is perfect for $r,n\in\mathbb{N}$. This implies, graph $G_{-r,n}$ $\cong$ $G_{o,n} * G_{-r}$ is perfect for $r,n\in\mathbb{N}$.
	\end{proof}		
	
	\begin{illu}  [{\cite[Illustration 4.2]{js24}}] \rm
		Consider the integral sum graph $G_{-5,7}$ as given in Figure 43. Vertices $v_{1}, v_{2}, v_{3}, v_{4}, v_{0}, v_{-1}, v_{-2}, v_{-3}$ (Vertices joined by the edges with magenta color) in $G_{-5,7}$ form a clique. Thus, $\omega(G_{-5,7})$ = 8. These vertices can be colored by 8 distinct colors, say $c_1,c_2,c_3,c_4,c_5,c_6,c_7, c_8$, respectively. Vertices $v_{-4}$ and $v_{-5}$ are nonadjacent to $v_{-3}$ and can be colored by $c_1$, the color of $v_1$. Vertices $v_{7},v_{6}$ and $v_{5}$ are nonadjacent to $v_{4}$ and can be colored by $c_4$. Therefore, $\chi(G_{-5,7})$ = 8. This implies, $G_{-5,7}$ is 1-perfect. Likewise we can show that every induced subgraph of $G_{-5,7}$ is 1-perfect. This implies, graph $G_{-5,7}$ is perfect. Edge coloring of $G_{-5,7}$ with 8 colors is shown in Figure 44.
	\end{illu}

	\begin{figure}
		\centering
		\begin{minipage}{0.41\textwidth}
			\centering
			\resizebox{1\textwidth}{!}{%
				\begin{tikzpicture}[scale =0.8]
					\node (c10) at (7,1.8) [circle,draw,scale=0.8,fill=magenta!80] {-3};
					\node (c9) at (7.5,4.1) [circle,draw,scale=0.8,fill=orange] {-2};
					\node (c8) at (9.2,5.9) [circle,draw,scale=0.8,fill=cyan] {-1};
					\node (c0) at (11.5,6.6) [circle,draw,scale=0.8,fill=pink]{0};
					\node (c1) at (13.8,5.9) [circle,draw,scale=0.8,fill=purple!50]{1};
					\node (c2) at (15.5,4.1) [circle,draw,scale=0.8,fill=lightgray] {2};
					\node (c3) at (16,1.8) [circle,draw,scale=0.8,fill=green] {3};
					\node (c4) at (15.6,0.1) [circle,draw,scale=0.8,fill=yellow] {4};
					\node (c5) at (14.5,-1.3) [circle,draw,scale=0.8,fill=yellow] {5};
					\node (c6) at (12.9,-2.1) [circle,draw,scale=0.8,fill=yellow] {6};
					\node (c7) at (11,-2.3) [circle,draw,scale=0.8,fill=yellow] {7};
					\node (c12) at (9.1,-1.7) [circle,draw,scale=0.8,fill=magenta!80] {-5};
					\node (c11) at (7.5,-0.1) [circle,draw,scale=0.8,fill=magenta!80] {-4};	
					
					\draw (c0)[magenta, thick] --node[][above]{ } (c1);
					\draw (c0)[magenta, thick] --node[near start] { } (c2);
					\draw (c0)[magenta, thick] --node[near start] { } (c3);
					\draw (c0)[magenta, thick] --node[near start] { } (c4);
					\draw (c0)[green, thick] --node[near start] {} (c5);
					\draw (c0)[green, thick]--node[near start] {} (c6);
					\draw (c0)[green, thick] --node[near start] {} (c7);
					\draw (c0)[magenta, thick] --node[near start] {} (c8);
					\draw (c0)[magenta, thick] --node[near start] {} (c9);
					\draw (c0)[magenta, thick] --node[near start] {} (c10);
					\draw (c0)[green, thick] --node[near start] {} (c11);
					\draw (c0)[green, thick] --node[near start] { } (c12);
					
					\draw (c1)[magenta, thick] --node[near start] { } (c8);
					\draw (c1)[magenta, thick] --node[near start] { } (c9);
					\draw (c1)[magenta, thick] --node[near start] { } (c10);
					\draw (c1)[green, thick] --node[near start] {}(c11);
					\draw (c1)[green, thick] --node[near start] {} (c12);
					
					\draw (c2)[magenta, thick] --node[near start] { } (c8);
					\draw (c2)[magenta, thick] --node[near start] { } (c9);
					\draw (c2)[magenta, thick] --node[near start] {}(c10);
					\draw (c2)[green, thick] --node[near start] {} (c11);
					\draw (c2)[green, thick] --node[near start] {} (c12);
					
					\draw (c3)[magenta, thick] --node[near start] { } (c8);
					\draw (c3)[magenta, thick] --node[near start] {} (c9);
					\draw (c3)[magenta, thick] --node[near start] {} (c10);
					\draw (c3)[green, thick] --node[near start] {}(c11);
					\draw (c3)[green, thick] --node[near start] {} (c12);
					
					\draw (c4)[magenta, thick] --node[near start] {} (c8);
					\draw (c4)[magenta, thick] --node[near start] {} (c9);
					\draw (c4)[magenta, thick] --node[near start] {} (c10);
					\draw (c4)[green, thick] --node[near start] {}(c11);
					\draw (c4)[green, thick] --node[near start] {}(c12);
					
					\draw (c5)[green, thick] --node[near start] {} (c8);
					\draw (c5)[green, thick] --node[near start] {} (c9);
					\draw (c5)[green, thick] --node[near start] {}(c10);
					\draw (c5)[black, thick] --node[near start] {} (c11);
					\draw (c5)[black, thick] --node[near start] {} (c12);
					
					\draw (c6)[green, thick] --node[near start] {} (c8);
					\draw (c6)[green, thick] --node[near start] {} (c9);
					\draw (c6)[green, thick] --node[near start] {} (c10);
					\draw (c6)[black, thick] --node[near start] { } (c11);
					\draw (c6)[black, thick] --node[near start] { } (c12);
					
					\draw(c7)[green, thick] --node[near start] { } (c8);
					\draw (c7)[green, thick] --node[near start] { } (c9);
					\draw (c7)[green, thick] --node[near start] { } (c10);
					\draw (c7)[black, thick] --node[near start] { } (c11);
					\draw (c7)[black, thick] --node[near start] { } (c12);
					
					\draw (c1)[magenta, thick] --node[above] { } (c2);
					\draw (c1)[green, thick] --node[near start] { }(c3);
					\draw (c1)[magenta, thick] --node[near start] {{ }} (c4);
					\draw (c1)[green, thick] --node[near start] { } (c5);
					\draw (c1)[green, thick] --node[near start] { }(c6);
					
					\draw (c2)[magenta, thick] --node[right] { } (c3);
					\draw (c2)[magenta, thick] --node[near end] { } (c4);
					\draw (c2)[green, thick] --node[near start] { }(c5);
					
					\draw (c3)[magenta, thick] --node[right] { }(c4);
					
					\draw (c8)[magenta, thick] --node[left] { }(c9);
					\draw (c8)[magenta, thick] --node[near end] { }(c10);
					\draw (c8)[green, thick] --node[near start] { } (c11);
					
					\draw (c9)[magenta, thick] --node[left] { } (c10);
			\end{tikzpicture} }%
			
			{\small	Fig. 43. $G_{-5, 7}:$ vertex coloring}
			\label{$G_{-5, 7}$ with vertex coloring.}
		\end{minipage}\hfill
		\begin{minipage}{0.41\textwidth}
			\centering
			\resizebox{1\textwidth}{!}{%
				\begin{tikzpicture}[scale =0.8]
					\node (c10) at (7,1.8) [circle,draw,scale=0.8,fill = yellow!30] {-3};
					\node (c9) at (7.5,4.1) [circle,draw,scale=0.8, fill = yellow!30] {-2};
					\node (c8) at (9.2,5.9) [circle,draw,scale=0.8, fill = yellow!30] {-1};
					\node (c0) at (11.5,6.6) [circle,draw,scale=0.8, fill = yellow!30]{0};
					\node (c1) at (13.8,5.9) [circle,draw,scale=0.8, fill = yellow!30]{1};
					\node (c2) at (15.5,4.1) [circle,draw,scale=0.8, fill = yellow!30] {2};
					\node (c3) at (16,1.8) [circle,draw,scale=0.8, fill = yellow!30] {3};
					\node (c4) at (15.6,0.1) [circle,draw,scale=0.8, fill = yellow!30] {4};
					\node (c5) at (14.5,-1.3) [circle,draw,scale=0.8, fill = yellow!30] {5};
					\node (c6) at (12.9,-2.1) [circle,draw,scale=0.8, fill = yellow!30] {6};
					\node (c7) at (11,-2.3) [circle,draw,scale=0.8, fill = yellow!30] {7};
					\node (c12) at (9,-1.7) [circle,draw,scale=0.8, fill = yellow!30] {-5};
					\node (c11) at (7.5,-0.1) [circle,draw,scale=0.8, fill = yellow!30] {-4};	
					
					\draw (c0)[gray, thick] --node[][above]{} (c1);
					\draw (c0)[orange, thick] --node[near start] {} (c2);
					\draw (c0)[cyan, thick] --node[near start] {} (c3);
					\draw (c0)[purple, thick] --node[near start] { } (c4);
					\draw (c0)[blue, thick] --node[near start] {} (c5);
					\draw (c0)[green, thick]--node[near start] {} (c6);
					\draw (c0)[brown, thick] --node[near start] {} (c7);
					\draw (c0)[violet, thick] --node[near start] {} (c8);
					\draw (c0)[magenta, thick] --node[near start] {} (c9);
					\draw (c0)[red, thick] --node[near start] {} (c10);
					\draw (c0)[olive, thick] --node[near start] {} (c11);
					\draw (c0)[black, thick] --node[near start] {} (c12);
					
					\draw (c1)[orange, thick] --node[near start] {} (c8);
					\draw (c1)[cyan, thick] --node[near start] {} (c9);
					\draw (c1)[purple, thick] --node[near start] {} (c10);
					\draw (c1)[blue, thick] --node[near start] {}(c11);
					\draw (c1)[green, thick] --node[near start] {} (c12);
					
					\draw (c2)[cyan, thick] --node[near start] {} (c8);
					\draw (c2)[purple, thick] --node[near start] {} (c9);
					\draw (c2)[blue, thick] --node[near start] {}(c10);
					\draw (c2)[green, thick] --node[near start] { } (c11);
					\draw (c2)[brown, thick] --node[near start] {} (c12);
					
					\draw (c3)[purple, thick] --node[near start] {} (c8);
					\draw (c3)[blue, thick] --node[near start] {} (c9);
					\draw (c3)[green, thick] --node[near start] {} (c10);
					\draw (c3)[brown, thick] --node[near start] {}(c11);
					\draw (c3)[violet, thick] --node[near start] {} (c12);
					
					\draw (c4)[blue, thick] --node[near start] {} (c8);
					\draw (c4)[green, thick] --node[near start] {} (c9);
					\draw (c4)[brown, thick] --node[near start] {} (c10);
					\draw (c4)[violet, thick] --node[near start] {}(c11);
					\draw (c4)[magenta, thick] --node[near start] {}(c12);
					
					\draw (c5)[green, thick] --node[near start] {} (c8);
					\draw (c5)[brown, thick] --node[near start] {} (c9);
					\draw (c5)[violet, thick] --node[near start] {}(c10);
					\draw (c5)[magenta, thick] --node[near start] {{ }} (c11);
					\draw (c5)[red, thick] --node[near start] {} (c12);
					
					\draw (c6)[brown, thick] --node[near start] {} (c8);
					\draw (c6)[violet, thick] --node[near start] {} (c9);
					\draw (c6)[magenta, thick] --node[near start] {} (c10);
					\draw (c6)[red, thick] --node[near start] {} (c11);
					\draw (c6)[olive, thick] --node[near start] {} (c12);
					
					\draw(c7)[magenta, thick] --node[near start] {} (c8);
					\draw (c7)[olive, thick] --node[near start] {} (c9);
					
					\draw (c7)[gray, thick] --node[near start] {} (c10);
					\draw (c7)[orange, thick] --node[near start] {} (c11);
					\draw (c7)[cyan, thick] --node[near start] {} (c12);
					
					\draw (c1)[violet, thick] --node[above] { } (c2);
					\draw (c1)[magenta, thick] --node[near start] {}(c3);
					\draw (c1)[red, thick] --node[near start] {} (c4);
					\draw (c1)[olive, thick] --node[near start] {} (c5);
					\draw (c1)[black, thick] --node[near start] { }(c6);
					
					\draw (c2)[red, thick] --node[right] {{ }} (c3);
					\draw (c2)[olive, thick] --node[near end] { } (c4);
					\draw (c2)[black, thick] --node[near start] { }(c5);
					
					\draw (c3)[black, thick] --node[right] {}(c4);
					
					\draw (c8)[red, thick] --node[left] { }(c9);
					\draw (c8)[olive, thick] --node[near end] {}(c10);
					\draw (c8)[black, thick] --node[near start] {} (c11);
					
					\draw (c9)[black, thick] --node{ } (c10);
				\end{tikzpicture}					
			}%
			
			{\small	Fig. 44. $G_{-5,7}$ with edge coloring}
			\label{$G_{-5,7}$ with edge coloring.}
		\end{minipage}
	\end{figure}
	
\section{$(a, d)$-CMS Decomposition of $K_{n}$ and $G_{0,m}$ using integral sum labeling}

In this section, we discuss {\it Ascending Subgraph Decomposition (ASD)}, $(a, d)$-ASD, {\em $(a, d)$-Continuous Monotonic Decomposition} or $(a, d)$-CMD, and $(a, d)$-Continuous Monotonic Subgraph Decomposition (CMSD) of graphs and present results obtained in \cite{vs15} on $(a,d)$-CMSD of $K_{n}$ and $G_{0,m}$ using integral sum labeling. 

Vilfred and Suryakala \cite{vs15} defined $(a, d)$-Continuous Monotonic Subgraph Decomposition ($(a, d)$-CMSD) of a graph. They used integral sum labeling to study $(a, d)$-CMSD of $K_{n}$ and $G_{0,m}$ and obtained the following results. (i)~ For $n \geq 3$, $K_n$ admits $(a, d)$-CMSD into triangular books for some $a$ and $d$, $a,d\in\mathbb{N}$; ~(ii)~ For $n\in\mathbb{N}$, $G_{0,2n}$, $G_{0,4n+2}$ and $G_{0,4n+3}$ admit $(a, d)$-CMSD into triangular books with bookmark for some $a$ and $d$, $a,d\in\mathbb{N}$; ~(iii)~ $G_{0,4n+1}$ admits Ascending Subgraph Decomposition (ASD) but doesn't admit $(a, d)$-ASD and $(a, d)$-CMD into triangular books with bookmark for any $a,d\in\mathbb{N}$;  ~(iv) ~ For $n\in\mathbb{N}$, $G_{0,4n+2}$, $G_{0,4n}$ and $G_{0,4n-1}$ admit $(a, d)$-CMSD into Fans with a handle for some $a$ and $d$, $a,d\in \mathbb{N}$ and ~(v)~ $G_{0,4n+1}$ admits ASD into Fans with a handle and one $P_2$ but doesn't admit $(a, d)$-ASD and $(a,d)$-CMD into Fans with a handle for any $a,d \in \mathbb{N}$. These results are presented in this section. 

Throughout this section, vertices of $K_n$ and  $G_{0,{n-1}}$  are considered as the vertices of an $n$-gon ordered in the anti-clockwise direction. Definitions of triangular book with a bookmark and fan graph with a handle are given in subsection 5.5. 

\subsection{Ascending Subgraph Decomposition (ASD) and $(a, d)$-ASD}

Alavi \cite{ab87} introduced the concept of {\it Ascending Subgraph Decomposition (ASD)} of a graph $G$ with size $(n$+1$)C_2$ as the decomposition of $G$ into $n$ subgraphs $G_1,G_2,\dots,G_n$ without isolated vertices such that each $G_i$ is isomorphic to a proper subgraph of $G_{i+1}$  and  $|E(G_i)|  = i$ for $1 \leq i\leq n$. Nagarajan \cite{nn06} generalized ASD to $(a, d)$-ASD of graph $G$ with size $\frac{(2a+(n-1)d)n}{2}$ as the decomposition of $G$ into $n$ subgraphs $G_1,G_2,\dots,G_n$ without isolated vertices such that each $G_i$ is isomorphic to a proper subgraph of $ G_{i+1}$ and $|E(G_i)|  = a+(i-1)d$ for $1 \leq i \leq n$. Clearly, $ASD$ of a graph $G$ and its $(1,1)$-ASD are the same.

\subsection{Continuous Monotonic Decomposition (CMD) and $(a, d)$-CMD}
 Gnana Dhas \cite{gp00} defined {\em $(a, d)$-Continuous Monotonic Decomposition} or $(a, d)$-CMD of a graph $G$ of size $\frac{(2a+(n-1)d)n}{2}$ as the decomposition of $G$ into $n$ subgraphs $G_1$, $G_2$, $\dots$, $G_n$ such that each $G_i$ is connected and  $|E(G_i)| = a+(i-1)d$ for $i = 1,2,\dots,n$. 
Clearly, CMD of a graph $G$ and its $(1,1)$-CMD are the same.

Among the family of graphs some graphs may have $(a,d)$-ASD, some may have $(a,d)$-CMD, some may have both $(a,d)$-ASD and $(a,d)$-CMD and the others have neither $(a,d)$-ASD nor $(a,d)$-CMD. Huaitand \cite{hk98} studied ASD of regular graphs and proved that every regular bipartite graph as ASD. Nagarajan \cite{nn06} studied $(a,d)$-ASD of wheels. Finding graphs having either $(a,d)$-ASD or $(a,d)$-CMD seems to be difficult and finding graphs having both $(a,d)$-ASD and $(a,d)$-CMD seems to be more difficult, $a,d\in\mathbb{N}$. 

\subsection{On $(a, d)$-Continuous Monotonic Subgraph Decomposition ($(a, d)-CMSD$)}

While studying decomposition of integral sum graphs, Vilfred and Suryakala \cite{vs15} came across graphs having both $(a, d)$-ASD and $(a, d)$-CMD. This motivated them to define Continuous Monotonic Subgraph Decomposition (CMSD) and $(a, d)$-CMSD of graphs as follows.

\begin{dfn} \cite{vs15}	A decomposition of graph $G$ that is both $(a,d)$-ASD and $(a,d)$-CMD is called an {\em $(a,d)$-Continuous Monotonic Subgraph Decomposition} or $(a,d)$-CMSD of $G$, $a,d\in\mathbb{N}$. Thus, $(a,d)$-CMSD of graph $G$ with size $\frac{(2a+(n-1)d)n}{2}$ is the decomposition of $G$ into $n$ subgraphs $G_1,G_2,\dots,G_n$ without isolated vertices such that each $G_i$ is {\it connected and isomorphic to a proper subgraph} of $G_{i+1}$ and  $|E(G_i)| = a+(i-1)d$ for $1\leq  i\leq n$ and $a,d,n\in\mathbb{N}$.
\end{dfn}

\subsection{On $(a, d)$-CMSD of $K_n$ into triangular books}

Results on CMSD and $(a, d)$-CMSD of $K_n$ into triangular books are presented in this subsection, $a,d,n\in\mathbb{N}$. It is also proved that for $n\geq 3$, $K_n$ admits $(1, 1)$-CMSD into stars.

\begin{theorem} [{\cite[Theorem 2.1]{vs15}}] \label{9.2}	{\rm
	For $m\in\mathbb{N}$ and $m\geq 2$, 
	
	(i) $K_{2m}$ admits $(1, 4)$-CMSD into triangular books and 
	
	(ii) $K_{2m-1}$ admits $(3, 4)$-CMSD into triangular books.}
\end{theorem}

\begin{proof}\quad Let $V(K_n)$ = $\{0,1,\dots,n-1\}$. We have $|E(K_n)|$ = $nC_2$. 
	\\
	(i) Let $n$ = $2m$, $m\in\mathbb{N}$. Then, $K_n$ = $K_{2m}$ and an $(1, 4)$-CMSD of $K_{2m}$ into triangular books is obtained as follows. 
	
	$K_{2m} = TB_{2m-2}(0,1)\cup TB_{2m-4}(2,3)\cup\dots\cup TB_2(2m-4,2m-3) \cup TB_0(2m-2$, $2m-1)$ where $TB_{2m-2j}(2j-2,2j-1)$ in $K_{2m}$ represents triangular book with spine $(2j-2,2j-1)$ and $(2j-2,2j-1,2j)$, $(2j-2,2j-1,2j+1)$, $\dots$, $(2j-2,2j-1,2m-1)$ as the $(2m-2j)$ number of triangular pages and is a connected subgraph, $j$ = $1,2,\dots,m$. In $K_{2m}$, $(0,1)$ is the spine for $TB_{2m-2}(0,1)$, both the vertices 0 and 1 are adjacent to the remaining 2m-2 vertices, $2,3,\dots,2m-1$ and each one is of degree $2m-1$ in $TB_{2m-2}(0,1);$ $(2,3)$ is the spine for $TB_{2m-4}(2,3)$, both the vertices 2 and 3 are adjacent to the $2m-4$ vertices, $4,5,\dots,2m-1$ and each one is of degree $2m-1$ in $TB_{2m-2}(0,1)\cup TB_{2m-4}(2,3);$ $(4,5)$ is the spine for $TB_{2m-6}(4,5)$, both the vertices 4 and 5 are adjacent to the $2m-6$ vertices, $6,7,\dots, 2m-1$ and each one is of degree $2m-1$ in $TB_{2m-2}(0,1)\cup TB_{2m-4}(2,3)\cup TB_{2m-6}(4,5);$ . . . ; $(2m-4,2m-3)$ is the spine for $TB_2(2m-4,2m-3)$, both the vertices $2m-4$ and $2m-3$ are adjacent to the 2 vertices, $2m-2$ and $2m-1$ and each one is of degree $2m-1$ in $TB_{2m-2}(0,1)$ $\cup$ $TB_{2m-4}(2,3)$ $\cup$ $TB_{2m-6}(4,5)$ $\cup$ $\dots$ $\cup$ $TB_2(2m-4,2m-3)$; $(2m-2,2m-1)$ is the spine for $TB_0(2m-2,2m-1)$ which is a triangular book without pages and each one of the vertices $2m-2$ and $2m-1$ is of degree $2m-1$ in $TB_{2m-2}(0,1)$ $\cup$ $TB_{2m-4}(2,3)$ $\cup$ $TB_{2m-6}(4,5)$ $\cup$ $\dots$ $\cup$ $TB_2(2m-4,2m-3)$ $\cup$ $TB_0(2m-2,2m-1)$ = $K_{2m}$. Also, $|E(TB_0(2m-2,2m-1))|$ = $1 < |E(TB_2(2m-4$, $2m-3))|$ = $5 < |E(TB_4(2m-6,2m-5))|$ = $9 < \dots < |E(TB_{2m -4}(2,3))|$ = $4m-7 < |E(TB_{2m-2}(0,1))|$ = $4m-3$. And clearly, $TB_0(2m-2,2m-1)$ is a connected subgraph of $TB_2(2m-4,2m-3)$ which is a connected subgraph of $TB_4(2m-6,2m-5)$ which is a connected subgraph of $\dots$ which is a connected subgraph of $TB_{2m-4}(2,3)$ which is a connected subgraph of $TB_{2m-2}(0,1)$, without vertex labels. Thus, $K_{2m}$ admits $(1,4)-{\it CMSD}$ into triangular books for $m \geq 2$. 
	
	In different colors, $(1, 4)$-CMSD of $K_4$, $K_6$ and $K_8$ are shown in  Figures 45, 46 and 47, respectively. Here, $K_4 = TB_2(0,1)\cup TB_0(2,3)$, $K_6 = TB_4(0,1)\cup  TB_2(2,3)\cup TB_0(4,5)$ and $K_8 = TB_6(0,1)\cup TB_4(2,3)\cup TB_2(4,5) \cup TB_0(6,7)$.
		
	\begin{center}
	\begin{tikzpicture}[scale =0.7]
	\node (a3) at (5,1.5)  [circle,draw,scale=0.6] {3};
	\node (a2) at (7,1.5)  [circle,draw,scale=0.6]{2};
	\node (a1) at (7,-.5)  [circle,draw,scale=0.6] {1};
	\node (a0) at (5,-.5)  [circle,draw,scale=0.6] {0};

	\draw (a0) -- (a1);
	\draw (a0) -- (a2);
	\draw (a0) -- (a3);

	\draw (a1) -- (a2);
	\draw (a1) -- (a3);
	
{\color{red}	
	\draw (a2) -- (a3);}
	
	\node (b4) at (10,2)  [circle,draw,scale=0.6] {4};
	\node (b3) at (12,2)  [circle,draw,scale=0.6]{3};
	\node (b2) at (13,.75)  [circle,draw,scale=0.6] {2};
	\node (b1) at (12,-.5)  [circle,draw,scale=0.6] {1};
	\node (b0) at (10,-.5)  [circle,draw,scale=0.6] {0};
	\node (b5) at (9,.75)  [circle,draw,scale=0.6]{5};
	
		\draw (b0) -- (b1);
		\draw (b0) -- (b2);
		\draw (b0) -- (b3);
		\draw (b0) -- (b4);
		\draw (b0) -- (b5);
		
		\draw (b1) -- (b2);
		\draw (b1) -- (b3);
		\draw (b1) -- (b4);
		\draw (b1) -- (b5);
		
	{\color{red}
		\draw (b2) -- (b3);
		\draw (b2) -- (b4);
		\draw (b2) -- (b5);
		
		\draw (b3) -- (b4);
        \draw (b3) -- (b5);}

	{\color{green}
		\draw (b4) -- (b5);}
	
	\node (c5) at (16,2.75)  [circle,draw,scale=0.6] {5};
	\node (c6) at (15,1.75)  [circle,draw,scale=0.6] {6};
	\node (c4) at (18,2.75)  [circle,draw,scale=0.6]{4};
	\node (c3) at (19,1.75)  [circle,draw,scale=0.6]{3};
	\node (c2) at (19,.5)  [circle,draw,scale=0.6] {2};
	\node (c1) at (18,-.5)  [circle,draw,scale=0.6] {1};
	\node (c0) at (16,-.5)  [circle,draw,scale=0.6] {0};
	\node (c7) at (15,.5)  [circle,draw,scale=0.6]{7};
	
		\draw (c0) -- (c1);
		\draw (c0) -- (c2);
		\draw (c0) -- (c3);
		\draw (c0) -- (c4);
		\draw (c0) -- (c5);
        \draw (c0) -- (c6);
        \draw (c0) -- (c7);
				
        \draw (c1) -- (c2);
        \draw (c1) -- (c3);
        \draw (c1) -- (c4);
        \draw (c1) -- (c5);
        \draw (c1) -- (c6);
        \draw (c1) -- (c7);		

	{\color{red}
		\draw (c2) -- (c3);
		\draw (c2) -- (c4);
		\draw (c2) -- (c5);
		\draw (c2) -- (c6);
        \draw (c2) -- (c7);		

		\draw (c3) -- (c4);
		\draw (c3) -- (c5);
		\draw (c3) -- (c6);
		\draw (c3) -- (c7);}
		
	{\color{green}
		\draw (c4) -- (c5);
		\draw (c4) -- (c6);
		\draw (c4) -- (c7);
	
	\draw (c5) -- (c6);
	\draw (c5) -- (c7);}

	{\color{violet}
	\draw (c6) -- (c7);}

	\end{tikzpicture}
	
	\vspace{.1cm}	
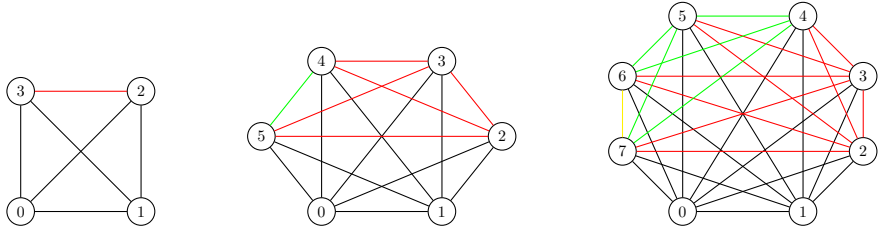
	{\footnotesize	Fig. 45. (1,4)-CMSD of $K_{4}$ \hspace{.2cm} Fig. 46. (1,4)-CMSD of $K_{6}$
	\hspace{.2cm} Fig. 47. (1,4)-CMSD of $K_8$ }
	\end{center}
	
	\noindent
	(ii) Consider $K_{2m-1}$. Decomposition of $K_{2m-1}$ into triangular books of $(3, 4)$-CMSD is obtained as follow. 
	
	$K_{2m+1} = TB_{2m-1}(0,1) \cup TB_{2m-3}(2,3) \cup \cdots \cup TB_3(2m-4,2m-3) \cup TB_1(2m-2$, $2m-1)$ where $TB_{2m+1-2j}(2j-2,2j-1)$ in $K_{2m+1}$ represents triangular book with spine $(2j-2,2j-1)$ and $(2j-2,2j-1,2j)$, $(2j-2,2j-1,2j+1)$, $\dots$, $(2j-2,2j-1,2m)$ as the $(2m+1-2j)$ number of triangular pages and is a connected subgraph, $j$ = $1,2,\dots,m$. The above decomposition of $K_{2m+1}$ is similar to the decomposition given in case $(i)$~ except $K_{2m+1}$ admits $(3, 4)-{\it CMSD}$ into triangular books since $|E(TB_1(2m-2, 2m-1))|$ = $3 < |E(TB_3(2m-4,2m-3))|$ = $7 < |E(TB_5(2m-6$, $2m-5))|$ = $11 < \dots < |E(TB_{2m-3}(2,3))|$ = $4m-5 < |E(TB_{2m-1}(0,1))|$ = $4m-1$ and $TB_1(2m-2,2m-1)$ is a connected subgraph of $TB_3(2m-4,2m-3)$ which is a connected subgraph of $TB_5(2m-6,2m-5)$ which is a connected subgraph of $\dots$ which is a connected subgraph of $TB_{2m-3}(2,3)$ which is a connected subgraph of $TB_{2m-1}(0,1)$, without vertex labels. 
			
\begin{center}
	\begin{tikzpicture}[scale = 0.6]
	\node (a2) at (7,1.5)  [circle,draw,scale=0.6]{2};
	\node (a1) at (8,-.5)  [circle,draw,scale=0.6] {1};
	\node (a0) at (6,-.5)  [circle,draw,scale=0.6] {0};
	
	\draw (a0) -- (a1);
	\draw (a0) -- (a2);
	
	\draw (a1) -- (a2);
		
	\node (b3) at (13,2)  [circle,draw,scale=0.6]{3};
	\node (b2) at (14.5,.75)  [circle,draw,scale=0.6] {2};
	\node (b1) at (14,-.5)  [circle,draw,scale=0.6] {1};
	\node (b0) at (12,-.5)  [circle,draw,scale=0.6] {0};
	\node (b4) at (11.5,.75)  [circle,draw,scale=0.6]{4};
	
	\draw (b0) -- (b1);
	\draw (b0) -- (b2);
	\draw (b0) -- (b3);
	\draw (b0) -- (b4);
	
	\draw (b1) -- (b2);
	\draw (b1) -- (b3);
	\draw (b1) -- (b4);
	
	{\color{red}
		\draw (b2) -- (b3);
		\draw (b2) -- (b4);
		
		\draw (b3) -- (b4);}
		
	\node (c6) at (17,.5)  [circle,draw,scale=0.6]{6};
	\node (c5) at (17,1.75)  [circle,draw,scale=0.6] {5};
	\node (c4) at (19,2.75)  [circle,draw,scale=0.6]{4};
	\node (c3) at (21,1.75)  [circle,draw,scale=0.6]{3};
	\node (c2) at (21,.5)  [circle,draw,scale=0.6] {2};
	\node (c1) at (20,-.5)  [circle,draw,scale=0.6] {1};
	\node (c0) at (18,-.5)  [circle,draw,scale=0.6] {0};
	
	\draw (c0) -- (c1);
	\draw (c0) -- (c2);
	\draw (c0) -- (c3);
	\draw (c0) -- (c4);
	\draw (c0) -- (c5);
	\draw (c0) -- (c6);
	
	\draw (c1) -- (c2);
	\draw (c1) -- (c3);
	\draw (c1) -- (c4);
	\draw (c1) -- (c5);
	\draw (c1) -- (c6);
	
	{\color{red}
		\draw (c2) -- (c3);
		\draw (c2) -- (c4);
		\draw (c2) -- (c5);
		\draw (c2) -- (c6);
		
		\draw (c3) -- (c4);
		\draw (c3) -- (c5);
		\draw (c3) -- (c6);}
	
	{\color{green}
		\draw (c4) -- (c5);
		\draw (c4) -- (c6);
		
		\draw (c5) -- (c6);}
		
	\end{tikzpicture}
	
	\vspace{.1cm}	
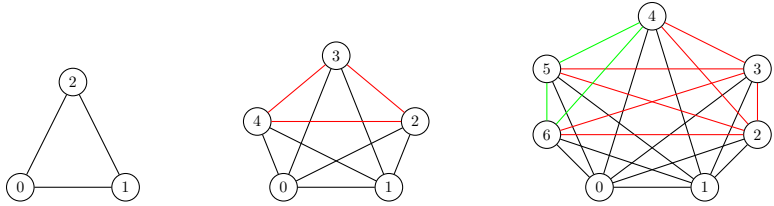
{\footnotesize	Fig. 48. (3, 4)-CMSD of $K_{3}$ \hspace{.15cm} Fig. 49. (3, 4)-CMSD of $K_{5}$
	\hspace{.15cm} Fig. 50. (3, 4)-CMSD of $K_7$}
\end{center}

	(3, 4)-CMSD of $K_3, K_5$ and $K_7$ are shown in  Figures 48, 49 and 50, respectively, with different colors. Here, (3,4)-CMSDs are $K_3$ = $TB_1(0,1)$, $K_5$ = $TB_3(0,1)$ $\cup$ $TB_1(2,3)$ and $K_7$ = $TB_5(0,1)$ $\cup$ $TB_3(2,3)$ $\cup$ $TB_1(4,5)$. Hence the result.  
\end{proof}

\begin{cor} [{\cite[Corollary 2.2]{vs15}}] \label{9.3} {\rm
	$K_n$ admits $(a, d)$-CMSD into triangular books for some $a$ and $d$, $a,d\in\mathbb{N}$. }\hfill  $\Box$
\end{cor}

\begin{theorem} [{\cite[Theorem 2.3]{vs15}}] \label{9.4} {\rm
			For $n\geq 3$, $K_n$ admits $(1, 1)$-CMSD into stars.}
\end{theorem}

\begin{proof} Let $V(K_n)$ = $\{v_0,v_1,v_2,\dot,v_{n-1}\}$ and $i$ be the label of vertex $v_i$ in $K_n$ for $i$ = 0,1,2,...,$n-1$ and $n \geq 3$. Consider the following decomposition of $K_n$.

 $K_n$ = $K_{1,1}(0; 1)$ $\cup$ $K_{1,2}(2; 0,1)$ $\cup$ $K_{1,3}(3; 0,1,2)$ $\cup$ $\dots$ $\cup$ $K_{1,n-1}(n-1; 0,1,2,\dots,n-2)$ where $K_{1,j}(j; 0,1,\dots,j-1)$ is the star $K_{1,j}$ with internal vertex $j$ and leaves, $0,1,\dots,j-1$, $1\leq j\leq n-1$ and $K_{1,1}(0; 1)$ is the edge with end points 0 and 1. Clearly, the above decomposition is a $(1,1)$-CMSD of $K_n$ into stars. 
 
	(1, 1)-CMSD of $K_6$, $K_7$ and $K_8$ are shown,  as given in the proof, in Figures 51, 52 and 53 with different colors. Here, (1, 1)-CMSDs are $K_6$ = $K_{1,1}(0; 1)$ $\cup$ $K_{1,2}(2; 0,1)$ $\cup$ $K_{1,3}(3; 0,1,2)$ $\cup$ $K_{1,4}(4; 0,1,2,3)$ $\cup$ $K_{1,5}(5; 0,1,2,3,4)$, $K_7$ = $K_{1,1}(0; 1)$ $\cup$ $K_{1,2}(2; 0,1)$ $\cup$ $K_{1,3}(3; 0,1,2)$ $\cup$ $\dots$ $\cup$ $K_{1,6}(6; 0,1,2,\dots,5)$ and $K_8$ = $K_{1,1}(0; 1)$ $\cup$ $K_{1,2}(2; 0,1)$ $\cup$ $K_{1,3}(3; 0,1,2)$ $\cup$ $\dots$ $\cup$ $K_{1,7}(7; 0,1,2,\dots,6)$. 	
\end{proof}

	\begin{center}
	\begin{tikzpicture}[scale =0.8]
\node (a4) at (5,2)  [circle,draw,scale=0.6] {4};
\node (a3) at (7,2)  [circle,draw,scale=0.6]{3};
\node (a2) at (8,.75)  [circle,draw,scale=0.6] {2};
\node (a1) at (7,-.5)  [circle,draw,scale=0.6] {1};
\node (a0) at (5,-.5)  [circle,draw,scale=0.6] {0};
\node (a5) at (4,.75)  [circle,draw,scale=0.6]{5};

{\color{brown}
\draw (a0) -- (a1);}

{\color{red}
	\draw (a2) -- (a0);
	\draw (a2) -- (a1);}

{\color{green}
	\draw (a3) -- (a0);
	\draw (a3) -- (a1);
    \draw (a3) -- (a2);}

{\color{blue}
	\draw (a4) -- (a0);
	\draw (a4) -- (a1);
	\draw (a4) -- (a2);
	\draw (a4) -- (a3);}

{\color{orange}
	\draw (a5) -- (a0);
	\draw (a5) -- (a1);
	\draw (a5) -- (a2);
	\draw (a5) -- (a3);
	\draw (a5) -- (a4);}

\node (b6) at (10,.5)  [circle,draw,scale=0.6]{6};
\node (b5) at (10,1.75)  [circle,draw,scale=0.6] {5};
\node (b4) at (12,2.75)  [circle,draw,scale=0.6]{4};
\node (b3) at (14,1.75)  [circle,draw,scale=0.6]{3};
\node (b2) at (14,.5)  [circle,draw,scale=0.6] {2};
\node (b1) at (13,-.5)  [circle,draw,scale=0.6] {1};
\node (b0) at (11,-.5)  [circle,draw,scale=0.6] {0};

{\color{brown}
	\draw (b0) -- (b1);}

{\color{red}
	\draw (b2) -- (b0);
	\draw (b2) -- (b1);}

{\color{green}
	\draw (b3) -- (b0);
	\draw (b3) -- (b1);
	\draw (b3) -- (b2);}

{\color{blue}
	\draw (b4) -- (b0);
	\draw (b4) -- (b1);
	\draw (b4) -- (b2);
	\draw (b4) -- (b3);}

{\color{orange}
\draw (b5) -- (b0);
\draw (b5) -- (b1);
\draw (b5) -- (b2);
\draw (b5) -- (b3);
\draw (b5) -- (b4);}

\draw (b6) -- (b0);
\draw (b6) -- (b1);
\draw (b6) -- (b2);
\draw (b6) -- (b3);
\draw (b6) -- (b4);
\draw (b6) -- (b5);

	\node (c5) at (17,2.75)  [circle,draw,scale=0.6] {5};
	\node (c6) at (16,1.75)  [circle,draw,scale=0.6] {6};
	\node (c4) at (19,2.75)  [circle,draw,scale=0.6]{4};
	\node (c3) at (20,1.75)  [circle,draw,scale=0.6]{3};
	\node (c2) at (20,.5)  [circle,draw,scale=0.6] {2};
	\node (c1) at (19,-.5)  [circle,draw,scale=0.6] {1};
	\node (c0) at (17,-.5)  [circle,draw,scale=0.6] {0};
	\node (c7) at (16,.5)  [circle,draw,scale=0.6]{7};
	
{\color{brown}
	\draw (c0) -- (c1);}

{\color{red}
	\draw (c2) -- (c0);
	\draw (c2) -- (c1);}

{\color{green}
	\draw (c3) -- (c0);
	\draw (c3) -- (c1);
	\draw (c3) -- (c2);}

{\color{blue}
	\draw (c4) -- (c0);
	\draw (c4) -- (c1);
	\draw (c4) -- (c2);
	\draw (c4) -- (c3);}

{\color{orange}
\draw (c5) -- (c0);
\draw (c5) -- (c1);
\draw (c5) -- (c2);
\draw (c5) -- (c3);
\draw (c5) -- (c4);}

\draw (c6) -- (c0);
\draw (c6) -- (c1);
\draw (c6) -- (c2);
\draw (c6) -- (c3);
\draw (c6) -- (c4);
\draw (c6) -- (c5);

{\color{purple}
\draw (c7) -- (c0);
\draw (c7) -- (c1);
\draw (c7) -- (c2);
\draw (c7) -- (c3);
\draw (c7) -- (c4);
\draw (c7) -- (c5);
\draw (c7) -- (c6);}
\end{tikzpicture}
	
	\vspace{.1cm}	
{\small	Fig. 51. $K_6$ = $K_{1,1}(0; 1) $ \hspace{1.2cm} Fig. 52. $K_7$ = $K_{1,1}(0; 1)$ \hspace{1.2cm} Fig. 53. $K_8$ = $K_{1,1}(0; 1)$
	
\hspace{.5cm} $\cup K_{1,2}(2; 0,1) \cup \dots$ \hspace{2cm}  $\cup K_{1,2}(2; 0,1) \cup \dots$ \hspace{2cm}  $ \cup K_{1,2}(2; 0,1) \cup \dots$
	
\hspace{.5cm}	$\cup K_{1,5}(5; 0,1,2,3,4)$ \hspace{2cm}  $\cup K_{1,6}(6; 0,1,...,5)$ \hspace{2cm}  $\cup K_{1,7}(7; 0,1,...,6)$ }
\end{center}

\subsection{On $(a, d)$-CMSD of $G_{0,m}$ into triangular books with bookmark}

Here, we present results obtained in \cite{vs15} on $(a, d)$-CMSD of $G_{0,m}$ into triangular books with bookmark, $a,d,m\in\mathbb{N}$. 

\begin{theorem} [{\cite[Theorem 2.4]{vs15}}] \label{9.5} {\rm
	For $m\in\mathbb{N}$, $G_{0,2m}$ admits $(2, 2)$-CMSD into triangular books with bookmark.}
\end{theorem}

\begin{proof} For $n \geq 2$, in the sum graph $G_n$, $|E(G_n)| = \frac{1}{2}(\frac{n(n-1)}{2} -  \left\lfloor\frac{n}{2}\right\rfloor)$, $d(v_j) = n-1-j$ if $1 \leq j\leq  \left\lfloor\frac{n+1}{2}\right\rfloor$  and $d(v_j) = n-j$ if  $ \left\lfloor\frac{n+1}{2}\right\rfloor +1\leq j\leq n$ where $\left\lfloor x\right\rfloor$ is the floor of $x$. Therefore $|E(G_{0,2m})| = 2m +|E(G_{2m})| = 2m+ \frac{1}{2}\big(\frac{2m(2m-1)}{2} - \left\lfloor\frac{2m}{2}\right\rfloor \big) = m(m+1)$ where $G_{0,2m} = K_1*G_{2m}$. The proof of the theorem is similar to the proof given to Theorem \ref{9.2}. Let $V(G_{0,2m}) = \{v_0,v_1,v_2,\dots,v_{2m}\}$ where $j$ is the integral sum label of vertex $v_j$ in the integral sum graph $G_{0,2m}$, $0 \leq j\leq 2m$. $(2, 2)$-CMSD of $G_{0,2m}$ into triangular books with bookmark is obtained as follows.
	
	$G_{0,2m}$ = $TB_0(0,2m-1)(0,2m)$ $\cup$ $TB_1(1,2m-2; 0)(1,2m-1)$ $\cup$ $TB_2(2,2m-3; 0,1)(2,2m-2)$ $\cup$ $TB_3(3,2m-4; 0,1,2)(3,2m-3)$ $\cup$ $\dots$ $\cup$ $TB_{m-1}(m-1,m; 0,1,2, . . . , m-2)(m-1,m+1)$ where $TB_j(j,2m-(j+1); 0,1,2,\dots,j-1)(j,2m-j)$ represents triangular book with spine $(j,2m-(j+1))$, book mark $(j,2m-j)$ and leaves $0,1,2,\dots,j-1$ for $j$ = $1,2,\dots,m-1$ and $TB_0(0,2m-1)(0,2m)$ is the triangular book with spine $(0,2m-1)$, bookmark $(0,2m)$ and without any leaf. This implies, $G_{0,2m}$ admits $(2, 2)$-CMSD into triangular books with bookmark since $|E(TB_0(0, 2m-1)(0,2m))|$ = 2 $< |E(TB_1(1,2m-2; 0)(1,2m-1))|$ = 4 $< |E(TB_2(2,2m-3; 0,1)(2,2m-2))|$ = 6 $<$ $\dots$ $< |E(TB_{m-2}(m-2,m+1)(m-2,m+2))|$ = $2m-2$ $< |E(TB_{m-1}(m-1,m; 0,1,2,\dots,m-1)(m-1,m+1))|$ = $2m$ and $TB_0(0,2m-1)(0,2m)$ is a connected subgraph of $TB_1(1,2m-2; 0)(1,2m-1)$ which is a connected subgraph of $TB_2(2,2m-3; 0,1)(2,2m-2)$ which is a connected subgraph of $TB_3(3,2m-4; 0,1,2)(3,2m-3)$ which is a connected subgraph of . . . which is a connected subgraph of $TB_{m-1}(m-1,m; 0,1,2,...$, $m-2)(m-1,m+1)$. 	

\begin{center}
\begin{tikzpicture}[scale =.8]
	\node (c6) at (18,3)  [circle,draw,scale=0.6]{6};
	\node (c5) at (20,2)  [circle,draw,scale=0.6]{5};
	\node (c4) at (20,.5)  [circle,draw,scale=0.6] {4};
	\node (c3) at (19,-.75)  [circle,draw,scale=0.6] {3};
	\node (c2) at (17,-.75)  [circle,draw,scale=0.6] {2};
	\node (c1) at (16,0.5)  [circle,draw,scale=0.6]{1};
	\node (c0) at (16,2)  [circle,draw,scale=0.6] {0};
			
			\draw (c2) -- (c0);
			\draw (c2) -- (c1);
			\draw (c2) -- (c3);
			\draw (c2) -- (c4);
			\draw (c3) -- (c0);
			\draw (c3) -- (c1);

     {\color{red}
			\draw (c1) -- (c0);
			\draw (c1) -- (c4);
			\draw (c1) -- (c5);
			\draw (c4) -- (c0);}

   {\color{green}			
			\draw (c0) -- (c5);
			\draw (c0) -- (c6);}
			

\node (c8) at (26.5,3.5)  [circle,draw,scale=0.6] {8};
\node (c7) at (28,2.75)  [circle,draw,scale=0.6]{7};
\node (c6) at (28.5,1.5)  [circle,draw,scale=0.6]{6};
\node (c5) at (28,.25)  [circle,draw,scale=0.6] {5};
\node (c4) at (27,-.75)  [circle,draw,scale=0.6] {4};
\node (c3) at (25,-.75)  [circle,draw,scale=0.6] {3};
\node (c2) at (23.75,.25)  [circle,draw,scale=0.6]{2};
\node (c1) at (23.5,1.5)  [circle,draw,scale=0.6] {1};
\node (c0) at (24.5,2.75)  [circle,draw,scale=0.6] {0};

\draw (c3) -- (c0);
\draw (c3) -- (c1);
\draw (c3) -- (c2);
\draw (c3) -- (c4);
\draw (c3) -- (c5);
\draw (c4) -- (c0);
\draw (c4) -- (c1);
\draw (c4) -- (c2);

\draw (c2) -- (c0);
\draw (c2) -- (c1);
\draw (c2) -- (c5);
\draw (c6) -- (c1);

{\color{red}
	\draw (c2) -- (c0);
	\draw (c2) -- (c1);
	\draw (c2) -- (c5);
	\draw (c2) -- (c6);

    \draw (c5) -- (c0);
    \draw (c5) -- (c1);}

{\color{green}			
	\draw (c1) -- (c0);
	\draw (c1) -- (c6);
	\draw (c0) -- (c6);
    \draw (c1) -- (c7);}

{\color{blue}			
	\draw (c0) -- (c7);
	\draw (c0) -- (c8);}

	\end{tikzpicture}				
	
\medskip
\small{	Fig. 54. $G_{0,6}$ = $TB_0(0,5)(0,6)$  \hspace{3cm}  Fig. 55. $G_{0,8}$ = $TB_0(0,7)(0,8)$
	
\hspace{1cm} $\cup TB_1(1,4; 0)(1,5)$  \hspace{3cm}  $\cup$ $TB_1(1,6; 0)(1,7)$  

\hspace{2.5cm}  $\cup$ $TB_2(2,3;0,1)(2,4)$ \hspace{2.5cm}  $\cup$ $TB_2(2,5; 0,1)(2,6)$	$\cup TB_3(3,4; 0, 1,2)(3,5)$ }
\end{center}

\medskip

\begin{center}
	\begin{tikzpicture}[scale = 0.7]
		
	\node (c10) at (25,5.7)  [circle,draw,scale=0.5] {10};
	\node (c9) at (26,4.8)  [circle,draw,scale=0.6] {9};
	\node (c8) at (27,3.5)  [circle,draw,scale=0.6]{8};
	\node (c7) at (27.5,2)  [circle,draw,scale=0.6]{7};
	\node (c6) at (28,.25)  [circle,draw,scale=0.6] {6};
	\node (c5) at (26,-.75)  [circle,draw,scale=0.6] {5};
	\node (c4) at (23,-.75)  [circle,draw,scale=0.6] {4};
	\node (c3) at (20.5,.25)  [circle,draw,scale=0.6]{3};
	\node (c2) at (19.5,1.75)  [circle,draw,scale=0.6] {2};
	\node (c1) at (19.75,3.5)  [circle,draw,scale=0.6] {1};
	\node (c0) at (21,4.75)  [circle,draw,scale=0.6] {0};

\draw (c4) -- (c0);
\draw (c4) -- (c1);
\draw (c4) -- (c2);
\draw (c4) -- (c3);
\draw (c4) -- (c5);
\draw (c4) -- (c6);

\draw (c5) -- (c0);
\draw (c5) -- (c1);
\draw (c5) -- (c2);
\draw (c5) -- (c3);
		
	{\color{red}
		\draw (c3) -- (c0);
		\draw (c3) -- (c1);
		\draw (c3) -- (c2);
		\draw (c3) -- (c6);
		\draw (c3) -- (c7);
		
		\draw (c6) -- (c0);
		\draw (c6) -- (c1);
		\draw (c6) -- (c2);}
	
	{\color{green}			
\draw (c2) -- (c0);
\draw (c2) -- (c1);
\draw (c2) -- (c7);
\draw (c2) -- (c8);

\draw (c7) -- (c0);
\draw (c7) -- (c1);}
	
	{\color{blue}			
		\draw (c1) -- (c0);
		\draw (c1) -- (c8);
		\draw (c1) -- (c9);

		\draw (c8) -- (c0);}
	
	{\color{orange}			
		\draw (c0) -- (c9);
		\draw (c0) -- (c10);}
	
	\end{tikzpicture}				
	
	\vspace{.1cm}	
\small{	Fig. 56. $G_{0,10}$ = $TB_0(0,9)(0,10)$ $\cup$ $TB_1(1,8;0)(1,9) \cup  TB_2(2,7;0,1)(2,8) $ \\ $\cup TB_3(3,6;$ $0,1,2)(3,7)$ $\cup$ $TB_4(4,5; 0,1,2,3)(4,6)$}	
\end{center}

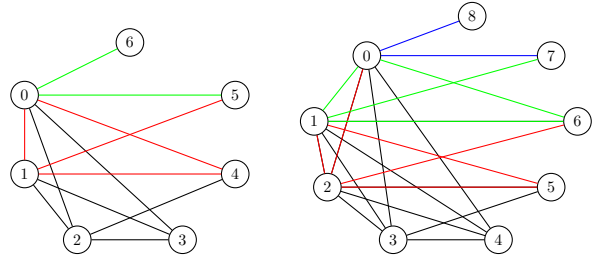
Graphs $G_{0,6}, G_{0,8}$ and $G_{0,10}$ are shown in their $(2,2)$-CMSD in different colors in Figures 54, 55 and 56,  respectively. In these $(2,2)$-CMSDs, $G_{0,6}$ = $TB_0(0,5)(0,6)\cup TB_1(1,4; 0)(1,5)$ $\cup$ $TB_2(2,3;0,1)$ $(2,4)$, $G_{0,8}$ = $TB_0(0,7)(0,8)$ $\cup$ $TB_1(1,6; 0)(1,7)$ $\cup$ $TB_2(2,5; 0,1)(2,6)$ $\cup$ $TB_3(3,4;$ $0, 1,2)(3,5)$ and $G_{0,10} = TB_0(0,9)(0,10) \cup TB_1(1$, $8;0)(1,9) \cup  TB_2(2,7;0,1)(2,8) \cup TB_3(3,6; 0,1,2)(3,7) \cup TB_4(4,5; 0,1,2,3)(4,6)$. 

Hence the result is proved. 
\end{proof}

\begin{theorem} [{\cite[Theorem 2.5]{vs15}}] \label{9.6} {\rm
	For $n\in\mathbb{N}$, $G_{0,4n+1}$ doesn't admit $(a, d)$-ASD and $(a,d)$-CMD into triangular books with bookmark for any $a,d\in\mathbb{N}$.}
\end{theorem}

\begin{proof} For $n\in\mathbb{N}$ and $k\in \mathbb{N}_0$, we have $|E(G_{0,4n+1})|$ = $\frac{1}{2}\big((4n+1)(4n+4)/2 -  \left\lfloor (4n+1)/2 \right\rfloor \big)$ = $(n+1)(4n+1) - n$ = $(2n+1)^2$ and $|E(TB_k(u,v)(u,w))|$ = $2k+2$. For $m,n\in\mathbb{N}$, if $G_{0,4n+1}$ admits $(a, d)$-ASD or $(a,d)$-CMD into triangular books with bookmark for any $a,d\in\mathbb{N}$, then let $G_{0,4n+1}$ = $TB_{k_1}^*(u_1,v_1)$ $\cup$ $TB_{k_2}^*(u_2,v_2)$ $\cup$ $\dots$ $\cup$ $TB_{k_m}^*(u_m,v_m)$ where $TB_{k_1}^*(u_1,v_1)$, $TB_{k_2}^*(u_2,v_2)$, $\dots$, $TB_{k_m}^*(u_m,v_m)$ are edge disjoint triangular books with bookmark in $G_{0,4n+1}$, $u_i,v_i\in V(G_{0,4n+1})$, $0 \leq  k_1 < k_2 < \dots < k_m$, $1 \leq i \leq m$ and $k_1,k_2,\dots,k_m \in\mathbb{N}_0$. Then $|E(G_{0,4n+1})|$ = $|E(TB_{k_1}^*(u_1,v_1))|$ + $|E(TB_{k_2}^*(u_2,v_2))|$ + $\dots$ + $|E(TB_{k_m}^*(u_m,v_m))|$ which implies, $(2n+1)^2$ = $(2k_1+2)$ + $(2k_2+2)$ + $\dots$ + $(2k_m+2)$ which is not possible since the L.H.S. is an odd number whereas the R.H.S. is an even number. Hence the result is true by the method of contradiction.  
\end{proof}

\begin{cor} [{\cite[Corollary 2.6]{vs15}}] \label{9.7} {\rm
	For $n\in\mathbb{N}$, $G_{0,4n+1}$ doesn't admit $(a, d)$-CMSD into triangular books with bookmark for any $a,d\in\mathbb{N}$. \hfill $\Box$}
\end{cor}

\begin{theorem} [{\cite[Theorem 2.7]{vs15}}] \label{9.8}  {\rm 	For $n\in\mathbb{N}$, 
\\
(i)~~~ $G_{0,4n}$ admits $(6, 8)$-CMSD into triangular books with bookmark;
\\
(ii)~~ $G_{0,4n+1}$ can be decomposed into triangular books with bookmark;
\\
(iii)~ $G_{0,4n+2}$ admits $(2, 8)$-CMSD into triangular books with bookmark and 
\\
(iv)~ $G_{0,4n+3}$ admits $(4, 8)$-CMSD into triangular books with bookmark.}
\end{theorem}

\begin{proof}\quad 	Let $V(G_{0,n})$ = $\{v_0, v_1, v_2, . . . , n\}$ and vertices of $G_{0,n}$ be subscript-labeled.  In the sum graph $G_n$, $|E(G_n)|$ = $\frac{1}{2}\big(n(n-1)/2 - \left\lfloor n/2 \right\rfloor\big)$, $d(v_j)$ = $n-1-j$ if $1\leq j \leq \left\lfloor (n+1)/2 \right\rfloor$ and $d(v_j)$ = $n-j$ if $\left\lfloor (n+1)/2 \right\rfloor + 1 \leq j \leq n$ where $\left\lfloor x \right\rfloor$ is the floor of $x$ and $n\in\mathbb{N}$. Therefore $|E(G_{0,n})|$ = $n +|E(G_n)|$ = $\frac{1}{2}\big(n(n+3)/2 - \left\lfloor n/2 \right\rfloor\big)$. Consider the following four cases of $n$ and the proof is similar to the proof given to Theorem \ref{9.2}. 
	
\medskip
\noindent
\textbf{Case (i)} $n$ = $4m$, $m\in\mathbb{N}$.\\
In this case, $(6, 8)$-CMSD of $G_{0,n} = G_{0,4m}$ into triangular books with bookmark is obtained as follows. 
	
$G_{0,4m}$ = $TB_{4m-2}(0,1; 2,3,\dots,4m-1)(0,4m)$ $\cup$ $TB_{4m-6}(2,3; 4,5,\dots,4m-3)(2,4m-2)$ $\cup$  $TB_{4m-10}(4$, $5; 6,7,...,4m-5)(4,4m-4)$ $\cup$ $TB_{4m-14}(6,7; 8,9,..., 4m-7)(6,4m-6)$ $\cup$ $\dots$ $\cup$ $TB_6(2m-4,2m-3; 2m-2, 2m-1,\dots,2m+3)(2m-4,2m+4)$ $\cup$ $TB_2(2m-2, 2m-1; 2m,2m+1)(2m-2, 2m+2)$ where  $TB_{4m-(2+4j)}(2j,2j+1$; $2j+2,2j+3$, $\dots$, $4m-2j-1)(2j,4m-2j)$ represents triangular book in $G_{0,4m}$ with spine $(2j,2j+1)$, {\it pendant vertex} with label $4m-2j$ and leaves $2j+2$, $2j+3$, $\dots$, $4m-2j-1$ and is a connected subgraph for $j$ = $0,1,2,...,m-1$. In this decomposition all the edges of $G_{0,4m}$ are partitioned into the edges of triangular books with bookmark and  $|E(TB_2(2m-2,2m-1; 2m,2m+1)$ $(2m-2,2m+2))|$ = 6 $<$ $|E(TB_6(2m-4, 2m-3; 2m-2,2m-1,\dots,2m+3)(2m-4, 2m+4))|$ = 14 $<$ $|E(TB_{10}(2m-6, 2m-5; 2m-4,2m-3,\dots,2m+5)(2m-6,2m+6))|$ = 22 $<$ $\dots$ $<$ $|E(TB_{4m-6}(2,3; 4,5,\dots,4m-3)(2,4m-2))|$ = $8m-10$ $<$ $|E(TB_{4m-2}(0,1; 2,3,...,4m-1)(0, 4m))|$ = $8m-2$ and $TB_{4m-2}(0,1; 2,3,...$, $4m-1)(0,4m)$ is a connected subgraph of $TB_{4m-6}(2,3; 4,5,\dots,4m-3)$ $(2,4m-2)$ which is a connected subgraph of $TB_{4m-10}(4,5; 6,7,\dots,4m-5)(4,4m-4)$ which is a connected subgraph of $TB_{4m-14}(6,7; 8,9,\dots,4m-7)(6,4m-6)$ which is a connected subgraph of $\dots$ which is a connected subgraph of $TB_6(2m-4,2m-3;$ $2m-2,2m-1,\dots,2m+3)(2m-4,2m+4)$ which is a connected subgraph of $TB_2(2m-2,2m-1; 2m,2m+1)(2m-2,2m+2)$. Thus, $G_{0,4m}$ admits $(6, 8)$-CMSD into triangular books with bookmark.  
	
	$(6,8)$-CMSD of $G_{0,4}, G_{0,8}$ and $G_{0,12}$ are shown in different colors in Figures 57, 58 and 59. Here, (6,8)-CMSD of the graphs are $G_{0,4}$ = $TB_2(0,1; 2,3)(0,4)$, $G_{0,8}$ = $TB_6(0,1; 2,3,4,5,6,7)(0,8)$ $\cup$ $TB_2(2,3; 4,5)$ $(2,6)$ and $G_{0,12}$ = $TB_{10}(0,1; 2,3,\dots, 11)(0,12)$ $\cup$ $TB_6(2,3; 4,5,\dots,9)(2,10)$ $\cup$ $TB_2(4,5; 6, 7)(4,8)$.	
	
	\begin{center}
		\begin{tikzpicture}[scale =0.8]
		\node (a4) at (9,1.25)  [circle,draw,scale=0.6] {4};
		\node (a3) at (10.5,2.5)  [circle,draw,scale=0.6] {3};
		\node (a2) at (12,1.25)  [circle,draw,scale=0.6]{2};
		\node (a1) at (11.5,-.5)  [circle,draw,scale=0.6] {1};
		\node (a0) at (9.5,-.5)  [circle,draw,scale=0.6] {0};
		
		\draw (a0) -- (a1);
		\draw (a0) -- (a2);
		\draw (a0) -- (a3);
		\draw (a0) -- (a4);
		
		\draw (a1) -- (a2);
		\draw (a1) -- (a3);
				
\node (c5) at (18.5,3.5)  [circle,draw,scale=0.6] {5};
\node (c4) at (20.5,2.75)  [circle,draw,scale=0.6]{4};
\node (c3) at (21,1.5)  [circle,draw,scale=0.6]{3};
\node (c2) at (20.5,0)  [circle,draw,scale=0.6] {2};
\node (c1) at (19.25,-.75)  [circle,draw,scale=0.6] {1};
\node (c0) at (17.75,-.75)  [circle,draw,scale=0.6] {0};
\node (c8) at (16.5,0)  [circle,draw,scale=0.6]{8};
\node (c7) at (16,1.5)  [circle,draw,scale=0.6] {7};
\node (c6) at (16.5,2.75)  [circle,draw,scale=0.6] {6};

\draw (c0) -- (c1);
\draw (c0) -- (c2);
\draw (c0) -- (c3);
\draw (c0) -- (c4);
\draw (c0) -- (c5);
\draw (c0) -- (c6);
\draw (c0) -- (c7);
\draw (c0) -- (c8);

\draw (c1) -- (c2);
\draw (c1) -- (c3);
\draw (c1) -- (c4);
\draw (c1) -- (c5);
\draw (c1) -- (c6);
\draw (c1) -- (c7);

{\color{red}
	\draw (c2) -- (c3);
	\draw (c2) -- (c4);
	\draw (c2) -- (c5);
	\draw (c2) -- (c6);
	
	\draw (c3) -- (c4);
    \draw (c3) -- (c5);}

\end{tikzpicture}				
	
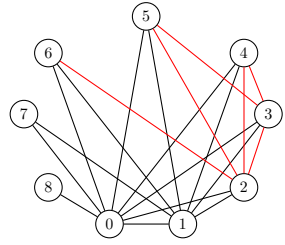
\small{	Fig. 57. $G_{0,4}$ = $TB_2(0,1:2,3)(0,4)$  \hspace{.5cm}  Fig. 58. $G_{0,8}$ = $TB_6(0,1:2,3,...,7)(0,8)$  
	
	\hspace{8cm}  $\cup$  $TB_2(2,3;4,5)(2,6)$}
\end{center}

\begin{center}
	\begin{tikzpicture}[scale = 0.6]
	
	\node (c8) at (22,5.5)  [circle,draw,scale=0.6] {8};
	\node (c7) at (24,6)  [circle,draw,scale=0.6] {7};
	\node (c6) at (26,6)  [circle,draw,scale=0.6] {6};
	\node (c5) at (28,5)  [circle,draw,scale=0.6] {5};
	\node (c4) at (29,3.6)  [circle,draw,scale=0.6]{4};
	\node (c3) at (28.7,1.5)  [circle,draw,scale=0.6]{3};
	\node (c2) at (27.5,.25)  [circle,draw,scale=0.6] {2};
	\node (c1) at (25,-.75)  [circle,draw,scale=0.6] {1};
	\node (c0) at (23,-.75)  [circle,draw,scale=0.6] {0};
	\node (c12) at (21,0)  [circle,draw,scale=0.5]{12};
	\node (c11) at (20,1.5)  [circle,draw,scale=0.5] {11};
	\node (c10) at (20,3)  [circle,draw,scale=0.5] {10};
	\node (c9) at (20.5,4.5)  [circle,draw,scale=0.6] {9};
	
\draw (c0) -- (c1);
\draw (c0) -- (c2);
\draw (c0) -- (c3);
\draw (c0) -- (c4);
\draw (c0) -- (c5);
\draw (c0) -- (c6);
\draw (c0) -- (c7);
\draw (c0) -- (c8);
\draw (c0) -- (c9);
\draw (c0) -- (c10);
\draw (c0) -- (c11);
\draw (c0) -- (c12);

\draw (c1) -- (c2);
\draw (c1) -- (c3);
\draw (c1) -- (c4);
\draw (c1) -- (c5);
\draw (c1) -- (c6);
\draw (c1) -- (c7);
\draw (c1) -- (c8);
\draw (c1) -- (c9);
\draw (c1) -- (c10);
\draw (c1) -- (c11);
	
{\color{red}
	\draw (c2) -- (c3);
	\draw (c2) -- (c4);
	\draw (c2) -- (c5);
	\draw (c2) -- (c6);
	\draw (c2) -- (c7);
	\draw (c2) -- (c8);
	\draw (c2) -- (c9);
	\draw (c2) -- (c10);
	
	\draw (c3) -- (c4);
	\draw (c3) -- (c5);
	\draw (c3) -- (c6);
	\draw (c3) -- (c7);
	\draw (c3) -- (c8);
	\draw (c3) -- (c9);}
	
	{\color{green}			
		\draw (c4) -- (c5);
		\draw (c4) -- (c6);
		\draw (c4) -- (c7);
		\draw (c4) -- (c8);
		
		\draw (c5) -- (c6);
		\draw (c5) -- (c7);}
		
	\end{tikzpicture}				
	
	\vspace{.2cm}	
	\small{	Fig. 59. $G_{0,12} = TB_{10}(0,1; 2,3,\dots$, $11)(0,12)$ 
		
	\hspace{4cm}	$\cup$ $TB_6(2,3; 4$, $5,\dots,9)(2,10)$ $\cup TB_2(4,5; 6, 7)(4,8)$}	
\end{center}

\noindent
\textbf{Case (ii):}~ $n$ = $4m+1$, $m\in\mathbb{N}$.\\
In this case, decomposition of $G_{0,4m+1}$ into triangular books with bookmark is obtained as follows.

$G_{0,4m+1}$ = $TB_{4m-1}(0,1; 2,3,\dots,4m)(0,4m+1)$ $\cup$ $TB_{4m-5}(2,3; 4,5,\dots,4m-2)(2,4m-1)$ $\cup$ $TB_{4m-9}$ $(4,5; 6,7,\dots,4m-4)(4,4m-3)$ $\cup$ $TB_{4m-13}(6,7; 8,9,\dots, 4m-6)(6,4m-5)$ $\cup$ $\dots$ $\cup$ $TB_7(2m-4,2m-3; 2m-2,2m-1,\dots,2m+4)(2m-4,2m+5)$ $\cup$ $TB_3(2m-2,2m-1; 2m,2m+1,2m+2)(2m-2,2m+3)$ $\cup$ $TB_0(2m,2m+1)$ where $TB_{4m+1-(2+4j)}(2j,2j+1; 2j+2,2j+3,\dots,4m-2j)(2j, 4m+1-2j)$ represents triangular book in $G_{0,4m+1}$ with spine $(2j,2j+1)$, pendant vertex with label $4m+1-2j$ and leaves $2j+2,2j+3,\dots,4m-2j$ and is a connected subgraph for $j$ = $0,1,2,\dots,m-1$ and $TB_0(2m,2m+1)$ is a triangular book with spine $(2m,2m+1)$ and without any leaf. All the edges of $G_{0,4m+1}$ are covered under this decomposition and $|E(TB_0(2m,2m+1))|$ = 1 $<$ $|E(TB_3(2m-2,2m-1;2m, 2m+1$, $2m+2)(2m-2,2m+3))|$ = 8 $<$ $|E(TB_7(2m-4,2m-3;$ $ 2m-2,2m-1,...,2m+4)(2m-4,2m+5))|$ = 16 $< \dots <$ $|E(TB_{4m-5}(2,3; 4,5,\dots,4m-2)(2,4m-1))|$ = $8m-8$ $<$ $|E(TB_{4m-1}(0,1; 2,3, \dots,4m)(0,4m+1))|$ = $8m$ and $TB_0(2m,2m+1)$ is a connected subgraph of $TB_3(2m-2,2m-1; 2m,2m+1,2m+2)(2m-2,2m+3)$ which is a connected subgraph of $TB_7(2m-4,2m-3; 2m-2,2m-1,...,2m+4) (2m-4,2m+5)$ which is a connected subgraph of ... which is a connected subgraph of $TB_{4m-5}(2,3; 4,5,...,4m-2)(2,4m-1)$ which is a connected subgraph of $TB_{4m-1}(0,1; 2,3,\dots,4m)(0,4m+1)$, without vertex labels. Thus, $G_{0,4m+2}$ is decomposed into triangular books with bookmark. 
	
Decomposition of graphs $G_{0,5}$, $G_{0,9}$ and $G_{0,13}$ into triangular books with bookmark are shown in different colors in Figures 60, 61 and 62. Here, the decomposition is $G_{0,5}$ = $TB_3(0,1; 2,3,4)(0,5)$ $\cup$ $TB_0(2,3)$, $G_{0,9}$ = $TB_7(0,1; 2,3,\dots,8)(0,9 )$ $\cup$ $TB_3(2,3; 4,5,6)(2,7)$ $\cup$ $TB_0(4,5)$ and $G_{0,13}$ = $TB_{11}(0,1$; $2,3,\dots,12)(0,13)$ $\cup$ $TB_7(2,3; 4,5,\dots,10)(2,11)$ $\cup$ $TB_3(4,5; 6,7,8)(4,9)$ $\cup$ $TB_0(6,7)$.

\begin{center}
	\begin{tikzpicture}[scale =.8]
	\node (a5) at (9,1)  [circle,draw,scale=0.6] {5};
	\node (a4) at (10,2.5)  [circle,draw,scale=0.6] {4};
	\node (a3) at (12,2.5)  [circle,draw,scale=0.6] {3};
	\node (a2) at (13,1)  [circle,draw,scale=0.6]{2};
	\node (a1) at (12,-.5)  [circle,draw,scale=0.6] {1};
	\node (a0) at (10,-.5)  [circle,draw,scale=0.6] {0};
	
	\draw (a0) -- (a1);
	\draw (a0) -- (a2);
	\draw (a0) -- (a3);
	\draw (a0) -- (a4);
	\draw (a0) -- (a5);
	
	\draw (a1) -- (a2);
	\draw (a1) -- (a3);
	\draw (a1) -- (a4);
{\color{red}
	\draw (a2) -- (a3);}
	
	\node (c9) at (16,0)  [circle,draw,scale=0.6]{9};
	\node (c8) at (16,1.25)  [circle,draw,scale=0.6]{8};
    \node (c7) at (16.5,2.5)  [circle,draw,scale=0.6] {7};
    \node (c6) at (17.5,3.25)  [circle,draw,scale=0.6] {6};
	\node (c5) at (18.5,3.5)  [circle,draw,scale=0.6] {5};
	\node (c4) at (20,2.75)  [circle,draw,scale=0.6]{4};
	\node (c3) at (20.5,1.5)  [circle,draw,scale=0.6]{3};
	\node (c2) at (20,0)  [circle,draw,scale=0.6] {2};
	\node (c1) at (18.75,-.75)  [circle,draw,scale=0.6]{1};
	\node (c0) at (17.25,-.75)  [circle,draw,scale=0.6] {0};
	
	\draw (c0) -- (c1);
	\draw (c0) -- (c2);
	\draw (c0) -- (c3);
	\draw (c0) -- (c4);
	\draw (c0) -- (c5);
	\draw (c0) -- (c6);
	\draw (c0) -- (c7);
	\draw (c0) -- (c8);
	\draw (c0) -- (c9);
	
	\draw (c1) -- (c2);
	\draw (c1) -- (c3);
	\draw (c1) -- (c4);
	\draw (c1) -- (c5);
	\draw (c1) -- (c6);
	\draw (c1) -- (c7);
	\draw (c1) -- (c8);
	
	{\color{red}
		\draw (c2) -- (c3);
		\draw (c2) -- (c4);
		\draw (c2) -- (c5);
		\draw (c2) -- (c6);		
		\draw (c2) -- (c7);
		
		\draw (c3) -- (c4);
		\draw (c3) -- (c5);
		\draw (c3) -- (c6);}

	{\color{green}	
        \draw (c4) -- (c5);}
	\end{tikzpicture}				
	
	\small{	\hspace{2.5cm} Fig. 60.  \hspace{3cm}  Fig. 61. $G_{0,9}$ = $TB_7(0,1;2,3,...,8)(0,9)$   
		
\hspace{.5cm}$G_{0,5}$ 	= $TB_3(0,1;2,3,4)(0,5)$ $\cup$ $TB_0(2,3)$  \hspace{1cm}	 $\cup$ $TB_3(2,3; 4,5,6)(2,7)$ $\cup$ $TB_0(4,5)$}
\end{center}

\noindent
\textbf{Case (iii):}~ $n = 4m+2$, $m\in\mathbb{N}.$

In this case, $(2, 8)$-CMSD of $G_{0,4m+2}$ into triangular books with bookmark is obtained as follows. 

$G_{0,4m+2}$ = $TB_{4m}(0,1; 2,3,\dots,4m+1)(0,4m+2)$ $\cup$ $TB_{4m-4}(2,3; 4,5,\dots,4m-1)$ $(2,4m)$ $\cup$ $TB_{4m-8}(4$, $5; 6,7,\dots,4m-3)(4,4m-2)$ $\cup$ $TB_{4m-12}(6,7; 8,9,\dots,4m-5)(6,4m-4)$ $\cup$ $\dots$ $\cup$ $TB_8(2m-4,2m-3$; $2m-2, 2m-1,\dots,2m+5)(2m-4,2m+6)$ $\cup$ $TB_4(2m-2$, $2m-1; 2m,2m+1,2m+2, 2m+3)(2m-2,2m+4)$ $\cup$ $TB_0(2m,2m+1)$ $(2m,2m+2)$. Here $TB_{4m-4j}(2j,2j+1; 2j+2,2j+3, \dots, 4m-2j+1)(2j,4m-2j+2)$ represents triangular book in $G_{0,4m+2}$ with spine $(2j,2j+1)$, pendant vertex $4m-2j+2$ and leaves $2j+2,2j+3,\dots,4m-2j+1$ and is a connected subgraph for $j$ = $0,1,2,...,m-1$ and $TB_0(2m,2m+ 1)(2m,2m+2)$ is a triangular book with spine $(2m,2m+1)$, pendant vertex with label $2m+2$ and without any leaf. In this decomposition, all the edges of $G_{0,4m+2}$ are partitioned into edges of triangular books with bookmark and $|E(TB_0(2m,2m+1)$ $(2m,2m+2))|$ = 2 $<$ $|E(TB_4(2m-2,2m-1; 2m,2m+1,2m+2,2m+3)(2m-2,2m+4))|$ = 10 $<$ $|E(TB_8(2m-4,2m-3$; $2m-2, 2m-1,\dots,2m+5)(2m-4,2m+6))|$ = 18 $< \dots <$ $|E(TB_{4m-4}(2,3; 4,5,\dots,4m-1)(2, 4m))|$ = $8m-6$ $<$ $|E(TB_{4m}(0,1; 2,3,\dots,4m+1)(0,4m+2))|$ = $8m+2$ and $TB_0(2m,2m+1)(2m,2m+2)$ is a connected subgraph of $TB_4(2m-2,2m-1; 2m,2m+1,2m+2, 2m+3)(2m-2,2m+4)$ which is a connected subgraph of $TB_8(2m-4,2m-3; 2m-2,2m-1,\dots,2m+5)(2m-4,2m+6)$ which is a connected subgraph of $\dots$ which is a connected subgraph of $TB_{4m-4}(2,3; 4,5,\dots,4m-1)(2,4m)$ which is a connected subgraph of $TB_{4m}(0,1; 2,3,\dots,4m+1)(0,4m+2)$, without vertex labels.
Thus, $G_{0,4m+2}$ admits $(2, 8)$-CMSD into triangular books with bookmark. 

\begin{center}
	\begin{tikzpicture}[scale = 0.5]	
	
\node (c13) at (21,-.25)  [circle,draw,scale=0.5]{13};
\node (c12) at (20.5,1)  [circle,draw,scale=0.5]{12};
\node (c11) at (20,3)  [circle,draw,scale=0.5] {11};
\node (c10) at (20.5,4.5)  [circle,draw,scale=0.5] {10};
\node (c9) at (21.5,6)  [circle,draw,scale=0.6] {9};	
\node (c8) at (23,7.25)  [circle,draw,scale=0.6] {8};
\node (c7) at (25.25,8)  [circle,draw,scale=0.6] {7};
\node (c6) at (27.75,7.25)  [circle,draw,scale=0.6] {6};
\node (c5) at (29.5,5.75)  [circle,draw,scale=0.6] {5};
\node (c4) at (30,4)  [circle,draw,scale=0.6]{4};
\node (c3) at (29.5,1.8)  [circle,draw,scale=0.6]{3};
\node (c2) at (28,0)  [circle,draw,scale=0.6] {2};
\node (c1) at (25.5,-1)  [circle,draw,scale=0.6] {1};
\node (c0) at (23,-1)  [circle,draw,scale=0.6] {0};
	
	\draw (c0) -- (c1);
	\draw (c0) -- (c2);
	\draw (c0) -- (c3);
	\draw (c0) -- (c4);
	\draw (c0) -- (c5);
	\draw (c0) -- (c6);
	\draw (c0) -- (c7);
	\draw (c0) -- (c8);
	\draw (c0) -- (c9);
	\draw (c0) -- (c10);
	\draw (c0) -- (c11);
	\draw (c0) -- (c12);
	\draw (c0) -- (c13);
	
	\draw (c1) -- (c2);
	\draw (c1) -- (c3);
	\draw (c1) -- (c4);
	\draw (c1) -- (c5);
	\draw (c1) -- (c6);
	\draw (c1) -- (c7);
	\draw (c1) -- (c8);
	\draw (c1) -- (c9);
	\draw (c1) -- (c10);
	\draw (c1) -- (c11);
	\draw (c1) -- (c12);
	
	{\color{red}
		\draw (c2) -- (c3);
		\draw (c2) -- (c4);
		\draw (c2) -- (c5);
		\draw (c2) -- (c6);
		\draw (c2) -- (c7);
		\draw (c2) -- (c8);
		\draw (c2) -- (c9);
		\draw (c2) -- (c10);
		\draw (c2) -- (c11);
		
		\draw (c3) -- (c4);
		\draw (c3) -- (c5);
		\draw (c3) -- (c6);
		\draw (c3) -- (c7);
		\draw (c3) -- (c8);
		\draw (c3) -- (c9);
		\draw (c3) -- (c10);}
	
	{\color{green}			
		\draw (c4) -- (c5);
		\draw (c4) -- (c6);
		\draw (c4) -- (c7);
		\draw (c4) -- (c8);
		\draw (c4) -- (c9);
		
		\draw (c5) -- (c6);
		\draw (c5) -- (c7);
		\draw (c5) -- (c8);}

{\color{orange}			
	\draw (c6) -- (c7);}
	
	
	\node (d14) at (36,-.5)  [circle,draw,scale=0.5]{14};
	\node (d13) at (35,.5)  [circle,draw,scale=0.5]{13};
	\node (d12) at (35,1.75)  [circle,draw,scale=0.5]{12};
	\node (d11) at (35,3.5)  [circle,draw,scale=0.5] {11};
	\node (d10) at (35.5,4.75) [circle,draw,scale=0.5]{10};
	\node (d9) at (36.5,6.5)  [circle,draw,scale=0.6] {9};	
	\node (d8) at (38,7.5)  [circle,draw,scale=0.6] {8};
	\node (d7) at (40,8)  [circle,draw,scale=0.6] {7};
	\node (d6) at (42.5,7.5)  [circle,draw,scale=0.6] {6};
	\node (d5) at (44,6.5)  [circle,draw,scale=0.6] {5};
	\node (d4) at (45,4)  [circle,draw,scale=0.6]{4};
	\node (d3) at (44.5,2)  [circle,draw,scale=0.6]{3};
	\node (d2) at (42.5,0)  [circle,draw,scale=0.6] {2};
	\node (d1) at (40,-1)  [circle,draw,scale=0.6] {1};
	\node (d0) at (38,-1)  [circle,draw,scale=0.6] {0};
	
	\draw (d0) -- (d1);
	\draw (d0) -- (d2);
	\draw (d0) -- (d3);
	\draw (d0) -- (d4);
	\draw (d0) -- (d5);
	\draw (d0) -- (d6);
	\draw (d0) -- (d7);
	\draw (d0) -- (d8);
	\draw (d0) -- (d9);
	\draw (d0) -- (d10);
	\draw (d0) -- (d11);
	\draw (d0) -- (d12);
	\draw (d0) -- (d13);
	\draw (d0) -- (d14);
	
	\draw (d1) -- (d2);
	\draw (d1) -- (d3);
	\draw (d1) -- (d4);
	\draw (d1) -- (d5);
	\draw (d1) -- (d6);
	\draw (d1) -- (d7);
	\draw (d1) -- (d8);
	\draw (d1) -- (d9);
	\draw (d1) -- (d10);
	\draw (d1) -- (d11);
	\draw (d1) -- (d12);
	\draw (d1) -- (d13);
	
	{\color{red}
		\draw (d2) -- (d3);
		\draw (d2) -- (d4);
		\draw (d2) -- (d5);
		\draw (d2) -- (d6);
		\draw (d2) -- (d7);
		\draw (d2) -- (d8);
		\draw (d2) -- (d9);
		\draw (d2) -- (d10);
		\draw (d2) -- (d11);
    	\draw (d2) -- (d12);
		
		\draw (d3) -- (d4);
		\draw (d3) -- (d5);
		\draw (d3) -- (d6);
		\draw (d3) -- (d7);
		\draw (d3) -- (d8);
		\draw (d3) -- (d9);
		\draw (d3) -- (d10);
		\draw (d3) -- (d11);}
	
	{\color{green}			
		\draw (d4) -- (d5);
		\draw (d4) -- (d6);
		\draw (d4) -- (d7);
		\draw (d4) -- (d8);
		\draw (d4) -- (d9);
    	\draw (d4) -- (d10);
		
		\draw (d5) -- (d6);
		\draw (d5) -- (d7);
		\draw (d5) -- (d8);
    	\draw (d5) -- (d9);}
	
	{\color{orange}			
		\draw (d6) -- (d7);
		\draw (d6) -- (d8);}
	
\end{tikzpicture}				
	
\medskip	
\small{	Fig. 62. $G_{0,13} = TB_{11}(0,1; 2,3,\dots$, $12)(0,13)$ \hspace{1cm}  Fig. 63. $G_{0,14} = TB_{12}(0,1; 2,3,\dots,13)(0,14)$ 
		
\hspace{2cm}   $\cup$ $TB_7(2,3; 4, 5,\dots,10)(2,11)$	\hspace{3.5cm} $\cup$ $TB_8(2,3; 4, 5,\dots,11)(2,12)$  
  
 \hspace{1.5cm}   $\cup$ $TB_3(4,5; 6, 7, 8)(4,9)$ $\cup TB_0(6, 7)$ \hspace{1.5cm} $\cup$ $TB_4(4,5; 6, 7, 8, 9)(4,10)$ $\cup$ $TB_0(6, 7)(6,8)$  }
	\end{center}
	
\medskip
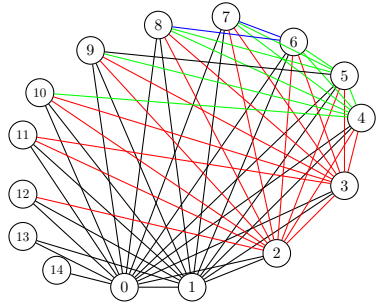
In Figure 63, $(2,8)$-CMSD of $G_{0,14}$  = $TB_{12}(0,1; 2,3,\dots, 13)(0,14)$ $\cup$ $TB_8(2,3; 4,5,\dots,11)(2, 12)$ $\cup$ $TB_4(4,5; 6,7,8,9)(4,10)$ $\cup$ $TB_0(6,7)(6,8)$ is presented. 

\medskip
\noindent	
\textbf{Case (iv):}~$n = 4m+3$, $m\in\mathbb{N}$.

	In this case, (4,8)-CMSD of $G_{0,4m+3}$ into triangular books with bookmark is obtained as follows. 
	
	$G_{0,4m+3}$ = $TB_{4m+1}(0,1; 2,3,\dots,4m+2)(0,4m+3)$ $\cup$ $TB_{4m-3}(2,3; 4,5,\dots,4m)(2,4m+1)$ $\cup$ $TB_{4m-7}$ $(4,5; 6,7,\dots, 4m-2)(4,4m-1)$ $\cup$ $TB_{4m-11}(6,7; 8,9,\dots,4m-4)(6,4m-3)$ $\cup$ $\dots$ $\cup$ $TB_9(2m-4,2m-3$; $2m-2,2m-1,\dots,2m+6)(2m-4, 2m+7)$ $\cup$ $TB_5(2m-2,2m-1; 2m,2m+1,\dots, 2m+4)(2m-2,2m+5)$ $\cup$ $TB_1(2m,2m+1;2m+2)(2m,2m+3)$ where $TB_{4m+1-4j}(2j,2j+1; 2j+2,2j+3,\dots,4m-2j+2)(2j,4m-2j+3)$ represents triangular book in $G_{0,4m+3}$ with spine $(2j,2j+1)$, pendant vertex $4m-2j+3$ and leaves $2j+2,2j+3,\dots,4m-2j+2$ and is a connected subgraph for $j$ = $0,1,2,\dots,m$. In this decomposition all the edges of $G_{0,4m+3}$ are partitioned into edges of triangular books with bookmark and $|E(TB_1(2m,2m+1; 2m+2)(2m,2m+3))|$ = $4 < |E(TB_5(2m-2,2m-1; 2m,2m+1,\dots,2m+4)$ $(2m-2,2m+5))|$ = $12 < |E(TB_9(2m-4,2m-3; 2m-2,2m-1,\dots,2m+6)(2m-4, 2m+7))|$ = $20 < \dots <$ $|E(TB_{4m-3}(2,3; 4,5,\dots,4m)(2,4m+1))|$ = $8m-4 <  |E(TB_{4m+1}(0,1;2,3,\dots,4m+2)(0,4m+3))|$ = $8m+4$
	and $TB_1(2m,2m+1; 2m+2)$ $(2m,2m+3)$ is a connected subgraph of $TB_5(2m-2$, $2m-1; 2m,2m+1,\dots, 2m+4)(2m-2,2m+5)$ which is a connected subgraph of $TB_9(2m-4,2m-3; 2m-2,2m-1,\dots, 2m+6)(2m-4,2m+7)$ which is a connected subgraph of $\dots$ which is a connected subgraph of $TB_{4m-3}(2,3; 4,5,\dots,4m)(2,4m+1)$ which is a connected subgraph of $TB_{4m+1}(0,1; 2,3,\dots,4m+2)(0,4m+3)$, without vertex labels. Thus, $G_{0,4m+3}$ admits $(4, 8)$-CMSD into triangular books with bookmark. Hence we get the result. 
	
	$(4,8)$-CMSD of $G_{0,15}$ = $TB_{13}(0,1; 2,3,\dots,14)(0,15)$ $\cup$ $TB_9(2,3; 4,5,\dots,12)(2, 13)$ $\cup$ $TB_5(4,5; 6,7,8$, $9,10)(4,11)$ $\cup$ $TB_1(6,7;8)(6,9)$ is shown in Figure 64.
\end{proof}

	\begin{center}
		\begin{tikzpicture}[scale = 0.6]	
		
	\node (c15) at (21,-2.5)  [circle,draw,scale=0.5]{15};
	\node (c14) at (19.5,-1)  [circle,draw,scale=0.5]{14};
	\node (c13) at (18.5,.5)  [circle,draw,scale=0.5]{13};
	\node (c12) at (18,2)  [circle,draw,scale=0.5]{12};
	\node (c11) at (18.5,3.5)  [circle,draw,scale=0.5] {11};
	\node (c10) at (19.5,4.75)[circle,draw,scale=0.5] {10};
	\node (c9) at (20.5,6) [circle,draw,scale=0.6] {9};	
	\node (c8) at (22,7)  [circle,draw,scale=0.6] {8};
	\node (c7) at (24.5,7.5)  [circle,draw,scale=0.6] {7};
	\node (c6) at (27,7)  [circle,draw,scale=0.6] {6};
	\node (c5) at (29,5.5)  [circle,draw,scale=0.6] {5};
	\node (c4) at (30,3.5)  [circle,draw,scale=0.6]{4};
	\node (c3) at (30,1.5)  [circle,draw,scale=0.6]{3};
	\node (c2) at (29,-.5)  [circle,draw,scale=0.6] {2};
	\node (c1) at (26.5,-2.5)  [circle,draw,scale=0.6] {1};
	\node (c0) at (24,-2.5)  [circle,draw,scale=0.6] {0};
		
		\draw (c0) -- (c1);
		\draw (c0) -- (c2);
		\draw (c0) -- (c3);
		\draw (c0) -- (c4);
		\draw (c0) -- (c5);
		\draw (c0) -- (c6);
		\draw (c0) -- (c7);
		\draw (c0) -- (c8);
		\draw (c0) -- (c9);
		\draw (c0) -- (c10);
		\draw (c0) -- (c11);
		\draw (c0) -- (c12);
		\draw (c0) -- (c13);
		\draw (c0) -- (c14);
		\draw (c0) -- (c15);
		
		\draw (c1) -- (c2);
		\draw (c1) -- (c3);
		\draw (c1) -- (c4);
		\draw (c1) -- (c5);
		\draw (c1) -- (c6);
		\draw (c1) -- (c7);
		\draw (c1) -- (c8);
		\draw (c1) -- (c9);
		\draw (c1) -- (c10);
		\draw (c1) -- (c11);
		\draw (c1) -- (c12);
		\draw (c1) -- (c13);
		\draw (c1) -- (c14);
		
		{\color{red}
			\draw (c2) -- (c3);
			\draw (c2) -- (c4);
			\draw (c2) -- (c5);
			\draw (c2) -- (c6);
			\draw (c2) -- (c7);
			\draw (c2) -- (c8);
			\draw (c2) -- (c9);
			\draw (c2) -- (c10);
			\draw (c2) -- (c11);
			\draw (c2) -- (c12);
     		\draw (c2) -- (c13);
			
			\draw (c3) -- (c4);
			\draw (c3) -- (c5);
			\draw (c3) -- (c6);
			\draw (c3) -- (c7);
			\draw (c3) -- (c8);
			\draw (c3) -- (c9);
			\draw (c3) -- (c10);
			\draw (c3) -- (c11);
			\draw (c3) -- (c12);}
		
		{\color{green}			
			\draw (c4) -- (c5);
			\draw (c4) -- (c6);
			\draw (c4) -- (c7);
			\draw (c4) -- (c8);
			\draw (c4) -- (c9);
			\draw (c4) -- (c10);
			\draw (c4) -- (c11);
			
			\draw (c5) -- (c6);
			\draw (c5) -- (c7);
			\draw (c5) -- (c8);
	    	\draw (c5) -- (c9);
    		\draw (c5) -- (c10);}

		{\color{orange}			
			\draw (c6) -- (c7);
			\draw (c6) -- (c8);
			\draw (c6) -- (c9);
		
			\draw (c7) -- (c8);	}
		
		\end{tikzpicture}				
		
		\vspace{.1cm}	
		\small{	Fig. 64. $G_{0,15} = TB_{13}(0,1; 2,3,\dots,14)(0,15) \cup TB_9(2,3; 4, 5,\dots,12)(2,13)$
			
			\hspace{3cm}	 $\cup TB_5(4,5; 6, 7, 8, 9, 10)(4,11)$ $\cup TB_1(6, 7; 8)(6,9)$}
\end{center}
	
\subsection{On $(a, d)$-CMD and $(a,d)$-CMSD of $G_{0,m}$ into fan with a handle}
	
	In this section, we present results on $(a,d)$-Ascending subgraph decomposition ($(a, d)$-ASD), $(a,d)$-continuous monotonic decomposition ($(a, d)$-CMD) and $(a,d)$-continuous monotonic subgraph decomposition ($(a, d)$-CMSD) of $G_{0,m}$ into fan with a handle, $a,d,m\in\mathbb{N}$. Fan with a handle is defined in section 5.2.

\begin{theorem} [{\cite[Theorem 2.8]{vs15}}] \label{9.9}  {\rm
	For $m\in\mathbb{N}$, $G_{0,4m+1}$ does not admit $(a,d)$-ASD and $(a,d)$-CMD into Fans with a handle for any $a,d\in\mathbb{N}$.}
\end{theorem}

\begin{proof}
	If possible, let $G_{0,4m+1}$ admit $(a,d)$-ASD into Fans with a handle for some $a,d\in\mathbb{N}$. Then, let $G_{0,4m+1}\cong F_{n_1}^*~\cup F_{n_2}^*~\cup\dots~\cup  F_{n_k}^*$ where $F_{n_1}^*, F_{n_2}^*, \dots, F_{n_k}^*$ are edge disjoint fans with handle for some $n_1,n_2,\dots,n_k \in\mathbb{N}$ and $2\leq n_1 < n_2 <\dots < n_k$. Then, $|E(G_{0,4m+1})| = |E(F_{n_1}^*)| + | E(F_{n_2}^*)| + \dots + |E(F_{n_k}^*)|$ which implies, $(2m+1)^2 = 2n_1+2n_2+\dots +2n_k$ which is a contradiction since the L.H.S. is an odd number whereas the R.H.S. is an even number. Hence the result.     
\end{proof}

\begin{cor} [{\cite[Corollary 2.9]{vs15}}] \label{9.10}  {\rm
	For $m\in\mathbb{N}$, $G_{0,4m+1}$ does not admit $(a,d)$-CMSD into Fans with a handle for any $a,d\in\mathbb{N}$. \hfill $\Box$}
\end{cor}

\begin{theorem} [{\cite[Theorem 2.10]{vs15}}] \label{9.11}  {\rm 	For $n\in\mathbb{N}$,
	\begin{enumerate}
		\item[(i)] 	$G_{0,4n+1}$ can be decomposed into Fans with a handle and one $P_2$;
		\item[(ii)]	$G_{0,4n+2}$ admits $(2,8)$-CMSD into Fans with a handle; 
		\item[(iii)]	$G_{0,4n-1}$ admits $(4,8)$-CMSD into Fans with a handle; and
		\item[(iv)]	$G_{0,4n}$ admits  $(6,8)$-CMSD into Fans with a handle.
	\end{enumerate}}
\end{theorem}

\begin{proof} For $n\geq 3$, $F_{n-1}^*$, fan with a handle has $n+1$ vertices and $2(n-1)$ edges. Let $V(G_{0,n}) = \{v_0, v_1, v_2,\dots , v_n\}$ where $v_j$ is the vertex with integral sum label $j$ in $G_{0,n}$, $0 \leq j \leq n$. In the sum graph $G_n$, $|E(G_n)| = \frac{1}{2}\big(\frac{n(n-1)}{2} - \left\lfloor \frac{n}{2}\right\rfloor\big)$, $d(v_j) = n-1-j$ if $1\leq j\leq  \left\lfloor \frac{n+1}{2} \right\rfloor$ and $d(v_j) = n-j$ if   $\left\lfloor \frac{n+1}{2}\right\rfloor +1\leq j\leq n$ where $\left\lfloor x\right\rfloor $ is the floor of $x$. Now, consider decomposition of $G_{0,n}$ into Fans with a handle for different values of $n$ separately.
	
	In $G_{0,n}$, the subset $\{v_i v_j : i+j = n~\text{or}~ n-1,~ 0\leq i,j\leq n\}\cup \{v_0 v_i : i = 1,2,\dots,n-2\}$ of $E(G_{0,n})$ forms $F_{n-1}^*$, fan graph with cycle $(v_0  v_{n-1}  v_1  v_{n-2}  \dots  v_{\left\lfloor\frac{n}{2}\right\rfloor})$, pendant edge $v_0v_n$ attached at the apex vertex $v_0$ and $n-3$ concurrent edges, $v_0 v_js$ for $j = 1,2,\dots, \left\lfloor \frac{n}{2}\right\rfloor -1,  \left\lfloor \frac{n}{2}\right\rfloor +1,  \left\lfloor \frac{n}{2}\right\rfloor  +2,\dots,n-2$. Using the definition of integral sum labeling, $G_n-(\{v_n, v_{n-1}\}~\cup\{v_i v_j : i+j = n ~\text{or}~ n-1,$ $1\leq i,j \leq n-2\})$ = $G_n - \{n, n-1, [n], [n-1]\} = G_{n-2}$. Also using Theorem \ref{4.9}, $G_{n-2}-\{v_1, v_{n-2}\}$ is isomorphic to unlabeled graph $G_{n-4}$. Therefore $G_{0,n}-(\{v_0, v_n, v_{n-1}, v_{n-2}\}$ $\cup$ $\{v_i v_0 : i+j = n ~\text{or}~ n-1$, $1\leq i,j\leq n-2\})$ is isomorphic to unlabeled graph $G_{0,n-4}$. This also follows from Theorem \ref{2.11}. Relabeling the vertices $v_1, v_2,\dots, v_{n-3}$ in the resultant graph $G_{0,n} - (\{v_0, v_n, v_{n-1}, v_{n-2}\}$ $\cup$ $\{v_i v_j : i+j = n~\text{or}~ n-1, 1\leq  i,j\leq n-2\})$ as $v_0, v_1,\dots, v_{n-4}$ using the bijection $i\rightarrow i-1$ among the vertex labels and continuing the same technique of choosing the vertex subset $\{v_0, v_{n-4}, v_{n-5}, v_{n-6}\}$ and the relabeled edge subset $\{v_i v_j : i+j = n-4 ~\text{or}~ n- 5$ for $0\leq i,j\leq n-4\}$ $\cup$ $\{v_0 v_i : i = 1,2,\dots,n-6\}$ (in the relabeled graph) which form a fan $F_{n-5}^*$. The underlying graph of the subgraph $G_{0,n-4}-(\{v_0, v_{n-4}, v_{n-5}, v_{n-6}\}$ $\cup$ $\{v_i v_j : i+j = n-4 ~\text{or}~ n- 5, 1\leq  i,j\leq n-6\})$ of the relabeled graph $G_{0,n-4}$ is isomorphic to the underlying graph of $G_{0,n-8}$. Continue the above process. And to complete the proof, we consider the following four cases of $n$.
	
	\vspace{.1cm}
	\noindent
	\textbf{Case (i):}~ $n = 4m+1$, $m\in\mathbb{N}$. 
	
	In this case, $G_{0,n} = G_{0,4m+1} = F_{4m}^* \cup F_{4(m-1)}^* \cup F_{4(m-2)}^* \cup \dots \cup F_8^* \cup F_4^* \cup P_2(m$, $m+1)$ = $P_2(m,m+1) \cup  (\bigcup\limits_{j=0}^{m-1}F_{4(m-j)}^*)$ where $F_{4m}^*$, $F_{4(m-1)}^*$, $F_{4(m-2)}^*$, $\dots$, $F_8^*$, $F_4^*$, $P_2(m$, $m+1)$ are edge disjoint subgraphs of $G_{0,4m+1}$; here $F_{4m}^*$ is the Fan with the handle $(v_0, v_{4m+1})$, apex vertex $v_0$ and $P_{4m}$ = $v_{4m} v_1 v_{4m-1} v_2 \dots v_{2m+2} v_{2m-1} v_{2m+1} v_{2m}$; $|E(G_{0,4m+1})|$ = $4m+1 +|E(G_{4m+1})|$ = $4m+1 + (4m+1)(4m)/4 - (2m)/2$ = $(2m+1)^2$; $| E(P_2)|$  = 1 $<$ $|E(F_4^*)|$ = 8 $<$ $| E(F_8^*)|$  = 16 $< \dots <$ $| E(F_{n-5}^*)|$  = $2(n-5)$ $<$ $|E(F_{n-1}^*)|$ = $2(n-1)$ = $8m$; $|E(P_2)|+|E(F_4^*)|$ + $|E(F_8^*)|$ + $\dots$ + $| E(F_{n-1}^*)|$ = $1+8+16+\dots+8m$ = $4m^2$ + $4m+1$ = $(2m+1)^2$ and $v_j$ is the vertex with integral sum label $j$ in $G_{0,4m+1}$, $j\in [0,4m+1]$.
	Moreover, $P_2$ is a connected subgraph of $F_4^*$ which is a connected subgraph of $F_8^*$ which is a connected subgraph of $\dots$ which is a connected subgraph of $F_{n-5}^*$ which is a connected subgraph of $F_{n-1}^*$, without vertex labels.
	Thus, $G_{0,4m+1}$ can be decomposed into Fans with a handle and one $P_2$. To see the structure of these Fans with a handle and one $P_2$ in the above decomposition, see Figures 65 and 65.1 to 65.4.
	
	Decomposition of $G_{0,13}$ into Fans with a handle and one $P_2$ as given in the proof of the theorem is given in Figure 65 and its subgraph decomposition is shown separately in Figures 65.1 to 65.4. Here, $G_{0,13}$ = $F_{12}^* \cup F_{8}^* \cup F_{4}^* \cup P_2(3,4)$.

	\begin{center}
		\begin{tikzpicture}[scale = 0.6]	
		
\node (c13) at (21,-.25)  [circle,draw,scale=0.5]{13};
\node (c12) at (20.5,1)  [circle,draw,scale=0.5]{12};
\node (c11) at (20,3)  [circle,draw,scale=0.5] {11};
\node (c10) at (20.5,4.5)  [circle,draw,scale=0.5] {10};
\node (c9) at (21.5,6)  [circle,draw,scale=0.6] {9};	
\node (c8) at (23,7.25)  [circle,draw,scale=0.6] {8};
\node (c7) at (25.25,8)  [circle,draw,scale=0.6] {7};
\node (c6) at (27.75,7.25)  [circle,draw,scale=0.6] {6};
\node (c5) at (29.5,5.75)  [circle,draw,scale=0.6] {5};
\node (c4) at (30,4)  [circle,draw,scale=0.6]{4};
\node (c3) at (29.5,2)  [circle,draw,scale=0.6]{3};
\node (c2) at (28,0)  [circle,draw,scale=0.6] {2};
\node (c1) at (25.5,-1)  [circle,draw,scale=0.6] {1};
\node (c0) at (23,-1)  [circle,draw,scale=0.6] {0};
		
		\draw (c0) -- (c1);
		\draw (c0) -- (c2);
		\draw (c0) -- (c3);
		\draw (c0) -- (c4);
		\draw (c0) -- (c5);
		\draw (c0) -- (c6);
		\draw (c0) -- (c7);
		\draw (c0) -- (c8);
		\draw (c0) -- (c9);
		\draw (c0) -- (c10);
		\draw (c0) -- (c11);
		\draw (c0) -- (c12);
		\draw (c0) -- (c13);
		
		\draw (c12) -- (c1);
		\draw (c1) -- (c11);
		\draw (c11) -- (c2);
		\draw (c2) -- (c10);
		\draw (c10) -- (c3);
		\draw (c3) -- (c9);
		\draw (c9) -- (c4);
		\draw (c4) -- (c8);
		\draw (c8) -- (c5);
		\draw (c5) -- (c7);
		\draw (c7) -- (c6);
		
		{\color{red}
		\draw (c1) -- (c2);
\draw (c1) -- (c3);
\draw (c1) -- (c4);
\draw (c1) -- (c5);
\draw (c1) -- (c6);
\draw (c1) -- (c7);
\draw (c1) -- (c8);
\draw (c1) -- (c9);
\draw (c1) -- (c10);
			
			\draw (c9) -- (c2);
			\draw (c2) -- (c8);
			\draw (c8) -- (c3);
			\draw (c3) -- (c7);
			\draw (c7) -- (c4);
			\draw (c4) -- (c6);
			\draw (c6) -- (c5);}
		
		{\color{green}			
			\draw (c2) -- (c3);
			\draw (c2) -- (c4);
			\draw (c2) -- (c5);
			\draw (c2) -- (c6);
			\draw (c2) -- (c7);

			\draw (c6) -- (c3);			
			\draw (c3) -- (c5);
			\draw (c5) -- (c4);}
		
		{\color{orange}			
			\draw (c3) -- (c4);}
		
		\end{tikzpicture}				
		
		\vspace{.1cm}	
		\small{	Fig. 65. $G_{0,13}$ = $F_{12}^* \cup F_{8}^* \cup F_{4}^* \cup P_2(3,4)$}
\end{center}

\medskip
	\begin{center}
		\begin{tikzpicture}[scale = 0.5]	
				
\node (d13) at (21,-.25)  [circle,draw,scale=0.5]{13};
\node (d12) at (20.5,1)  [circle,draw,scale=0.5]{12};
\node (d11) at (20,3)  [circle,draw,scale=0.5] {11};
\node (d10) at (20.5,4.5)  [circle,draw,scale=0.5] {10};
\node (d9) at (21.5,6)  [circle,draw,scale=0.6] {9};	
\node (d8) at (23,7.25)  [circle,draw,scale=0.6] {8};
\node (d7) at (25.25,8)  [circle,draw,scale=0.6] {7};
\node (d6) at (27.75,7.5)  [circle,draw,scale=0.6] {6};
\node (d5) at (29.5,5.75)  [circle,draw,scale=0.6] {5};
\node (d4) at (30,4)  [circle,draw,scale=0.6]{4};
\node (d3) at (29.5,2)  [circle,draw,scale=0.6]{3};
\node (d2) at (28,0.25)  [circle,draw,scale=0.6] {2};
\node (d1) at (26,-1)  [circle,draw,scale=0.6] {1};
\node (d0) at (23,-1)  [circle,draw,scale=0.6] {0};

\draw (d0) -- (d1);
\draw (d0) -- (d2);
\draw (d0) -- (d3);
\draw (d0) -- (d4);
\draw (d0) -- (d5);
\draw (d0) -- (d6);
\draw (d0) -- (d7);
\draw (d0) -- (d8);
\draw (d0) -- (d9);
\draw (d0) -- (d10);
\draw (d0) -- (d11);
\draw (d0) -- (d12);
\draw (d0) -- (d13);

\draw (d12) -- (d1);
\draw (d1) -- (d11);
\draw (d11) -- (d2);
\draw (d2) -- (d10);
\draw (d10) -- (d3);
\draw (d3) -- (d9);
\draw (d9) -- (d4);
\draw (d4) -- (d8);
\draw (d8) -- (d5);
\draw (d5) -- (d7);
\draw (d7) -- (d6);

\node (c13) at (34,-.25)  [circle,draw,scale=0.5]{13};
\node (c12) at (33.5,1)  [circle,draw,scale=0.5]{12};
\node (c11) at (33,3)  [circle,draw,scale=0.5] {11};
\node (c10) at (33.5,4.5)  [circle,draw,scale=0.5] {10};
\node (c9) at (34.5,6)  [circle,draw,scale=0.6] {9};	
\node (c8) at (36,7.25)  [circle,draw,scale=0.6] {8};
\node (c7) at (38.25,8)  [circle,draw,scale=0.6] {7};
\node (c6) at (40.75,7.25)  [circle,draw,scale=0.6] {6};
\node (c5) at (42.5,5.75)  [circle,draw,scale=0.6] {5};
\node (c4) at (43,4)  [circle,draw,scale=0.6]{4};
\node (c3) at (42.5,2)  [circle,draw,scale=0.6]{3};
\node (c2) at (41,0)  [circle,draw,scale=0.6] {2};
\node (c1) at (38.5,-1)  [circle,draw,scale=0.6] {1};
\node (c0) at (36,-1)  [circle,draw,scale=0.6] {0};
		
		{\color{red}
			\draw (c1) -- (c2);
			\draw (c1) -- (c3);
			\draw (c1) -- (c4);
			\draw (c1) -- (c5);
			\draw (c1) -- (c6);
			\draw (c1) -- (c7);
			\draw (c1) -- (c8);
			\draw (c1) -- (c9);
			\draw (c1) -- (c10);
			
			\draw (c9) -- (c2);
			\draw (c2) -- (c8);
			\draw (c8) -- (c3);
			\draw (c3) -- (c7);
			\draw (c7) -- (c4);
			\draw (c4) -- (c6);
			\draw (c6) -- (c5);}
				
		\end{tikzpicture}				
		
		\vspace{.1cm}	
		\small{	Fig. 65.1. $F^*_{12}$ in $G_{0,13}$ \hspace{2cm} Fig. 65.2. $F^*_{8}$ in $G_{0,13}$ }
\end{center}

\medskip
	\begin{center}
		\begin{tikzpicture}[scale = 0.5]	
		
\node (d13) at (21,-.25)  [circle,draw,scale=0.5]{13};
\node (d12) at (20.5,1)  [circle,draw,scale=0.5]{12};
\node (d11) at (20,3)  [circle,draw,scale=0.5] {11};
\node (d10) at (20.5,4.5)  [circle,draw,scale=0.5] {10};
\node (d9) at (21.5,6)  [circle,draw,scale=0.6] {9};	
\node (d8) at (23,7.25)  [circle,draw,scale=0.6] {8};
\node (d7) at (25.25,8)  [circle,draw,scale=0.6] {7};
\node (d6) at (27,7.5)  [circle,draw,scale=0.6] {6};
\node (d5) at (29,6.75)  [circle,draw,scale=0.6] {5};
\node (d4) at (30,4.5)  [circle,draw,scale=0.6]{4};
\node (d3) at (29.75,2)  [circle,draw,scale=0.6]{3};
\node (d2) at (28,0)  [circle,draw,scale=0.6] {2};
\node (d1) at (25.5,-1)  [circle,draw,scale=0.6] {1};
\node (d0) at (23,-1)  [circle,draw,scale=0.6] {0};
		
		{\color{green}			
			\draw (d2) -- (d3);
			\draw (d2) -- (d4);
			\draw (d2) -- (d5);
			\draw (d2) -- (d6);
			\draw (d2) -- (d7);
			
			\draw (d6) -- (d3);			
			\draw (d3) -- (d5);
			\draw (d5) -- (d4);}
		

\node (c13) at (34,-.25)  [circle,draw,scale=0.5]{13};
\node (c12) at (33.5,1)  [circle,draw,scale=0.5]{12};
\node (c11) at (33,3)  [circle,draw,scale=0.5] {11};
\node (c10) at (33.5,4.5)  [circle,draw,scale=0.5] {10};
\node (c9) at (34.5,6)  [circle,draw,scale=0.6] {9};	
\node (c8) at (36,7.25)  [circle,draw,scale=0.6] {8};
\node (c7) at (38.25,8)  [circle,draw,scale=0.6] {7};
\node (c6) at (40.75,7.25)  [circle,draw,scale=0.6] {6};
\node (c5) at (42.5,5.75)  [circle,draw,scale=0.6] {5};
\node (c4) at (43,4)  [circle,draw,scale=0.6]{4};
\node (c3) at (42.5,2)  [circle,draw,scale=0.6]{3};
\node (c2) at (41,0)  [circle,draw,scale=0.6] {2};
\node (c1) at (38.5,-1)  [circle,draw,scale=0.6] {1};
\node (c0) at (36,-1)  [circle,draw,scale=0.6] {0};
		
		{\color{orange}			
			\draw (c3) -- (c4);}
		
		\end{tikzpicture}				
		
		\vspace{.1cm}	
		\small{	Fig. 65.3 $F^*_{4}$ in $G_{0,13}$  \hspace{2cm} Fig. 65.4 $P_2$ in $G_{0,13}$}
\end{center}

	\noindent
	\textbf{Case (ii):}~ $n$ = $4m+2$, $m\in\mathbb{N}$.
	 
	In this case, $G_{0,n}$ = $G_{0,4m+2}$ = $F_{4m+1}^*$ $\cup$ $F_{4(m-1)+1}^*$ $\cup$ $F_{4(m-2)+1}^*$ $\cup$ $\dots$ $\cup$ $F_5^*$ $\cup$ $F_1^*$ = $ \bigcup\limits_{j=0}^{m} F_{4(m-j)+1}^*$ where $F_{4m+1}^*$, $F_{4(m-1)+1}^*$, $F_{4(m-2)+1}^*$, $\dots$, $F_5^*$, $F_1^*$ = $P_3(m+1,m$, $m+2)$ are edge disjoint subgraphs of $G_{0,4m+2}$; $F_{4m+1}^*$ is the Fan with the handle $(v_0, v_{4m+2})$, apex vertex $v_0$ and $P_{4m+1}$ = $v_{4m+1}$ $v_1$ $v_{4m}$ $v_2$ $\dots$ $v_{2m-1}$ $v_{2m+2}$ $v_{2m}$ $v_{2m+1}$; $F_1^*$ = $P_3(v_{m+1},v_m,v_{m+2})$ is the path $v_{m+1} v_m v_{m+2}$ in $G_{0,4m+2}$ and $v_j$ is the vertex with integral sum label $j$ in $G_{0,4m+2}$, $j\in [0,4m+2]$. Also  $|E(G_{0,4m+2})|$ = $4m+2+|E(G_{4m+2})|$ = $4m+2 + (4m+2)(4m+1)/4 - (2m+1)/2$ = $4m^2+6m+2$ = $2(2m+1)(m+1)$, $|E(F_1^*)|$ = 2 $<$ $| E(F_5^*)|$ = 10 $<$ $|E(F_9^*)|$ = 18 $<\dots$ $<$ $|E(F_{n-5}^*)|$ = $2(n-5)$ $<$ $|E(F_{n-1}^*)|$ = $2(n-1)$ = $2(4m+1)$ = $8m+2$ and $2+10+18+ \dots +(2+8m)$ = $2(2m+1)(m+1)$. Thus, $G_{0,4m+2}$ admits $(2, 8)$-CMSD into Fans with a handle. Here $F_1^*$ is the trivial fan with a handle. 
	
\medskip	
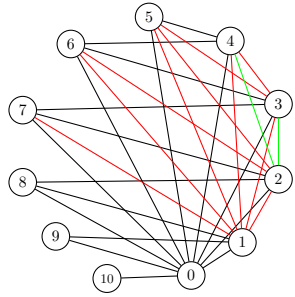
	$(2, 8)$-CMSD of $G_{0,10}$ into Fans with a handle is shown in Figure 66 and its subgraph decomposition is shown separately in Figures 66.1 to 66.3. Here, $G_{0,10}$ = $F_{9}^* \cup F_{5}^* \cup F_1^*$.
	
\begin{figure}	
	\centering
	\begin{tikzpicture}[scale = 0.5]
	
\node (c0) at (12,-2.1)  [circle,draw,scale=.6,  ] {0};	
\node (c1) at (14.6,-1) [circle,draw,scale=.6,  ] {1};
\node (c2) at (15.8,1.1) [circle,draw,scale=.6,  ] {2};
\node (c3) at (15.8,3.6) [circle,draw,scale=.6,  ] {3};
\node (c4) at (14.2,5.7) [circle,draw,scale=.6,  ]{4};
\node (c5) at (11.5,6.5)  [circle,draw,scale=.6,  ]{5};	
\node (c6) at (8.9,5.6) [circle,draw,scale=.6,  ] {6};
\node (c7) at (7.3,3.4) [circle,draw,scale=.6,  ] {7};
\node (c8) at (7.3,1) [circle,draw,scale=.6, ] {8};
\node (c9) at (8 ,-0.8) [circle,draw,scale=.6,  ] {9};
\node (c10) at (9.5,-2.2) [circle,draw,scale=.5,  ] {10};
	
\draw (c0) -- (c1);
\draw (c0) -- (c2);
\draw (c0) -- (c3);
\draw (c0) -- (c4);
\draw (c0) -- (c5);
\draw (c0) -- (c6);
\draw (c0) -- (c7);
\draw (c0) -- (c8);
\draw (c0) -- (c9);
\draw (c0) -- (c10);

\draw (c9) -- (c1);
\draw (c1) -- (c8);
\draw (c8) -- (c2);
\draw (c2) -- (c7);
\draw (c7) -- (c3);
\draw (c3) -- (c6);
\draw (c6) -- (c4);
\draw (c4) -- (c5);

{\color{red}
	\draw (c1) -- (c2);
	\draw (c1) -- (c3);
	\draw (c1) -- (c4);
	\draw (c1) -- (c5);
	\draw (c1) -- (c6);
	\draw (c1) -- (c7);
	
	\draw (c6) -- (c2);
	\draw (c2) -- (c5);
	\draw (c5) -- (c3);
	\draw (c3) -- (c4);}

{\color{green}			
	\draw (c2) -- (c3);
	\draw (c2) -- (c4);}

\node (d0) at (23.5,-2.1)  [circle,draw,scale=.6,  ] {0};	
\node (d1) at (26.5,-1) [circle,draw,scale=.6,  ] {1};
\node (d2) at (27.4,1.1) [circle,draw,scale=.6,  ] {2};
\node (d3) at (27.4,3.6) [circle,draw,scale=.6,  ] {3};
\node (d4) at (25.8,5.7) [circle,draw,scale=.6,  ]{4};
\node (d5) at (23.1,6.5)  [circle,draw,scale=.6,  ]{5};	
\node (d6) at (20.5,5.6) [circle,draw,scale=.6,  ] {6};
\node (d7) at (18.9,3.4) [circle,draw,scale=.6,  ] {7};
\node (d8) at (18.9,1) [circle,draw,scale=.6, ] {8};
\node (d9) at (20,-0.8) [circle,draw,scale=.6,  ] {9};
\node (d10) at (20.7,-2.2) [circle,draw,scale=.5,  ] {10};

\draw (d0) -- (d1);
\draw (d0) -- (d2);
\draw (d0) -- (d3);
\draw (d0) -- (d4);
\draw (d0) -- (d5);
\draw (d0) -- (d6);
\draw (d0) -- (d7);
\draw (d0) -- (d8);
\draw (d0) -- (d9);
\draw (d0) -- (d10);

\draw (d9) -- (d1);
\draw (d1) -- (d8);
\draw (d8) -- (d2);
\draw (d2) -- (d7);
\draw (d7) -- (d3);
\draw (d3) -- (d6);
\draw (d6) -- (d4);
\draw (d4) -- (d5);

\end{tikzpicture}
			
\medskip	
\small{ Fig. 66. $G_{0,10}$ = $F_{9}^* \cup F_{5}^* \cup P_3(4,2,3)$  \hspace{.5cm} Fig. 66.1. $F_{9}^*$ in $G_{0,10}$ \hspace{1cm}} \label{G_{0,10}}
\end{figure}	

\medskip	
 \begin{figure}	
 	\centering
 	\begin{tikzpicture}[scale = 0.5]
 	
 	\node (c0) at (12.9,-2.1)  [circle,draw,scale=.6,  ] {0};	
 	\node (c1) at (14.6,-1) [circle,draw,scale=.6,  ] {1};
 	\node (c2) at (15.8,1.1) [circle,draw,scale=.6,  ] {2};
 	\node (c3) at (15.8,3.6) [circle,draw,scale=.6,  ] {3};
 	\node (c4) at (14.2,5.7) [circle,draw,scale=.6,  ]{4};
 	\node (c5) at (11.5,6.5)  [circle,draw,scale=.6,  ]{5};	
 	\node (c6) at (8.9,5.6) [circle,draw,scale=.6,  ] {6};
 	\node (c7) at (7.3,3.4) [circle,draw,scale=.6,  ] {7};
 	\node (c8) at (7.3,1) [circle,draw,scale=.6, ] {8};
 	\node (c9) at (8.4 ,-0.8) [circle,draw,scale=.6,  ] {9};
 	\node (c10) at (10.1,-2.2) [circle,draw,scale=.5,  ] {10};
 	 	
 	{\color{red}
 		\draw (c1) -- (c2);
 		\draw (c1) -- (c3);
 		\draw (c1) -- (c4);
 		\draw (c1) -- (c5);
 		\draw (c1) -- (c6);
 		\draw (c1) -- (c7);
 		
 		\draw (c6) -- (c2);
 		\draw (c2) -- (c5);
 		\draw (c5) -- (c3);
 		\draw (c3) -- (c4);}
 		
 	\node (d0) at  (23.5,-2.1)[circle,draw,scale=.6,]{0};	
 	\node (d1) at (26.5,-1) [circle,draw,scale=.6,  ] {1};
 	\node (d2) at (27.4,1.1) [circle,draw,scale=.6,  ] {2};
 	\node (d3) at (27.4,3.6) [circle,draw,scale=.6,  ] {3};
 	\node (d4) at (25.8,5.7) [circle,draw,scale=.6,  ]{4};
 	\node (d5) at (23.1,6.5)  [circle,draw,scale=.6,  ]{5};	
 	\node (d6) at (20.5,5.6) [circle,draw,scale=.6,  ] {6};
 	\node (d7) at (18.9,3.4) [circle,draw,scale=.6,  ] {7};
 	\node (d8) at (18.9,1) [circle,draw,scale=.6, ] {8};
 	\node (d9) at (20,-0.8) [circle,draw,scale=.6,  ] {9};
 	\node (d10) at (20.7,-2.2) [circle,draw,scale=.5,  ] {10};
 	
	{\color{green}			
	\draw (d2) -- (d3);
	\draw (d2) -- (d4);}
 	
 	\end{tikzpicture}
 	
 \medskip	
 	\small{ 	Fig. 66.2 $F_{5}^*$ in $G_{0,10}$ \hspace{1.5cm} Fig. 66.3. $F_1^*$ in $G_{0,10}$} \label{G_{0,10}}
\end{figure}	

\medskip
\noindent
\textbf{Case (iii) :}~ $n = 4m-1$, $m\in\mathbb{N}$.\\
	In this case, $G_{0,n}$ = $G_{0,4m-1}$ = $F_{4m-2}^*~\cup F_{4m-6}^*~\cup   F_{4m-10}^*~\cup\dots~\cup F_6^*~\cup F_2^*$ = $\bigcup\limits_{j=0}^{m-1}F_{4(m-j)-2}^*$, $|E(G_{0,4m-1})|$ = $4m-1+|E(G_{4m-1})|$ = $4m-1 + (4m-1)(4m-2)/4 - (2m-1)/2$ = $4m^2$, $|E(F_2^*)|$ = 4 $<$ $|E(F_6^*)|$ = 12 $<$ $| E(F_{10}^*)|$ = 20 $< \dots <$ $|E(F_{4m-6}^*)|$ = $2(4m-6)$ $<$ $|E(F_{4m-2}^*)|$ = $2(4m-2)$ where $F_{4m-2}^*$, $F_{4(m-1)-2}^*$, $F_{4(m-2)-2}^*$, $\dots$, $F_6^*$, $F_2^*$ are edge disjoint subgraphs of $G_{0,4m-1}$; $F_{4m-2}^*$ is the Fan with handle $(v_0, v_{4m-1})$ and $v_j$ is the vertex with integral sum label $j$ in $G_{0,4m-1}$, $j\in [0,4m-1]$. This implies, $G_{0,4m-1}$ admits $(4, 8)$-CMSD into Fans with a handle. 
	
\medskip	
	$(4, 8)$-CMSD of $G_{0,11}$ into Fans with a handle is shown in Figure 67 and its subgraph decomposition is shown separately in Figures 67.1 to 67.3. Here, $G_{0,11}$ = $F_{10}^* \cup F_{6}^* \cup F_{2}^*$.

\begin{figure}	
	\centering
\begin{tikzpicture}[scale = 0.5]
		
\node (c0) at (11.5,-2.25) [circle,draw,scale=.6,  ] {0};
\node (c1) at (14.2,-1.8)  [circle,draw,scale=.6,  ] {1};
\node (c2) at (15.8,-0.3) [circle,draw,scale=.6,  ] {2};
\node (c3) at (16.5,2) [circle,draw,scale=.6,  ] {3};
\node (c4) at (16,4.3) [circle,draw,scale=.6,  ] {4};
\node (c5) at (14.3,6.2) [circle,draw,scale=.66,  ]{5};
\node (c6) at (11.5,6.6)  [circle,draw,scale=.6,  ]{6};
\node (c7) at (9.2,5.9) [circle,draw,scale=.6,  ] {7};
\node (c8) at (7.6,4.2) [circle,draw,scale=.6,  ] {8};
\node (c9) at (7,2) [circle,draw,scale=.6, ] {9};
\node (c10) at (7,-0.3) [circle,draw,scale=.5,  ] {10};
\node (c11) at (8,-2.25) [circle,draw,scale=.5,  ] {11};
		
		\draw (c0) -- (c1);
		\draw (c0) -- (c2);
		\draw (c0) -- (c3);
		\draw (c0) -- (c4);
		\draw (c0) -- (c5);
		\draw (c0) -- (c6);
		\draw (c0) -- (c7);
		\draw (c0) -- (c8);
		\draw (c0) -- (c9);
		\draw (c0) -- (c10);
		\draw (c0) -- (c11);
		
		\draw (c10) -- (c1);
		\draw (c1) -- (c9);
		\draw (c9) -- (c2);
		\draw (c2) -- (c8);
		\draw (c8) -- (c3);
		\draw (c3) -- (c7);
		\draw (c7) -- (c4);
		\draw (c4) -- (c6);
		\draw (c6) -- (c5);
		
		{\color{red}
			\draw (c1) -- (c2);
			\draw (c1) -- (c3);
			\draw (c1) -- (c4);
			\draw (c1) -- (c5);
			\draw (c1) -- (c6);
			\draw (c1) -- (c7);
			\draw (c1) -- (c8);
			
			\draw (c7) -- (c2);
			\draw (c2) -- (c6);
			\draw (c6) -- (c3);
			\draw (c3) -- (c5);
			\draw (c5) -- (c4);}
		
		{\color{green}			
			\draw (c2) -- (c3);
			\draw (c2) -- (c5);
			\draw (c2) -- (c4);
			\draw (c4) -- (c3); }
		
\node (d0) at (24,-2.25) [circle,draw,scale=.6,  ] {0};
\node (d1) at (26.5,-1.8)  [circle,draw,scale=.6,  ] {1};
\node (d2) at (27.8,-0.3) [circle,draw,scale=.6,  ] {2};
\node (d3) at (28.5,2) [circle,draw,scale=.6,  ] {3};
\node (d4) at (28,4.3) [circle,draw,scale=.6,  ] {4};
\node (d5) at (26.3,6.2) [circle,draw,scale=.6,  ]{5};
\node (d6) at (23.5,6.6)  [circle,draw,scale=.6,  ]{6};
\node (d7) at (21.2,5.9) [circle,draw,scale=.6,  ] {7};
\node (d8) at (19.6,4.2) [circle,draw,scale=.6,  ] {8};
\node (d9) at (19,2) [circle,draw,scale=.6, ] {9};
\node (d10) at (19,-0.3) [circle,draw,scale=.5,  ] {10};
\node (d11) at (20,-2.25) [circle,draw,scale=.5,  ] {11};
		
		\draw (d0) -- (d1);
		\draw (d0) -- (d2);
		\draw (d0) -- (d3);
		\draw (d0) -- (d4);
		\draw (d0) -- (d5);
		\draw (d0) -- (d6);
		\draw (d0) -- (d7);
		\draw (d0) -- (d8);
		\draw (d0) -- (d9);
		\draw (d0) -- (d10);
		\draw (d0) -- (d11);
		
		\draw (d10) -- (d1);
		\draw (d1) -- (d9);
		\draw (d9) -- (d2);
		\draw (d2) -- (d8);
		\draw (d8) -- (d3);
		\draw (d3) -- (d7);
		\draw (d7) -- (d4);
		\draw (d4) -- (d6);
		\draw (d6) -- (d5);
		
		\end{tikzpicture}
		
		\vspace{.1cm}	
		\small{\hspace{.5cm} Fig. 67. $G_{0,11}$ = $F_{10}^* \cup F_{6}^* \cup F_{2}^*$  \hspace{2cm} Fig. 67.1. $F_{10}^*$ in $G_{0,11}$} \label{G_{0,10}}
	\end{figure}	
\begin{figure}	
	\centering
	\begin{tikzpicture}[scale = 0.5]
	
	\node (c0) at (11.5,-2.25) [circle,draw,scale=.6,  ] {0};
	\node (c1) at (14.2,-1.8)  [circle,draw,scale=.6,  ] {1};
	\node (c2) at (16,-0.3) [circle,draw,scale=.6,  ] {2};
	\node (c3) at (16.5,2) [circle,draw,scale=.6,  ] {3};
	\node (c4) at (16,4.3) [circle,draw,scale=.6,  ] {4};
	\node (c5) at (14,6.4) [circle,draw,scale=.6,  ]{5};
	\node (c6) at (11.5,6.6)  [circle,draw,scale=.6,  ]{6};
	\node (c7) at (9.2,5.9) [circle,draw,scale=.6,  ] {7};
	\node (c8) at (7.6,4.2) [circle,draw,scale=.6,  ] {8};
	\node (c9) at (7,2) [circle,draw,scale=.6, ] {9};
	\node (c10) at (7.8,-0.3) [circle,draw,scale=.5,  ] {10};
	\node (c11) at (9,-2.25) [circle,draw,scale=.5,  ] {11};
		
	{\color{red}
		\draw (c1) -- (c2);
		\draw (c1) -- (c3);
		\draw (c1) -- (c4);
		\draw (c1) -- (c5);
		\draw (c1) -- (c6);
		\draw (c1) -- (c7);
		\draw (c1) -- (c8);
		
		\draw (c7) -- (c2);
		\draw (c2) -- (c6);
		\draw (c6) -- (c3);
		\draw (c3) -- (c5);
		\draw (c5) -- (c4);}
		
	\node (d0) at (24,-2.25) [circle,draw,scale=.6,  ] {0};
	\node (d1) at (26.2,-1.8)  [circle,draw,scale=.6,  ] {1};
	\node (d2) at (27.8,-0.3) [circle,draw,scale=.6,  ] {2};
	\node (d3) at (28.5,2) [circle,draw,scale=.6,  ] {3};
	\node (d4) at (28,4.3) [circle,draw,scale=.6,  ] {4};
	\node (d5) at (26.3,6.2) [circle,draw,scale=.6,  ]{5};
	\node (d6) at (23.5,6.6)  [circle,draw,scale=.6,  ]{6};
	\node (d7) at (21.2,5.9) [circle,draw,scale=.6,  ] {7};
	\node (d8) at (19.6,4.2) [circle,draw,scale=.6,  ] {8};
	\node (d9) at (19,2) [circle,draw,scale=.6, ] {9};
	\node (d10) at (19.8,-0.3) [circle,draw,scale=.5,  ] {10};
	\node (d11) at (21,-2.25) [circle,draw,scale=.5,  ] {11};
	
	{\color{green}			
	\draw (d2) -- (d3);
	\draw (d2) -- (d5);
	\draw (d2) -- (d4);
	\draw (d4) -- (d3);}
	\end{tikzpicture}
	
	\vspace{.1cm}	
	\small{ Fig. 67.2. $F_{6}^*$ in $G_{0,11}$\hspace{2cm} Fig. 67.3. $F_{2}^*$ in $G_{0,11}$} \label{G_{0,11}}
\end{figure}

\medskip	
	\noindent
	\textbf{Case (iv):}~$n = 4m$, $m\in\mathbb{N}$. \\
	In this case, $G_{0,n}$ = $G_{0,4m}$ = $F_{4m-1}^*~\cup F_{4m-5}^*~\cup F_{4m-9}^*~\cup \dots~\cup   F_7^* ~\cup F_3^*$ = $\bigcup\limits_{j=0}^{m-1} F_{4(m-j)-1}^*$, $|E(G_{0,4m})|$ = $4m +|E(G_{4m})|$ = $4m + 4m(4m-1)/4 - (2m)/2$ = $2m(2m+1)$, $|E(F_3^*)|$ = 6 $<$ $|E(F_7^*)|$  = 14 $<$ $|E(F_{11}^*)|$ = 22 $< \dots <$ $|E(F_{4m-5}^*)|$ = $2(4m-5)$ $<$ $|E(F_{4m-1}^*)$ = $2(4m-1)$ and $6+14+ \dots +(6+8(m-1))$ = $2m(2m+1)$ where $F_{4m-1}^*$, $F_{4(m-1)-1}^*$, $F_{4(m-2)-1}^*$, $\dots$, $F_7^*$, $F_3^*$ are edge disjoint subgraphs of $G_{0,4m}$, $F_{4m-1}^*$ is the Fan with the handle $(v_0, v_{4m})$ and $v_j$ is the vertex with integral sum label $j$ in $G_{0,4m}$, $j\in [0,4m]$. Thus, $G_{0,4m}$ admits $(6, 8)$-CMSD into Fans with a handle.
	
	\medskip	
		$(6,8)$-CMSD of $G_{0,12}$ into Fans with a handle is shown in Figure 68 and its subgraph decomposition is shown separately in Figures 68.1 to 68.3. Here, $G_{0,12}$ = $F_{11}^* \cup F_{7}^* \cup F_{3}^*$.
	
	Thus, in all the above four cases we could prove the result. 
	\end{proof} 
	
\begin{center}
	\begin{tikzpicture}[scale = 0.5]
	
	\node (c12) at (19.5,-1)  [circle,draw,scale=0.5]{12};
    \node (c11) at (19,1)  [circle,draw,scale=0.5] {11};
    \node (c10) at (19,3)  [circle,draw,scale=0.5] {10};
    \node (c9) at (19.5,4.5)  [circle,draw,scale=0.6] {9};
	\node (c8) at (20.75,5.75) [circle,draw,scale=0.6] {8};
	\node (c7) at (22.75,6.5)  [circle,draw,scale=0.6] {7};
	\node (c6) at (24.5,7)  [circle,draw,scale=0.6] {6};
	\node (c5) at (26.75,6.25)  [circle,draw,scale=0.6] {5};
	\node (c4) at (28.5,4.5)  [circle,draw,scale=0.6]{4};
	\node (c3) at (28.75,2)  [circle,draw,scale=0.6]{3};
	\node (c2) at (27.5,-0.5)  [circle,draw,scale=0.6] {2};
	\node (c1) at (25.5,-2)  [circle,draw,scale=0.6] {1};
	\node (c0) at (22.5,-2.25) [circle,draw,scale=0.6] {0};
	
	\draw (c0) -- (c1);
	\draw (c0) -- (c2);
	\draw (c0) -- (c3);
	\draw (c0) -- (c4);
	\draw (c0) -- (c5);
	\draw (c0) -- (c6);
	\draw (c0) -- (c7);
	\draw (c0) -- (c8);
	\draw (c0) -- (c9);
	\draw (c0) -- (c10);
	\draw (c0) -- (c11);
	\draw (c0) -- (c12);
	
	\draw (c11) -- (c1);
	\draw (c1) -- (c10);
	\draw (c10) -- (c2);
	\draw (c2) -- (c9);
	\draw (c9) -- (c3);
	\draw (c3) -- (c8);
	\draw (c8) -- (c4);
	\draw (c4) -- (c7);
	\draw (c7) -- (c5);
	\draw (c5) -- (c6);
	
	{\color{red}
		\draw (c1) -- (c2);
		\draw (c1) -- (c9);
		\draw (c1) -- (c8);
		\draw (c8) -- (c2);
		\draw (c2) -- (c7);
		\draw (c7) -- (c3);
		\draw (c3) -- (c6);
		\draw (c6) -- (c4);
		\draw (c4) -- (c5);
		
		\draw (c1) -- (c7);
		\draw (c1) -- (c6);
		\draw (c1) -- (c5);
		\draw (c1) -- (c4);
		\draw (c1) -- (c3);}
	
	{\color{green}			
		\draw (c2) -- (c3);
		\draw (c2) -- (c6);
		\draw (c2) -- (c5);
		\draw (c5) -- (c3);
		\draw (c3) -- (c4);

		\draw (c2) -- (c4);}

\node (d8) at (34.75,5.75)  [circle,draw,scale=0.6] {8};
\node (d7) at (36.75,6.5)  [circle,draw,scale=0.6] {7};
\node (d6) at (38.5,7)  [circle,draw,scale=0.6] {6};
\node (d5) at (40.5,6)  [circle,draw,scale=0.6] {5};
\node (d4) at (42,4.5)  [circle,draw,scale=0.6]{4};
\node (d3) at (42.5,2.5)  [circle,draw,scale=0.6]{3};
\node (d2) at (41.5,-0.5)  [circle,draw,scale=0.6] {2};
\node (d1) at (39.5,-2)  [circle,draw,scale=0.6] {1};
\node (d0) at (36.5,-2.25)  [circle,draw,scale=0.6] {0};
\node (d12) at (33.5,-1)  [circle,draw,scale=0.5]{12};
\node (d11) at (33,1)  [circle,draw,scale=0.5] {11};
\node (d10) at (33,3)  [circle,draw,scale=0.5] {10};
\node (d9) at (33.5,4.5)  [circle,draw,scale=0.6] {9};
	
	\draw (d0) -- (d1);
	\draw (d0) -- (d2);
	\draw (d0) -- (d3);
	\draw (d0) -- (d4);
	\draw (d0) -- (d5);
	\draw (d0) -- (d6);
	\draw (d0) -- (d7);
	\draw (d0) -- (d8);
	\draw (d0) -- (d9);
	\draw (d0) -- (d10);
	\draw (d0) -- (d11);
	\draw (d0) -- (d12);
	
	\draw (d11) -- (d1);
	\draw (d1) -- (d10);
	\draw (d10) -- (d2);
	\draw (d2) -- (d9);
	\draw (d9) -- (d3);
	\draw (d3) -- (d8);
	\draw (d8) -- (d4);
	\draw (d4) -- (d7);
	\draw (d7) -- (d5);
	\draw (d5) -- (d6);
	
	\end{tikzpicture}				
		
\vspace{.1cm}	
\small{ Fig. 68. $G_{0,12}$ = $F_{11}^* \cup F_{7}^* \cup F_{3}^*$ \hspace{2cm} Fig. 68.1. $F_{11}^*$ in $G_{0,12}$ \hspace{.5cm} } \label{G_{0,12}}
\end{center}

\vspace{.2cm}	
\begin{center}
	\begin{tikzpicture}[scale = 0.5]
	
	\node (c8) at (20.75,5.75)  [circle,draw,scale=0.6] {8};
	\node (c7) at (22.75,6.5)  [circle,draw,scale=0.6] {7};
	\node (c6) at (24.5,7)  [circle,draw,scale=0.6] {6};
	\node (c5) at (27,6.5)  [circle,draw,scale=0.6] {5};
	\node (c4) at (28.5,4.5)  [circle,draw,scale=0.6]{4};
	\node (c3) at (28.5,2)  [circle,draw,scale=0.6]{3};
	\node (c2) at (27.5,-.5)  [circle,draw,scale=0.6] {2};
	\node (c1) at (25.5,-2)  [circle,draw,scale=0.6] {1};
	\node (c0) at (22.5,-2.25)  [circle,draw,scale=0.6] {0};
	\node (c12) at (19.5,-1)  [circle,draw,scale=0.5]{12};
	\node (c11) at (19,1)  [circle,draw,scale=0.5] {11};
	\node (c10) at (19,3)  [circle,draw,scale=0.5] {10};
	\node (c9) at (19.5,4.5)  [circle,draw,scale=0.6] {9};
		
	{\color{red}
    	\draw (c1) -- (c2);
    	\draw (c1) -- (c9);
    	\draw (c1) -- (c8);
		\draw (c8) -- (c2);
		\draw (c2) -- (c7);
		\draw (c7) -- (c3);
		\draw (c3) -- (c6);
		\draw (c6) -- (c4);
		\draw (c4) -- (c5);
		
		\draw (c1) -- (c7);
		\draw (c1) -- (c6);
		\draw (c1) -- (c5);
		\draw (c1) -- (c4);
		\draw (c1) -- (c3);
	   	\draw (c1) -- (c2);}
		
	\node (d12) at (33.5,-1)  [circle,draw,scale=0.5]{12};
    \node (d11) at (33,1)  [circle,draw,scale=0.5] {11};
    \node (d10) at (33,3)  [circle,draw,scale=0.5] {10};
    \node (d9) at (33.5,4.5)  [circle,draw,scale=0.6] {9};
	\node (d8) at (34.75,5.75)  [circle,draw,scale=0.6]{8};
	\node (d7) at (36.75,6.5)  [circle,draw,scale=0.6] {7};
	\node (d6) at (38.5,7)  [circle,draw,scale=0.6] {6};
	\node (d5) at (40.5,6.5)  [circle,draw,scale=0.6] {5};
	\node (d4) at (42,5)  [circle,draw,scale=0.6]{4};
	\node (d3) at (42.5,1.5)  [circle,draw,scale=0.6]{3};
	\node (d2) at (41.5,-.5)  [circle,draw,scale=0.6] {2};
	\node (d1) at (39.5,-2)  [circle,draw,scale=0.6] {1};
	\node (d0) at (36.5,-2.25)  [circle,draw,scale=0.6]{0};
	
	{\color{green}			
	\draw (d2) -- (d3);
	\draw (d2) -- (d6);
	\draw (d2) -- (d5);
	\draw (d5) -- (d3);	
	\draw (d3) -- (d4);
	\draw (d4) -- (d2);}
	
	\end{tikzpicture}				
	
	\vspace{.1cm}	
	\small{ Fig. 68.2. $F_{7}^*$ in $G_{0,12}$ \hspace{3cm} Fig. 68.3. $F_{3}^*$ in $G_{0,12}$  } \label{G_{0,12}}	
\end{center}

\subsection{Conjecture on $(a, d)$-CMSD of $K_n$ into fan with a handle}

Vilfred and Suryakala \cite{vs15} proved that $n\equiv ~0,1~ mod (4)$ is a necessary condition for the existence of $(a, d)$-CMSD of $K_n$ into families of Fans with a handle. We improve this condition to Theorem \ref{9.12}. Based on these, a conjecture is proposed on the existence of $(a, d)$-CMSD of $K_n$ into families of Fans with a handle. 

\begin{theorem} [{\cite[Theorem 2.11]{vs15}}] \label{9.12a}  {\rm The necessary condition for the existence of $(a, d)$-CMSD of $K_n$ into families of Fans with a handle is $n(n-1)$ = $2k(2a + (k-1)d)$ and $n\equiv ~0,1~ mod (4)$, $a,d,k\in\mathbb{N}$. \hfill $\Box$}
\end{theorem}

\begin{theorem}  \label{9.12}  {\rm The necessary condition for the existence of $(a, d)$-CMSD of $K_n$ into families of Fans with a handle is $n(n-1)$ = $2k(2a + (k-1)d)$ and $n\equiv ~0,1~ mod (4)$, $a,d,k\in\mathbb{N}$.}
\end{theorem}

\begin{proof} Let $K_n$ admit $(a, d)$-CMSD into families of Fans with a handle for some $a,d\in\mathbb{N}$. Let $K_n$ = $F_{a}^*$ $\cup$ $F_{a+d}^*$ $\cup$ $F^*_{a+2d}$ $\cup$ $\dots$ $\cup$ $F_{a+(k-1)d}^*$ where $F_{a}^*$, $F_{a+d}^*$, $\dots$, $F_{a+(k-1)d}^*$ are edge disjoint Fans, each with a handle for some $a,d,k \in \mathbb{N}$. This implies, $|E(K_n)|$ = $|E(F_{a}^*)|$ + $|E(F_{a+d}^*)|$ + $\dots$ + $|E(F_{a+(k-1)d}^*)|$ which implies, $\binom{n}{2}$ = $a + a+d + a+2d$ + $\dots$ + $a + (k-1)d$ = $\frac{k}{2}(2a + (k-1)d)$. This implies, $n(n-1)$ = $2k(2a + (k-1)d)$. This also implies, $n \equiv 0,1~ mod (4)$ since $2k(2a + (k-1)d)$ is a multiple of 4 for any $a,d,k\in\mathbb{N}$. Hence we get the result.
\end{proof}

\medskip
The following conjecture is proposed in \cite{vs15} on the existence of $(a, d)$-CMSD of $K_n$ into families of Fans with a handle for some $a,d\in\mathbb{N}$.

\begin{conj} [{\cite[Conjecure]{vs15}}]  \label{9.13}  {\rm 	
		For every $n\in\mathbb{N}$, the following decompositions of $K_{4n}$ and $K_{4n+1}$  into Fan with a handle exist. 
		\begin{enumerate}
			\item [\rm (i)]  $K_{4n}$ = $F^*_{3n-1} \cup F^*_{3n-2} \cup F^*_{3n-3} \cup \dots \cup F^*_{n+1} \cup F^*_n$ which is a $(2n, 2)$-CMSD of $K_{4n}$ into Fan with a handle and 
			
			\item [\rm (ii)] $K_{4n+1}$ = $F^*_{3n} \cup F^*_{3n-1} \cup F^*_{3n-2} \cup \dots \cup F^*_{n+2} \cup F^*_{n+1}$ which is a $(2n+2, 2)$-CMSD of $K_{4n+1}$  into Fan with a handle. \hfill    $\Box$
	\end{enumerate} }
\end{conj}

The conjecture is verified on $K_{4n}$ and $K_{4n+1}$ for $n$ = 1, 2 in problem \ref{p8}.

\medskip
Note that in the above conjecture, $|E(F^*_{3n-1})| + |E(F^*_{3n-2})| + |E(F^*_{3n-3})| + \dots + |E(F^*_{n+1})| + |E(F^*_n)|$ = $2(3n-1) + 2(3n-2)+ 2(3n-3)+ \dots + 2(n+1) + 2n$ = $2(2n(n+3n-1)/2)$ = $4n(4n-1)/2$ = $|E(K_{4n})|$ and $|E(F^*_{3n})| + |E(F^*_{3n-1})|$ + $|E(F^*_{3n-2})| + \dots + |E(F^*_{n+2})| + |E(F^*_{n+1})|$ = $2(3n) + 2(3n-1)+ 2(3n-2)+ \dots$ + $2(n+2) + 2(n+1)$ = $2(2n(n+1+3n)/2)$ = $4n(4n+1)/2$ = $|E(K_{4n+1})|$.  

\begin{prm} [{\cite[Problem 2]{vs15}}]  \label{p8} {\rm Show that the following hold.
 \begin{enumerate}
 	\item [\rm (a1)] $K_4$ = $F^*_2 \cup F^*_1$ is a (2, 2)-CMSD of $K_4$ into Fan with a handle.

 	\item [\rm (a2)] $K_8$ = $F^*_5 \cup F^*_4 \cup F^*_3 \cup F^*_2$ is a $(4, 2)$-CMSD of $K_8$ into Fan with a handle.
 	 	
	\item [\rm (b1)] $K_5$ = $F^*_3 \cup F^*_2$ is a (4, 2)-CMSD of $K_8$ into Fan with a handle. 	

 	\item [\rm (b2)] $K_9$ = $F^*_6 \cup F^*_5 \cup F^*_4 \cup F^*_3$ is a (6, 2)-CMSD of $K_9$ into Fan with a handle.
 \end{enumerate}     }
\end{prm}

\noindent
{\bf Solution}\quad For $d$ = 2 and $n$ = 4, 5, 8 and 9,  $(a, d)$-CMSD of $K_n$ into families of Fans with a handle are given in Figures 65 to 68.4 as follows.

\medskip
\noindent
$(a1)$ $(2, 2)$-CMSD of $K_4$ into families of Fans with a handle are given in Figures 69, 69.1, 69.2 as follows.		

 	\begin{center}
 		\begin{tikzpicture}[scale = 0.9]
 		\node (a3) at (6,1.5)  [circle,draw,scale=0.6] {3};
 		\node (a2) at (8,1.5)  [circle,draw,scale=0.6]{2};
 		\node (a1) at (8,-.5)  [circle,draw,scale=0.6] {1};
 		\node (a0) at (6,-.5)  [circle,draw,scale=0.6] {0};
 		
 		\draw (a0) -- (a3);
 		\draw (a0) -- (a1);
 		
 		\draw (a0) -- (a2);
 		
 		\draw (a2) -- (a1);
 		
 {\color{red}	
	\draw (a3) -- (a1);
	\draw (a3) -- (a2);}

 		\node (b3) at (11,1.5)  [circle,draw,scale=0.6] {3};
 		\node (b2) at (13,1.5)  [circle,draw,scale=0.6]{2};
 		\node (b1) at (13,-.5)  [circle,draw,scale=0.6] {1};
 		\node (b0) at (11,-.5)  [circle,draw,scale=0.6] {0};
 		
 			\draw (b0) -- (b3);
            \draw (b0) -- (b1);

            \draw (b0) -- (b2);

            \draw (b2) -- (b1);
 		
 		\node (c3) at (16,1.5)  [circle,draw,scale=0.6] {3};
 		\node (c2) at (18,1.5)  [circle,draw,scale=0.6]{2};
 		\node (c1) at (18,-.5)  [circle,draw,scale=0.6] {1};
 		\node (c0) at (16,-.5)  [circle,draw,scale=0.6] {0};
 		
 		{\color{red}	
 			\draw (c3) -- (c1);
 			\draw (c3) -- (c2);}
\end{tikzpicture}
 		
 		\vspace{.1cm}
 		\noindent
 		\small{Fig. 69. $K_4$ = $F^*_2 \cup F^*_1$  \hspace{1cm} Fig. 69.1. $F^*_2$ in $K_4$ \hspace{1cm} Fig. 69.2. $F^*_1$ in $K_4$}
 	\end{center}

\medskip
\noindent
$(a2)$ $(4, 2)$-CMSD of $K_8$ into families of Fans with a handle are given in Figures 70, 70.1 to 70.4 as follows.		

\begin{center}
\begin{tikzpicture}[scale = 0.9]

\node (a7) at (13.5,0.5)  [circle,draw,scale=0.6]{7};
\node (a6) at (13.5,1.75)  [circle,draw,scale=0.6] {6};
\node (a5) at (15,2.75)  [circle,draw,scale=0.6] {5};
\node (a4) at (17,2.75)  [circle,draw,scale=0.6]{4};
\node (a3) at (18,1.75)  [circle,draw,scale=0.6]{3};
\node (a2) at (18,.5)  [circle,draw,scale=0.6] {2};
\node (a1) at (17,-.5)  [circle,draw,scale=0.6] {1};
\node (a0) at (15,-.5)  [circle,draw,scale=0.6] {0};

\draw (a0)  -- (a7);
\draw (a0)  -- (a6);
\draw (a6)  -- (a5);
\draw (a5)  -- (a4);
\draw (a4)  -- (a3);
\draw (a3)  -- (a2);
\draw (a2)  -- (a0);

\draw (a0)  -- (a5);
\draw (a0)  -- (a4);
\draw (a0)  -- (a3);

{\color{red}
\draw (a1)  -- (a0);
\draw (a1)  -- (a3);
\draw (a3)  -- (a5);
\draw (a5)  -- (a7);
\draw (a7)  -- (a4);
\draw (a4)  -- (a1);

\draw (a1)  -- (a5);
\draw (a1)  -- (a7);}

{\color{green}
\draw (a2)  -- (a5);
\draw (a2)  -- (a4);
\draw (a4)  -- (a6);
\draw (a6)  -- (a1);
\draw (a1)  -- (a2);

\draw (a6)  -- (a2);}

{\color{blue}
	\draw (a7)  -- (a2);
    \draw (a7)  -- (a6);
    \draw (a6)  -- (a3);
    \draw (a3)  -- (a7);}

\node (b7) at (21,.5)  [circle,draw,scale=0.6]{7};
\node (b6) at (21,1.75)  [circle,draw,scale=0.6] {6};
\node (b5) at (22,2.75)  [circle,draw,scale=0.6] {5};
\node (b4) at (24,2.75)  [circle,draw,scale=0.6]{4};
\node (b3) at (25,1.75)  [circle,draw,scale=0.6]{3};
\node (b2) at (25,.5)  [circle,draw,scale=0.6] {2};
\node (b1) at (24,-.5)  [circle,draw,scale=0.6] {1};
\node (b0) at (22,-.5)  [circle,draw,scale=0.6] {0};

\draw (b0)  -- (b7);
\draw (b0)  -- (b6);
\draw (b6)  -- (b5);
\draw (b5)  -- (b4);
\draw (b4)  -- (b3);
\draw (b3)  -- (b2);
\draw (b2)  -- (b0);

\draw (b0)  -- (b5);
\draw (b0)  -- (b4);
\draw (b0)  -- (b3);
\end{tikzpicture}

{\small Fig. 70. $K_8$ = $F^*_5 \cup F^*_4 \cup F^*_3 \cup F^*_2$  \hspace{2cm} Fig. 70.1. $F^*_5$ in $K_8$ \hspace{1cm}} 
\end{center}
	
\medskip
\begin{center}
	\begin{tikzpicture}[scale = 0.85]
	
\node (a7) at (9,0.5)  [circle,draw,scale=0.6]{7};
\node (a6) at (9,1.75)  [circle,draw,scale=0.6] {6};
\node (a5) at (10,2.75)  [circle,draw,scale=0.6] {5};
\node (a4) at (11.5,2.75)  [circle,draw,scale=0.6]{4};
\node (a3) at (12.5,1.75)  [circle,draw,scale=0.6]{3};
\node (a2) at (12.5,.5)  [circle,draw,scale=0.6] {2};
\node (a1) at (11.5,-.5)  [circle,draw,scale=0.6] {1};
\node (a0) at (10,-.5)  [circle,draw,scale=0.6] {0};
	
{\color{red}
	\draw (a1)  -- (a0);
	\draw (a1)  -- (a3);
	\draw (a3)  -- (a5);
	\draw (a5)  -- (a7);
	\draw (a7)  -- (a4);
	\draw (a4)  -- (a1);
	
	\draw (a1)  -- (a5);
	\draw (a1)  -- (a7);}
	
	\node (b7) at (14,.5)  [circle,draw,scale=0.6]{7};
	\node (b6) at (14,1.75)  [circle,draw,scale=0.6] {6};
	\node (b5) at (15,2.75)  [circle,draw,scale=0.6] {5};
	\node (b4) at (17,2.75)  [circle,draw,scale=0.6]{4};
	\node (b3) at (18,1.75)  [circle,draw,scale=0.6]{3};
	\node (b2) at (18,.5)  [circle,draw,scale=0.6] {2};
	\node (b1) at (17,-.5)  [circle,draw,scale=0.6] {1};
	\node (b0) at (15,-.5)  [circle,draw,scale=0.6] {0};
	
{\color{green}
\draw (b2)  -- (b5);
\draw (b2)  -- (b4);
\draw (b4)  -- (b6);
\draw (b6)  -- (b1);
\draw (b1)  -- (b2);

\draw (b6)  -- (b2);}

\node (c7) at (19,.5)  [circle,draw,scale=0.6]{7};
\node (c6) at (19,1.75)  [circle,draw,scale=0.6] {6};
\node (c5) at (20,2.75)  [circle,draw,scale=0.6] {5};
\node (c4) at (22,2.75)  [circle,draw,scale=0.6]{4};
\node (c3) at (23,1.75)  [circle,draw,scale=0.6]{3};
\node (c2) at (23,.5)  [circle,draw,scale=0.6] {2};
\node (c1) at (22,-.5)  [circle,draw,scale=0.6] {1};
\node (c0) at (20,-.5)  [circle,draw,scale=0.6] {0};

{\color{blue}
	\draw (c7)  -- (c2);
	\draw (c7)  -- (c6);
	\draw (c6)  -- (c3);
    \draw (c3)  -- (c7);}
\end{tikzpicture}

\vspace{.2cm}
{\small Fig. 70.2. $F^*_4$ in $K_8$ \hspace{1cm} Fig. 70.3. $F^*_3$ in $K_8$ \hspace{1cm}  Fig. 70.4. $F^*_2$ in $K_8$  } 
\end{center}

\medskip
\noindent
$(b1)$ $(4, 2)$-CMSD of $K_5$ into families of Fans with a handle are given in Figures 71, 71.1, 71.2 as follows.		

\begin{center}
	\begin{tikzpicture}[scale = 0.9]
	
	\node (a4) at (7.5,1.5)  [circle,draw,scale=0.6] {4};
	\node (a3) at (9,2.5)  [circle,draw,scale=0.6]{3};
	\node (a2) at (10.5,1.5)  [circle,draw,scale=0.6]{2};
	\node (a1) at (10,-.5)  [circle,draw,scale=0.6] {1};
	\node (a0) at (8,-.5)  [circle,draw,scale=0.6]{0};
	
	\draw (a0)  -- (a4);
	\draw (a0)  -- (a1);
	\draw (a0)  -- (a3);
	\draw (a3)  -- (a1);
	\draw (a1)  -- (a2);
	
	\draw (a2)  -- (a0);
	
	{\color{red} 
		\draw (a4)  -- (a1);
		\draw (a4)  -- (a2);
		\draw (a2)  -- (a3);
		\draw (a3)  -- (a4);}
	
	\node (u4) at (12.5,1.5)  [circle,draw,scale=0.6] {4};
	\node (u3) at (14,2.5)  [circle,draw,scale=0.6]{3};
	\node (u2) at (15.5,1.5)  [circle,draw,scale=0.6]{2};
	\node (u1) at (15,-.5)  [circle,draw,scale=0.6] {1};
	\node (u0) at (13,-.5)  [circle,draw,scale=0.6]{0};
	
	\draw (u0)  -- (u4);
	\draw (u0)  -- (u1);
	\draw (u0)  -- (u3);
	\draw (u3)  -- (u1);
	\draw (u1)  -- (u2);
	
	\draw (u2)  -- (u0);
	
	\node (v4) at (17.5,1.5)  [circle,draw,scale=0.6] {4};
	\node (v3) at (19,2.5)  [circle,draw,scale=0.6]{3};
	\node (v2) at (20.5,1.5)  [circle,draw,scale=0.6]{2};
	\node (v1) at (20,-.5)  [circle,draw,scale=0.6] {1};
	\node (v0) at (18,-.5)  [circle,draw,scale=0.6]{0};
	
	{\color{red} 
		\draw (v4)  -- (v1);
		\draw (v4)  -- (v2);
		\draw (v2)  -- (v3);
		\draw (v3)  -- (v4);}
	\end{tikzpicture}
	
	\vspace{.1cm}
	\noindent
	\small{Fig. 71. $K_5$ = $F^*_3 \cup F^*_2$  \hspace{1cm} Fig. 71.1. $F^*_3$ in $K_5$ \hspace{1cm} Fig. 71.2. $F^*_2$ in $K_5$}
\end{center}

\medskip
\noindent
$(b2)$ $(6, 2)$-CMSD of $K_9$ into families of Fans with a handle are given in Figures 72, 72.1 to 72.4 as follows.		

\begin{center}
	\begin{tikzpicture}[scale = 0.8]
	
	\node (a8) at (13,.5)  [circle,draw,scale=0.6]{8};
	\node (a7) at (12.5,2)  [circle,draw,scale=0.6] {7};
	\node (a6) at (13.5,3.5)  [circle,draw,scale=0.6] {6};
	\node (a5) at (15.5,4.5)  [circle,draw,scale=0.6] {5};
	\node (a4) at (17.5,3.5)  [circle,draw,scale=0.6]{4};
	\node (a3) at (18.5,2)  [circle,draw,scale=0.6]{3};
	\node (a2) at (18,.5)  [circle,draw,scale=0.6] {2};
	\node (a1) at (16.5,-.75)  [circle,draw,scale=0.6] {1};
	\node (a0) at (14.5,-.75)  [circle,draw,scale=0.6] {0};
	
\draw (a0)  -- (a8);
\draw (a0)  -- (a7);
\draw (a7)  -- (a6);
\draw (a6)  -- (a4);
\draw (a4)  -- (a3);
\draw (a3)  -- (a2);
\draw (a2)  -- (a1);
\draw (a1)  -- (a0);

\draw (a0)  -- (a6);
\draw (a0)  -- (a4);
\draw (a0)  -- (a3);
\draw (a0)  -- (a2);    
	
{\color{red}
\draw (a7)  -- (a1);
\draw (a7)  -- (a3);
\draw (a3)  -- (a5);
\draw (a5)  -- (a2);
\draw (a2)  -- (a4);
\draw (a4)  -- (a8);
\draw (a8)  -- (a7);

\draw (a7)  -- (a5);
\draw (a7)  -- (a4);	
\draw (a7)  -- (a2);}

{\color{green}
\draw (a8)  -- (a5);
\draw (a8)  -- (a1);
\draw (a1)  -- (a3);
\draw (a3)  -- (a6);
\draw (a6)  -- (a2);		
\draw (a2)  -- (a8);

\draw (a8)  -- (a3);
\draw (a8)  -- (a6); }

{\color{brown}
\draw (a5)  -- (a0);
\draw (a5)  -- (a4);
\draw (a4)  -- (a1);
\draw (a1)  -- (a6);
\draw (a6)  -- (a5);
\draw (a5)  -- (a1);}

	\node (b8) at (22,.5)  [circle,draw,scale=0.6]{8};
	\node (b7) at (21.5,2)  [circle,draw,scale=0.6] {7};
	\node (b6) at (22.5,3.5)  [circle,draw,scale=0.6] {6};
	\node (b5) at (24.5,4.5)  [circle,draw,scale=0.6] {5};
	\node (b4) at (26.5,3.5)  [circle,draw,scale=0.6]{4};
	\node (b3) at (27.5,2)  [circle,draw,scale=0.6]{3};
	\node (b2) at (27,.5)  [circle,draw,scale=0.6] {2};
	\node (b1) at (25.5,-.75)  [circle,draw,scale=0.6] {1};
	\node (b0) at (23.5,-.75)  [circle,draw,scale=0.6] {0};
	
	\draw (b0)  -- (b8);
    \draw (b0)  -- (b7);
    \draw (b7)  -- (b6);
    \draw (b6)  -- (b4);
    \draw (b4)  -- (b3);
    \draw (b3)  -- (b2);
    \draw (b2)  -- (b1);
    \draw (b1)  -- (b0);
    
    \draw (b0)  -- (b6);
    \draw (b0)  -- (b4);
    \draw (b0)  -- (b3);
    \draw (b0)  -- (b2);    
\end{tikzpicture}
	
	\vspace{.2cm}
	{\small Fig. 72. $K_9$ = $F^*_6 \cup F^*_5 \cup F^*_4 \cup F^*_3$  \hspace{2cm} Fig. 72.1. $F^*_6$ in $K_9$ \hspace{3cm}} 
\end{center}

\medskip
\begin{center}
	\begin{tikzpicture}[scale = 0.57]
	
	\node (a8) at (14.5,.5)  [circle,draw,scale=0.6]{8};
	\node (a7) at (14,2)  [circle,draw,scale=0.6] {7};
	\node (a6) at (15,3.5)  [circle,draw,scale=0.6] {6};
	\node (a5) at (17,4.5)  [circle,draw,scale=0.6] {5};
	\node (a4) at (19,3.5)  [circle,draw,scale=0.6]{4};
	\node (a3) at (20,2)  [circle,draw,scale=0.6]{3};
	\node (a2) at (19.5,.5)  [circle,draw,scale=0.6] {2};
	\node (a1) at (18,-.75)  [circle,draw,scale=0.6] {1};
	\node (a0) at (16,-.75)  [circle,draw,scale=0.6] {0};
	
{\color{red}
		\draw (a7)  -- (a1);
		\draw (a7)  -- (a3);
		\draw (a3)  -- (a5);
		\draw (a5)  -- (a2);
		\draw (a2)  -- (a4);
		\draw (a4)  -- (a8);
		\draw (a8)  -- (a7);
		
		\draw (a7)  -- (a5);
		\draw (a7)  -- (a4);	
		\draw (a7)  -- (a2); }
		
	\node (b8) at (22,.5)  [circle,draw,scale=0.6]{8};
	\node (b7) at (21.5,2)  [circle,draw,scale=0.6] {7};
	\node (b6) at (22.5,3.5)  [circle,draw,scale=0.6] {6};
	\node (b5) at (24.5,4.5)  [circle,draw,scale=0.6] {5};
	\node (b4) at (26.5,3.5)  [circle,draw,scale=0.6]{4};
	\node (b3) at (27.5,2)  [circle,draw,scale=0.6]{3};
	\node (b2) at (27,.5)  [circle,draw,scale=0.6] {2};
	\node (b1) at (25.5,-.75)  [circle,draw,scale=0.6] {1};
	\node (b0) at (23.5,-.75)  [circle,draw,scale=0.6] {0};
	
{\color{green}
	\draw (b8)  -- (b5);
	\draw (b8)  -- (b1);
	\draw (b1)  -- (b3);
	\draw (b3)  -- (b6);
	\draw (b6)  -- (b2);		
	\draw (b2)  -- (b8);
	
	\draw (b8)  -- (b3);
	\draw (b8)  -- (b6); }
	
\node (c8) at (29.5,.5)  [circle,draw,scale=0.6]{8};
\node (c7) at (29,2)  [circle,draw,scale=0.6] {7};
\node (c6) at (30,3.5)  [circle,draw,scale=0.6] {6};
\node (c5) at (32,4.5)  [circle,draw,scale=0.6] {5};
\node (c4) at (34,3.5)  [circle,draw,scale=0.6]{4};
\node (c3) at (35,2)  [circle,draw,scale=0.6]{3};
\node (c2) at (34.5,.5)  [circle,draw,scale=0.6] {2};
\node (c1) at (33,-.75)  [circle,draw,scale=0.6] {1};
\node (c0) at (31,-.75)  [circle,draw,scale=0.6] {0};

{\color{brown}
	\draw (c5)  -- (c0);
	\draw (c5)  -- (c4);
	\draw (c4)  -- (c1);
	\draw (c1)  -- (c6);
	\draw (c6)  -- (c5);
	\draw (c5)  -- (c1);
}
\end{tikzpicture}
	
	\vspace{.2cm}
	{\small Fig. 72.2. $F^*_5$ in $K_9$ \hspace{1cm} Fig. 72.3. $F^*_4$ in $K_9$ \hspace{1cm} Fig. 72.4. $F^*_3$ in $K_9$ } 
\end{center}

\section{Decomposition of $G_n$, $G^c_n$ and $K_n$ using sum and anti-$\mathbb{Z}$-sum labelings}

In this section, we present results on decomposition of graphs $G_n$, $G^c_n$ and $K_n$ using sum and anti-integral-sum labelings. 

The complement of an integral sum graph $G$ satisfies the property that $e$ = $uv$ is an edge of $G^c$ if and only if the sum of the labels on vertices $u$ and $v$ is not a vertex label in $G^c$. From this idea, anti-sum labeling and anti-integral sum labeling are defined \cite{vm13ca}. 

The concepts of sum, anti-sum, integral sum and anti-integral sum labeling are used to decompose graphs $G_n$, $G^c$ and $K_n$ and these results are presented in this section. For more on anti-sum and anti-integral sum labelings see section 2.3.

Throughout this section, we use the following notations: $K_1(i)$ represents a vertex with label $i$ and $P_2(i,j)$ represents an edge whose vertex labels are $i$ and $j$, $i,j\in\mathbb{Z}$. And for simplification, vertices of any sum or integral sum graph $G$ are named by their corresponding sum or integral sum label.

\subsection{Decomposition of $G^+_n$ using integral sum labeling}

In \cite{vm13ca}, a few recurrence relations on $G^+_n$ and ${(G_n^+)}^c$ are derived using sum and anti-sum labelings. These relations are used in decomposition of complete graph $K_n$ = $G^+_n$ $\cup$ ${(G_n^+)}^c$ and we present these results here.

\begin{dfn}	 {\rm \cite{h69}}\quad A graph $G$ is {\em decomposable} into the subgraphs $G_1$, $G_2$, $\dots$, $G_n$ of $G$, if no $G_i$ has isolated vertices and the edge set of $G$ can be partitioned into the subsets $E(G_1)$, $E(G_2)$, $\dots$, $E(G_n)$, $i$ = $1,2,\dots, n$. 
	
Graph $G$ is said to be {\em $H$-decomposable}, if $G_i$ $\cong$ $H$ for every $i$, $i$ = $1,2,\dots, n$. If $G$ is $H$-decomposable, then we say that $H$ divides $G$ and we write $H/G$. 
\end{dfn} 

\begin{theorem} [{\cite[Theorem 3.1]{vm13ca}}] \label{9.1.2} \quad {\rm For $n \geq 2$, 
		\noindent
		\begin{enumerate} 
			\item[{\rm (i)}] $G_{2n}^+$ $\cong$ $(G_{2n-1}^+$ $\cup$ $G^+( \{2n \}))$ $\bigcup$ \ $( {\cup}_{i=1}^{n-1}$ $(i,2n-i))$;\\
			\item[{\rm (ii)}] $G_{2n}^+$ $\cong$ $(G^+( [1,n-1] \cup [n+1,2n])$ $\cup$ $G^+( \{n \}))$  $\bigcup$ $( {\cup}_{i=1}^{n-1} (i,n))$ $\bigcup$ $( {\cup_{i=1}^{\left\lfloor \frac{n-1}{2}\right\rfloor}} (i,n-i))$;\\
			\item[{\rm (iii)}] $G_{2n+1}^+$ $\cong$ $(G_{2n}^+$ $\cup$ $G^+( \{2n+1 \}))$ $\bigcup$ $( {\cup}_{i=1}^{n}$ $(i,2n+1-i))$ and\\
			\item[{\rm (iv)}] $G_{2n+1}^+$ $\cong$ $(G^+( [1,n] \cup [n+2,2n+1])$ $\cup$ $G^+( \{n+1 \}))$
			
			\hspace{6cm} $\bigcup$ $( \cup_{i=1}^{n}$ $(i,n+1))$ $\bigcup$ $( \cup_{i=1}^{\left\lfloor \frac{n}{2}\right\rfloor}$ $(i,n+1-i))$.  
	\end{enumerate} } 
\end{theorem}

\begin{proof}\quad Proof is based on the Principle of Mathematical Induction on n. 
\begin{enumerate} 
   \item[{\rm (i)}] Partition the set of vertex labels of $G^+_{2n}$ into 2, one with $1,2,\dots,2n-1$ and the other with $2n$. Then, $G^+_{2n}$ can be considered as the union of disjoint sum graphs $G^+_{2n-1}$ = $G^+({1,2,...,2n-1})$ and $G^+(\{2n\})$ with additional edges (edges other than the edges of $G^+_{2n-1}$) $(1,2n-1)$, $(2,2n-2)$, $\dots$, $(n-1,2n-(n-1))$ = $(n-1,n+1)$. Now, result $(i)$ follows from the definition of sum graph.  
\item[{\rm (ii)}] In this case, consider two disjoint sum graphs, one with vertex labels $1,2,\dots,n-1,n+1$, $n+2,\dots,2n$ and the other with $n$. Construct $G^+_{2n}$ from the above two disjoint graphs and with additional edges (other than the edges of $G^+(\{1,2,...,n-1,n+1,n+2,...,2n\})$) $(1,n-1)$, $(2,n-2)$, $\dots$, $(\left\lfloor (n-1)/2 \right\rfloor, n- \left\lfloor (n-1)/2 \right\rfloor)$ and $(n,1),(n,2),\dots,(n,n-1)$. Then, result (ii) follows from the definition of sum graph. 
	\item[{\rm (iii)}] Similar to the proof as given to result (i).
	\item[{\rm (iv)}] Similar to the proof as given to result (ii).
\end{enumerate} 
\end{proof}

\begin{theorem} [{\cite[Theorem 3.2]{vm13ca}}] \label{9.1.3} \quad  {\rm For $n \geq 2$, 
		\noindent
		\begin{enumerate} 
			\item[{\rm (i)}] $G_{2n}^+$ $\cong$ $(G^+( [1, n-1] \cup [n+2, 2n])$ $\cup$ $G^+( \{n,n+1 \})$ 
			
			\hspace{1cm} $\bigcup$ $( \cup_{i=1}^{n-1}$ $((n,i)$ $\cup$ $(n+1,i)))$ $\bigcup$ $( \cup_{i=1}^{\left\lfloor \frac{n-1}{2}\right\rfloor}$ $(i,n-i))$ 
			
			\hspace{5cm} $\bigcup$ $( \cup_{i=1}^{\left\lfloor \frac{n}{2}\right\rfloor}$ $(i,n+1-i))$ and 
			\item[{\rm (ii)}] $G_{2n+1}^+$ $\cong$ $(G^+( [1, n-1] \cup [n+2, 2n+1])$ $\cup$ $({G^+( \{n,n+1 \}))}^c$ 
			
			\hspace{2cm} $\bigcup$ $( \cup_{i=1}^{n-1}((n,i)$ $\cup$ $(n+1,i)))$ $\bigcup$ $( \cup_{i=1}^{\left\lfloor \frac{n-1}{2}\right\rfloor}(i,n-i))$
			
			\hspace{6cm} $\bigcup$ $( \cup_{i=1}^{\left\lfloor \frac{n}{2}\right\rfloor}(i,n+1-i))$. 
	\end{enumerate}  } 
\end{theorem}

\begin{proof}\quad The proof is similar to the proof of Theorem \ref{9.1.2}, except the following vertex labeling. Here, ${(G^+(\{n,n+1\}))}^c$ = $P_2(n,n+1)$. i.e., ${(G^+(\{n,n+1\}))}^c$ is an edge whose vertex labels are $n$ and $n+1$. 
\end{proof}

\begin{theorem} [{\cite[Theorem 3.3]{vm13ca}}] \label{9.1.4} \quad {\rm For $n \geq 2$, 
		\noindent
		\begin{enumerate} 
			\item[{\rm (i)}] $G_{2n}^+$ $\cong$ $K_1(n)$ $\cup$ $K_1(n+1)$ 
			
			\hspace{1cm} 	$\bigcup$ $( \cup_{j=2}^{n} (K_1(n-j+1)$ $\cup K_1(n+j)$  
			
			\hfill $\cup$ $( \cup_{i=1}^{j-1}((n-j+1,n-i+1)$ $\cup (n-j+1,n+i)))))$ and
			\item[{\rm (ii)}] $G_{2n+1}^+$ $\cong$ $\bigcup_{j=1}^{n} (K_1(n-j+1)$ $\cup$ $K_1(n+j)$
			
			\hfill $\bigcup$ $( \cup_{i=1}^{j}((n-j+1,n-i+1)$ $\cup$ $(n-j+1,n+i))))$  $\bigcup$ $K_1(2n+1)$. 
	\end{enumerate} } 
\end{theorem}

\begin{proof}
	\begin{enumerate} 
		\item[{\rm (i)}]  In this case, start with the pair of vertices $K_1(n)$, $K_1(n+1)$ and after wards at each time consider one additional pair of vertices with these vertices in the order of $K_1(n-1)$, $K_1(n+2)$; $K_1(n-2)$, $K_1(n+3)$; $\dots$; $K_1(1)$, $K_1(2n)$ and corresponding additional edges so that at each time the resultant graph is an induced subgraph of $G^+{2n}$. 
		\item[{\rm (ii)}]   In this case, start with the pair of vertices $K_1(n)$, $K_1(n+1)$ and the edge $(n,n+1)$ and after wards at each time we consider one pair of additional vertices with these in the order of $K_1(n-1)$, $K_1(n+2)$; $K_1(n-2)$, $K_1(n+3)$; $\dots$; $K_1(1)$, $K_1(2n)$ and corresponding edges so that at each time the resultant graph is an induced subgraph of $G^+_{2n+1}$. And with $K_1(2n+1)$ at the end of the above process, we obtain graph $G^+_{2n+1}$. Hence the result.
	\end{enumerate} 
\end{proof}

\subsection{Decomposition of $G^c_n$ using anti-integral-sum labeling}

In \cite{vm13ca}, a few recurrence relations on ${G^+_n}^c$ similar to the recurrence relations derived in the previous subsection with respect to $G^+_n$ are obtained and from these recurrence relations, a few results on decomposition of graph ${G^+_n}^c$ are derived and presented here. 

\begin{theorem} [{\cite[Theorem 4.1]{vm13ca}}] \label{9.2.1} \quad {\rm For $n \geq 2$, 
		\noindent
		\begin{enumerate} 
			\item[{\rm (i)}] ${(G_{2n}^+)}^c$ $\cong$ $({(G_{2n-1}^+)}^c$ $*$ $G^+( \{2n \}))$ - $( \cup_{i=1}^{n-1} (i,2n-i));$
			\item[{\rm (ii)}] ${(G_{2n}^+)}^c$ $\cong$ ${((G^+( [1, n-1] \cup [n+1, 2n]))}^c$ $*$ $G^+( \{n \}))$  
			
			\hspace{3cm} - $((( \cup_{i=1}^{n-1} (i,n))$ $\bigcup$ $( \cup_{i=1}^{\left\lfloor \frac{n-1}{2}\right\rfloor}(i,n-i))));$
			\item[{\rm (iii)}] ${(G_{2n+1}^+)}^c \cong ({(G_{2n}^+)}^c * G^+( \{2n+1 \})) - ( \cup_{i=1}^{n} (i,2n+1-i))$ and 
			\item[{\rm (iv)}] ${(G_{2n+1}^+)}^c$  $\cong$ $({(G^+( [1, n] \cup [n+2, 2n+1]))}^c$ $*$ $G^+( \{n+1 \}))$  
			
			\hspace{2cm} - $(( \cup_{i=1}^{n} (i,n+1))$ $\bigcup$ $( \cup_{i=1}^{\left\lfloor \frac{n}{2}\right\rfloor}(i,n+1-i)))$. \hfill $\Box$
	\end{enumerate} } 
\end{theorem}

\begin{theorem} [{\cite[Theorem 4.2]{vm13ca}}] \label{9.2.2} \quad {\rm For $n \geq 2$, 
		\noindent
		\begin{enumerate} 
			\item[{\rm (i)}] ${(G_{2n}^+)}^c$  $\cong$ $({(G^+( [1, n-1] \cup [n+2, 2n]))}^c * (G^+( [n,n+1]))^c)$ 
			
			\hspace{1.5cm} -  $((\bigcup_{i=1}^{n-1}((n,i)$ $\cup$ $(n+1,i)))$ 
			
			\hspace{2cm} $\bigcup$ $(\cup_{i=1}^{\left\lfloor \frac{n-1}{2}\right\rfloor}(i,n-i))$
			$\bigcup$ $(\cup_{i=2}^{\left\lfloor \frac{n}{2}\right\rfloor}(i,n+1-i)));$ 
			
			\item[{\rm (ii)}] ${(G_{2n+1}^+)}^c$  $\cong$ $({(G^+( [1, n-1] \cup [n+2, 2n+1]))}^c$ $*$ $G^+([n,n+1]))$ 
			
			\hspace{1.5cm} - $((\bigcup_{i=1}^{n-1} ((n,i)$ $\cup$ $(n+1,i)))$ 
			
			\hspace{2cm}$\bigcup$ $(\cup_{i=1}^{\left\lfloor \frac{n-1}{2}\right\rfloor}(i,n-i))$ $\bigcup$ $(\cup_{i=2}^{\left\lfloor \frac{n}{2}\right\rfloor}(i,n+1-i)))$. \hfill $\Box$
	\end{enumerate} } 
\end{theorem}

\begin{theorem} [{\cite[Theorem 4.3]{vm13}}] \label{9.2.3} \quad {\rm  For $n \geq 2$, 
		\noindent
		\begin{enumerate} 
			\item[{\rm (i)}] ${(G_{2n}^+)}^c  \cong P_2(n,n+1)$ $\bigcup ( \cup_{j=2}^{n} (P_2(n-j+1,n+j) \bigcup ( \cup_{i=1}^{2j-2}(n-j+i+1,n+j))))$
			
			\hspace{.6cm} = $P_2(n,n+1)$ $\bigcup$ $(\cup_{j=2}^{n} (P_2(n-j+1,n+j)$ 
			
			\hfill $\bigcup$ $(\cup_{i=1}^{j}((n+j,n-i+1)$ $\cup$ $(n+j,n+j-i+1))))$ and 
			
			\item[{\rm (ii)}] ${(G_{2n+1}^+)}^c$  $\cong$ $K_1(n+1)$ $\bigcup$ $(\cup_{j=2}^{n+1}K_1(n+j)$ $\cup$ $K_1(n-j+2)$  
			
			\hspace{4cm} $\bigcup$ $(\cup_{i=1}^{2j-2}(n-j+i+1,n+j)))$
			
			\hspace{1cm} = $K_1(n+1)$ $\bigcup$ $(\cup_{j=2}^{n+1}K_1(n+j)$ $\cup$ $K_1(n-j+2)$  
			
			\hspace{2.5cm} $\bigcup$ $(\cup_{i=1}^{j}((n+j,n-i+1)$ $\cup$ $(n+j,n+j-i+1))))$. \hfill $\Box$
	\end{enumerate} } 
\end{theorem}

\begin{proof}\quad The result follows from the recurrence relations of Theorem \ref{9.2.2} and from the structure of graphs ${(G^+_{2n})}^c$ and ${(G^+_{2n+1})}^c$.   
\end{proof}

\subsection{Decomposition of $K_n$ using sum and anti-integral-sum labelings}

Using decomposition of graphs $G_n^+$ and ${(G_n^+)}^c$ as given in subsections 9.1 and 9.2 and the relation, $K_n$ = $G_n^+ \cup {(G_n^+)}^c$, the following decomposition results on $K_n$ are obtained in \cite{vm13ca} and presented here. 

\begin{theorem} [{\cite[Theorem 5.1]{vm13ca}}] \quad \label{d1} {\rm  For $n \geq 2$, 
		\begin{enumerate} 
			\item[{\rm (i)}] $K_{2n}$ $\cong$ $G_{2n-1}^+$ $\cup$ ${(G_{2n-1}^+)}^c$ $\cup$ $K_1(2n)$~ 
			
			\hspace{3cm} $\bigcup$ $(\cup_{i=1}^{n-1} ((n-i,2n)$ $\cup$ $(n+i,2n)))$ $\cup$ $(n,2n)$; \\
			\item[{\rm (ii)}] $K_{2n}$ $\cong$ $G^+( [1, n-1] \cup [n+1, 2n])$  $\cup$ ${(G^+( [1, n-1] \cup [n+1, 2n]))}^c$ $\cup$ $K_1(n)$
			
			\hspace{3cm}   $\bigcup$ $(\cup_{i=1}^{n-1}$ $((n,n-i)$  $\cup$ $(n,n+i)))$ $\cup$ $(n,2n)$;\\
			\item[{\rm (iii)}] $K_{2n+1}$ $\cong$ $G_{2n}^+$ $\cup$ ${(G_{2n}^+)}^c$ $\cup$ $K_1(2n+1)$ 
			
			\hspace{3cm} $\bigcup$ $(\cup_{i=1}^{n} ((n+1-i,2n+1)$ $\cup$ $(n+i,2n+1)))$ and \\
			\item[{\rm (iv)}] $K_{2n+1}$ $\cong$ $G^+( [1, n] \cup [n+2, 2n+1])$  $\cup$ ${(G^+ ([1, n] \cup [n+2, 2n+1]))}^c$ $\cup$ $K_1(n+1)$  
			
			\hfill $\bigcup$ $(\cup_{i=1}^{n}$ $(n+1-i,n+1)$ $\cup$ $(n+1,n+1+i)))$. 
	\end{enumerate} } 
\end{theorem}

\begin{proof}\quad Result follows from the relation $G_n \cup G^c_n$ = $K_n$ and using Theorems \ref{9.1.2} and \ref{9.2.1}. 
\end{proof}

\begin{theorem} [{\cite[Theorem 5.2]{vm13ca}}] \quad \label{d2} {\rm  For $n \geq 2$, 
		\begin{enumerate} 
			\item[{\rm (i)}] $K_{2n}$ $\cong$ $G^+([1, n-1]$ $\cup$ $[n+2, 2n])$
			
			\hspace{1cm}  $\bigcup$ ${(G^+( [1, n-1] \cup [n+2, 2n]))}^c$ $\bigcup$ $P_2(n,n+1)$ 
			
			\hfill $\bigcup$ $(\cup_{i=1}^{n-1}$ $((n,n-i)$ $\cup$ $(n,n+1+i)$ $\cup$ $(n+1,n-i)$ $\cup$ $(n+1,n+1+i)));$ 
			\item[{\rm (ii)}] $K_{2n+1}$ $\cong$ $G^+([1, n-1] \cup [n+2, 2n+1])$ 
			
			\hspace{1cm} $\cup$ ${(G^+([1, n-1] \cup [n+2, 2n+1]))}^c$ 
			
			\hspace{1cm} $\cup$ $P_2(n,n+1)$ $\cup$ $(n,n+2)$ $\cup$ $(n+1,n+2)$ 
			
			\hfill $\bigcup$ $(\cup_{i=1}^{n-1}$ $((n,i)$ $\cup$ $(n,2n+2-i)$ $\cup$ $(n+1,i)$ $\cup$ $(n+1,2n+2-i)))$.  
	\end{enumerate} } 
\end{theorem}
\begin{proof}\quad Result follows from Theorems \ref{9.1.3} and \ref{9.2.2}. 
\end{proof}

\begin{theorem} [{\cite[Theorem 5.3]{vm13ca}}] \quad \label{d3} {\rm  For $n \geq 2$, 
		\noindent
		\begin{enumerate} 
			\item[{\rm (i)}] $K_{2n} \cong P_2(n,n+1)$ 
			
			\hspace{1cm} $\bigcup$ $(\cup_{j=2}^{n} (P_2(n-j+1,n+j)$ 
			
			\hspace{1cm} $\bigcup$ $(\cup_{i=1}^{j-1}((n-j+1,n-i+1) \cup (n-j+1,n+i)$
			
			\hfill  $\cup$ $(n+j,n-i+1)$ $\cup$ $(n+j,n+j-i)))))$ and 
			\item[{\rm (ii)}] $K_{2n+1}$ $\cong$ $K_1(n+1)$ $\bigcup$ $(\cup_{j=1}^{n} (P_2(n+1-j,n+1+j)$  
			
			\hspace{4cm} $\bigcup$ $(\cup_{i=1}^{2j-1}(n+1-j+i,n+1+j))))$.
\end{enumerate} } 
\end{theorem}
\begin{proof}\quad Result follows from Theorems \ref{9.1.4} and \ref{9.2.3}. 
\end{proof}

\begin{rem} The authors feel that the above type of decomposition of $K_{2n+1}$ will help to settle open problems like Ringel's Tree Packing conjecture \cite{g26}. 
\end{rem}

\section{Edge coloring of integral sum graphs using integral sum labeling}

Vilfred \cite{vm13c} introduced {\em edge-sum class} and {\em edge-sum chromatic number} of an integral sum graph and studied their properties in \cite{js24, vm13c,vl22} similar to edge coloring of graphs. These are presented in this section.

\subsection{On edge-sum coloring, and edge-sum chromatic number of $G^+(S)$}  

An assignment of colors to the vertices of a graph such that adjacent vertices have the distinct colors is called a {\em proper coloring} of the graph. A {\em color class} is the set of all vertices with any one color and that is an independent set. The {\em chromatic number $\chi{(G)}$} of a graph $G$ is the minimum number of colors required to color the vertices of $G$ for which $G$ has a proper coloring \cite{h69}. Similarly, an assignment of colors to the edges of a graph $G$ in such a way that adjacent edges have distinct colors is termed as {\em proper edge coloring}. An {\em edge color class} is the set of all edges of $G$ with any single color and it is an independent set of edges. The {\em edge chromatic number $\chi^{'}{(G)}$} is defined as the minimum number of colors required to color the edges of $G$ for which $G$ has a proper edge coloring \cite{h69}. 
	
\begin{dfn}  \cite{vm13c}  Let $G^+(S)$ be an integral sum graph. In $G^+(S)$, the set of all edges, each with same edge-sum number, say $i$, is called the {\em edge-sum class} and is denoted by $E_i(G^+(S))$ or simply by $E_i$. 
\end{dfn}

\begin{dfn} \cite{vm13c} An integral sum graph $G^+(S)$ with vertex labeling function $f$ is said be an {\em edge-sum color graph} if any two edges of $G$ have the same color if and only if their edge-sum numbers are same. 
\end{dfn}

Thus, in an edge-sum color graph, edges of an edge-sum class have same color and edges of different edge-sum classes have different colors. 

\begin{dfn} \cite{vm13c} In an integral sum graph $G^+(S)$, the number of distinct non-empty edge-sum classes is called the {\em edge-sum chromatic  number} of $G^+(S)$ and is denoted by $\chi'_{\mathbb{Z}-sum}(G^+(S))$.
\end{dfn}

\subsection{Edge-sum classes and partition of edges in an $\mathbb{Z}$-sum graph} 

Given an integral sum graph $ G^+(S)$, we assume that vertex $u_j$ has label $j$. Every edge of an integral sum graph has an {\em induced edge-sum number} which is the sum of labels of its end vertices and adjacent edges have different edge-sum numbers. 

\begin{dfn} Let $G^+(S)$ be an integral sum graph. In $G^+(S)$, the set of all edges, each with same edge-sum number, say $i$, is called the {\em edge-sum class} and is denoted by $E_i(G^+(S))$ or simply by $E_i$. 
\end{dfn}

For example, consider $G^+(S)$ = $G_{-2,4}$. The edge-sum classes of $G^+(S)$ are

\begin{enumerate}
\item [\rm 1.] $E_0$ = $\{ (-1, 1), (-2, 2)\}$ = {\color{cyan}  $\{$edges with cyan color$\}$, }

\item [\rm 2.] $E_1$ = $\{ (0, 1), (-1, 2), (-2, 3) \}$ = {\color{black}  $\{$edges with black color$\}$, }

\item [\rm 3.] $E_2$ = $\{(0, 2), (-1, 3), (-2, 4) \}$ = {\color{red}  $\{$edges with red color$\}$, }

\item [\rm 4.] $E_3$ = $\{(0, 3), (1, 2), (-1, 4) \}$ = {\color{green}  $\{$edges with green color$\}$,  }

\item [\rm 5.] $E_4$ = $\{(0, 4), (1, 3)\}$ = {\color{orange}  $\{$edges with orange color$\}$, }

\item [\rm 6.] $E_{-1}$ = $\{(0, -1), (1, -2)\}$ = {\color{blue}  $\{$edges with blue color$\}$ } and

\item [\rm 7.] $E_{-2}$ = $\{(0, -2)\}$ = {\color{violet} $\{$edges with violet color$\}$.} 
\end{enumerate}
 
\medskip
\begin{center}
	\begin{tikzpicture} [scale = 1.8] 
		
		\node (a0) at (14,2)  [circle,draw,scale=0.6]{0};
		\node (a1) at (15,1.5)  [circle,draw,scale=0.6]{1};
		\node (a2) at (15.5,.5)  [circle,draw,scale=0.6] {2};
		\node (a3) at (15,-.25)  [circle,draw,scale=0.6] {3};
		\node (a4) at (14,-.25)  [circle,draw,scale=0.6] {4};
		\node (a5) at (13,.5)  [circle,draw,scale=0.5]{-2};
		\node (a6) at (13,1.5)  [circle,draw,scale=0.5] {-1};
		
		\draw (a0)[black]  -- (a1);
		\draw (a6)[black]  -- (a2);
		\draw (a5)[black]  -- (a3);	
		
		\draw (a0)[red, thick] -- (a2);
		\draw (a6)[red, thick]  -- (a3);
		\draw (a5)[red, thick]  -- (a4); 
		
		\draw (a0)[green, thick] -- (a3);
		\draw (a1)[green, thick] -- (a2);
		\draw (a6)[green, thick] -- (a4); 
		
		\draw (a0)[orange, thick]  -- (a4);
		\draw (a1)[orange, thick]  -- (a3);	
		
		\draw (a0)[violet, thick]  -- (a5); 
		
		\draw (a0)[blue, thick]  -- (a6);
		\draw (a5)[blue, thick]  -- (a1); 
		
		\draw (a5)[cyan, thick]  -- (a2);
		\draw (a6)[cyan, thick]  -- (a1);
		
		\node (b0) at (17,2)  [circle,draw,scale=0.6]{0};
		\node (b1) at (18,1.5)  [circle,draw,scale=0.6]{1};
		\node (b2) at (18.5,.5)  [circle,draw,scale=0.6] {2};
		\node (b3) at (18,-.25)  [circle,draw,scale=0.6] {3};
		\node (b4) at (17,-.25)  [circle,draw,scale=0.6] {4};
		\node (b5) at (16,.5)  [circle,draw,scale=0.5]{-2};
		\node (b6) at (16,1.5)  [circle,draw,scale=0.5] {-1};
		
		\draw (b0)[black]  -- (b1);
		\draw (b6)[black]  -- (b2);
		\draw (b5)[black]  -- (b3);	
		
		\draw (b0)[red, thick] -- (b2);
		\draw (b6)[red, thick]  -- (b3);
		\draw (b5)[red, thick]  -- (b4); 
		
		\draw (b1)[green, thick] -- (b2);
		\draw (b0)[green, thick] -- (b3);
		\draw (b6)[green, thick] -- (b4); 
		
		\draw (b1)[cyan, thick]  -- (b3);	
		\draw (b0)[cyan, thick]  -- (b4);
		\draw (b5)[cyan, thick]  -- (b2);
		
		\draw (b6)[blue, thick]  -- (b1);
		\draw (b0)[blue, thick]  -- (b5); 
		
		\draw (b0)[orange, thick]  -- (b6);
		\draw (b5)[orange, thick]  -- (b1); 		
	\end{tikzpicture}
	
   \vspace{.2cm}
	Fig. 73. $G_{-2, 4}$ with   \hspace{2cm} Fig. 74. $G_{-2, 4}$ with   
	
	\hspace{1cm}  $\chi'_{\mathbb{Z}-sum}(G_{-2, 4})$ = 7  \hspace{2cm}    $\chi'(G_{-2, 4})$ = 6 = $\Delta(G_{-2, 4})$ 
\end{center}

\medskip
Integral sum graph $G_{-2,4}$ with its edge-sum coloring and 
edge coloring are given in Figures 73 and 74, respectively. One set of edge coloring classes of $G_{-2, 4}$, correspondiing to Figure 74, is as follows. 

\medskip
$\{ \{(0,1), (-1,2), (-2,3)\}, \{(0,2), (-1,3), (-2,4)\}, \{(1,2), (0,3), (-1,4)\}$,

\hspace{4cm} $\{(1,3), (0,4), (-2,2)\}, \{(0,-2), (-1,1)\}, \{(0,-1), (1,-2)\} \}$.

\medskip
The reason for calling $E_i$ as the edge-sum class is based on Theorems \ref{p1} and \ref{p2}. Now, we prove that the non-empty sets of all edge-sum classes of an integral sum graph partition its edge set. 

\begin{theorem} [{\cite[Theorem 2.1]{vl22}}] \label{p1} {\rm \quad Let $G^+(S)$ be an integral sum graph of order $n$, $n \in \mathbb{N}$. Then, the following holds.
	\noindent
	\begin{enumerate}
		\item[{\rm (i)}] $E(G^{+}(S))$ = $\bigcup_{i \in S}{E_i}$. 

		\item[{\rm (ii)}] Two edge-sum classes are either equal or disjoint. That is, for every $i,j \in S$, either $E_i$ = $E_j$ or $E_i \cap E_j$ = $\emptyset$.  
	
		\item[{\rm (iii)}] For $i,j \in S$ if $E_i \neq \emptyset$, $E_j \neq \emptyset$ and $i \neq j$, then $E_i \neq E_j$. 
	
		\item[{\rm (iv)}] A non-empty edge-sum class is an independent set of edges of the integral sum graph. That is, no two elements (edges) of an edge-sum class have common vertex.
	
		\item[{\rm (v)}] Number of distinct non-empty edge-sum classes of $G^{+}(S)$ is less than or equal to $n$, the order of the graph.
	
		\item[{\rm (vi)}] Number of distinct non-empty edge-sum classes of $G^{+}(S)$ is equal to $n$ if and only if every vertex label occurs as the induced sum of at least one edge in $G^+(S).$
	\end{enumerate} }   
\end{theorem}

\begin{proof} 
	\begin{enumerate}
		\item[{\rm (i)}]  Every edge of $G^+(S)$ has an induced edge sum number, say $i$, which is the sum of the labels of its end vertices and $E_i$ is the set of all edges, each with induced edge sum number $i$ in $G^+(S)$, $i\in S$. Also, either $E_i$ = $\emptyset$ when there is no edge in $G^+(S)$ with edge sum number $i$ or $E_i \neq \emptyset$ when $G^+(S)$ has at least one edge with induced edge-sum number $i$, $i\in S$. Hence,  $E(G^{+}(S))$ = $\bigcup_{i \in S}{E_i}$.
	
		\item[{\rm (ii)}] The result is true when $E_i$ = $\emptyset$ or $E_j$ = $\emptyset$ or $E_i$ = $E_j$ = $\emptyset$, $i,j\in S$.
		
		Let $E_i, E_j \neq \emptyset$, $i,j\in S$. Then, $i$ = $j$ if and only if $E_i$ = $E_j$ since $E_i$ is the set of all edges of $G^+(S)$ with induced edge-sum number $i$, $i,j\in S$. Hence, we get the result $(ii)$.

		\item[{\rm (iii)}] 	Follows from result $(ii)$.
		
		\item[{\rm (iv)}] Let $e_i$ = $u v$ and $e_j$ = $u w$ be two adjacent edges in $G^+(S)$, $i \neq j$ and $i,j\in S$. By the definition of integral sum labeling, vertices $v$ and $w$ can not have same integral sum labeling in $G^+(S)$ and thereby induced edge-sum numbers of edges $u v$ and $u w$ are different. And so $e_i$ and $e_j$ belong to different edge-sum classes. Hence, we get result $(iv)$.
		
		\item[{\rm (v)}] Follows from the fact that every $i\in S$ need not be induced edge-sum number of an edge in $G^+(S)$ but every induced edge-sum number belongs to $S$.
		
		\item[{\rm (vi)}] 	This follows from the definition of integral sum labeling and from the fact that $n$ = $|S|$ = $|V(G^+)|$ and each edge belongs to an edge-sum class, say, $E_i$, $i \in S$.
	\end{enumerate} 
\end{proof}

\begin{theorem} [{\cite[Theorem 2.2]{vl22}}] \label{p2} \quad {\rm  The set of all non-empty edge-sum classes of an integral sum graph $G^+(S)$ partitions the set of all edges of the graph.} 
\end{theorem}

\begin{proof}\quad  In an integral sum graph $G^+(S)$, every edge has an induced edge-sum number which belongs to $S$ and also, $E(G^+(S))$ = $\sum_{i\in S}{E_i}$ using Theorem \ref{p1}(i). By (ii) and (iii) of Theorem \ref{p1}, any two non-empty edge-sum classes in an integral sum graph are either equal or disjoint. Thus, we get the result. 
\end{proof}

We have seen that the set of all non-empty edge-sum classes of an integral sum graph partition the edge set of the graph. For a given integral sum graph $G^{+}(S)$, the set of all edge-sum classes is unique. Property $(iv)$ helps us to consider an integral sum graph as {\em edge sum color graph} by applying same color to all edges in an edge-sum class and different colors to different edge-sum classes. 

\subsection{Chromatic and edge chromatic numbers of $G_{-r, n}$ }

In this subsection, we present results on the chromatic and the edge chromatic numbers of $G_{-r, n}$,  $r \in \mathbb{N}_0$ and $n \in \mathbb{N}$. In \cite{vm13c} (also see \cite{vl22}), Vilfred proposed Theorem \ref{thm Grs 2 chrom} (ii) as a conjecture and it was completely settled in \cite{js24}.  

\begin{theorem} [{\cite[Theorem 3.5]{vm13c}}] \label{thm Grs 1 chrom}\quad {\rm  
		Let $n \in \mathbb{N}$.
		\begin{enumerate}
			\item [\rm (i)] $\chi^{'}\left(G_{-1,1} \right)$ = $3$.
			
			\item [\rm (ii)] $\chi^{'}\left(G_{-1,n} \right)$ = $n+1$ for $n\geq 2$. 
			
			\item [\rm (iii)] $\chi^{'}\left(G_{-n,n} \right)$ = $2n$ for $ 2\leq n \leq 6$.
	\end{enumerate}}
\end{theorem}

\begin{proof} See the proof of Theorem \ref{thm Grs 2 chrom}(ii). 
\end{proof}

In \cite{vm13c} (also see \cite{vl22}), Vilfred proposed Theorem \ref{thm Grs 2 chrom} (ii) as a conjecture. The conjecture was settled completely in \cite{js24} and the proof is presented here.  

\begin{theorem} [{\cite[Theorem 4.4]{js24}}] \label{thm Grs 2 chrom} \quad {\rm 	Let $r,n \in \mathbb{N}$.
		\begin{enumerate}
			\item [\rm (i)] the chromatic number of $G_{-r,n}$ is $\chi\left(G_{-r,n} \right)$ = $ 1 + \left \lceil{\frac{r}{2}} \right\rceil + \left \lceil{\frac{n}{2}} \right\rceil$.
			
			\item [\rm (ii)] the edge chromatic number of $G_{-r,n}$ is $\chi^{'}\left(G_{-r,n} \right)$ = $r+n$.
	\end{enumerate}}
\end{theorem}	

\begin{proof}\quad For $-r \leq i \leq n$ and $r,n \in \mathbb{N}$, let $f : V(G_{- r,n}) \rightarrow \mathbb{Z}$ defined by $f(v_i) = i$ be the integral sum labeling of the graph $G_{-r,n}$.	

\medskip
\begin{enumerate}
	\item [\rm (i)] For $r,n \in \mathbb{N}$, integral sum graph $G_{-r,n}$ is perfect by Theorem \ref{thm Grs pf}. Therefore  $\chi\left(G_{-r,n} \right)$ = $\omega(G_{-r,n} )$ = $ 1 + \left \lceil{\frac{r}{2}} \right\rceil + \left \lceil{\frac{n}{2}} \right\rceil$ follows from Lemma \ref{resGrs}. Hence the result is true in this case. Figure 43 shows the vertex coloring and clique in $G_{-5,7}$.

\medskip
	\item [\rm (ii)] For $r,n \in \mathbb{N}$, the maximum degree of integral sum graph $G_{-r,n}$ is $\Delta (G_{-r,n})$ = $r+n$. And using Theorem \ref{1.1}, we get, $\Delta (G_{-r,n})$ $\leq$ ${\chi}^{'}(G_{-r,n}) \leq \Delta (G_{-r,n})+1$. Here we present a proper edge coloring of the integral sum graph $G_{-r,n}$ with $r+n$ colors so ${\chi}^{'}(G_{-r,n}) = r + n$ for $r,n \in \mathbb{N}$.
	
	\medskip
	\indent Let $c_k$ denote $k^{th}$ color assigned to an edge and $C_k$ denote the color class of edges, each with color $c_k$ in $G_{-r,n}$, $r,n \in \mathbb{N}$. 
	For $r,n \in \mathbb{N}$, color the edges of $G_{-r,n}$ as follows.
	
\medskip
	$v_0 v_j$ $\mapsto$ $c_j$, ~$1 \leq j \leq n$; ~ i.e., edge $v_0v_j$ is taking color $c_j$, $1 \leq j \leq n$;
	
\medskip
	$v_0 v_{-i}$ $\mapsto$ $c_{i+n}$, ~$1 \leq i \leq r$;
	
\medskip
	$v_{-i} v_j$ $\mapsto$ $c_{i+j}$, ~ $1 \leq i \leq r$ and $1 \leq j \leq n-1$;
	
\medskip
	$v_{-i} v_n$ $\mapsto$ $c_{2i+n}$,~ $1 \leq i \leq \left\lfloor \frac{r}{2}\right\rfloor$;
	
\medskip
	$v_{-i} v_n$ $\mapsto$ $c_{i-\left\lfloor \frac{r}{2}\right\rfloor}$,~ $\left\lfloor \frac{r}{2}\right\rfloor < i \leq r$, $ n > \left\lceil \frac{r}{2}\right\rceil$;
	
\medskip
	$v_{-i} v_n$ $\mapsto$ $c_{i-\left\lfloor \frac{r}{2}\right\rfloor}$, ~$\left\lfloor \frac{r}{2}\right\rfloor < i \leq \left\lfloor \frac{r}{2}\right\rfloor + n-1$, $ n \leq \left\lceil \frac{r}{2}\right\rceil$;
	
\medskip
	$v_{-i} v_n$ $\mapsto$ $c_{2\left( i-\left\lfloor \frac{r}{2}\right\rfloor\right) -n+1}$, ~$\left\lfloor \frac{r}{2}\right\rfloor + n \leq i \leq r$, $ n \leq \left\lceil \frac{r}{2}\right\rceil$;
	
\medskip
	$v_i v_j$ $\mapsto$ $c_{i+j+r}$,~ $1 \leq i,j,i+j \leq n$ and $i < j$ and
	
\medskip
	$v_{-i} v_{-j}$ $\mapsto$ $c_{i+j+n}$,~ $1 \leq i,j,i+j \leq r$ and $i < j$.
	
\medskip
	\indent It is clear from the above edge coloring that colors $c_1$ to $c_{n+r}$ are assigned to the edges of $G_{-r,n}$ and no more edge colors are required. And also colors of edges at each vertex of $G_{-r,n}$ are all distinct by the following. We have $G_{-r,n}$ $\cong$ $K_1* (G_{-r}*G_n)$ and $G_{-r}*G_n$ $\cong$ $G_{-r} \cup G_n \cup K_{r, n}$, $r,n\in\mathbb{N}$.
	
	In $G_{-r,n}$, colors $c_1$, $c_2$, \dots, $c_{n+r}$ are assigned to the $n+r$ edges incident at the vertex $v_0$.
	
\medskip
	\indent In $K_{-r,n}$,  we get the following possible colors taken by its edges incident at each of its vertices under the coloring already assigned to the edges of $G_{-r,n}$.  
	
\medskip
\begin{enumerate}
 \item [\rm	(a)] For $1 \leq j \leq n-1$, distinct colors $c_{j+1}$, $c_{j+2}$, \dots, $c_{j+r}$ are assigned to the $r$ edges incident at $v_j$. And these colors are also different from color $c_j$ of the edge $v_0 v_j$ at $ v_0 $.  
	
\item [\rm	(b)] For $n> \left\lceil \frac{r}{2}\right\rceil $,  distinct colors  $c_{2+n}$, $c_{4+n}$, ..., $c_{2\left( \left\lfloor \frac{r}{2}\right\rfloor-1\right) +n}$, $c_{2\left\lfloor\frac{r}{2}\right\rfloor+n}$, $c_1, c_2, \dots, $ $c_{r-\left\lfloor \frac{r}{2}\right\rfloor}$ are assigned to the $r$ edges incident at $v_n$ and these colors are different from color $c_n$ of the edge $v_0 v_n$ at $v_n$. Clearly, $2\left\lfloor \frac{r}{2}\right\rfloor+n$ $\leq$ $r+n$ and thus the colors assigned are from the $r+n$ colors only.
	
\item [\rm	(c)] For $n\leq \left\lceil \frac{r}{2}\right\rceil $, distinct colors $c_{2+n}$, $c_{4+n}$, ...., $c_{2\left( \left\lfloor \frac{r}{2}\right\rfloor-1\right) +n}$, $c_{2\left\lfloor\frac{r}{2}\right\rfloor+n}$, $c_1, c_2, \dots, $ $c_{\left\lfloor \frac{r}{2}\right\rfloor+n-1}$, $c_{n+1}$, $c_{n+3}$, $\dots$, $c_{2\left\lceil \frac{r}{2}\right\rceil -n+1}$ are assigned to the $r$ edges incident at $v_n$ and these colors are different from color $c_n$ of the edge $v_0 v_n$ at $v_n$. Clearly, $2\left\lfloor \frac{r}{2}\right\rfloor+n$, $2\left\lceil \frac{r}{2}\right\rceil -n+1$ $\leq$ $r+n$ and atmost $r+n$ colors are assigned.
	
\item [\rm	(d)] For $1 \leq i \leq \left\lfloor \frac{r}{2}\right\rfloor$, distinct colors $c_{i+1}$, $c_{i+2}$, \dots, $c_{i+n-1}$, $c_{2i+n}$ are assigned to the $n$ edges incident at $v_{-i}$. And all these colors are different from color $c_{i+n}$ of the edge $v_0 v_{-i}$ at $v_{-i}$.
	
\item [\rm	(e)] For $n> \left\lceil \frac{r}{2}\right\rceil $ and $\left\lfloor \frac{r}{2}\right\rfloor < i \leq r$, distinct colors $c_{i+1}$, $c_{i+2}$, \dots, $c_{i+n-1}$, $c_{i-\left\lfloor \frac{r}{2}\right\rfloor}$ are assigned to the $n$ edges incident at $v_{-i}$ and all these colors are different from color $c_{i+n}$ of the edge $v_0 v_{-i}$ at $v_{-i}$.
	
\item [\rm	(f)] For $n\leq \left\lceil \frac{r}{2}\right\rceil $ and $\left\lfloor \frac{r}{2}\right\rfloor < i \leq \left\lfloor \frac{r}{2}\right\rfloor +n-1$, distinct colors $c_{i+1}$, $c_{i+2}$, \dots, $c_{i+n-1}$, $c_{i-\left\lfloor \frac{r}{2}\right\rfloor}$ are assigned to the $n$ edges incident at $v_{-i}$ and all these colors are different from color $c_{i+n}$ of the edge $v_0 v_{-i}$ at $v_{-i}$. 
	
\item [\rm	(g)] For $n\leq \left\lceil \frac{r}{2}\right\rceil $ and $\left\lfloor \frac{r}{2}\right\rfloor +n \leq i \leq r$, distinct colors $c_{i+1}$, $c_{i+2}$, \dots, $c_{i+n-1}$, $c_{2\left( i-\left\lfloor \frac{r}{2}\right\rfloor\right) -n+1}$ are assigned to the $n$ edges incident at $v_{-i}$ and all these colors are different from color $c_{i+n}$ of the edge $v_0 v_{-i}$ at $v_{-i}$.
	\end{enumerate}
\end{enumerate}	
	Thus, in the above edge coloring of $G_{-r,n}$, edges of $K_{-r,n}$ take colors from the $r+n$ colors and colors of edges incident at each vertex of $K_{-r,n}$ are all distinct. 
	
	\indent In $G_{n}$, for $1 \leq i,j \leq n$, $i \neq j$ and $3 \leq i+j \leq n$, color of edge $v_i v_j$ is $c_{i+j+r}$ and thereby edges of $G_n$ are taking colors only from the colors $c_{3+r}$, $c_{4+r}$, \dots, $c_{n+r}$; edge colors $c_{i+j+r}$ at $v_i$ are all distinct and also different from $c_i$ of the edge $v_0 v_{i}$ at $v_{i}$ as well as that of other possible colors of edges which incident at $v_{i}$ in $G_{-r,n}$ where $i \neq j$ and $3 \leq i+j \leq n$. The same arguments also hold when $i$ and $j$ are interchanged. 
	
	\indent In $G_{-r}$, for $1 \leq i,j \leq r$, $i \neq j$ and $3 \leq i+j \leq r$, color of edge $v_{-i} v_{-j}$ is $c_{i+j+n}$ and thereby edges of $G_{-r}$ are taking colors only from the colors $c_{3+n}$, $c_{4+n}$, \dots, $c_{r+n}$; edge colors $c_{i+j+n}$ at $v_{-i}$ are all distinct and also different from $c_{i+n}$ of the edge $v_0 v_{-i}$ at $v_{-i}$ as well as that of other possible colors of edges which incident at $v_{-i}$ in $G_{-r,n}$ where $i \neq j$ and $3 \leq i+j \leq r$. The same arguments also hold when ${i}$ and $j$ are interchanged in this case. 
	
	\indent Thus, in the above edge coloring, edges of $G_{-r,n}$ takes only $r+n$ number of colors and colors of edges incident at each vertex of $G_{-r,n}$ are all distinct and thereby ${\chi}^{'}(G_{-r,n})$ = $r+n~$ since $~\Delta (G_{-r,n})$ $\leq$ ${\chi}^{'}(G_{-r,n}) \leq \Delta (G_{-r,n})+1$ by Vizing's Theorem \ref{1.1}. 
\end{proof}	

\begin{figure}
	\begin{center}
		\centering
		\resizebox{.95\textwidth}{!}{%
			\begin{tikzpicture}[scale=.8]
				
				\node (c0) at (11.75,6.6) [circle,draw,scale=0.8, fill = yellow!30]{0};
				\node (c1) at (14.3,5.9) [circle,draw,scale=0.8, fill = yellow!30]{1};
				\node (c2) at (16,4.1) [circle,draw,scale=0.8, fill = yellow!30] {2};
				\node (c3) at (16.5,1.8) [circle,draw,scale=0.8, fill = yellow!30] {3};
				\node (c4) at (16.1,0.1) [circle,draw,scale=0.8, fill = yellow!30] {4};
				\node (c5) at (15,-1.3) [circle,draw,scale=0.8, fill = yellow!30] {5};
				\node (c6) at (13.4,-2.1) [circle,draw,scale=0.8, fill = yellow!30] {6};
				\node (c7) at (11.5,-2.3) [circle,draw,scale=0.8, fill = yellow!30] {7};
				\node (c8) at (9.72,5.9) [circle,draw,scale=0.8, fill = yellow!30] {-1};
				\node (c9) at (8,4.1) [circle,draw,scale=0.8, fill = yellow!30] {-2};
				\node (c10) at (7.5,1.8) [circle,draw,scale=0.8,fill = yellow!30] {-3};
				\node (c11) at (8,-0.1) [circle,draw,scale=0.8, fill = yellow!30] {-4};
				\node (c12) at (9.5,-1.7) [circle,draw,scale=0.8, fill = yellow!30] {-5};
				
				\draw (c0)[gray, thick] --node[][above]{} (c1);
				\draw (c0)[orange, thick] --node[near start] {} (c2);
				\draw (c0)[cyan, thick] --node[near start] {} (c3);
				\draw (c0)[purple, thick] --node[near start] { } (c4);
				\draw (c0)[blue, thick] --node[near start] {} (c5);
				\draw (c0)[green, thick]--node[near start] {} (c6);
				\draw (c0)[brown, thick] --node[near start] {} (c7);
				\draw (c0)[violet, thick] --node[near start] {} (c8);
				\draw (c0)[magenta, thick] --node[near start] {} (c9);
				\draw (c0)[red, thick] --node[near start] {} (c10);
				\draw (c0)[olive, thick] --node[near start] {} (c11);
				\draw (c0)[black, thick] --node[near start] {} (c12);
				
				\draw (c1)[orange, thick] --node[near start] {} (c8);
				\draw (c1)[cyan, thick] --node[near start] {} (c9);
				\draw (c1)[purple, thick] --node[near start] {} (c10);
				\draw (c1)[blue, thick] --node[near start] {}(c11);
				\draw (c1)[green, thick] --node[near start] {} (c12);
				
				\draw (c2)[cyan, thick] --node[near start] {} (c8);
				\draw (c2)[purple, thick] --node[near start] {} (c9);
				\draw (c2)[blue, thick] --node[near start] {}(c10);
				\draw (c2)[green, thick] --node[near start] { } (c11);
				\draw (c2)[brown, thick] --node[near start] {} (c12);
				
				\draw (c3)[purple, thick] --node[near start] {} (c8);
				\draw (c3)[blue, thick] --node[near start] {} (c9);
				\draw (c3)[green, thick] --node[near start] {} (c10);
				\draw (c3)[brown, thick] --node[near start] {}(c11);
				\draw (c3)[violet, thick] --node[near start] {} (c12);
				
				\draw (c4)[blue, thick] --node[near start] {} (c8);
				\draw (c4)[green, thick] --node[near start] {} (c9);
				\draw (c4)[brown, thick] --node[near start] {} (c10);
				\draw (c4)[violet, thick] --node[near start] {}(c11);
				\draw (c4)[magenta, thick] --node[near start] {}(c12);
				
				\draw (c5)[green, thick] --node[near start] {} (c8);
				\draw (c5)[brown, thick] --node[near start] {} (c9);
				\draw (c5)[violet, thick] --node[near start] {}(c10);
				\draw (c5)[magenta, thick] --node[near start] {{ }} (c11);
				\draw (c5)[red, thick] --node[near start] {} (c12);
				
				\draw (c6)[brown, thick] --node[near start] {} (c8);
				\draw (c6)[violet, thick] --node[near start] {} (c9);
				\draw (c6)[magenta, thick] --node[near start] {} (c10);
				\draw (c6)[red, thick] --node[near start] {} (c11);
				\draw (c6)[olive, thick] --node[near start] {} (c12);
				
				\draw(c7)[magenta, thick] --node[near start] {} (c8);
				\draw (c7)[olive, thick] --node[near start] {} (c9);
				
				\draw (c7)[gray, thick] --node[near start] {} (c10);
				\draw (c7)[orange, thick] --node[near start] {} (c11);
				\draw (c7)[cyan, thick] --node[near start] {} (c12);
				
				\draw (c1)[violet, thick] --node[above] { } (c2);
				\draw (c1)[magenta, thick] --node[near start] {}(c3);
				\draw (c1)[red, thick] --node[near start] {} (c4);
				\draw (c1)[olive, thick] --node[near start] {} (c5);
				\draw (c1)[black, thick] --node[near start] { }(c6);
				
				\draw (c2)[red, thick] --node[right] {{ }} (c3);
				\draw (c2)[olive, thick] --node[near end] { } (c4);
				\draw (c2)[black, thick] --node[near start] { }(c5);
				
				\draw (c3)[black, thick] --node[right] {}(c4);
				
				\draw (c8)[red, thick] --node[left] { }(c9);
				\draw (c8)[olive, thick] --node[near end] {}(c10);
				\draw (c8)[black, thick] --node[near start] {} (c11);
				
				\draw (c9)[black, thick] --node{ } (c10);
				
				\node (a0) at (22.5,6.6) [circle,draw,scale=0.8,fill=yellow] {0};
				\node (a1) at (24.8,5.9)[circle,draw,scale=0.8,fill=yellow]{1};
				\node (a2) at (26.5,4.1) [circle,draw,scale=0.8,fill=yellow] {2};
				\node (a3) at (27,1.8) [circle,draw,scale=0.8,fill=yellow] {3};	
				\node (a4) at (20.2,5.9) [circle,draw,scale=0.8,fill=yellow] {-1};
				\node (a5) at (18.5,4.1) [circle,draw,scale=0.8,fill=yellow] {-2};
				\node (a6) at (18,1.8) [circle,draw,scale=0.8,fill=yellow] {-3};	
				\node (a7) at (18.5,-0.1) [circle,draw,scale=0.8,fill=yellow] {-4};
				\node (a8) at (20.1,-1.7) [circle,draw,scale=0.8,fill=yellow] {-5};
				\node (a9) at (22.3,-2.3) [circle,draw,scale=0.8,fill=yellow] {-6};
				\node (a10) at (24.6,-1.9) [circle,draw,scale=0.8,fill=yellow] {-7};
				\node (a11) at (26.3,-.4) [circle,draw,scale=0.8,fill=yellow] {-8};
				
				\draw (a0)[ orange, thick] --node[]{ } (a1);
				\draw (a0)[ cyan, thick] --node[near start] { } (a2);
				\draw (a0)[ purple, thick] --node[near start] { } (a3);
				\draw (a0)[ blue, thick] --node[near start] { } (a4);
				\draw (a0)[ green, thick] --node[near start] {} (a5);
				\draw (a0)[ brown, thick]--node[near start] {} (a6);
				\draw (a0)[ violet, thick] --node[near start] {} (a7);
				\draw (a0)[ magenta, thick] --node[near start] {} (a8);
				\draw (a0)[ red, thick] --node[near start] {} (a9);
				\draw (a0)[ olive, thick] --node[near start] {} (a10);
				\draw (a0)[ gray, thick] --node[near start] {} (a11);
				
				\draw (a1)[ cyan, thick] --node[near start] { } (a4);
				\draw (a1)[ purple, thick] --node[near start] { } (a5);
				\draw (a1)[ blue, thick] --node[near start] { } (a6);
				\draw (a1)[ green, thick] --node[near start] { } (a7);
				\draw (a1)[ brown, thick] --node[near start] { } (a8);
				\draw (a1)[ violet, thick] --node[near start] { } (a9);
				\draw (a1)[ magenta, thick] --node[near start] { } (a10);
				\draw (a1)[ red, thick] --node[near start] {}(a11);
				
				\draw (a2)[ purple, thick] --node[near end] { } (a4);
				\draw (a2)[ blue, thick] --node[near start] { }(a5);
				\draw (a2)[ green, thick] --node[near end] { } (a6);
				\draw (a2)[ brown, thick] --node[near start] { }(a7);
				\draw (a2)[ violet, thick] --node[near start] { } (a8);
				\draw (a2)[ magenta, thick] --node[near start] { } (a9);
				\draw (a2)[ red, thick] --node[near start] {}(a10);
				\draw (a2)[ olive, thick] --node[near start] {} (a11);
				
				\draw (a3)[ green, thick] --node[right] { }(a4);
				\draw (a3)[ violet, thick] --node[right] { }(a5);
				\draw (a3)[ red, thick] --node[right] { }(a6);
				\draw (a3)[ gray, thick] --node[right] { }(a7);
				\draw (a3)[ orange, thick] --node[near start] { } (a8);
				\draw (a3)[ cyan, thick] --node[near start] {} (a9);
				\draw (a3)[ blue, thick] --node[near start] {} (a10);
				\draw (a3)[ brown, thick] --node[near start] {}(a11);
				
				\draw (a4)[ brown, thick] --node[near start] {}(a5);
				\draw (a4)[ violet, thick] --node[near start] {} (a6);
				\draw (a4)[ magenta, thick] --node[near start] {}(a7);
				\draw (a4)[ red, thick] --node[near start] {} (a8);
				\draw (a4)[ olive, thick] --node[near start] {} (a9);
				\draw (a4)[ gray, thick] --node[near start] {} (a10);
				
				\draw (a5)[ magenta, thick] --node[near start] {}(a6);
				\draw (a5)[ red, thick] --node[near start] {} (a7);
				\draw (a5)[ olive, thick] --node[near start] {} (a8);
				\draw (a5)[ gray, thick] --node[near start] {} (a9);
				
				\draw (a6)[ olive, thick] --node[near start] {} (a7);
				\draw (a6)[ gray, thick] --node[near start] {} (a8);
				
				\draw (a1)[ olive, thick] --node[above] { } (a2);
				
		\end{tikzpicture} }%
		
		Fig. 75. $G_{-5, 7}$ with edge coloring \hspace{2cm} Fig. 76. $G_{-8, 3}$ with edge coloring
		\label{$G_{-8, 3}$ with edge coloring.}
	\end{center}
\end{figure}

\begin{illu} [{\cite[Illustration 4.2]{js24}}] \label{iluG-57}\nonumber\rm
	Consider $G_{-5,7}$, the case where $n > \left\lceil \frac{m}{2}\right\rceil$. We have $\Delta(G_{-5,7})$ = 12. Using the algorithm given in the above proof, we color the edges of $G_{-5,7}$ as follows.  \\
	${\color{gray}c_1} \rightarrow v_0v_1, v_7v_{-3}$;\\
	${\color{orange}c_2} \rightarrow v_0v_2,  v_1v_{-1}, v_7v_{-4}$;\\
	${\color{cyan}c_3} \rightarrow v_0v_3, v_1v_{-2}, v_2v_{-1}, v_7v_{-5}$;\\
	${\color{purple}c_4} \rightarrow v_0v_4, v_1v_{-3}, v_2v_{-2}, v_3v_{-1}$;\\
	${\color{blue}c_5} \rightarrow v_0v_5, v_1v_{-4}, v_2v_{-3}, v_3v_{-2},  v_4v_{-1}$;\\
	${\color{green}c_6} \rightarrow v_0v_6, v_1v_{-5}, v_2v_{-4}, v_3v_{-3}, v_4v_{-2}, v_5v_{-1}$;\\
	${\color{brown}c_7} \rightarrow v_0v_7, v_2v_{-5}, v_3v_{-4}, v_4v_{-3}, v_5v_{-2}, v_6v_{-1}$;\\
	${\color{violet}c_8} \rightarrow v_0v_{-1}, v_3v_{-5}, v_4v_{-4}, v_5v_{-3}, v_6v_{-2}, v_1v_2$;\\
	${\color{magenta}c_{9}} \rightarrow v_0v_{-2}, v_4v_{-5}, v_5v_{-4}, v_6v_{-3}, v_7v_{-1}, v_1v_3$\\
	${\color{red}c_{10}} \rightarrow v_0v_{-3}, v_5v_{-5}, v_6v_{-4}, v_1v_4, v_2v_3, v_{-1}v_{-2}$;\\
	${\color{olive}c_{11}} \rightarrow v_0v_{-4}, v_6v_{-5}, v_7v_{-2}, v_1v_5, v_2v_4, v_{-1}v_{-3}$;\\
	${\color{black}c_{12}} \rightarrow v_0v_{-5}, v_1v_6, v_2v_5, v_3v_4, v_{-1}v_{-4}, v_{-2}v_{-3}$.\\
	
	The color assignments of the integral sum graph $G_{-5,7}$ are displayed in the Table \ref{abcde} and the edge coloring of $G_{-5,7}$ is shown in Figure 75. Here, 12 colors are used to color the edges of the integral sum graph $G_{-5,7}$. This implies, $\chi'(G_{-5,7})$ = 12 = $\Delta(G_{-5,7})$ using Vizing theorem \cite{v65}. 
\end{illu}

Similarly, by the following illustration, it is shown that edges of the integral sum graph $G_{-8,3}$, the case where $n < \left\lceil \frac{m}{2}\right\rceil$, can be colored using 11 colors by the algorithm given in the proof of Theorem \ref{thm Grs 2 chrom}. 

\begin{illu} [{\cite[Illustration 4.5]{js24}}] \label{iluG-8,3}\nonumber\rm
	Consider $G_{-8,3}$, the case where $n < \left\lceil \frac{m}{2}\right\rceil$.   We have $\Delta(G_{-8,3})$ = 11. Using the algorithm given in the above proof, we can color the edges of $G_{-8,3}$ as follows.  \\
	${\color{orange}c_1} \rightarrow v_0v_1, v_{-5}v_3$;\\
	${\color{cyan}c_2} \rightarrow v_0v_2,  v_{-1}v_1, v_{-6}v_3$;\\
	${\color{purple}c_3} \rightarrow v_0v_3, v_{-1}v_2, v_{-2}v_1$;\\
	${\color{blue}c_4} \rightarrow v_0v_{-1}, v_
	{-2}v_2, v_{-3}v_1, v_{-7}v_3$;\\
	${\color{green}c_5} \rightarrow v_0v_{-2}, v_{-1}v_3, v_{-3}v_2, v_{-4}v_1$;\\
	${\color{brown}c_6} \rightarrow v_0v_{-3}, v_{-1}v_{-2}, v_{-4}v_2, v_{-5}v_1, v_{-8}v_3$;\\
	${\color{violet}c_7} \rightarrow v_0v_{-4},v_{-1}v_{-3}, v_{-2}v_{3} ,v_{-5}v_2, v_{-6}v_1$;\\
	${\color{magenta}c_8} \rightarrow v_0v_{-5}, v_{-1}v_{-4}, v_{-2}v_{-3}, v_{-6}v_2, v_{-7}v_1$;\\
	${\color{red}c_9} \rightarrow v_0v_{-6}, v_{-1}v_{-5}, v_{-2}v_{-4}, v_{-3}v_3, v_{-7}v_2, v_{-8}v_1$;\\
	${\color{olive}c_{10}} \rightarrow v_0v_{-7}, v_{-1}v_{-6}, v_{-2}v_{-5}, v_{-3}v_{-4}, v_{-8}v_2$;\\
	${\color{gray}c_{11}} \rightarrow v_0v_{-8}, v_{-1}v_{-7}, v_{-2}v_{-6}, v_{-3}v_{-5}, v_{-4}v_3, v_{1}v_2$.\\
	
	The color assignments of edges of the integral sum graph $G_{-8,3}$ are displayed in the Table \ref{abcd} and the edge coloring of $G_{-8,3}$ is shown in Figure 76. Here, 11 colors are used to color the edges of the integral sum graph $G_{-8,3}$. This implies, $\chi'(G_{-8,3})$ = 11 = $\Delta(G_{-8,3})$ using Vizing Theorem \ref{1.1}. 
\end{illu}

\begin{table} 
	\caption{{Colors assigned to edges of integral sum graph $G_{-5,7}$}.}
	\label{abcde}
	\begin{center}
		\scalebox{0.8}{
			{\bf	\scriptsize
				\begin{tabular}{||c|c|ccccccccc c||} \hline \hline
					&  &  \multicolumn{10}{c||}{ } \\ 
					\tiny	{Sl. No.:} & Color &  \multicolumn{10}{c||}{Edges}\\ 
					&  &  \multicolumn{10}{c||}{ }\\  \hline \hline
					&  &  &  &  &  &  &  &  &  &  & \\ 
					\tiny	$1$ &	$C_1$&	$(0,1)$ &  & & & & 	& & $(7,-3)$ &	&	\\
					&  &  &  &  &  &  &  &  &  &  & \\  \hline
					&  &  &  &  &  &  &  &  &  &  & \\ 
					\tiny	$2$ & $C_2$ & { $(0,2)$} & $(1,-1)$ &	  & &	&	& &	$(7,-4)$	 &	&	 \\ 
					&  &  &  &  &  &  &  &  &  &  & \\  \hline
					&  &  &  &  &  &  &  &  &  &  & \\ 
					\tiny	$3$ & $C_3$ & $(0,3)$ &	$(1,-2)$ & $(2,-1)$ &	&	& & & $(7,-5)$ &	&		 \\ 
					&  &  &  &  &  &  &  &  &  &  & \\  \hline
					&  &  &  &  &  &  &  &  &  &  & \\ 		
					\tiny	$4$ & $C_4$ & $(0,4)$ &	$(1,-3)$ & $(2,-2)$ & $(3,-1)$ &	& & & & 	&		 \\ 
					&  &  &  &  &  &  &  &  &  &  & \\  \hline
					&  &  &  &  &  &  &  &  &  &  & \\ 		
					\tiny	$5$ & $C_5$ & $(0,5)$ &	$(1,-4)$ & $(2,-3)$ &	$(3,-2)$ & $(4,-1)$ & & & & 	&		 \\ 
					&  &  &  &  &  &  &  &  &  &  & \\  \hline
					&  &  &  &  &  &  &  &  &  &  & \\ 		
					\tiny	$6$ & $C_6$ & $(0,6)$ &	$(1,-5)$ & $(2,-4)$  & $(3,-3)$ & $(4,-2)$	& $(5,-1)$ & & & 	&		 \\ 
					&  &  &  &  &  &  &  &  &  &  & \\  \hline
					&  &  &  &  &  &  &  &  &  &  & \\ 		
					\tiny	$7$ & $C_7$ & $(0,7)$ &	& $(2,-5)$ & $(3,-4)$ &	$(4,-3)$ & $(5,-2)$ & $(6,-1)$ & & 	&		 \\ 
					&  &  &  &  &  &  &  &  &  &  & \\  \hline
					&  &  &  &  &  &  &  &  &  &  & \\ 		
					\tiny	$8$ & $C_8$ & $(0,-1)$ & $(1,2)$	& & $(3,-5)$ & $(4,-4)$ & $(5,-3)$	& $(6,-2)$ &	&  &\\ 
					&  &  &  &  &  &  &  &  &  &  & \\  \hline
					&  &  &  &  &  &  &  &  &  &  & \\ 		
					\tiny	$9$ & $C_9$ & $(0,-2)$ & $(1,3)$ &  & &	$(4,-5)$ & $(5,-4)$ & $(6,-3)$ & $(7,-1)$ & 	&	\\ 
					&  &  &  &  &  &  &  &  &  &  & \\  \hline
					&  &  &  &  &  &  &  &  &  &  & \\ 		
					\tiny	$10$ & $C_{10}$ & $(0,-3)$ & $(1,4)$ & $(2,3)$ & & &	$(5,-5)$ & $(6,-4)$ & 	&	$(-1,-2)$  & \\ 
					&  &  &  &  &  &  &  &  &  &  & \\  \hline
					&  &  &  &  &  &  &  &  &  &  & \\ 
					\tiny	$11$ & $C_{11}$ & $(0,-4)$ & $(1,5)$ & $(2,4)$ & & &	& $(6,-5)$ & $(7,-2)$ & $(-1,-3)$  &	\\ 
					&  &  &  &  &  &  &  &  &  &  & \\  \hline
					&  &  &  &  &  &  &  &  &  &  & \\ 
					\tiny	$12$ & $C_{12}$ & $(0,-5)$ & $(1,6)$	& $(2,5)$	& $(3,4)$ &  & & & & $(-1,-4)$ & $(-2,-3)$  \\  
					&  &  &  &  &  &  &  &  &  &  & \\  \hline\hline
		\end{tabular}} }
	\end{center}
\end{table} 

\begin{table} 
	\caption{{Colors assigned to edges of integral sum graph $G_{-8,3}$}.}
	\label{abcd}
	\begin{center}
		\scalebox{0.8}{
			{\bf	\scriptsize
				\begin{tabular}{||c|c|ccccccccc c||} \hline \hline
					&  &  \multicolumn{10}{c||}{ } \\ 
					\tiny	{Sl. No.:} & Color &  \multicolumn{10}{c||}{Edges}\\ 
					&  &  \multicolumn{10}{c||}{ }\\  \hline \hline
					&  &  &  &  &  &  &  &  &  &  & \\ 
					\tiny	$1$&	$C_1$&	$(0,1)$&		  &	 & & & 	$(-5,3)$		& & &	&	\\
					&  &  &  &  &  &  &  &  &  &  & \\  \hline
					&  &  &  &  &  &  &  &  &  &  & \\ 
					\tiny	$2$&	$C_2$&	{ $(0,2)$}&$	(-1,1)$&	  & &	&	&	$(-6,3)$	& &	&	 \\ 
					&  &  &  &  &  &  &  &  &  &  & \\  \hline
					&  &  &  &  &  &  &  &  &  &  & \\ 
					\tiny	$3$&$	C_3$&$	(0,3)$&	$(-1,2)$&$(-2,1)$&	&	& & & & 	&		 \\ 
					&  &  &  &  &  &  &  &  &  &  & \\  \hline
					&  &  &  &  &  &  &  &  &  &  & \\ 		
					\tiny	$4$&$	C_4$&$	(0,-1)$&	      &$(-2,2)$&$	(-3,1)$&	&	& &	$(-7,3)$&  &\\ 
					&  &  &  &  &  &  &  &  &  &  & \\  \hline
					&  &  &  &  &  &  &  &  &  &  & \\ 		
					\tiny	$5$&$	C_5$&$	(0,-2)$&$	(-1,3)$&	   &$	(-3,2)$&	$(-4,1)$	& & & & 	&	\\ 
					&  &  &  &  &  &  &  &  &  &  & \\  \hline
					&  &  &  &  &  &  &  &  &  &  & \\ 		
					\tiny	$6$&$	C_6$&$	(0,-3)$&	$(-1,-2)$&	   &      	&	$(-4,2)$&	$(-5,1)$ & & 	&	$(-8,3)$  & \\ 
					&  &  &  &  &  &  &  &  &  &  & \\  \hline
					&  &  &  &  &  &  &  &  &  &  & \\ 
					\tiny	$7$&$	C_7$&	$(0,-4)$&$	(-1,-3)$& $	(-2,3)$	& & &	$	(-5,2)$ &	$(-6,1)$	&  	& &	\\ 
					&  &  &  &  &  &  &  &  &  &  & \\  \hline
					&  &  &  &  &  &  &  &  &  &  & \\ 
					\tiny	$8$&$	C_8$&	$(0,-5)$&$	(-1,-4)$&$	(-2,-3)$	& & & & $(-6,2)$&$	(-7,1)$& &  \\ 
					&  &  &  &  &  &  &  &  &  &  & \\  \hline
					&  &  &  &  &  &  &  &  &  &  & \\ 
					\tiny	$9$&$	C_9$&	$(0,-6)$&$(-1,-5)$& $	(-2,-4)$ &$	(-3,3)$	&  &		& & 	$(-7,2)$ &$	(-8,1)$	& \\ 
					&  &  &  &  &  &  &  &  &  &  & \\  \hline
					&  &  &  &  &  &  &  &  &  &  & \\ 
					\tiny	$10$&$	C_{10}$&$(0,-7)$&$	(-1,-6)$& $	(-2,-5)$	& $(-3,-4)$ & & & & & 					$(-8,2)$	& \\ 
					&  &  &  &  &  &  &  &  &  &  & \\  \hline
					&  &  &  &  &  &  &  &  &  &  & \\ 
					\tiny	$11$&$	C_{11}$&$(0,-8)$& $(-1,-7)$& $(-2,-6)$ &	$(-3,-5)$	&$ (-4,3)$ & & & & 		&			$(1,2)$ \\ 
					&  &  &  &  &  &  &  &  &  &  & \\  \hline\hline
		\end{tabular}} }
	\end{center}
\end{table} 

\subsection{Integral sum graphs of class 1 and class 2} 

For any graph $G$, Vizing's Theorem \ref{1.1} gives a bound for its edge chromatic number.	Based on the bounds, graphs are clssified into class 1 and 2 as follows \cite{v65}.

\medskip
A simple graph $G$ is of {\em class 1} if $\chi^{'}\left(G \right) $ = $\Delta(G)$. It is of {\em class 2} if $\chi^{'}\left(G \right) $ = $\Delta(G)+1$. Similar to this classification,  in \cite{vm13c}, integral sum graphs are classified into {\em integral sum graph of class 1} and {\em integral sum graph of class 2}. Here, we present a few results on these integral sum graphs.

For a given integral sum graph $G$, the edge coloring of $G$ obtained by considering each edge-sum class as an edge color class need not be a minimal edge coloring of $G$. Also, for a given integral sum graph $G^+(S)$, $\chi'_{\mathbb{Z}-sum}(G^+(S))$ need not be equal to the order of $S$ and $\chi'_{\mathbb{Z}-sum}(G^+(S))$ = $|S|$ if and only if $|E_i| \neq \emptyset$  for every $i \in S$.

\begin{dfn} \cite{vm13c}  Integral sum graph $G^+(S)$ is said to be of {\em class 1} if the edge-sum chromatic number and edge chromatic number of $G^+(S)$ are the same. Otherwise, $G^+(S)$ is of {\em class 2}. That is, an integral sum graph $G^+(S)$ is of {\em class 1} if $\chi'_{\mathbb{Z}-sum}(G^+(S))$ = $\chi'(G^+(S))$ and is of {\em class 2} if $\chi'_{\mathbb{Z}-sum}(G^+(S))$ $\neq$ $\chi'(G^+(S))$.
\end{dfn} 

 Note that in \cite{vm13c}, we called integral sum graph of class 1 and integral sum graph of class 2 as {\em edge sum-perfect color graph} and {\em edge sum-non-perfect color graph}, respectively. Hereafter, we use the terms integral sum graph of class 1 and 2 only. 

In Section 10.1, we obtained all the edge-sum classes of $G_{-2,4}$ and from the corresponding edge sum coloring, we get, $\chi'_{\mathbb{Z}-sum}(G_{-2, 4})$ = 7. Edge-sum coloring of $G_{-2, 4}$ with  $\chi'_{\mathbb{Z}-sum}(G_{-2, 4})$ = 7 is given in Figure 73 and a minimum edge coloring of $G_{-2, 4}$ with $\chi'(G_{-2, 4})$ = 6 = $\Delta(G_{-2, 4})$, using Vizing Theorem \ref{1.1}, is given in Figure 74. This implies, $G_{-2, 4}$ is a class 2 integral sum graph.  


\medskip
 The following example illustrates the edge-sum coloring of $G_{-1,5}$.

\begin{enumerate}
	\item [\rm 1.] $E_1$ = $\{(0, 1), (-1, 2)\}$ = {\color{black}  $\{$edges with black color$\}$, }
	
	\item [\rm 2.] $E_2$ = $\{(0, 2), (-1, 3)\}$ = {\color{red}  $\{$edges with red color$\}$, }
	
	\item [\rm 3.] $E_3$ = $\{(0, 3), (1, 2), (-1, 4)\}$  = {\color{green}  $\{$edges with green color$\}$, } 
	
	\item [\rm 4.] $E_4$ = $\{(0, 4), (1, 3), (-1, 5)\}$ = {\color{orange}  $\{$edges with orange color$\}$, } 
	
	\item [\rm 5.] $E_5$ = $\{(0, 5), (1, 4), (2, 3)\}$ = {\color{violet}  $\{$edges with violet color$\}$, }

	\item [\rm 6.] $E_0$ = $\{(-1, 1)\}$ = {\color{cyan}  $\{$edge with cyan color$\}$ } and

	\item [\rm 7.] $E_{-1}$ = $\{(0, -1) \}$ = {\color{blue}  $\{$edge with blue color$\}$ }. 
\end{enumerate}

This implies, $\chi'_{\mathbb{Z}-sum}(G_{-1, 5})$ = 7. An edge-sum coloring of $G_{-1, 5}$ with  $\chi'_{\mathbb{Z}-sum}(G_{-1, 5})$ = 7 is given in Figure 77. On the other hand, one set of edge color classes of $G_{-1,5}$ with $\chi'(G_{-1, 5})$ = 6 is given in Figure 78 and is 

$\{\{ (0,1), (-1,2)\}$ = {\color{black}  $\{$edges with black color$\}$, }  

~$\{(0,2), (-1, 3)\}$ = {\color{red}  $\{$edges with red color$\}$, }  

~$\{(0,3), (1,2), (-1,4)\}$ = {\color{green}  $\{$edges with green color$\}$, }  

~$\{(0,4), (1,3), (-1,5)\}$ = {\color{cyan}  $\{$edges with cyan color$\}$, }  

~$\{(0,5), (-1, 1)\}$ = {\color{blue}  $\{$edges with blue color$\}$, }  

~$\{(0, -1), (1,4), (2,3)\}$ = {\color{orange}  $\{$edges with orange color$\}$}  $\}$ which is a minimum edge coloring of $G_{-1, 5}$, follows from Vizing Theorem \ref{1.1} since $\Delta(G_{-1,5})$ = 6. This implies, $\chi'(G_{-1, 5})$ = 6. 

\medskip
Thus, $\chi'_{\mathbb{Z}-sum}(G_{-1, 5})$ = 7 $\neq$ $\chi'(G_{-1, 5})$ = 6. This implies, $G_{-1, 5}$ is an integral sum graph of class 2.  

\medskip
Whereas Star graph $S_n$ is an integral sum graph of class 1, see Theorem \ref{a} which also proves that for $n \in \mathbb{N}$, $G_{-1,1}$ and $G_{0,n}$ are integral sum graphs of class 1 and for $n \geq 2$, $G_{-1,n}$ is an integral sum graph of class 2.

\medskip

\begin{center}
	\begin{tikzpicture} [scale = 1.3] 
		
		\node (a0) at (13,2)  [circle,draw,scale=0.6]{0};
		\node (a1) at (14.25,1.75)  [circle,draw,scale=0.6]{1};
		\node (a2) at (15,.5)  [circle,draw,scale=0.6] {2};
		\node (a3) at (14.25,-0.5)  [circle,draw,scale=0.6] {3};
		\node (a4) at (13,-0.5)  [circle,draw,scale=0.6] {4};
		\node (a5) at (12,-.25)  [circle,draw,scale=0.6]{5};
		\node (a6) at (12,1)  [circle,draw,scale=0.5] {-1};
		
		\draw (a0)[black]  -- (a1);
		\draw (a6)[black]  -- (a2);
		
		\draw (a0)[red, thick] -- (a2);
		\draw (a6)[red, thick]  -- (a3);
		
		\draw (a0)[green, thick] -- (a3);
		\draw (a1)[green, thick] -- (a2);
		\draw (a6)[green, thick] -- (a4); 
		
		\draw (a0)[orange, thick]  -- (a4);
		\draw (a1)[orange, thick]  -- (a3);	
		
		\draw (a0)[violet, thick]  -- (a5); 
		\draw (a1)[violet, thick]  -- (a4); 
		\draw (a2)[violet, thick]  -- (a3); 
		
		\draw (a0)[blue, thick]  -- (a6);
		
		\draw (a6)[cyan, thick]  -- (a1);
		
		\draw (a5)[orange, thick] -- (a6);
		
		\node (b0) at (18.5,2)  [circle,draw,scale=0.6]{0};
		\node (b1) at (19.75,1.75)  [circle,draw,scale=0.6]{1};
		\node (b2) at (20.5,.5)  [circle,draw,scale=0.6] {2};
		\node (b3) at (19.75,-0.5)  [circle,draw,scale=0.6] {3};
		\node (b4) at (18.5,-0.5)  [circle,draw,scale=0.6] {4};
		\node (b5) at (17.5,-.25)  [circle,draw,scale=0.6]{5};
		\node (b6) at (17.5,1)  [circle,draw,scale=0.5] {-1};
		
		\draw (b0)[black]  -- (b1);
		\draw (b6)[black]  -- (b2);
		
		\draw (b0)[red, thick] -- (b2);
		\draw (b6)[red, thick]  -- (b3);
		
		\draw (b1)[green, thick] -- (b2);
		\draw (b0)[green, thick] -- (b3);
		\draw (b6)[green, thick] -- (b4); 
		
		\draw (b1)[cyan, thick]  -- (b3);	
		\draw (b0)[cyan, thick]  -- (b4);
		\draw (b6)[cyan, thick]  -- (b5);
		
		\draw (b6)[blue, thick]  -- (b1);
		\draw (b0)[blue, thick]  -- (b5); 
		
		\draw (b0)[orange, thick]  -- (b6);
		\draw (b1)[orange, thick]  -- (b4);
		\draw (b2)[orange, thick]  -- (b3);
	\end{tikzpicture}
	
	\medskip
	Fig. 77. $G_{-1, 5}$ with $\chi'_{\mathbb{Z}-sum}(G_{-1, 5})$ = 7  \hspace{1cm} Fig. 78. $G_{-1, 5}$ with $\chi'(G_{-1, 5})$ = 6 = $\Delta(G_{-1, 5})$ \hspace{.5cm} 
	
\end{center}

\begin{theorem} [{\cite[Theorem 2.5]{vl22}}] \label{a} {\rm For $m,n \in \mathbb{N}$, the following hold:
	\begin{enumerate}
		\item [\rm (a)]  $\chi'_{\mathbb{Z}-sum}(G_{-m,n})$ = $m+n+1$; 
		\item [\rm (b)]  $G_{-1,1}$ and $G_{0,n}$ are class 1 integral sum graphs, 
		
		$\chi'_{\mathbb{Z}-sum}(G_{-1,1})$ = 3 = $\chi'(G_{-1,1})$ and 
		
		$\chi'_{\mathbb{Z}-sum}(G_{0, n})$ = $n$ = $\chi'(G_{0, n})$; 
		\item [\rm (c)]  Star graph $S_{n}$ is a class 1 integral sum graph and 
		
		$\chi'_{\mathbb{Z}-sum}(S_{n})$ = $n - 1$ = $\chi'(S_{n})$; 
		\item [\rm (d)] For $n \geq 2$, $G_{-1,n}$ is a class 2 integral sum graph,
		
		$\chi'_{\mathbb{Z}-sum}(G_{-1,n})$ = $n+2$ and $\chi'(G_{-1,n})$ = $1+n$; and  
		\item [\rm (e)] For $2 \leq n \leq 6$, $G_{-n,n}$ is a class 2 integral sum graph,
		
		$\chi'_{\mathbb{Z}-sum}(G_{-n,n})$ = $2n+1$ and $\chi'(G_{-n,n})$ = $2n$.
	\end{enumerate}
}
\end{theorem}

\begin{proof}
\begin{enumerate}
	\item [\rm (a)] The edge-sum classes of the integral sum graph $G_{-m,n}$ are $E_{-m}$, $E_{-m+1}$, $\dots$, $E_{-1}$, $E_0$, $E_1$, $\dots$, $E_n$ and each one is non-empty. Hence, $\chi'_{\mathbb{Z}-sum}(G_{-m,n})$ = $m+n+1$ for $m,n \in \mathbb{N}$. 

	\item [\rm (b)] By observation, it is clear that $\chi'_{\mathbb{Z}-sum}(G_{-1,1})$ = 3 = $\chi'(G_{-1,1})$. 
	
	For $n \in \mathbb{N}$, the edge-sum classes of integral sum graph $G_{0,n}$ are $E_0$, $E_1$, $\dots$, $E_n$ and each one, except $E_0$, is non-empty and hence $\chi'_{\mathbb{Z}-sum}(G_{0,n})$ = $n$. Using Vizing’s Theorem \ref{1.1}, we get $n \leq \chi'(G_{0,n}) \leq n+1$ since $\Delta(G_{0,n})$ = $n$ = $deg(u_0)$. Now, $\{\{(0,1)\}$, $\{(0,2)\}$, $\{(0,3), (1,2)\}$, $\{(0,4), (1,3)\}$, $\{(0,5)$, $(1,4)$, $(2,3)\}$, $\{(0,6)$, $(1, 5)$, $(2,4)\}$, $\dots$, $\{(0,n)$, $(1,n-1)$, $(2,n-2)$, $\dots$, $(\left\lfloor \frac{n-1}{2}\right\rfloor, \left\lfloor \frac{n+2}{2}\right\rfloor)\}\}$ is a set of edge color classes of $G_{0,n}$ and it is of order $n$. This implies, $\chi'(G_{0,n})$ = $n$ and thereby, $\chi'(G_{0,n})$ = $n$ =  $\chi'_{\mathbb{Z}-sum}(G_{0,n})$. 

	\item [\rm (c)] We have Star graph $S_{n}$ = $K_1 * ((n - 1)K_1)$ = $K_{1,n-1}$, $n \geq 2$. Let $n \geq 2$, $V(S_n)$ = $\{u_0, v_1, v_2, \dots, v_{n-1}\}$, $d(u_0)$ = $n-1$ and $d(v_1)$ = 1 = $d(v_2)$ = $\dots$ = $d(v_{n-1})$. Define labeling $f:$ $V(S_n) \to \mathbb{Z}$ $\ni$ 	$f(u_0)$ = 0, $f(v_1)$ = $t$, $f(v_{j+1})$ = 	$d\big( \sum^{j}_{i=1} f(v_i) \big) + t$, $j$ = 1, 2, $\dots$, $n - 2$.
						
	Then the sequence of vertex labeling of $S_n$ is $\{f(u_0) = 0\}$ $\cup$ $\{f(v_i)$ = $t{(d + 1)}^{i-1}\}^{n-1}_{i=1}$, $d,t\in\mathbb{N}$. Also, 
		
	$(i)$ for $i$ = 1, 2, ..., $n-1$, $f(u_0)$ + $f(v_i)$ = 0 + $t{(d + 1)}^{i-1}$  = $f(v_i)$ and
		
	$(ii)$ for $n \geq 2$, $1 \leq i < j < k \leq n-1$ and $i,j,k,d,t\in\mathbb{N}$,  
	$f(v_i)$ + $f(v_j)$ $\neq$ $f(v_k)$ since 
		$t{(d + 1)}^i$ + $t{(d + 1)}^j$ = $t{(d + 1)}^i(1+{(d + 1)}^{j-i})$ $\neq$ $t{(d + 1)}^k$, 	$1 \leq j - i < j < k$. This implies, $u_0$ and $v_i$ are adjacent for every $i$ whereas $v_j$ and $v_k$ are non-adjacent for every $j$ and $k$, $j \neq k$, $1 \leq j,k \leq n-1$. Thus, $f$ is an integral sum labeling of $S_n$, $n \geq 2$. Hence $S_n$ is an integral sum graph.
	
	$S_n$ has $n$-1 edge-sum classes, each a singleton set and thereby
	$\chi'_{\mathbb{Z}-sum}(S_{n})$ = $n - 1$ and $\Delta(S_n)$ = $d(u_0)$ = $n-1$ = $\chi'(S_{n})$ since all the $n-1$ edges at $u_0$ take different colors. This implies, $S_n$ is a class 1 integral sum graph for $n \geq 2$.
	
	\item [\rm (d)]  For $n \geq 2$, the edge-sum classes of integral sum graph $G_{-1,n}$ are $E_0$, $E_1$, $\dots$, $E_n$, $E_{-1}$, each one is non-empty and hence $\chi'_{\mathbb{Z}-sum}(G_{-1,n})$ = $n+2$. $\Delta (G_{-1,n})$ = $1+n$ = degree of the vertex whose label is 0 in $G_{-1,n}$ and 
	
	$\{ \{(0,1), (-1,2)\}$, 
	
	$\{(0,2), (-1,3)\}$, 
	
	$\{(0,3), (-1,4), (1,2)\}$, 
	
	$\{(0,4)$, $(-1,5)$, $(1,3)\}$, 
	
	$\{(0,5), (-1,6), (1,4), (2,3)\}$, 
	
	$\{(0,6), (-1,7), (1,5)$, $(2,4)\}$, 
	
	$\{(0,7)$, $(-1,8), (1,6), (2,5), (3,4)\}$, 
	
	$\{(0,8), (-1,9), (1,7), (2,6), (3,5)\}$,
	
	$\dots$,
	
	$\{(0,n-1)$, $(-1,n)$, $(1,n-2)$, $(2,n-3)$, $\dots$, $(\left\lfloor \frac{n-3}{2}\right\rfloor, \left\lfloor \frac{n+2}{2}\right\rfloor )\}$, 
	
	$\{(0,n)$, $(-1,1)\}$, 
	
	$\{(0,-1)$, $(1,n-1)$, $(2,n-2), (3,n-3), \dots, ( \left\lfloor \frac{n-1}{2}\right\rfloor, \left\lfloor \frac{n+2}{2}\right\rfloor )\}\}$ 
	\\
	is a set of edge color classes of $G_{-1,n}$ which is of order $n+1$. This implies, the edge chromatic number of $G_{-1,n}$ is $n+1$, using Vizing's theorem \cite{v65}. This implies, for $n \geq 2$, $\chi'(G_{-1,n})$ = $1+n$ $\neq$ $\chi'_{\mathbb{Z}-sum}(G_{-1,n})$ = $n+2$. Hence $G_{-1,n}$ is a class 2 integral sum graph. 

	\item [\rm (e)] At first, for $n$ = 2,3,4,5,6, a set of edge color class of order $2n$ of $G_{-n,n}$ is obtained as follows:
	\begin{enumerate}
		\item [\rm (e1)]  A set of edge color classes of order 4 of $G_{-2,2}$ is 
		
		$\{\{(0,1), (-2,2)\}$, 
		
		~$\{(0,2) (-1,1)\}$, 
		
		~$\{(0,-1), (-2,1)\}$, 
		
		~$\{(0,-2), (-1,2)\}\}$. 

		\item [\rm (e2)]  A set of edge color classes of order 6 of $G_{-3,3}$ is 
		
		$\{\{(0,1)$, $(-1,2)$, $(-2,3)\}$, 
		
		~$\{(0,2)$, $(-1,3)$, $(1,-3)\}$, 
		
		~$\{(0,3)$, $(1,2)$, $(-1,-2)\}$, 
		
		~$\{(0,-1)$, $(1,-2)$, $(2,-3)\}$, 
		
		~$\{(0,-2)$, $(1,-1)$, $(3,-3)\}$, 
		
		~$\{(0,-3)$, $(2,-2)\}\}$. 

		\item [\rm (e3)]  A set of edge color classes of order 8 of $G_{-4,4}$ is 
		
		$\{\{(0,1)$, $(-3,3)$, $(-2,4)$, $(-1,2)\}$, 
		
		~$\{(0,2)$, $(-4,1)$, $(-3,4)$, $(-1,3)\}$, 
		
		~$\{(0,3)$, $(-4,4)$, $(1,2)\}$, 
		
		~$\{(0,4), (-4,3), (-2,2)$, $(-1,1)\}$, 
		
		~$\{(0,-1), (-3,2), (-2,3)\}$, 
		
		~$\{(0,-2)$, $(-3,1)$, $(-1, 4)\}$, 
		
		~$\{(0,-3), (-4,2), (-2,-1), (1,3)\}$, 
		
		~$\{(0,-4)$, $(-3,-1)$, $(-2,1)\}\}$. 

		\item [\rm (e4)]  A set of edge color classes of order 10 of $G_{-5,5}$ is 
		
		$\{\{(0,1)$, $(-4,4)$, $(-3,5)$, $(-2,3), (-1,2)\}$, 
		
		~$\{(0,2), (-3,4), (-2,1), (-1,3)\}$, 
		
		~$\{(0,3), (-1,4), (1,2)\}$, 
		
		~$\{(0,4), (-5,3), (-4,1)$, $(-3,2)$, $(-2,5)\}$, 
		
		~$\{(0,5)$, $(-5,1)$, $(-4,-1)$, $(-3,-2)$, $(2,3)\}$, 
		
		~$\{(0,-1)$, $(-5,5)$, $(-4,3)$, $(-3,1)$, $(-2,2)\}$, 
		
		~$\{(0,-2)$, $(-5,4)$, $(-4,2), (-3,3), (-1,1)\}$, 
		
		~$\{(0,-3), (-2,-1), (1,3)\}$, 
		
		~$\{(0,-4),(-5,2), (-2,4), (-1,5)\}$, 
		
		~$\{(0,-5), (-4,5), (-3,-1), (1,4)\}\}$. 

		\item [\rm (e5)]  A set of edge color classes of order 12 of $G_{-6,6}$ is 
		
		$\{\{(0,1)$, $(-5,6)$, $(-4,4)$, $(-3,5)$, $(-2,3)$, $(-1,2)\}$, 
		
		~$\{(0,2)$, $(-6,4)$, $(-5, 1)$, $(-3,3)$, $(-2,5)\}$, 
		
		~$\{(0,3)$, $(-6,6)$, $(-4,5)$, $(-3,1)$, $(-2,2)$, $(-1,4)\}$, 
		
		~$\{(0,4), (-6,1), (-4,2)$, $(-2,6)$, $(-1, 5)\}$, 
		
		~$\{(0,5), (-1,6), (1,4), (2,3)\}$, 
		
		~$\{(0,6), (-4,3)$, $(2,4)$, $(-1,1)\}$, 
		
		~$\{(0, -1), (-5,5), (-3,2), (-2,1)\}$, 
		
		~$\{(0,-2)$, $(-5,4)$, $(-4$, $1)$, $(-3,6)$, $(-1,3)\}$, 
		
		~$\{(0,-3)$, $(-6,3)$, $(-5,-1)$, $(-4,-2)$, $(1,5)\}$, 
		
		~$\{(0,-4)$, $(-6,5)$, $(-5,2)$, $(-3,-1)$, $(-2,4)$, $(1,3)\}$, 
		
		~$\{(0,-5)$, $(-6,2), (-4$, $-1)$, $(-3,-2)\}$, 
		
		~$\{(0,-6)$, $(-5,3)$, $(-4,6)$, $(-3,5)$, $(-2,-1)$,  $(1,2)\}\}$.
	\end{enumerate} 

\medskip
	Using $(a)$, we get $\chi'_{\mathbb{Z}-sum}(G_{-m,n}) = m+n+1$ for $m,n \in \mathbb{N}$ and thereby $\chi'_{\mathbb{Z}-sum}(G_{-n,n})$ = $2n+1$ for $n \in \mathbb{N}$. For $n \geq 2$, $\Delta (G_{-n,n})$ = $2n$ = degree of the vertex whose label is $0$ in $G_{-n,n}$. This implies, for $n$ = 2,3,4,5,6, using Vizing's Theorem \ref{1.1}, the edge chromatic number of $G_{-n,n}$ is $2n$. Thus, for $2 \leq n \leq 6$, $\chi'(G_{-n,n})$ = $2n \neq \chi'_{\mathbb{Z}-sum}(G_{-n,n})$ = $2n+1$ and thereby, $G_{-n,n}$ is a class 2 integral sum graph for $n$ = 2,3,4,5,6. 
\end{enumerate}
Hence we get the result.
\end{proof}

\begin{theorem} [{\cite[Theorem 2.6]{vl22}}] \label{b} {\rm Let $n \geq 2$ and $m,n\in\mathbb{N}$. 
		
$\chi'_{\mathbb{Z}-sum}(G_{-m, n})$ =
$m+n+1$ and $\chi'(G_{-m,n})$ = $m + n$. And in particular, 

\medskip
$\chi'_{\mathbb{Z}-sum}(G_{-n, n})$ = $2n+1$ and $\chi'(G_{-n,n})$ = $2n$. }
\end{theorem}

\begin{proof}\quad 	In Theorem \ref{a}, it is proved that for $m,n\in\mathbb{N}$, $\chi'_{\mathbb{Z}-sum}(G_{-m,n})$ = $m+n+1$. Moreover, $\Delta(G_{-m,n})$ = $m + n$ and therefore by Vizig’s Theorem \ref{1.1}, $m+n \leq \chi'(G_{-m,n}) \leq m+n + 1$. A proper edge coloring of the integral sum graph $G_{-m, n}$ with $m+n$ colors is presented in Theorem \ref{thm Grs 2 chrom}(ii) so that $\chi'(G_{-m, n}) = m + n$, $m,n \in \mathbb{N}$.	Hence we get the result.
\end{proof}	

 In \cite{js24}, it is proved that $G_{n}$, $G_{0, n}$ $G_{-m, n}$ are of class 1 graphs, $m,n \in \mathbb{N}$. Here, we present the result.
  
\begin{theorem} [{\cite[Theorem 4.6]{js24}}] {\rm For positive integers $r$ and $n$, graphs $G_n$, $ G_{0,n} $ and $G_{-r,n}$ are of class 1.}
\end{theorem}	

\begin{proof}
	Using Theorems \ref{thm Gn chrom}, \ref{thm G0n chrom} and \ref{thm Grs 2 chrom}, we get ${\chi}^{'}(G_n) $ = $n-2$ = $\Delta(G_n)$ = degree of the vertex with integral sum labeling 1 in $G_{n}$, ${\chi}^{'}(G_{0,n}) $ = $n$ = $\Delta(G_{0,n})$ = degree of the vertex with integral sum labeling 0 in $G_{0,n}$ and ${\chi}^{'}(G_{-r,n}) $ = $r+n$ = $\Delta(G_{-r,n})$ = degree of the vertex with integral sum labeling 0 in $G_{-r,n}$, $r,n\in\mathbb{N}$. 
\end{proof}	

In Theorem \ref{a}, it is proved that for $2 \leq n \leq 6$, $\chi'_{\mathbb{Z}-sum}(G_{-n,n})$ = $2n+1$ $\neq$ $\chi'(G_{-n,n})$ = $2n$ and thereby, $G_{-n,n}$ is an integral sum graph of class 2 for $n$ = 2,3,4,5,6. For $n\in\mathbb{N}$, $\chi'(G_{0,n})$ = $n$ = $\chi'_{\mathbb{Z}-sum}(G_{0,n})$ and thereby $G_{0, n}$ is an integral sum graph of class 1. And for $n \geq 2$, $\chi'(G_{-1, n})$ = $n+1$ $\neq$ $n+2$ = $\chi'_{\mathbb{Z}-sum}(G_{-1, n})$ and thereby $G_{-1,n}$  is an integral sum graph of class 2 for $n \geq 2$ and $n\in\mathbb{N}$.

\subsection{Open problems on integral sum graphs of class 1 and class 2}

Followings are some open problems for further research. 

\begin{oprm}\cite{vl22} {\rm  Characterize integral sum graphs $G^+(S)$ such that 
\begin{enumerate}
	\item [\rm (1)] $\chi'_{\mathbb{Z}-sum}(G^+(S))$ = order of the graph $G^+(S)$;
	
	\item [\rm (2)]  $G^+(S)$ is a class 1 integral sum graph. i.e., $\chi'_{\mathbb{Z}-sum}(G^+(S))$ = $\chi'(G^+(S))$; 
\begin{enumerate}
	\item [\rm (2a)] $\chi'_{\mathbb{Z}-sum}(G^+(S))$ = $\chi'(G^+(S))$ = $\Delta(G^+(S))$; 
	\item [\rm (2b)] $\chi'_{\mathbb{Z}-sum}(G^+(S))$ = $\chi'(G^+(S))$ = $\Delta(G^+(S))+1$;
\end{enumerate}  	
	\item [\rm (3)] $G^+(S)$ is a class 2 integral sum graph. That is, $\chi'_{\mathbb{Z}-sum}(G^+(S))$ $\neq$ $\chi'(G^+(S))$;
\begin{enumerate}
	\item [\rm (3a)] $\chi'_{\mathbb{Z}-sum}(G^+(S))$ $<$ $\chi'(G^+(S))$ and
	\item [\rm (3b)] $\chi'_{\mathbb{Z}-sum}(G^+(S))$ $>$ $\chi'(G^+(S))$.
	\hfill $\Box$
\end{enumerate}  	
\end{enumerate}  }  
\end{oprm}

\section{Conclusion}

Even though integral sum labeling of graphs seems to be elimentary, we could derive interesting general results, especially different properties of integral sum graphs $G_{n}$, $G_{0,n}$ and $G_{-m,n}$, $m,n\in\mathbb{N}$. We also obtained some interesting elementary results on positive integers and open problems on properties of integral sum graphs which are presented below.

While deriving Theorems \ref{8.7}, \ref{8.8}, \ref{8.9}, \ref{8.11} and Corollary \ref{8.14} in section 4, the following simple properties of natural numbers are obtained.

\begin{theorem} {\rm (Theorem \ref{8.7})} \quad \label{150}  {\rm Let $m,n\in \mathbb{N}$. 
		\begin{enumerate}
			\item [\rm (i)] $(m+n)(m^2+5mn+n^2-1)$ is divisible by 3. 
			
			\item [\rm (ii)] $2(m^3+n^3)$ + $3m^2+4m+n$ is divisible by 6. 
			
			\item [\rm (iii)] $(m+n)(2{(m+n)}^2 + 9(m+n)+13)$ is divisible by 6.    \hfill $\Box$
	\end{enumerate} }
\end{theorem}

\begin{theorem}{\rm (Theorem \ref{8.15})} \label{151} \quad {\rm Let $n\in\mathbb{N}$. \begin{enumerate}
			\item [\rm (i)(a)]  $n(n+1)(7n-4)$ = $n(7n^2+3n-4)$ is divisible by 6.
			\item [\rm (b)]  $n(n+1)(7n+8)$ = $n(7n^2+15n+8)$ is divisible by 6.
			
			\item [\rm (ii)(a)] $n(n+1)(n+2)(7n-3)$ is divisible by 24.
			\item [\rm (b)] $n(n+1)(n+2)(7n+13)$ is divisible by 24.
			\item [\rm (c)] $n(n+1)(7n^2+15n-10)$ is divisible by 24. 
			\item [\rm (d)] $n(n+1)(7n^2+31n+22)$ is divisible by 24. \hfill $\Box$	
	\end{enumerate} }
\end{theorem}

In subsection 5.1, while considering integral sum labeling of Star graph $S_n$ with $V(S_n)$ = $\{u_0$, $v_1$, $v_2$, $\dots$, $v_{n-1}\}$, $d(u_0)$ = $n-1$, $d(v_1)$ = 1 = $d(v_2)$ = $\dots$ = $d(v_{n-1})$ and $f$ as its integral sum labeling, the following interesting sequence from its vertex labeling $f(v_i)$s are obtained.
 $$\{ f(v_i) = t{(d + 1)}^{i-1}\}^{n-1}_{i=1}.$$ 
 It is noted that, for each value of $t$ and $d$, the above sequence is a GP and is a strictly monotonic increasing sequence, $n \geq 2$ and $d,i,n,t\in\mathbb{N}$. 

In subsection 2.6, while studying existence of Hamiltonian cycles that exist in the integral sum graphs $G_{-m,n}$, the following open problem was proposed, $m\in \mathbb{N}_0$ and $n\in \mathbb{N}$.

\begin{oprm} {\rm (Open Problem 1) Find the number of distinct Hamiltonian cycles that exist in the integral sum graph $G_{-m,n}$, $m,n\in \mathbb{N}$. \hfill $\Box$}
\end{oprm}

In subsection 3.3, while studying maximal integral sum graphs of the form $G_{-r,s}$, having at least one of it’s vertices is of degree $r+s$, the following three open problems were proposed, $r,s\in\mathbb{N}$. 

\begin{oprm} {\rm (Open Problem 2)  
\quad Do there exist any other type of maximal integral sum graph of order $n+1$ that have no vertex of degree $n$? \hfill $\Box$}
\end{oprm}

\begin{oprm} 
{\rm (Open Problem 3) \quad Although no proper spanning super graph of a maximal integral sum graph is an integral sum graph, the converse is not known. i.e., Are there proper spanning subgraphs of a maximal integral sum graph that are integral sum graphs. \hfill $\Box$}
\end{oprm}

\begin{oprm} {\rm (Open Problem 4) \quad  For which graphs $G$ and $H$, are the graphs $K_1 * (G \cup H)$ and $K_1 * (G * H)$ are integral sum graphs, when $G$ and $H$ are either both sum graphs or are not both sum graphs? \hfill $\Box$}
\end{oprm}

\vspace{.2cm}
In subsection 5.6, while studying integral sum labeling of graphs $P_k * G_n$, the following conjecture was proposed, $k,n \geq 4$ and $k,n\in\mathbb{N}$.

\begin{conj} 
{\rm (Conjecture 1) \quad For $k,n \geq 4$, $P_k * G_n$ is not an integral sum graph, $k,n\in\mathbb{N}$. \hfill $\Box$  }
\end{conj}

\vspace{.2cm}
In subsection 8.7, while studying decompositions of $K_{n}$  into Fan with a handle, the following conjecture was proposed. In Problem \ref{p8}, it was verified on $K_{4n}$ and $K_{4n+1}$ for $n$ = 1, 2.

\begin{conj} 
{\rm (Conjecture 2) \quad 	For every $n\in\mathbb{N}$, the following decompositions of $K_{4n}$ and $K_{4n+1}$  into Fan with a handle exist. 
		\begin{enumerate}
			\item [\rm (i)]  $K_{4n}$ = $F^*_{3n-1} \cup F^*_{3n-2} \cup F^*_{3n-3} \cup \dots \cup F^*_{n+1} \cup F^*_n$ which is a $(2n, 2)$-CMSD of $K_{4n}$ into Fan with a handle and 
			
			\item [\rm (ii)] $K_{4n+1}$ = $F^*_{3n} \cup F^*_{3n-1} \cup F^*_{3n-2} \cup \dots \cup F^*_{n+2} \cup F^*_{n+1}$ which is a $(2n+2, 2)$-CMSD of $K_{4n+1}$  into Fan with a handle. \hfill    $\Box$
	\end{enumerate} }
\end{conj}

\vspace{.2cm}
In subsection 10.4, while studying characterisation of class 1 and 2 integral sum graphs, the following open problem was proposed. 

\begin{oprm} 
{\rm (Open Problem 5) \quad  Characterize integral sum graphs $G^+(S)$ such that 
		\begin{enumerate}
			\item [\rm (1)] $\chi'_{\mathbb{Z}-sum}(G^+(S))$ = order of the graph $G^+(S)$;
			
			\item [\rm (2)]  $G^+(S)$ is a class 1 integral sum graph. 
			\\
			i.e., $\chi'_{\mathbb{Z}-sum}(G^+(S))$ = $\chi'(G^+(S))$; 
			\begin{enumerate}
				\item [\rm (2a)] $\chi'_{\mathbb{Z}-sum}(G^+(S))$ = $\chi'(G^+(S))$ = $\Delta(G^+(S))$; 
				\item [\rm (2b)] $\chi'_{\mathbb{Z}-sum}(G^+(S))$ = $\chi'(G^+(S))$ = $\Delta(G^+(S))+1$;
			\end{enumerate}  	
			\item [\rm (3)] $G^+(S)$ is a class 2 integral sum graph. 
			\\
			i.e., $\chi'_{\mathbb{Z}-sum}(G^+(S))$ $\neq$ $\chi'(G^+(S))$;
			\begin{enumerate}
				\item [\rm (3a)] $\chi'_{\mathbb{Z}-sum}(G^+(S))$ $<$ $\chi'(G^+(S))$ and
				\item [\rm (3b)] $\chi'_{\mathbb{Z}-sum}(G^+(S))$ $>$ $\chi'(G^+(S))$.
				\hfill $\Box$
			\end{enumerate}  	
	\end{enumerate}  }  
\end{oprm} 

The authors feel that one can do a lot of research in this area.

\vspace{.2cm}
\noindent
{\bf Acknowledgement} The first author expresses his sincere thanks to the Central University of Kerala, Periye - 671 316, Kasaragod, Kerala, India and St. Jude's College, Thoothoor - 629 176, Kanyakumari District, Tamil Nadu, India for providing facilities to do this research work.

\end{document}